\documentclass[10pt]{amsart}
\usepackage{amsmath}
\usepackage{amsxtra}
\usepackage{amscd}
\usepackage{amsthm}
\usepackage{amsfonts}
\usepackage{amssymb}
\usepackage{eucal}
\usepackage[all]{xy}
\usepackage{graphicx}
\usepackage{comment}
\usepackage{epsfig}
\usepackage{psfrag}
\usepackage{mathrsfs}
\usepackage{amscd}
\usepackage{rotating}
\usepackage{lscape}
\usepackage{amsbsy}
\usepackage{verbatim}
\usepackage{moreverb}
\usepackage{url}

\newtheorem{cor}[subsubsection]{Corollary}
\newtheorem{lem}[subsubsection]{Lemma}
\newtheorem{prop}[subsubsection]{Proposition}

\newtheorem{conj}[subsubsection]{Conjecture}
\newtheorem{thm}[subsubsection]{Theorem}

\newtheorem{defn}[subsubsection]{Definition}

\theoremstyle{remark}
\newtheorem{rem}[subsubsection]{Remark}

\theoremstyle{remark}

\newcommand{\thmref}[1]{Theorem~\ref{#1}}

\newcommand{\secref}[1]{Sect.~\ref{#1}}
\newcommand{\lemref}[1]{Lemma~\ref{#1}}
\newcommand{\propref}[1]{Proposition~\ref{#1}}
\newcommand{\corref}[1]{Corollary~\ref{#1}}
\newcommand{\conjref}[1]{Conjecture~\ref{#1}}
\newcommand{\remref}[1]{Remark~\ref{#1}}

\numberwithin{equation}{section}

\newcommand{\nc}{\newcommand}
\nc{\renc}{\renewcommand}
\nc{\ssec}{\subsection}
\nc{\sssec}{\subsubsection}
\nc{\on}{\operatorname}

\nc\ol{\overline}
\nc\wt{\widetilde}
\nc\tboxtimes{\wt{\boxtimes}}
\nc\tstar{\wt{\star}}
\nc{\alp}{\alpha}

\nc{\ZZ}{{\mathbb Z}}
\nc{\NN}{{\mathbb N}}
\nc{\OO}{{\mathbb O}}
\renc{\SS}{{\mathbb S}}
\nc{\DD}{{\mathbb D}}
\nc{\GG}{{\mathbb G}}

\nc{\Fq}{{\mathbb F}_q}
\nc{\Fqb}{\ol{{\mathbb F}_q}}
\nc{\Ql}{\ol{{\mathbb Q}_\ell}}
\nc{\id}{\text{id}}
\nc\X{\mathcal X}

\nc{\Hom}{\on{Hom}}
\nc{\Lie}{\on{Lie}}
\nc{\Loc}{\on{Loc}}
\nc{\Pic}{\on{Pic}}
\nc{\Bun}{\on{Bun}}
\nc{\IC}{\on{IC}}
\nc{\Aut}{\on{Aut}}
\nc{\rk}{\on{rk}}
\nc{\Sh}{\on{Sh}}
\nc{\Perv}{\on{Perv}}
\nc{\pos}{{\on{pos}}}
\nc{\Conv}{\on{Conv}}
\nc{\Sph}{\on{Sph}}
\nc{\Sym}{\on{Sym}}
\nc{\BunBb}{\overline{\Bun}_B}
\nc{\BunNb}{\overline{\Bun}_N}
\nc{\BunTb}{\overline{\Bun}_T}
\nc{\BunBbm}{\overline{\Bun}_{B^-}}
\nc{\BunBbel}{\overline{\Bun}_{B,el}}
\nc{\BunBbmel}{\overline{\Bun}_{B^-,el}}
\nc{\Buno}{\overset{o}{\Bun}}
\nc{\BunPb}{{\overline{\Bun}_P}}
\nc{\BunBM}{\Bun_{B(M)}}
\nc{\BunBMb}{\overline{\Bun}_{B(M)}}
\nc{\BunPbw}{{\widetilde{\Bun}_P}}
\nc{\BunBP}{\widetilde{\Bun}_{B,P}}
\nc{\GUb}{\overline{G/U}}
\nc{\GUPb}{\overline{G/U(P)}}

\nc{\Hhom}{\underline{\on{Hom}}}
\nc\syminfty{\on{Sym}^{\infty}}
\nc\lal{\ol{\lambda}}
\nc\xl{\ol{x}}
\nc\thl{\ol{\theta}}
\nc\nul{\ol{\nu}}
\nc\mul{\ol{\mu}}
\nc\Sum\Sigma
\nc{\oX}{\overset{\circ}{X}{}}
\nc{\hl}{\overset{\leftarrow}h{}}
\nc{\hr}{\overset{\rightarrow}h{}}
\nc{\M}{{\mathcal M}}
\nc{\N}{{\mathcal N}}
\nc{\F}{{\mathcal F}}
\nc{\D}{{\mathcal D}}
\nc{\Q}{{\mathcal Q}}
\nc{\Y}{{\mathcal Y}}
\nc{\G}{{\mathcal G}}
\nc{\E}{{\mathcal E}}
\nc{\CalC}{{\mathcal C}}
\nc\Dh{\widehat{\D}}

\nc{\C}{{\mathcal C}}
\nc{\K}{{\mathcal K}}
\renewcommand{\H}{{\mathcal H}}

\nc{\T}{{\mathcal T}}
\nc{\V}{{\mathcal V}}
\renc{\P}{{\mathcal P}}
\nc{\A}{{\mathcal A}}
\nc{\B}{{\mathcal B}}
\nc{\U}{{\mathcal U}}

\nc{\Gr}{{\on{Gr}}}

\nc{\qo}{{\mathfrak q}}
\nc{\po}{{\mathfrak p}}
\nc{\s}{{\mathfrak s}}
\nc\w{\text{w}}

\renewcommand{\mod}{{\on{-mod}}}

\nc\mathi\iota
\nc\Spec{\on{Spec}}
\nc\Mod{\on{Mod}}
\nc{\tw}{\widetilde{\mathfrak t}}
\nc{\pw}{\widetilde{\mathfrak p}}
\nc{\qw}{\widetilde{\mathfrak q}}
\nc{\jw}{\widetilde j}

\nc{\grb}{\overline{\Gr_{X^{\fset}}}}
\nc{\I}{\mathcal I}

\nc{\lambdach}{{\check\lambda}}
\nc{\Lambdach}{{\check\Lambda}{}}
\nc{\much}{{\check\mu}}
\nc{\omegach}{{\check\omega}}
\nc{\nuch}{{\check\nu}}
\nc{\etach}{{\check\eta}}
\nc{\alphach}{{\check\alpha}}

\nc{\rhoch}{{\check\rho}}
\nc{\ch}{{\check\fh}}

\nc{\Hb}{\overline{\H}}

\nc{\BA}{{\mathbb{A}}}
\nc{\BC}{{\mathbb{C}}}
\nc{\BE}{{\mathbb{E}}}
\nc{\BG}{{\mathbb{G}}}
\nc{\BM}{{\mathbb{M}}}
\nc{\BO}{{\mathbb{O}}}
\nc{\BD}{{\mathbb{D}}}
\nc{\BN}{{\mathbb{N}}}
\nc{\BP}{{\mathbb{P}}}
\nc{\BR}{{\mathbb{R}}}
\nc{\BZ}{{\mathbb{Z}}}
\nc{\BS}{{\mathbb{S}}}
\nc{\BV}{{\mathbb{V}}}

\nc{\CA}{{\mathcal{A}}}
\nc{\CB}{{\mathcal{B}}}

\nc{\CE}{{\mathcal{E}}}
\nc{\CF}{{\mathcal{F}}}
\nc{\CH}{{\mathcal{H}}}

\nc{\CL}{{\mathcal{L}}}
\nc{\CC}{{\mathcal{C}}}
\nc{\CG}{{\mathcal{G}}}
\nc{\CM}{{\mathcal{M}}}
\nc{\CN}{{\mathcal{N}}}
\nc{\CK}{{\mathcal{K}}}
\nc{\CO}{{\mathcal{O}}}
\nc{\CP}{{\mathcal{P}}}
\nc{\CQ}{{\mathcal{Q}}}
\nc{\CR}{{\mathcal{R}}}
\nc{\CS}{{\mathcal{S}}}
\nc{\CT}{{\mathcal{T}}}
\nc{\CU}{{\mathcal{U}}}
\nc{\CV}{{\mathcal{V}}}
\nc{\CW}{{\mathcal{W}}}
\nc{\CX}{{\mathcal{X}}}
\nc{\CY}{{\mathcal{Y}}}
\nc{\CZ}{{\mathcal{Z}}}
\nc{\CI}{{\mathcal{I}}}

\nc{\csM}{{\check{\mathcal A}}{}}
\nc{\oM}{{\overset{\circ}{\mathcal M}}{}}
\nc{\obM}{{\overset{\circ}{\mathbf M}}{}}
\nc{\oCA}{{\overset{\circ}{\mathcal A}}{}}
\nc{\obA}{{\overset{\circ}{\mathbf A}}{}}
\nc{\ooM}{{\overset{\circ}{M}}{}}
\nc{\osM}{{\overset{\circ}{\mathsf M}}{}}
\nc{\vM}{{\overset{\bullet}{\mathcal M}}{}}
\nc{\nM}{{\underset{\bullet}{\mathcal M}}{}}
\nc{\oD}{{\overset{\circ}{\mathcal D}}{}}
\nc{\obC}{{\overset{\circ}{\mathbf C}}{}}
\nc{\obD}{{\overset{\circ}{\mathbf D}}{}}
\nc{\oA}{{\overset{\circ}{\mathbb A}}{}}
\nc{\op}{{\overset{\bullet}{\mathbf p}}{}}
\nc{\oU}{{\overset{\bullet}{\mathcal U}}{}}
\nc{\oZ}{{\overset{\circ}{\mathcal Z}}{}}
\nc{\ofZ}{{\overset{\circ}{\mathfrak Z}}{}}
\nc{\oF}{{\overset{\circ}{\fF}}}

\nc{\fa}{{\mathfrak{a}}}
\nc{\fb}{{\mathfrak{b}}}
\nc{\fc}{{\mathfrak{c}}}
\nc{\fd}{{\mathfrak{d}}}
\nc{\ff}{{\mathfrak{f}}}
\nc{\fg}{{\mathfrak{g}}}
\nc{\fgl}{{\mathfrak{gl}}}
\nc{\fh}{{\mathfrak{h}}}
\nc{\fj}{{\mathfrak{j}}}
\nc{\fl}{{\mathfrak{l}}}
\nc{\fm}{{\mathfrak{m}}}
\nc{\fn}{{\mathfrak{n}}}
\nc{\fu}{{\mathfrak{u}}}
\nc{\fr}{{\mathfrak{r}}}
\nc{\fs}{{\mathfrak{s}}}
\nc{\ft}{{\mathfrak{t}}}
\nc{\fz}{{\mathfrak{z}}}
\nc{\fsl}{{\mathfrak{sl}}}
\nc{\hsl}{{\widehat{\mathfrak{sl}}}}
\nc{\hgl}{{\widehat{\mathfrak{gl}}}}
\nc{\hg}{{\widehat{\mathfrak{g}}}}
\nc{\chg}{{\widehat{\mathfrak{g}}}{}^\vee}
\nc{\hn}{{\widehat{\mathfrak{n}}}}
\nc{\chn}{{\widehat{\mathfrak{n}}}{}^\vee}

\nc{\fA}{{\mathfrak{A}}}
\nc{\fB}{{\mathfrak{B}}}
\nc{\fD}{{\mathfrak{D}}}
\nc{\fE}{{\mathfrak{E}}}
\nc{\fF}{{\mathfrak{F}}}
\nc{\fG}{{\mathfrak{G}}}
\nc{\fK}{{\mathfrak{K}}}
\nc{\fL}{{\mathfrak{L}}}
\nc{\fM}{{\mathfrak{M}}}
\nc{\fN}{{\mathfrak{N}}}
\nc{\fP}{{\mathfrak{P}}}
\nc{\fU}{{\mathfrak{U}}}
\nc{\fV}{{\mathfrak{V}}}
\nc{\fZ}{{\mathfrak{Z}}}

\nc{\bb}{{\mathbf{b}}}
\nc{\bc}{{\mathbf{c}}}
\nc{\bd}{{\mathbf{d}}}
\nc{\bbf}{{\mathbf{f}}}
\nc{\be}{{\mathbf{e}}}
\nc{\bg}{{\mathbf{g}}}
\nc{\bi}{{\mathbf{i}}}
\nc{\bj}{{\mathbf{j}}}
\nc{\bn}{{\mathbf{n}}}
\nc{\bo}{{\mathbf{o}}}
\nc{\bp}{{\mathbf{p}}}
\nc{\bq}{{\mathbf{q}}}
\nc{\bt}{{\mathbf{t}}}
\nc{\bu}{{\mathbf{u}}}
\nc{\bv}{{\mathbf{v}}}
\nc{\bx}{{\mathbf{x}}}
\nc{\bs}{{\mathbf{s}}}
\nc{\by}{{\mathbf{y}}}
\nc{\bw}{{\mathbf{w}}}
\nc{\bA}{{\mathbf{A}}}
\nc{\bK}{{\mathbf{K}}}
\nc{\bB}{{\mathbf{B}}}
\nc{\bC}{{\mathbf{C}}}
\nc{\bG}{{\mathbf{G}}}
\nc{\bD}{{\mathbf{D}}}
\nc{\bH}{{\mathbf{H}}}
\nc{\bM}{{\mathbf{M}}}
\nc{\bN}{{\mathbf{N}}}
\nc{\bO}{{\mathbf{O}}}
\nc{\bT}{{\mathbf{T}}}
\nc{\bV}{{\mathbf{V}}}
\nc{\bW}{{\mathbf{W}}}
\nc{\bX}{{\mathbf{X}}}
\nc{\bZ}{{\mathbf{Z}}}
\nc{\bS}{{\mathbf{S}}}

\nc{\sA}{{\mathsf{A}}}
\nc{\sB}{{\mathsf{B}}}
\nc{\sC}{{\mathsf{C}}}
\nc{\sD}{{\mathsf{D}}}
\nc{\sF}{{\mathsf{F}}}
\nc{\sG}{{\mathsf{G}}}
\nc{\sK}{{\mathsf{K}}}
\nc{\sM}{{\mathsf{M}}}
\nc{\sO}{{\mathsf{O}}}
\nc{\sW}{{\mathsf{W}}}
\nc{\sQ}{{\mathsf{Q}}}
\nc{\sP}{{\mathsf{P}}}
\nc{\sV}{{\mathsf{V}}}
\nc{\sS}{{\mathsf{S}}}
\nc{\sT}{{\mathsf{T}}}
\nc{\sZ}{{\mathsf{Z}}}
\nc{\sll}{{\mathsf{l}}}
\nc{\sr}{{\mathsf{r}}}
\nc{\bk}{{\mathsf{k}}}
\nc{\sg}{{\mathsf{g}}}

\nc{\sff}{{\mathsf{f}}}
\nc{\sfb}{{\mathsf{b}}}
\nc{\sfc}{{\mathsf{c}}}
\nc{\sd}{{\mathsf{d}}}
\nc{\se}{{\mathsf{e}}}

\nc{\BK}{{\bar{K}}}

\nc{\tA}{{\widetilde{\mathbf{A}}}}
\nc{\tB}{{\widetilde{\mathcal{B}}}}
\nc{\tg}{{\widetilde{\mathfrak{g}}}}
\nc{\tG}{{\widetilde{G}}}
\nc{\TM}{{\widetilde{\mathbb{M}}}{}}
\nc{\tO}{{\widetilde{\mathsf{O}}}{}}
\nc{\tU}{{\widetilde{\mathfrak{U}}}{}}
\nc{\TZ}{{\tilde{Z}}}
\nc{\tx}{{\tilde{x}}}
\nc{\tbv}{{\tilde{\bv}}}
\nc{\tfP}{{\widetilde{\mathfrak{P}}}{}}
\nc{\tz}{{\tilde{\zeta}}}
\nc{\tmu}{{\tilde{\mu}}}

\nc{\urho}{\underline{\rho}}
\nc{\uB}{\underline{B}}
\nc{\uC}{{\underline{\mathbb{C}}}}
\nc{\ui}{\underline{i}}
\nc{\uj}{\underline{j}}
\nc{\ofP}{{\overline{\mathfrak{P}}}}
\nc{\oB}{{\overline{\mathcal{B}}}}
\nc{\og}{{\overline{\mathfrak{g}}}}
\nc{\oI}{{\overline{I}}}

\nc{\eps}{\varepsilon}
\nc{\hrho}{{\hat{\rho}}}

\nc{\one}{{\mathbf{1}}}
\nc{\two}{{\mathbf{t}}}

\nc{\Rep}{{\mathop{\operatorname{\rm Rep}}}}
\nc{\Tot}{{\mathop{\operatorname{\rm Tot}}}}
\nc{\Ker}{{\mathop{\operatorname{\rm Ker}}}}
\nc{\Hilb}{{\mathop{\operatorname{\rm Hilb}}}}
\nc{\End}{{\mathop{\operatorname{\rm End}}}}
\nc{\Ext}{{\mathop{\operatorname{\rm Ext}}}}
\nc{\CHom}{{\mathop{\operatorname{{\mathcal{H}}\it om}}}}
\nc{\GL}{{\mathop{\operatorname{\rm GL}}}}
\nc{\gr}{{\mathop{\operatorname{\rm gr}}}}
\nc{\Id}{{\mathop{\operatorname{\rm Id}}}}
\nc{\de}{{\mathop{\operatorname{\rm def}}}}
\nc{\length}{{\mathop{\operatorname{\rm length}}}}
\nc{\supp}{{\mathop{\operatorname{\rm supp}}}}

\nc{\Cliff}{{\mathsf{Cliff}}}
\nc{\Fl}{\on{Fl}}
\nc{\Fib}{{\mathsf{Fib}}}
\nc{\Coh}{{\on{Coh}}}
\nc{\QCoh}{{\on{QCoh}}}
\nc{\IndCoh}{{\on{IndCoh}}}
\nc{\FCoh}{{\mathsf{FCoh}}}

\nc{\reg}{{\text{\rm reg}}}

\nc{\cplus}{{\mathbf{C}_+}}
\nc{\cminus}{{\mathbf{C}_-}}
\nc{\cthree}{{\mathbf{C}_*}}
\nc{\Qbar}{{\bar{Q}}}
\nc\Eis{\on{Eis}}
\nc\Eisb{\ol\Eis{}}
\nc\Eisr{\on{Eis}^{rat}{}}
\nc\wh{\widehat}
\nc{\Def}{\on{Def_{\check{\fb}}(E)}}
\nc{\barZ}{\overline{Z}{}}
\nc{\barbarZ}{\overline{\barZ}{}}
\nc{\barpi}{\overline\pi}
\nc{\barbarpi}{\overline\barpi}
\nc{\barpip}{\overline\pi{}^+}
\nc{\barpim}{\overline\pi{}^-}

\nc{\hattimes}{\wh\otimes}

\nc{\bh}{{\bar{h}}}
\nc{\bOmega}{{\overline{\Omega(\check \fn)}}}

\nc{\seq}[1]{\stackrel{#1}{\sim}}

\nc{\cT}{{\check{T}}}
\nc{\cG}{{\check{G}}}
\nc{\cM}{{\check{M}}}
\nc{\cB}{{\check{B}}}
\nc{\cP}{{\check{P}}}

\nc{\fp}{{\mathfrak{p}}}
\nc{\fq}{{\mathfrak{q}}}

\nc{\ct}{{\check{\mathfrak t}}}
\nc{\cg}{{\check{\fg}}}
\nc{\cb}{{\check{\fb}}}
\nc{\cn}{{\check{\fn}}}
\nc{\cp}{{\check{\fp}}}
\nc{\cm}{{\check{\fm}}}

\nc{\sfp}{{\mathsf{p}}}
\nc{\sfq}{{\mathsf{q}}}

\nc{\cLambda}{{\check\Lambda}}

\nc{\cla}{{\check\lambda}}
\nc{\cmu}{{\check\mu}}
\nc{\cnu}{{\check\nu}}
\nc{\ceta}{{\check\eta}}

\nc{\DefbE}{{\on{Def}_{\cB}(E_\cT)}}

\nc{\imathb}{{\ol{\imath}}}
\nc{\rlr}{\overset{\longrightarrow}{\underset{\longrightarrow}\longleftarrow}}

\nc{\oBun}{\overset{\circ}\Bun}
\nc{\LocSys}{\on{LocSys}}
\nc{\BunBbb}{\ol{\ol{Bun}}_B}
\nc{\BunBr}{\Bun_B^{rat}}
\nc{\BunBrsg}{\Bun_B^{rat,\on{s.g.}}}
\nc{\BunBrp}{\Bun_B^{rat,polar}}
\nc{\BunBrpbg}{\Bun_B^{rat,polar,\on{b.g.}}}
\nc{\BunBrpsg}{\Bun_B^{rat,polar,\on{s.g.}}}
\nc{\BunTrp}{\Bun_T^{rat,polar}}
\nc{\BunTrpbg}{\Bun_T^{rat,polar,\on{b.g.}}}
\nc{\BunTrpsg}{\Bun_T^{rat,polar,\on{s.g.}}}
\nc{\BunNr}{\Bun_N^{rat}}
\nc{\BunNre}{\Bun_N^{enh,rat}}
\nc{\BunTr}{\Bun_T^{rat}}
\nc{\Vect}{\on{Vect}}
\nc{\Whit}{\on{Whit}}
\nc{\bTb}{\ol{\on{CT}}}
\nc{\Ran}{\on{Ran}}
\nc{\bTr}{\on{CT}^{rat}{}}
\nc\jmathr{\jmath^{rat}{}}
\nc{\ux}{\underline{x}}
\nc{\clambda}{{\check\lambda}}
\nc{\calpha}{{\check\alpha}}
\nc{\ind}{{\mathbf{ind}}}
\nc{\oblv}{{\mathbf{oblv}}}
\nc{\StinftyCat}{\on{DGCat}}
\nc{\inftygroup}{\infty\on{-Grpd}}

\nc{\fset}{\on{fSet}}
\nc{\LocSysG}{\LocSys_{\cG}}
\nc{\Sing}{\on{Sing}}
\nc{\dr}{{\on{dR}}}
\nc{\Dmod}{\on{D-mod}}
\nc{\Ind}{\on{Ind}}
\nc{\Sat}{\on{Sat}}
\nc{\Ho}{\on{Ho}}
\nc{\Res}{\on{Res}}
\nc{\sotimes}{\overset{!}\otimes}
\nc{\mmod}{{\on{-}}{\mathbf{mod}}}
\nc{\Maps}{\on{Maps}}
\nc{\CMaps}{{\mathcal Maps}}
\nc{\bMaps}{{\mathbf{Maps}}}

\nc{\dgSch}{\on{DGSch}}
\nc{\Sch}{\on{Sch}}
\nc{\affdgSch}{\on{DGSch}^{\on{aff}}}
\nc{\affSch}{\on{Sch}^{\on{aff}}}

\nc{\HC}{\CH\bC}

\nc{\csupp}{\supp}
\nc{\Arth}{\on{Arth}}
\nc{\ArthG}{{\on{Arth}_\cG}}
\nc{\ul}{\underline}
\nc{\Alg}{\on{-Alg}}

\begin{document}

\title[Singular support of coherent sheaves]{Singular support of coherent sheaves \\
and the geometric Langlands conjecture}

\author{D.~Arinkin and D.~Gaitsgory}

\email{arinkin@math.wisc.edu, gaitsgde@math.harvard.edu} 


\date{\today}

\date{\today}

\begin{abstract}
We define the notion of singular support of a coherent sheaf on a quasi-smooth derived (DG) scheme
or Artin stack, where ``quasi-smooth" means that it is a locally complete intersection in the derived sense.
This develops the idea of  ``cohomological" support of coherent sheaves on a locally complete intersection scheme introduced
by D.~Benson, S.~B.~Iyengar, and H.~Krause.
We study the behaviour of singular support under the direct and inverse image functors for coherent sheaves. 

We use the theory of singular support of coherent sheaves to formulate
the categorical Geometric Langlands Conjecture. We verify that it passes natural consistency tests:
it is compatible with the Geometric Satake equivalence and with the Eisenstein series functors. 
The latter compatibility is particularly important, as it fails in the original ``naive" form of the conjecture.
 
\end{abstract}

\maketitle

\tableofcontents

\section{Introduction}

\ssec{What we are trying to do}

\sssec{}

Let \footnote{MSC: 14F05, 14H60. Key words: cohomological support, geometric Langlands program.}
$G$ be a connected reductive group, and $X$ a smooth, connected and complete curve over a ground field
$k$, assumed algebraically closed and of characteristic $0$. Let $\Bun_G(X)$
be the moduli stack of $G$-bundles on $X$, and consider the (DG) category $\Dmod(\Bun_G(X))$. 

\medskip

The goal of the (classical, global and unramified) geometric Langlands program is to express the category $\Dmod(\Bun_G(X))$
in terms of the Langlands dual group $\cG$; more precisely, in terms of the (DG) category of quasi-coherent sheaves
on the (DG) stack $\LocSysG$ of local systems on $X$ with respect to $\cG$. 

\medskip

The naive guess, referred to by A.~Beilinson and V.~Drinfeld for a number of years as ``the best hope," says that the category
$\Dmod(\Bun_G(X))$ is simply equivalent to $\QCoh(\LocSysG)$. For example, this is indeed the case when $G$ is torus,
and the required equivalence is given by the Fourier transform of \cite{L1,L2} and \cite{Ro1,Ro2}. 

\medskip

However, the ``best hope" does not hold for groups other than the torus. For example, it fails in the simplest case
of $G=SL_2$ and $X=\BP^1$.  An explicit calculation showing this can be found in \cite{La}. 

\sssec{}

There is a heuristic reason for the failure of the ``best hope": 

\medskip

In the classical theory of automorphic forms one expects
that automorphic representations are parametrized not just by Galois representations, but by Arthur parameters, 
i.e., in addition to a homomorphism from the Galois group to $\cG$, one needs to specify a nilpotent element in $\cg$
centralized by the image of the Galois group. 

\medskip

In addition, there has been a general understanding that the presence of the commuting nilpotent element must
be ``cohomological in nature." As an incarnation of this, for automorphic representations realized in the cohomology
of Shimura varieties, the nilpotent element in question acts as the Lefschetz operator of multiplication by the Chern
class of the corresponding line bundle. 

\medskip

So, in the geometric theory one has been faced with the challenge of how to modify the Galois side, i.e.,
the category $\QCoh(\LocSysG)$, with the hint being that the solution should come from considering 
the stack of pairs $(\sigma,A)$, where $\sigma$ is a $\cG$-local system, and $A$ its endomorphism, i.e., a horizontal
section of the associated local system $\cg_\sigma$. 

\medskip

The general feeling, shared by many people who have looked at this problem, was that the sought-for modification
has to do with the fact that the DG stack $\LocSysG$ is not smooth. I.e., we need to modify the category $\QCoh(\LocSysG)$
by taking into account the singularities of $\LocSysG$.

\sssec{}

The goal of the present paper is to provide such a modification, and 
to formulate the appropriately modified version
of the ``best hope."

\medskip

In fact, one does not have to look very far for the possibilities to ``tweak" the category 
$\QCoh(\LocSysG)$. Recall that for any
reasonable algebraic DG stack $\CZ$, the category $\QCoh(\CZ)$ is compactly
generated by its subcategory $\QCoh(\CZ)^{\on{perf}}$ of perfect complexes. 

\medskip

Now, if $\CZ$ is non-smooth,
one can enlarge the category $\QCoh(\CZ)^{\on{perf}}$ to that of coherent complexes, denoted $\Coh(\CZ)$.
By passing to the ind-completion, one obtains the category $\IndCoh(\CZ)$, studied in \cite{DrGa0} and 
\cite{IndCoh}. 

\medskip

(We emphasize that the difference between $\QCoh(\CZ)$ and $\IndCoh(\CZ)$ is not
a ``stacky" phenomenon: it is caused not by automorphisms of points of $\CZ$, 
but rather by its singularities.)

\medskip

Thus, one can try $\IndCoh(\LocSysG)$ as a candidate for the Galois side of Geometric Langlands. However,
playing with the example of $X=\BP^1$ shows that, whereas $\QCoh(\LocSysG)$ was ``too small" to be 
equivalent to $\Dmod(\Bun_G)$, the category $\IndCoh(\LocSysG)$ is ``too large." 

\medskip

So, a natural guess for the category on the Galois side is that it should be a (full) subcategory of
$\IndCoh(\LocSysG)$ that contains $\QCoh(\LocSysG)$. This is indeed the shape of the answer 
that we will propose. Our goal is to describe the corresponding subcategory.

\sssec{}

Let us go back to the situation of a nice (=QCA, in the terminology of \cite{DrGa0}) algebraic DG stack $\CZ$. In 
general, it is not so clear how to describe categories that lie between $\QCoh(\CZ)$ and $\IndCoh(\CZ)$. 
However, the situation is more manageable when $\CZ$ is \emph{quasi-smooth}, that is, if its singularities are 
modelled, locally in the smooth topology, by a complete intersection. 

\medskip

For every point $z\in \CZ$ we can consider the derived cotangent space $T_z^*(\CZ)$, which is a complex of
vector spaces lying in cohomological degrees $\leq 1$. The assumption that $\CZ$ is quasi-smooth is 
equivalent to that the cohomologies vanish in degrees $<-1$ (smoothness is equivalent to the vanishing
of $H^{-1}$ as well). 

\medskip

It is easy to see that the assignment $z\mapsto H^{-1}\left(T_z^*(\CZ)\right)$ forms a well-defined 
\emph{classical}\footnote{``classical" as opposed to ``DG''} stack,
whose projection to $\CZ$ is affine, and which carries a canonical action of $\BG_m$ by dilations. 
We will denote this stack by $\Sing(\CZ)$. 

\medskip

We are going to show (see \secref{s:sing}) that to every $\CF\in \IndCoh(\CZ)$ one can assign its singular support,
denoted $\on{SingSupp}(\CF)$, which is a conical Zariski-closed subset in $\Sing(\CZ)$. 
It is easy to see that for $\CF\in \QCoh(\CZ)$, its singular support is contained in the zero-section of
$\Sing(\CZ)$. It is less obvious, but still true (see \thmref{t:zero sect}) that if $\CF$ is such that its singular support
is the zero-section, then $\CF$ belongs to $\QCoh(\CZ)\subset \IndCoh(\CZ)$. Thus, the singular support of
an object of $\IndCoh(\CZ)$ exactly measures the degree to which this object does not belong to
$\QCoh(\CZ)$. 

\medskip

For a fixed conical Zariski-closed subset $Y\subset \Sing(\CZ)$, we can consider the full subcategory 
$\IndCoh_Y(\CZ)\subset \IndCoh(\CZ)$ consisting of those objects, whose singular support lies in $\CY$. 
(The paper \cite{St} implies that the assignment
$Y\mapsto \IndCoh_Y(\CZ)$ establishes a bijection between subsets $Y$ as above and  
full subcategories of $\IndCoh(\CZ)$ satisfying certain
natural conditions.)  

\medskip

We should mention that the procedure of assigning the singular support to an object $\CF\in \IndCoh(\CZ)$
is cohomological in nature: we read if off (locally) from the action of the algebra of Hochschild cochains
(on smooth affine charts of $\CZ$) on our object. This loosely corresponds to the cohomological
nature of the Arthur parameter. To the best of our knowledge, the general idea of using cohomological operators to 
define support is due to D.~Benson, S.~B.~Iyengar, and H.~Krause \cite{Kr}. 

\medskip

The singular support for (ind)coherent complexes on a quasi-smooth DG scheme is similar to the 
singular support (also known as the ``characteristic variety") of a coherent $D$-module on a smooth variety, and also,
perhaps in a more remote way, to the singular support of constructible sheaves on a manifold.
Because
of this analogy, we use the name ``singular support" (rather than, say, ``cohomological support") for the support 
of ind-coherent sheaves. 

\sssec{}

Returning to the Galois side of Geometric Langlands, we note that the DG stack $\LocSysG$ is indeed
quasi-smooth. Moreover, we note that the corresponding stack $\Sing(\LocSysG)$ classifies
pairs $(\sigma,A)$, i.e., Arthur parameters. Thus, we rename
$$\ArthG:=\Sing(\LocSysG).$$

\medskip

Our candidate for a category lying between $\QCoh(\LocSysG)$ and $\IndCoh(\LocSysG)$
corresponds to a particular closed subset of $\ArthG$. Namely, let 
$$\on{Nilp}_{glob}\subset \ArthG,$$
be the subset of pairs $(\sigma,A)$ with nilpotent $A$. 

\medskip

So, we propose the following modified version of the ``best hope":

\begin{conj}  \label{c:main preview}
There exists an equivalence of categories
$$\Dmod(\Bun_G)\simeq \IndCoh_{\on{Nilp}_{glob}}(\LocSysG).$$
\end{conj}

Thus, given $\CM\in \Dmod(\Bun_G)$, one cannot really speak of ``the Arthur parameter" corresponding to $\CM$.
However, one can specify a conical closed subset of $\on{Nilp}_{glob}$ over which $\CM$ is supported.

\sssec{}  \label{sss:additional}

Geometric Langlands correspondence is more than simply \emph{an} equivalence of categories
as in Conjecture \ref{c:main preview}. Rather, the sought-for equivalence must satisfy  a number
of compatibility conditions. 

\medskip

Two of these conditions are discussed in the present paper: compatibility with the Geometric Satake
Equivalence and compatibility with the Eisenstein series functors.

\medskip

Two other conditions have to do with the description of the Whittaker D-module on the automorphic side, 
and of the construction of automorphic D-modules by localization from Kac-Moody representations. 

\medskip

On the Galois
side, the Whittaker D-module is supposed to correspond to the structure sheaf on $\LocSysG$.
The localization functor should correspond to the direct image functor with respect to the map to $\LocSysG$ from the scheme of opers. 

\medskip

Both of the latter procedures are insensitive to the singular aspects of $\LocSysG$, which is why we do not
discuss them in this paper. However, we plan to revisit these objects in a subsequent publication.

\sssec{}

Conjecture \ref{c:main preview} contains the following statement (which can actually be proved unconditionally):

\begin{conj}
The monoidal category $\QCoh(\LocSysG)$ acts on $\Dmod(\Bun_G)$.
\end{conj}

The above corollary allows one to take fibers of the category $\Dmod(\Bun_G)$ at 
points of $\LocSysG$, or more generally $S$-points for any test scheme $S$. Namely, 
for $\sigma:S\to \LocSysG$ we set
$$\Dmod(\Bun_G)_\sigma:=\QCoh(S)\underset{\QCoh(\LocSysG)}\otimes \Dmod(\Bun_G).$$

In particular, for a $k$-point $\sigma$ of $\LocSysG$, the corresponding category
$\Dmod(\Bun_G)_\sigma$ is that of Hecke eigensheaves with eigenvalue $\sigma$.

\medskip

However, we do not know (and have no reasons to believe) that 
this procedure can be refined to $\ArthG$. In other words, we do not expect that
the category $\QCoh(\ArthG)$ should act on $\Dmod(\Bun_G)$ and that it should be possible to
take the fiber of $\Dmod(\Bun_G)$ at a specified Arthur parameter. 

\sssec{}

Let $\CZ$ be a quasi-smooth DG stack and let $Y\subset \Sing(\CZ)$ be a conical Zariski-closed subset. If $Y$ contains the zero-section of $\Sing(\CZ)$, then 
the category $\IndCoh_Y(\CZ)$ contains $\QCoh(\CZ)$ as a full subcategory. 

\medskip

In particular, $\IndCoh_{\on{Nilp}_{glob}}(\LocSysG)$ contains $\QCoh(\LocSysG)$ as a full subcategory.

\medskip

Accepting Conjecture \ref{c:main preview}, we obtain that $\QCoh(\LocSysG)$ corresponds to
a certain full subcategory of $\Dmod(\Bun_G)$, consisting of objects whose support in $\ArthG$
is contained in the zero-section. In other words, its support only contains Arthur parameters $(\sigma,A)$ with $A=0$. 

\medskip

We denote this category by $\Dmod_{\on{temp}}(\Bun_G)$. In \corref{c:temp} we give an intrinsic
characterization of this subcategory in terms of the action of the Hecke functors.

\ssec{Results concerning Langlands correspondence}

The main theorems of this paper fall into two classes. On the one hand, we prove some general results
about the behavior of the categories $\IndCoh_Y(\CZ)$. On the other hand, we run some consistency
checks on Conjecture \ref{c:main preview}.

\medskip

We will begin with the review of the latter.

\sssec{}

Recall that the Hecke category $\on{Sph}(G,x)\simeq \Dmod(\Gr_{G,x})^{G(\wh\CO_x)}$ is a monoidal
category acting on $\Dmod(\Bun_G)$ by the Hecke functors. 

\medskip

The (derived) Geometric Satake Equivalence identifies $\on{Sph}(G,x)$ with a certain subcategory
of the category of ind-coherent sheaves on the DG stack 
$$\on{pt}/\cG\underset{\cg/\cG}\times \on{pt}/\cG,$$
which is a monoidal category under convolution.

\medskip

We prove (\thmref{t:Satake}) that the resulting subcategory of
$\IndCoh(\on{pt}/\cG\underset{\cg/\cG}\times \on{pt}/\cG)$
is determined by a singular support condition. 

\medskip

Namely, there is a natural isomorphism
$$\Sing(\on{pt}/\cG\underset{\cg/\cG}\times \on{pt}/\cG)\simeq \cg^*/\cG,$$
which allows us to view 
\[\on{Nilp}(\cg^*)/\cG\]
as a conical subset of $\Sing(\on{pt}/\cG\underset{\cg/\cG}\times \on{pt}/\cG)$. Here
$\on{Nilp}(\cg^*)\subset \cg^*$
is the cone of nilpotent elements. We show that the Geometric Satake Equivalence identifies $\on{Sph}(G,x)$ with the full subcategory
\[ \IndCoh_{\on{Nilp}(\cg^*)/\cG}(\on{pt}/\cG\underset{\cg/\cG}\times \on{pt}/\cG)\subset \IndCoh(\on{pt}/\cG\underset{\cg/\cG}\times \on{pt}/\cG)\]
corresponding to this conical subset.

\medskip

We also show (see \propref{p:Hecke action spec}) that the conjectural equivalence of ``modified best hope''
(Conjecture~\ref{c:main preview}) 
is consistent with the Geometric Satake Equivalence. Namely, we construct a natural action of the monoidal
category \[\IndCoh_{\on{Nilp}(\cg^*)/\cG}(\on{pt}/\cG\underset{\cg/\cG}\times \on{pt}/\cG)\] 
on the category $\IndCoh_{\on{Nilp}_{glob}}(\LocSysG)$. This action should correspond to the action of
the monoidal category $\on{Sph}(G,x)$ on $\Dmod(\Bun_G)$ under the equivalence of Conjecture~\ref{c:main preview}.

\sssec{}

The following fact was observed in \cite{La}, and independently, by R.~Bezrukavnikov. 

\medskip

Take 
$X=\BP^1$, and let $\delta_{\bf \one}\in \Dmod(\Bun_G)$ be the D-module of $\delta$-functions at
the trivial bundle $\bf\one\in\Bun_G$. Then
the Hecke action of $\on{Sph}(G,x)$ on $\delta_{\bf \one}$
defines an equivalence
$$\on{Sph}(G,x)\to \Dmod(\Bun_G).$$

It is easy to see that the action of $\IndCoh(\on{pt}/\cG\underset{\cg/\cG}\times \on{pt}/\cG)$ 
on the sky-scraper of the unique $k$-point of $\LocSysG$ defines an equivalence
$$\IndCoh(\on{pt}/\cG\underset{\cg/\cG}\times \on{pt}/\cG)\to \IndCoh(\LocSysG).$$

Thus, \thmref{t:Satake} implies the existence of \emph{an} equivalence as stated in 
\conjref{c:main preview} for $X=\BP^1$.

\medskip

To show that this equivalence is \emph{the} equivalence, one needs to verify the additional 
properties that one expects the equivalence of \conjref{c:main preview} to satisfy, see \secref{sss:additional}.
Only some of these properties have been verified so far. So, an interested reader is welcome to tackle
them.

\sssec{}

A fundamental construction in the classical theory of automorphic functions is that of Eisenstein series. In the geometric
theory, this takes the form of a functor
$$\Eis^P_!:\Dmod(\Bun_M)\to \Dmod(\Bun_G),$$
defined for every parabolic subgroup $P$ with Levi quotient $M$.

\medskip

The Eisenstein series functor $\Eis^P_!$ is defined as 
$$(\sfp^P)_!\circ (\sfq^P)^*,$$
where $\sfp^P$ and $\sfq^P$ are the maps in the diagram
\begin{gather} 
\xy
(-15,0)*+{\Bun_G}="X";
(15,0)*+{\Bun_M.}="Y";
(0,15)*+{\Bun_P}="Z";
{\ar@{->}_{\sfp^P} "Z";"X"};
{\ar@{->}^{\sfq^P} "Z";"Y"};
\endxy
\end{gather}

It is not obvious that the functor $\Eis^P_!$ makes sense, because
the functors $(\sfq^P)^*$ and $(\sfp^P)^!$ are being applied to D-modules 
that need not be holonomic. For some non-trivial reasons, the functor $\Eis^P_!$ 
is defined on the category $\Dmod(\Bun_M)$; see \secref{sss:Eis!} for a quick summary and
\cite[Proposition~1.2]{DrGa2} for a proof.

\medskip

On the spectral side one has an analogous functor
$$\Eis^{P}_{spec}:\IndCoh(\LocSys_{\cM})\to \IndCoh(\LocSys_\cG),$$
defined as 
$$(\sfp^P_{spec})^{\IndCoh}_*\circ (\sfq^P_{spec})^!$$
using the diagram 

\begin{gather}  
\xy
(-15,0)*+{\LocSysG}="X";
(15,0)*+{\LocSys_\cM.}="Y";
(0,15)*+{\LocSys_\cP}="Z";
{\ar@{->}_{\sfp^P_{spec}} "Z";"X"};
{\ar@{->}^{\sfq^P_{spec}} "Z";"Y"};
\endxy
\end{gather}

The Langlands correspondence for groups $G$ and $M$ is supposed to intertwine
the functors $\Eis^P_!$ and $\Eis^P_{spec}$ (up to tensoring by a line bundle). 

\medskip

Thus, a consistency check for \conjref{c:main preview} should imply:

\begin{thm}  \label{t:nilp to nilp}
The functor $\Eis^{P}_{spec}$ maps 
$$\IndCoh_{\on{Nilp}_{glob}}(\LocSys_{\cM})\to \IndCoh_{\on{Nilp}_{glob}}(\LocSys_{\cG}).$$
\end{thm}

We prove this theorem in \secref{s:eis} (see \propref{p:spectral Eis fine}).

\sssec{}

The main result of this paper is \thmref{t:generation}, which is a refinement of 
\thmref{t:nilp to nilp}. 

\medskip

Essentially, \thmref{t:generation} says that the choice
of $\IndCoh_{\on{Nilp}_{glob}}(\LocSysG)$ as a subcategory of $\IndCoh(\LocSysG)$
containing $\QCoh(\LocSysG)$ is the minimal one, if we want to have an equivalence
with $\Dmod(\Bun_G)$ compatible with the Eisenstein series functors.

\medskip

More precisely, \thmref{t:generation} says that the category $\IndCoh_{\on{Nilp}_{glob}}(\LocSysG)$ is
generated by the essential images of $\QCoh(\LocSys_{\cM})$ under the functors 
$\Eis^{P}_{spec}$ for all parabolics $P$ (including $P=G$).

\ssec{Results concerning the theory of singular support}

\sssec{}

As was already mentioned, the idea of support based on cohomological operations was pioneered by 
D.~Benson, S.~B.~Iyengar, and H.~Krause in \cite{Kr}. 

\medskip

Namely, let $\bT$ be a triangulated category (containing arbitrary direct sums), and let $A$ be an algebra
graded by non-negative even integers that acts on $\bT$. By this we mean that every homogeneous 
element $a\in A_{2n}$ defines a
natural transformation from the identity functor to the shift functor $\CF\mapsto \CF[2n]$. Given a
homogeneous element $a\in A$, one can attach to it the full subcategory $\bT_{\Spec(A)-Y_a}\subset \bT$
consisting of $a$-local objects, and its left orthogonal, denoted $\bT_{Y_a}$, to be thought of as 
consisting of objects ``set-theoretically supported on the set of zeroes of $a$." More generally, one
can attach the corresponding subcategories 
$$\bT_Y\subset \bT\supset \bT_{\Spec(A)-Y}$$
to any conical Zariski-closed subset $Y\subset \Spec(A)$.

\medskip

It is shown in {\it loc.cit.} that the categories $\bT_Y\subset \bT$ are very well-behaved. Namely, they
satisfy essentially the same properties as when $\bT=A\mod$, and we are talking about the usual notion
of support in commutative algebra. 

\sssec{}

Let us now take $\bT$ to be the homotopy category of the DG category $\IndCoh(Z)$, where $Z$ is an affine DG scheme. 
There is a universal choice of a graded algebra acting on $\bT$, namely $\on{HH}(Z)$, the Hochschild cohomology of 
$Z$.  

\medskip

We note that when $Z$ is quasi-smooth, there is a canonical map of graded algebras
$$\Gamma(\Sing(Z),\CO_{\Sing(Z)})\to \on{HH}(Z),$$
where the grading on $\CO_{\Sing(Z)}$ is obtained by scaling by $2$ the action of $\BG_m$
along the fibers of $\Sing(Z)\to Z$. 

\medskip

Thus, by \cite{Kr}, we obtain the desired assignment 
$$Y\subset \Sing(Z) \, \rightsquigarrow\, \IndCoh_Y(Z)\subset \IndCoh(Z).$$

For a given $\CF\in \IndCoh(Z)$, its singular support is by definition the smallest $Y$ such that
$$\CF\subset \IndCoh_Y(Z).$$

\begin{rem}
We chose the terminology ``singular support" by loose analogy with the theory of D-modules.
In the latter case, the singular support of a D-module is a conical subset of the (usual) cotangent
bundle, which measures the degree to which the D-module is not lisse.
\end{rem}

\begin{rem}
We do not presume to make a thorough review of the existing literature on the subject. However,
in Appendix \ref{s:rev} we will indicate how the notion of singular support developed in this
paper is related to several other approaches, due to D.~Orlov, L.~Positselski, G.~Stevenson, and M.~Umut Isik,
respectively.
\end{rem}

\sssec{}

If the above definition of singular support sounds a little too abstract, here is how it can be rewritten more explicitly.

\medskip

First, we consider the most basic example of a quasi-smooth (DG !) scheme. Namely, let $\CV$
be a smooth scheme, and let $\on{pt}\to \CV$ be a $k$-point. We consider the DG scheme
$$\CG_{\on{pt}/\CV}:=\on{pt}\underset{\CV}\times \on{pt}.$$

\medskip

Explicitly, let $V$ denote the tangent space to $\CV$ at $\on{pt}$. Then for every \emph{parallelization} of $\CV$
at $\on{pt}$, i.e., for an identification of the formal completion of $\CO_\CV$ at $\on{pt}$ with
$\wh{\Sym(V)}$, we obtain an isomorphism
$$\CG_{\on{pt}/\CV}\simeq \Spec(\Sym(V^*[1])).$$

\medskip

Now, Koszul duality defines an equivalence of DG categories
$$\on{KD}_{\on{pt}/\CV}:\IndCoh(\CG_{\on{pt}/\CV})\simeq \Sym(V[-2])\mod$$
(this equivalence does not depend on the choice of a parallelization). 

\medskip

It is easy to see that $\Sing(\CG_{\on{pt}/\CV})\simeq V^*$. In terms of this equivalence, the singular support of 
an object in $\IndCoh(\CG_{\on{pt}/\CV})$ becomes the usual support of the corresponding object in $\Sym(V[-2])\mod$.
(By definition, the support of a $\Sym(V[-2])$-module $M$ is the support of its cohomology $H^\bullet(M)$ viewed as a graded
$\Sym(V)$-module.)

\medskip

For example, for $\CF=\CO_{\CG_{\on{pt}/\CV}}$, its singular support is $\{0\}$. By contrast,
for $\CF$ being the skyscraper at the unique $k$-point of $\CG_{\on{pt}/\CV}$, its singular
support is all of $V^*$. 

\sssec{}   \label{sss:Koszul preview}

Suppose now that a quasi-smooth DG scheme $Z$ is given as a global complete intersection. By this we mean that
$$Z=\on{pt}\underset{\CV}\times \CU,$$
where $\CU$ and $\CV$ are smooth. (Any quasi-smooth DG scheme can be locally written in this form, see \corref{c:q-smooth lci}.)

\medskip

Then it is easy to see that there exists a canonical closed embedding
\begin{equation} \label{e:embed Sing preview}
\Sing(Z)\hookrightarrow V^*\times Z, 
\end{equation}
where $V^*$ is as above. 

\medskip

Now, we have a canonical isomorphism of DG schemes
$$Z\underset{\CU}\times Z\simeq \CG_{\on{pt}/\CV}\times Z$$
(see \secref{sss:action on Z}), so we have a map, denoted
$$\on{act}_{\on{pt}/\CV,Z}:\CG_{\on{pt}/\CV}\times Z\to Z.$$

We show in \corref{c:Koszul generation}(a), that an object $\CF\in \IndCoh(Z)$
has its singular support inside $Y\subset \Sing(Z)\subset V^*\times Z$ if and only
if the object
$$(\on{KD}_{\on{pt}/\CV}\otimes \on{Id})\circ \on{act}_{\on{pt}/\CV,Z}^!(\CF)\in 
\Sym(V[-2])\mod\otimes \IndCoh(Z)$$
is supported on $Y$ in the sense of commutative algebra. 

\sssec{}

Here is yet another, even more explicit, characterization of singular support of objects of $\Coh(Z)$,
suggested to us by V.~Drinfeld. 

\medskip

Let $(z,\xi)$ be an element of $\Sing(Z)$, where $z$ is a point of $Z$ and $0\neq \xi\in H^{-1}(T^*_z(Z))$.
We wish to know when this point belongs to $\on{SingSupp}(\CF)$ for a given $\CF\in \Coh(Z)$.

\medskip

Suppose that $Z$ is written as in \secref{sss:Koszul preview}. Then by \eqref{e:embed Sing preview}, 
$\xi$ corresponds to a cotangent vector to $\CV$ at $\on{pt}$. Let $f$ be a function on $\CV$ that
vanishes at $\on{pt}$, and whose differential equals $\xi$. 
Let $Z'\subset \CU$ be the hypersurface cut out by the pullback of $f$ to $U$. Note that $Z'$ 
\emph{is singular} at $\on{pt}$. Let $i$ denote the closed embedding $Z\hookrightarrow Z'$.

\medskip

We have:

\begin{thm} \label{t:drinfeld preview}
The element $(z,\xi)$ \emph{does not} belong to $\on{SingSupp}(\CF)$ if and only if
$i_*(\CF)\in \Coh(Z')$ is \emph{perfect} on a Zariski neighborhood of $z$.
\end{thm}

This theorem will be proved as \corref{c:Drinfeld}

\sssec{}  \label{sss:sing preview}

Here are some basic properties of the assignment
$$Y\mapsto \IndCoh_Y(Z):$$

\medskip

\noindent(a) As was mentioned before,
$$\CF\in \QCoh(\CF)\, \Leftrightarrow\, \on{SingSupp}(\CF)=\{0\};$$
this is \thmref{t:zero sect}.

\medskip

\noindent(b) The assignment $Y\rightsquigarrow \IndCoh_Y(Z)$ 
is Zariski-local (see \corref{c:on open +}).
In particular, this allows us to define singular support on non-affine DG schemes.

\medskip

\noindent(c) For $Z$ quasi-compact, the category $\IndCoh_Y(Z)$ is compactly
generated. This is easy for $Z$ affine (see \corref{c:with support comp gen})
and is a variant of the theorem of \cite{TT} in general
(see Appendix \ref{a:B}).

\medskip

\noindent(d) The category $\IndCoh_Y(Z)$ has a t-structure whose eventually
coconnective part identifies with that of $\QCoh(Z)$. Informally, the difference between 
$\IndCoh_Y(Z)$ and $\QCoh(Z)$ is ``at $-\infty$." This is the content of \secref{ss:t-structure}.

\medskip

\noindent(e) For $\CF\in \Coh(Z)$, 
$$\on{SingSupp}(\CF)=\on{SingSupp}(\BD_Z^{\on{Serre}}(\CF)),$$
where $\BD_Z^{\on{Serre}}$ is the Serre duality anti-involution of the
category $\Coh(Z)$. This is \propref{p:Serre}.

\medskip

\noindent(f) Singular support can be computed point-wise. Namely, for $\CF\in \IndCoh(Z)$ and 
a geometric point $\Spec(k')\overset{i_z}\longrightarrow Z$, the graded vector 
space $H^\bullet(i_z^!(\CF))$ is acted on by the algebra $\on{Sym}(H^1(T_z(Z))$, and
$$\on{SingSupp}(\CF)\subset Y\, \Leftrightarrow\, \forall\, z,\, \on{supp}_{\on{Sym}(H^1(T_z(Z))}\left(H^\bullet(i_z^!(\CF))\right)
\subset Y\cap H^{-1}(T^*_z(Z)).$$
This is \propref{p:pointwise}.

\medskip

\noindent(g) For $Z_1$ and $Z_2$ quasi-compact, and $Y_i\subset \Sing(Z_i)$, we have
$$\IndCoh_{Y_1}(Z_1)\otimes \IndCoh_{Y_2}(Z_2)=\IndCoh_{Y_1\times Y_2}(Z_1\times Z_2)$$
as subcategories of
$$\IndCoh(Z_1)\otimes \IndCoh(Z_2)\otimes \IndCoh(Z_1\times Z_2).$$
This is \lemref{l:products}.

\medskip

\noindent(h) An estimate on singular support ensures preservation of coherence. For example,
if $\CF',\CF''\in \Coh(Z)$ are such that the set-theoretic intersection of their respective singular supports is
contained in the zero-section of $\Sing(Z)$, then the tensor product $\CF'\otimes \CF''$
lives in finitely many cohomological degrees. This is proved in \propref{p:preserve coherence}.

\sssec{Functoriality}   \label{sss:funct preview}

Let $f:Z_1\to Z_2$ be a map between quasi-smooth (and quasi-compact) DG schemes. We have the functors
$$f^\IndCoh_*:\IndCoh(Z_1)\to \IndCoh(Z_2) \text{ and } f^!:\IndCoh(Z_2)\to \IndCoh(Z_1)$$
(they are adjoint if $f$ is proper). 

\medskip

We wish to understand how they behave in relation to the categories $\IndCoh_Y(Z)$. 

\medskip

First, we note that there exists a canonical map
$$Z_1\underset{Z_2}\times \Sing(Z_2)\to \Sing(Z_1);$$
we call this map ``the singular codifferential of $f$," and denote it by $\Sing(f)$. 

\medskip

The singular support behaves naturally under direct and inverse images: 

\begin{thm} \label{t:funct preview} Let $Y_i\subset \Sing(Z_i)$ be conical Zariski-closed subsets.

\smallskip

\noindent(a) If $\Sing(f)^{-1}(Y_1)\subset Y_2\underset{Z_2}\times Z_1$, then the functor
$f^\IndCoh_*$ maps $\IndCoh_{Y_1}(Z_1)$ to $\IndCoh_{Y_2}(Z_2)$.

\smallskip

\noindent(b) If $Y_2\underset{Z_2}\times Z_1\subset \Sing(f)^{-1}(Y_1)$, then $f^!$ maps 
 $\IndCoh_{Y_2}(Z_2)$ to $\IndCoh_{Y_1}(Z_1)$.

\end{thm}

This is proved in \propref{p:singsupp functoriality}.  

\medskip

Suppose now that the map $f$ is itself quasi-smooth. According to \lemref{l:codiff for q-smooth}, this
is equivalent to the condition that the singular codifferential map $\Sing(f)$ be a closed embedding.
For $Y_2\subset \Sing(Z_2)$, let
$$Y_1:=\Sing(f)(Y_2\underset{Z_2}\times Z_1)$$
be the corresponding subset in $\Sing(Z_1)$. 

\medskip

In \corref{c:quasi-smooth tensor up}, we will show:

\begin{thm}  \label{t:tensor up preview}
The functor $f^!$ defines an equivalence
$$\QCoh(Z_1)\underset{\QCoh(Z_2)}\otimes \IndCoh_{Y_2}(Z_2)\to \IndCoh_{Y_1}(Z_1).$$
\end{thm}

\medskip

Finally, we have the following crucial result (see \thmref{t:prop cons}):

\begin{thm} \label{t:proper preview} Assume that $f$ is proper, and let $Y_i\subset \Sing(Z_i)$ be such that
the composed map
$$\Sing(f)^{-1}(Y_1)\hookrightarrow  Z_1\underset{Z_2}\times \Sing(Z_2)\to \Sing(Z_2)$$
is surjective onto $Y_2$. Then the essential image of $\IndCoh_{Y_1}(Z_1)$ under $f^\IndCoh_*$
generates $\IndCoh_{Y_2}(Z_2)$.
\end{thm}

\begin{rem} Theorems~\ref{t:tensor up preview} and \ref{t:proper preview} are deeper than Theorem~\ref{t:funct preview}. Indeed, Theorem~\ref{t:funct preview} provides upper bounds on the essential images 
$f^\IndCoh_*(\IndCoh_{Y_2}(Z_2))$ and $f^!(\IndCoh_{Y_1}(Z_1))$; to a large extent, these bounds formally follow
from definitions. On the other hand, Theorems~\ref{t:tensor up preview} and \ref{t:proper preview} give a 
precise description of this image under some additional assumptions on
the morphism $f$; they are more ``geometric" in nature.
\end{rem}

\sssec{}

Finally, we remark that \thmref{t:tensor up preview} ensures that the assignment
$Y\mapsto \IndCoh_Y(Z)$ is local also in the smooth topology. This allows us to develop
the theory of singular support on DG Artin stacks:  

\medskip

For a quasi-smooth DG Artin stack $\CZ$ we introduce the classical Artin stack 
$\Sing(\CZ)$ by descending $\Sing(Z)$ over affine DG schemes $Z$ mapping smoothly
to $\CZ$ (any such $Z$ is automatically quasi-smooth).

\medskip

Given $Y\subset \Sing(\CZ)$ we define the category $\IndCoh_Y(\CZ)$ as the limit of
the categories 
$$\IndCoh_{Z\underset{\CZ}\times Y}(Z)$$
over $Z\to \CZ$ as above.

\medskip

One easily establishes the corresponding properties of the category $\IndCoh_Y(\CZ)$ 
by reducing to the case of schemes. The one exception is the question of compact
generation.

\medskip

At the moment we cannot show that $\IndCoh_Y(\CZ)$ is compactly generated for a general
nice (=QCA) algebraic DG stack. However, we can do it when $\CZ$ can be presented as a global
complete intersection (see \corref{c:category comp gen stacks}). Fortunately, this is the case for $\CZ=\LocSys_G$ for any
affine algebraic group $G$. 

\medskip

However, we do prove that the category $\IndCoh_Y(\CZ)$ is dualizable for a general $\CZ$
which is QCA (see \corref{c:dualizable}).

\ssec{How this paper is organized}

This paper is divided into three parts. Part I contains miscellaneous preliminaries, Part II develops
the theory of singular support for $\IndCoh$, and Part III discusses the applications to Geometric
Langlands.

\sssec{}

In \secref{s:q-smooth} we recall the notion of quasi-smooth DG scheme and morphism. We show that 
this condition is equivalent to that of locally complete intersection. We also introduce the classical scheme
$\Sing(Z)$ attached to a quasi-smooth DG scheme $Z$. 

\medskip

In \secref{s:sup} we review the theory of support in a triangulated category acted on by a graded algebra.
Most results of this section are contained in \cite{Kr}. (Note, however, that what we call support
is the closure of the support in the terminology of {\it loc.cit.})

\sssec{}

In \secref{s:sing} we introduce the notion of singular support of objects of $\IndCoh(Z)$, where
$Z$ is a quasi-smooth DG scheme, and establish the basic properties listed in \secref{sss:sing preview}.

\medskip

In \secref{s:gci} we study the case of a global complete intersection, and prove the description of
singular support via Koszul duality, mentioned in \secref{sss:Koszul preview}.

\medskip

In \secref{s:pointwise} we show how singular support of an object $\CF\in \IndCoh(Z)$ 
can be read off the behavior of the fibers of $\CF$ at geometric points of $Z$. 

\medskip

In \secref{s:funct} we establish the functoriality properties of categories $\IndCoh_Y(Z)$
mentioned in \secref{sss:funct preview}.

\medskip

In \secref{s:stacks} we develop the theory of singular support on Artin DG stacks. 

\medskip

In \secref{s:gci stacks} we prove that if a quasi-compact algebraic DG stack $\CZ$ is given as a global complete
intersection, then the subcategories of $\IndCoh(\CZ)$ defined by singular 
support are compactly generated.

\sssec{}

In \secref{s:LocSys} we recall the definition of the DG stack of $G$-local systems on a given 
scheme $X$, where $G$ is an algebraic group. The reason we decided to include this section  
rather than refer to some existing source is that, even for $X$ a smooth and complete curve,
the stack $\LocSys_G$ is an object of derived algebraic geometry, and the relevant definitions
do not seem to be present in the literature, although seem to be well-known in folklore. 

\medskip

In \secref{s:Langlands}, 
 we introduce the global nilpotent cone
$$\on{Nilp}_{glob}\subset \ArthG=\Sing(\LocSysG)$$
and formulate the Geometric Langlands Conjecture. 

\medskip

In \secref{s:Satake} we study how our proposed form of 
Geometric Langlands Conjecture interacts with the Geometric Satake Equivalence.

\medskip

In \secref{s:eis}, we study the functors of Eisenstein series on both the 
automorphic and the Galois side of the correspondence, and prove a consistency
result (\thmref{t:generation}) with Conjecture \ref{c:main preview}.

\sssec{}

In Appendix \ref{s:eq} we list several facts pertaining to the notion of action of
an algebraic group on a category, used in several places in the main body of
the paper. These facts are fully documented in \cite{GA}. 

\medskip

In Appendix \ref{s:Weil} we discuss the formation of mapping spaces in derived
algebraic geometry, and how it behaves under deformation theory.

\medskip

In Appendix \ref{a:B} we prove a version of the Thomason-Trobaugh theorem for categories
defined by singular support.

\medskip

In Appendix \ref{s:proof of finiteness} we prove a certain finite generation result for Exts between coherent
sheaves on a quasi-smooth DG scheme. 

\medskip

In Appendix \ref{s:E2} we make a brief review of the theory of $\BE_2$-algebras. 

\medskip

In Appendix \ref{s:Hoch} we collect some facts concerning the $\BE_2$-algebra of Hochschild cochains 
on an affine DG scheme, and its generalization that has to do with groupoids. 

\medskip

In Appendix \ref{s:Lie} we study the connection between the $\BE_2$-algebra of Hochschild cochains and Lie algebras,
and its generalization to the case of groupoids. 

\medskip

In Appendix \ref{s:rev} we review the relation of the theory of singular support developed
in this paper with several other approaches. 

\ssec{Conventions}

\sssec{}

Throughout the text we work with a ground field $k$, assumed algebraically closed 
and of \emph{characteristic $0$}. 

\sssec{$\infty$-categories and DG categories}

Our conventions follow completely those adopted in the paper \cite{DrGa1}, and we refer the reader
to Sect. 1,  where the latter are explained. 

\medskip

In particular: 

\smallskip

\noindent(i) When we say ``category" by default we mean ``$(\infty,1)$-category".  

\smallskip

\noindent(ii) For a category
$\bC$ and objects $\bc_1,\bc_2\in \bC$ we will denote by $\Maps_\bC(\bc_1,\bc_2)$
the $\infty$-groupoid of maps between them. We will denote by $\Hom_\bC(\bc_1,\bc_2)$
the \emph{set} $\pi_0(\Maps_\bC(\bc_1,\bc_2))$, i.e., $\Hom$ in the ordinary category
$\on{Ho}(\bC)$. 

\smallskip

\noindent(iii) All DG categories are assumed to be pretriangulated and, unless explicitly stated otherwise, 
cocomplete (that is, they contain arbitrary direct sums). 
All functors between DG categories are assumed to be exact and continuous (that is, commuting with arbitrary 
direct sums, or equivalently, with all colimits). In particular, all subcategories are by default
assumed to be closed under arbitrary direct sums.

Note that our terminology for functors follows \cite{Lu0} in that a colimit-preserving functor is called
continuous, rather than co-continuous.

\smallskip

\noindent(iv) For a DG category $\bC$ equipped with a t-structure, we will denote by
$\bC^{\leq 0}$ (resp. $\bC^{\geq 0}$) the corresponding subcategories
of connective (resp. coconnective) objects. We let $\bC^\heartsuit$ denote the heart of
the t-structure. We also let $\bC^+$ (resp. $\bC^-$) denote the subcategory of
eventually coconnective (resp. connective) objects.

\smallskip

\noindent(v) We let $\Vect$ denote the DG category of complexes of vector spaces;
thus, the usual category of $k$-vector spaces is denoted by $\Vect^\heartsuit$. 
The category of $\infty$-groupoids is denoted by $\inftygroup$.

\smallskip

\noindent(vi) If $\bC$ is a DG category, let $\CMaps_{\bC}(\bc_1,\bc_2)$ denote the 
corresponding object of $\Vect$. Sometimes, we view $\CMaps_{\bC}(\bc_1,\bc_2)$
as a spectrum. In particular, $\Maps_\bC(\bc_1,\bc_2)$ is recovered as the $0$-th
space of $\CMaps_{\bC}(\bc_1,\bc_2)$. We can also view
$$\Maps_\bC(\bc_1,\bc_2)\simeq \tau^{\leq 0}(\CMaps_{\bC}(\bc_1,\bc_2))\in \Vect^{\leq 0},$$
where $\Vect^{\leq 0}$ maps to $\inftygroup$ via the Dold-Kan correspondence. 

\smallskip

\noindent(vi') We will denote by $\Hom^\bullet_\bC(\bc_1,\bc_2)$ the graded vector space 
$H^\bullet(\CMaps_{\bC}(\bc_1,\bc_2))$. By definition
$$H^\bullet(\CMaps_{\bC}(\bc_1,\bc_2))=\underset{i}\oplus\, \Hom_{\on{Ho}(\bC)}(\bc_1,\bc_2[i]).$$
We sometimes use the notation 
$\Hom^\bullet(-,-)$ when instead of a DG category we have just a triangulated category $\bT$. That is,
we use $\Hom^\bullet_\bT(\bt_1,\bt_2)$ in place of the more common
$\Ext^\bullet_\bT(\bt_1,\bt_2)$. 

\smallskip

\noindent(vii)
Let $\bC'\subset\bC$ be a full subcategory of the DG category $\bC$. (Recall that $\bC$ is assumed to be cocomplete,
and $\bC'$ is assumed to be closed under direct sums.) Suppose that the inclusion $\bC'\to\bC$ admits a left
adjoint, which is automatically exact and continuous. We call this adjoint the \emph{localization functor}.

\smallskip

In this case, we let $\bC'':={}^{\perp}\!(\bC')$ be the left orthogonal complement of $\bC'$. The inclusion
$\bC''\to\bC$ admits a right adjoint, which is easily seen to be continuous (see \lemref{lem:exact of triangles},
which contains the version of this claim for triangulated categories). We call the right adjoint functor 
$\bC\to\bC''$ the \emph{colocalization functor}. The resulting diagram
\[\bC''\rightleftarrows\bC\rightleftarrows\bC'\]
is a \emph{short exact sequence of DG categories}.
 
\smallskip

\noindent(viii) For a monoidal DG category $\bC$, we use the terms ``a DG category with an action of $\bC$"
and ``a DG category tensored over $\bC$" interchangeably (see \cite[A.1.4]{Lu0} for definition).

\smallskip

\noindent(ix) For a (DG) associative algebra $A$, we denote by $A\mod$ the corresponding
DG category of $A$-modules.

\sssec{DG schemes and Artin stacks}

Conventions and notation regarding DG schemes and DG Artin stacks (and, more generally, prestacks)
follow \cite{Stacks}. In particular, $\affSch$, $\affdgSch$, $\Sch$, $\dgSch$, and $\on{PreStk}$ 
stand for the categories of classical affine schemes, affine DG schemes, classical schemes, DG schemes, and (DG) 
prestacks, respectively. A short review of the conventions can be found also in \cite[Sect. 0.6.4-0.6.5]{DrGa0}. 

\medskip

We denote $\on{pt}:=\Spec(k)$. 

\medskip

By default, all schemes/Artin stacks are derived. When they are classical, we will 
emphasize this explicitly. 

\medskip

All DG schemes and DG Artin stacks in this paper will be assumed locally almost of finite type over $k$
(see \cite[Sect. 1.3.9, 2.6.5, 3.3.1 and 4.9]{Stacks}), unless specified otherwise. We denote by 
$\dgSch_{\on{aft}}\subset\dgSch$ the full subcategory of DG schemes that are locally almost of finite type.

\medskip

We will also use the following convention: we will not distinguish between the notions
of classical scheme/Artin stack and that of $0$-coconnective DG scheme/Artin stack
(see \cite[Sect. 4.6.3]{Stacks} for the latter notion):

\medskip

By definition, classical schemes/Artin stacks form a full
subcategory among functors $(\affSch)^{\on{op}}\to \inftygroup$, while $0$-coconnective DG schemes/Artin stacks
form a full subcategory among functors $(\affdgSch)^{\on{op}}\to \inftygroup$. 
However, the two categories are equivalent: the equivalence is given by restriction 
along the fully faithful embedding $(\affSch)^{\on{op}}\subset (\affdgSch)^{\on{op}}$; the 
inverse procedure is left Kan extension, followed by sheafification. 

\medskip

A more detailed discussion of the notion of $n$-coconnectivity can be found in
\cite[Sect. 0.5]{IndSch}.

\medskip

For a given DG scheme/Artin stack $\CZ$, we will denote by $^{cl}\CZ$ the underlying
classical scheme/classical Artin stack. 

\sssec{}

Conventions regarding the categories $\QCoh(-)$ and $\IndCoh(-)$ on DG schemes/Artin stacks
follow those of \cite{QCoh} and \cite{IndCoh}, respectively. Conventions regarding the category
$\Dmod(-)$ follow \cite{Crys}.

\ssec{Acknowledgements}

The idea of this work emerged as a result of discussions with V.~Drinfeld. The very existence of
this paper expresses our debt to him.

\medskip

We are grateful to J.~Lurie for teaching us the foundations of derived algebraic geometry. We are grateful to
N.~Rozenblyum for clarifying various aspects of the theory of $\BE_2$-algebras. 

\medskip

We would like to thank A.~Beilinson, H.~Krause, V.~Lafforgue, A.~Neeman, T.~Pantev, G.~Stevenson
for helpful communications and discussions. 

\medskip

Special thanks are due to S.~Raskin, whose extensive comments on successive drafts of the paper led
to substantial improvements.

\medskip

We are very grateful to the anonymous referees for their numerous comments and suggestions.

\medskip

The work of D.A. is partially supported by NSF grants DMS-1101558 and DMS-0635607.
The work of D.G. is partially supported by NSF grant DMS-1063470.

\bigskip 

\bigskip

\centerline{\bf Part I: Preliminaries}

\section{Quasi-smooth DG schemes and the scheme of singularities}  \label{s:q-smooth}

We remind that all DG schemes in this paper are assumed locally almost of finite type
over the ground field $k$, unless explicitly stated otherwise.

\medskip

In this section we will recall the notions of quasi-smooth (DG) scheme and quasi-smooth morphism
between DG schemes. To a quasi-smooth DG scheme $Z$ we will attach a classical scheme $\Sing(Z)$
that ``controls" the singularities of $Z$. 

\medskip

Quasi-smoothness is the ``correct" DG version of being locally a complete intersection. 

\ssec{The notion of quasi-smoothness}

In this subsection we define quasi-smoothness in terms of the cotangent complex. We will
show that any quasi-smooth morphism is, locally on the source, a composition of a morphism
that can be obtained as a base change of $\on{pt}\to \BA^n$, followed by a smooth morphism.

\sssec{}  Recall the notion of smoothness for a map between DG schemes: 

\medskip

A map $Z_1\to Z_2$ is called smooth if the DG scheme $^{cl}\!Z_2\underset{Z_2}\times Z_1$ 
is classical, and the resulting map
$$^{cl}\!Z_2\underset{Z_2}\times Z_1\to {}^{cl}Z_2$$
is a smooth map of classical schemes. 

\medskip

In particular, a DG scheme is called smooth 
if its map to $\on{pt}:=\Spec(k)$ is smooth. We summarize the properties of smooth maps below. 

\begin{lem}  \hfill  \label{l:smooth sch}

\smallskip

\noindent{\em(a)}
A DG scheme $Z$ is smooth if and only if its cotangent
complex $T^*(Z)$ is a vector bundle, that is, $T^*(Z)$ is Zariski-locally isomorphic to $\CO_Z^n$.

\smallskip

\noindent{\em(b)} A smooth DG scheme is classical, and a smooth
classical scheme is smooth as a DG scheme (see, e.g., \cite[Sect. 8.4.2 and Proposition 9.1.4]{IndSch} for the proof).

\smallskip

\noindent{\em(c)}
A map $f:Z_1\to Z_2$ between DG schemes is smooth if and only if
its relative cotangent complex $T^*(Z_1/Z_2)$ is a vector bundle.
\end{lem}

\sssec{The definition of quasi-smoothness}

A DG scheme $Z$ is called quasi-smooth \footnote{A.k.a. ``l.c.i."=locally complete intersection.}
if its cotangent complex $T^*(Z)$ is perfect of Tor-amplitude
$[-1,0]$. 

\medskip

Equivalently, we require that, Zariski-locally on $Z$,  
the object $T^*(Z)\in \QCoh(Z)$ could be presented by a complex 
$$\CO_Z^{\oplus m}\to \CO_Z^{\oplus n}.$$
This is equivalent to the condition that all geometric fibers of $T^*(Z)$ are acyclic in degrees below $-1$.

\begin{rem}
We should emphasize that if $Z$ is a quasi-smooth DG scheme, the underlying classical scheme $^{cl}\!Z$ 
need \emph{not} be a locally complete intersection. In fact, \emph{any} classical affine scheme can be realized 
in this way for tautological reasons. 
\end{rem}

\sssec{}

The following is a particular case of \corref{c:shape of q-smooth}:

\begin{cor}  \label{c:q-smooth lci}
A DG scheme $Z$ is quasi-smooth if and only if it can be Zariski-locally presented as a fiber
product (in the category of DG schemes)
$$
\CD
Z  @>>>  \BA^n  \\
@VVV    @VVV   \\
\on{pt}  @>{\{0\}}>>  \BA^m.
\endCD
$$
\end{cor}

For future reference we record:

\begin{cor}  \label{c:shape of smooth}
Let $f:Z_1\to Z_2$ be a \emph{smooth} morphism between quasi-smooth DG schemes. Then, Zariski-locally on $Z_1$,
there exists a Cartesian diagram
$$
\CD
Z_1  @>>>  \CU_1 \\
@V{f}VV  @VV{f_\CU}V   \\
Z_2  @>>>  \CU_2 \\
@VVV    @VVV   \\
\on{pt}  @>>>  \CV,
\endCD
$$
where the DG schemes $\CU_1,\CU_2,\CV$ are smooth, and the map $f_\CU$
is smooth as well. 
\end{cor}

\begin{proof}
This easily follows from the fact that, given a smooth morphism $f:Z_1\to Z_2$ between DG schemes, and 
a closed embedding $Z_2\hookrightarrow \CU_2$ with $\CU_2$ smooth, we can, Zariski-locally on $Z_1$, 
complete it to a Cartesian square
$$
\CD
Z_1  @>>>  \CU_1 \\
@V{f}VV  @VV{f_\CU}V   \\
Z_2  @>>>  \CU_2 \\
\endCD
$$
with $f_\CU$ smooth. 
\end{proof}

\sssec{Quasi-smooth maps}

We say that a morphism of DG schemes $f:Z_1\to Z_2$ is quasi-smooth if the relative
cotangent complex $T^*(Z_1/Z_2)$ is perfect of Tor-amplitude $[-1,0]$.

\medskip

We note:

\begin{lem}  \label{l:factoring quasi-smooth}
A quasi-smooth morphism can be, Zariski-locally on the source, factored as a composition of a quasi-smooth closed
embedding, followed by a smooth morphism.
\end{lem}

\begin{proof}

With no restriction of generality, we can assume that both $Z_1$ and $Z_2$ are affine. Decompose $f:Z_1\to Z_2$ as 
$$Z_1 \overset{f'}\to  Z_2\times \BA^n \to  Z_2,$$
where the first arrow is a closed embedding. 

\medskip

It follows tautologically that the fact that $f$ is quasi-smooth implies that $f'$ is quasi-smooth.

\end{proof} 

The following gives an explicit description of quasi-smooth closed embeddings:

\begin{prop} \label{p:shape of q-smooth cl}
Let $f:Z_1\to Z_2$ be a quasi-smooth closed embedding. Then, Zariski-locally on $Z_2$, there
exists a Cartesian diagram (in the category of DG schemes)
$$
\CD
Z_1  @>{f}>> Z_2 \\
@VVV    @VVV   \\
\on{pt} @>\{0\}>>  \BA^m.
\endCD
$$
\end{prop}

\begin{proof}

Note that $T^*(Z_1/Z_2)[-1]$ is the derived conormal 
sheaf $\CN^*(Z_1/Z_2)$ to $Z_1$ inside $Z_2$. The conditions imply that 
$\CN^*(Z_1/Z_2)$ is a vector bundle. 

\medskip

Consider the restriction $$^{cl}\CN^*(Z_1/Z_2):=\CN^*(Z_1/Z_2)|_{^{cl}\!Z_1}.$$ This
is the classical conormal sheaf to $^{cl}\!Z_1$ in $^{cl}\!Z_2$. Locally on $Z_2$ we can choose 
sections 
$$\{{}^{cl}\!f_1,\dots,{}^{cl}\!f_m\} \in \on{ker}(\CO_{^{cl}\!Z_2}\to f_*(\CO_{^{cl}\!Z_1})),$$
whose differentials generate $^{cl}\CN^*(Z_1/Z_2)$. Lifting the above sections to 
sections $f_1,\dots,f_m$ of $\CO_{Z_2}$, we obtain a map 
$$Z_1\to \on{pt}\underset{\BA^m}\times Z_2.$$
We claim that the latter map is an isomorphism. Indeed, it is a closed embedding and induces an isomorphism
at the level of derived cotangent spaces. 

\end{proof}

Combining the statements of \lemref{l:factoring quasi-smooth} and \propref{p:shape of q-smooth cl}, we obtain:

\begin{cor} \label{c:shape of q-smooth}
A morphism $f:Z_1\to Z_2$ is quasi-smooth if and only if, Zariski-locally on the source, it can be included in a diagram
$$
\CD
Z_1  @>{f'}>>  Z_2\times \BA^n @>{\on{pr}}>>  Z_2 \\
@VVV    @VVV   \\
\on{pt} @>\{0\}>>  \BA^m,
\endCD
$$
in which the square is Cartesian (in the category of DG schemes), and $f'$ is a closed embedding. 
\end{cor}

\ssec{Cohomological properties of quasi-smooth maps}

\sssec{}

First, we note:

\begin{cor}  \label{c:f Tor dim}   Let $f:Z_1\to Z_2$ be quasi-smooth. Then it is
of bounded Tor dimension, locally in the Zariski topology on $Z_1$. 
\end{cor}

\begin{proof}
Follows from \corref{c:shape of q-smooth} by base change.
\end{proof}

\sssec{}  \label{sss:loc ev coconn}

Recall now the notion of \emph{eventually coconnective} morphism, see \cite[Definition 3.5.2]{IndCoh}.
Namely, a morphism $f:Z_1\to Z_2$ between quasi-compact DG schemes is eventually coconnective 
if $f^*$ sends $\Coh(Z_2)$ to $\QCoh(Z_1)^+$ (equivalently, to $\Coh(Z_1)$). 

\medskip

We will say a morphism $f$ is \emph{locally eventually coconnective} if it is eventually coconnective
locally in the Zariski topology on the source. 

\medskip

Evidently, a morphism of bounded Tor dimension is eventually coconnective\footnote{
In fact,
for maps between quasi-compact eventually coconnective schemes, the converse is also true, see \cite[Lemma 3.6.3]{IndCoh}}. 
Hence, from \corref{c:f Tor dim} we obtain:

\begin{cor}  \label{c:q-smooth eventualy coconn}
A quasi-smooth morphism is locally eventually coconnective.
\end{cor}

\begin{cor} \label{c:f^*} Suppose $f:Z_1\to Z_2$ is a quasi-smooth map between DG schemes.
Then the functor $$f^{\IndCoh,*}:\IndCoh(Z_2)\to \IndCoh(Z_1),$$
\emph{left} adjoint to $f_*^{\IndCoh}:\IndCoh(Z_1)\to \IndCoh(Z_2)$ 
is well-defined.
\end{cor}

\begin{proof} Follows from \corref{c:q-smooth eventualy coconn} by 
\cite[Proposition 3.5.4 and Lemma 3.5.8]{IndCoh}. \footnote{
Technically, the corresponding assertions in \cite{IndCoh} are made under the assumption 
that $Z_1$ and $Z_2$ be quasi-compact.
However, this restriction is immaterial because the question of existence of 
$f^{\IndCoh,*}$ is Zariski-local on $Z_1$.}
\end{proof}

\sssec{} 

Finally, let us recall the notion of Gorenstein morphism between DG schemes, see \cite[Definition 7.3.2]{IndCoh}. We have:

\begin{cor} \label{c:Gorenstein}
A quasi-smooth morphism $f:Z_1\to Z_2$ between DG schemes is Gorenstein.
\end{cor}

\begin{proof}
The claim is local in the Zariski topology on $Z_1$. By \corref{c:shape of q-smooth}, locally we can write $f$ as a composition of a smooth
morphism and a morphism which is obtained by base change from the embedding $\on{pt}\to\BA^n$. Now the statement follows
from the combination of the following facts: (i) the property of being Gorenstein survives compositions (obvious from the definition), 
(ii) smooth morphisms are Gorenstein (see \cite[Corollary 7.5.2]{IndCoh}); (iii) the map $\on{pt}\to\BA^n$
is Gorenstein (direct calculation); and (iv) the base change of a Gorenstein morphism
is Gorenstein (see \cite[Sect. 7.5.4]{IndCoh}). 
\end{proof}

As a particular case, we have:

\begin{cor}  \label{c:Gorenstein abs}
A quasi-smooth DG scheme is Gorenstein (that is, its dualizing complex $\omega_Z\in \IndCoh_Z$ is a cohomologically shifted line bundle).
\end{cor}

\begin{rem}
If $Z$ is a DG scheme and $n\in \BZ$ is such that $\omega_Z[-n]$ is a line bundle, one
can call $n$ the ``virtual dimension of $Z$."
\end{rem}

\ssec{The scheme of singularities}

In this subsection we will introduce one of the main characters of this paper. Namely, if $Z$ is a 
quasi-smooth DG scheme, we will attach to it a classical scheme $\Sing(Z)$ that measures how
far $Z$ is from being smooth. 

\sssec{}

Let $Z$ be a DG scheme such that $T^*(Z)\in \QCoh(Z)$ is perfect (as is the case for quasi-smooth DG schemes). 
We define the tangent complex $T(Z)\in\QCoh(Z)$ to be the dual of $T^*(Z)$. Since $T^*(Z)$ is perfect,
so is $T(Z)$, and the dual of $T(Z)$ identifies with $T^*(Z)$.
\begin{rem}
Of course, the tangent complex can be defined for a general DG scheme $Z$ locally almost of finite type. 
However, to avoid losing information, one has to use Serre's duality instead of the ``naive" dual. Thus
we obtain a variant of the tangent complex that belongs to the category $\IndCoh(Z)$. 
This will be addressed in more detail in \cite{Algebroids}.
\end{rem}

\sssec{}

Let $Z$ be quasi-smooth. Note that in this case $T(Z)$ is perfect of Tor-amplitude $[0,1]$. In particular,
it has cohomologies only in degrees $0$ and $1$; moreover $H^1(T(Z))$ measures the degree to which
$Z$ is non-smooth.

\medskip

For a quasi-smooth DG scheme $Z$ we define the classical scheme $\Sing(Z)$ which we will refer to
as ``the scheme of singularities of $Z$" as
$$^{cl}\!\left(\Spec_{Z}\left(\on{Sym}_{\CO_Z}(T(Z)[1])\right)\right).$$

Note that since we are passing to the underlying classical scheme, the above 
is the same as
$$\Spec_{^{cl}\!Z}\left(\on{Sym}_{\CO_{^{cl}\!Z}}\left(H^1(T(Z))\right)\right),$$
where $H^1(T(Z))$ is considered as a coherent sheaf on $^{cl}\!Z$.

\medskip

The scheme $\Sing(Z)$ carries a canonical $\BG_m$-action along the fibers of the projection
$\Sing(Z)\to{^{cl}\!Z}$; the action corresponds to the natural grading on the symmetric algebra 
$\on{Sym}_{\CO_{Z}}(T(Z)[1])$.
 
\sssec{}

By definition, $k$-points of $\Sing(Z)$ are pairs $(z,\xi)$, where $z\in Z(k)$
and $\xi\in H^{-1}(T_z^*(Z))$.

\sssec{}  \label{sss:expl presentations}

Suppose that $Z$ is presented as a fiber product
$$
\CD
Z  @>{\iota}>> \CU  \\ 
@VVV    @VVV \\
\on{pt}  @>>>  \CV,
\endCD
$$
where $\CU$ and $\CV$ are smooth. In this situation, we say that $Z$ is given as a global complete intersection.

\medskip

Let $V$ denote the tangent space to $\CV$ at $\on{pt}\in \CV$. We have 
$$T(Z)\simeq \on{Cone}(\iota^*(T(\CU))\to V\otimes \CO_Z)[-1].$$
Hence, we obtain a canonical map
$$V\otimes \CO_Z\to T(Z)[1],$$
which gives rise to a surjection of coherent sheaves
$$V\otimes \CO_{^{cl}\!Z}\to H^1(T(Z)),$$
and, hence, we obtain a $\BG_m$-equivariant closed embedding
\begin{equation} \label{e:expl presentation}
\Sing(Z)\hookrightarrow {}^{cl}(V^*\times Z)\hookrightarrow V^*\times Z,
\end{equation}
where $V^*$ is the scheme $\Spec(\Sym(V))$. 

\ssec{The singular codifferential}

In this subsection, to a map between quasi-smooth DG schemes we will attach its singular
codifferential, which can be thought as the ``$H^{-1}$-version"  of the usual codifferential,
the latter being the map between classical cotangent spaces. 

\sssec{}

Let $f:Z_1\to Z_2$ be a map between quasi-smooth DG schemes. 
Define
$$\Sing(Z_2)_{Z_1}:={}^{cl}\!\left(\Sing(Z_2)\underset{Z_2}\times Z_1\right)\simeq
{}^{cl}\!\left(\Spec_{Z_1}\left(\on{Sym}_{\CO_{Z_1}}\left(f^*(T(Z_2)[1])\right)\right)\right).$$

\medskip

Note that $f$ induces a
morphism $T(Z_1)\to f^*(T(Z_2))$. Taking the Zariski spectra of the corresponding symmetric 
algebras, we obtain a map 
\begin{equation}  \label{e:codifferential}
\Sing(f):\Sing(Z_2)_{Z_1}\to \Sing(Z_1).
\end{equation} 
We will refer to this map as the ``singular codifferential."

\sssec{} \label{sss:q-smooth maps}

We have the following characterization of quasi-smooth maps between quasi-smooth DG schemes:

\begin{lem}  \label{l:codiff for q-smooth}
Let $f:Z_1\to Z_2$ be a morphism between quasi-smooth
DG schemes. Then $f$ is quasi-smooth if and only if
the singular codifferential $\Sing(f)$ is a closed embedding.
\end{lem}

\begin{proof} The relative cotangent complex $T^*(Z_1/Z_2)$ is the cone of the 
codifferential $$f^*(T^*(Z_2))\to T^*(Z_1).$$
Thus, $f$ is quasi-smooth if and only if 
the induced map
\[(df_x)^*:H^{-1}(T^*(Z_2)_{f(x)})\to H^{-1}(T^*(Z_1)_{x})\]
is injective for every $x\in Z_1$. The latter condition is equivalent to surjectivity of the morphism 
$H^1(T(Z_1))\to H^1(f^*(T(Z_2)))$, which is equivalent to $\Sing(f)$ being 
a closed embedding.
\end{proof}

\begin{lem} \label{l:codiff for smooth}
Let $f:Z_1\to Z_2$ be a morphism between quasi-smooth
DG schemes. Then $f$ is smooth if and only if the following two conditions are satisfied:
\begin{itemize}
\item The (classical) differential $df_x:H^0(T(Z_1)_x)\to H^0(T(Z_2)_{f(x)})$ is surjective for all
$k$-points $x\in X$;
\item The singular codifferential $\Sing(f)$ is an isomorphism.
\end{itemize}
\end{lem}
\begin{proof} The argument is similar to the proof of Lemma~\ref{l:codiff for q-smooth}; we leave it to the reader.
\end{proof}

\section{Support in triangulated and DG categories}  \label{s:sup}

In this section we will review the following construction. Given a triangulated category $\bT$
and a commutative graded algebra $A$ mapping to its center, we will define full subcategories in $\bT$
corresponding to closed (resp. open) subsets of $\Spec(A)$ in the Zariski topology. 

\medskip

This construction is a variant of that from \cite{Kr}. Unlike \cite{Kr}, we do not assume that the 
categories are compactly generated. Also, we use a coarser notion of support; see \remref{r:suppkrause}.

\ssec{Localization with respect to homogeneous elements}

\sssec{} \label{sss:localization a}   

Let $\bT$ be a cocomplete triangulated category.  Let $A$ be a commutative algebra, graded by 
even integers, and equipped with a homomorphism to the graded center of $\bT$. That is, 
for every $\bt\in \bT$ we have a homomorphism of graded algebras
\begin{equation} \label{e:action on object}
A\to \underset{n}\bigoplus\, \Hom_{\bT}(\bt,\bt[2n]),
\end{equation}
and for every $\phi:\bt'\to \bt''[m]$, the diagram 
$$
\CD
A  @>>>  \underset{n\geq 0}\bigoplus\, \Hom_{\bT}(\bt',\bt'[2n])  \\
@VVV    @VV{\phi\circ -}V  \\
\underset{n\geq 0}\bigoplus\, \Hom_{\bT}(\bt'',\bt''[2n])  @>{-\circ \phi}>>  \underset{n\geq 0}\bigoplus\, \Hom_{\bT}(\bt',\bt''[2n+m])  
\endCD$$
commutes. In this situation, we say that $A$ acts on the triangulated category $\bT$.

\sssec{}  \label{sss:abs category on open}

Let $a\in A$ be a homogeneous element. We let $Y_a\subset \Spec(A)$ be the conical 
(i.e., $\BG_m$-invariant) Zariski-closed subset of $\Spec(A)$ cut out by $a$. Here
$\Spec(A)$ is the Zariski spectrum of $A$ (where $A$ is viewed as a plain commutative algebra).

\medskip

We define the full subcategory $\bT_{\Spec(A)-Y_a}\subset \bT$ to consist
of those objects $\bt\in \bT$ for which the map
$$a:\bt\to \bt[2k]$$
is an isomorphism, where $2k=\deg(a)$. 

\medskip

Clearly, the subcategory $\bT_{\Spec(A)-Y_a}$ is thick (i.e., triangulated and closed under direct summands), 
and closed under taking arbitrary direct sums. 

\medskip

The inclusion $\bT\hookleftarrow \bT_{\Spec(A)-Y_a}$ admits a left adjoint,
explicitly given by
\begin{equation} \label{e:loc}
\bt\mapsto hocolim\left(\bt\overset{a}\to \bt[2k]\overset{a}\to\dots\right).
\end{equation}

We denote the resulting endofunctor
$$\bT\to \bT_{\Spec(A)-Y_a}\to \bT$$
by $\on{Loc}_a$. 

\begin{rem}
Although taking a homotopy colimit in a triangulated category is an operation that is defined only up to a 
non-canonical isomorphism, the expression in \eqref{e:loc} is canonical by virtue of being
a left adjoint. 
\end{rem}

\sssec{} Recall that a full thick subcategory $\bT'\subset\bT$ is said to be \emph{left admissible} if the inclusion
$\bT'\hookrightarrow\bT$ admits a left adjoint. If this is the case, we let $\bT'':={}^{\perp}\!(\bT')$ be its left orthogonal;
the inclusion $\bT''\hookrightarrow\bT$ admits a right adjoint. We say that the resulting diagram
\[\bT''\rightleftarrows\bT\rightleftarrows\bT'\]
is a \emph{short exact sequence of triangulated categories} if $\bT'$ is closed under direct sums. Note that
$\bT''$, being a left orthogonal, is automatically closed under direct sums.

\begin{lem} \label{lem:exact of triangles}
Let 
\[\bT''\overset{F''}{\underset{G''}{\rightleftarrows}}\bT\overset{F'}{\underset{G'}{\rightleftarrows}}\bT'\]
be a short exact sequence of categories. Then all four functors $F'$, $F''$, $G'$, $G''$ are 
triangulated (preserve exact triangles and shifts)
and continuous (preserve arbitrary direct sums). 
\end{lem}

\begin{proof} The inclusion functors $F''$, $G'$ are triangulated and continuous for tautological reasons. 
The functor $F'$ is continuous because it is a left adjoint; it is triangulated because its right adjoint $G'$
is triangulated. 

\medskip

We now see that the composition $G'\circ F'$ is continuous. 
Therefore, the composition $F''\circ G''$ is continuous as well, being the cone of the adjunction map between
the identity functor and $G'\circ F'$. This implies that $G''$ is continuous. Finally, $G''$ is triangulated
because it is the right adjoint of a triangulated functor.
\end{proof}

\sssec{}
Let \[\bT_{Y_a}:={}^\perp\!(\bT_{\Spec(A)-Y_a})\subset \bT\] 
be the left orthogonal of $\bT_{\Spec(A)-Y_a}$. We obtain an exact sequence of triangulated categories
\[
\bT_{Y_a}\rightleftarrows \bT\rightleftarrows \bT_{\Spec(A)-Y_a}.
\]
Denote the composition
\[\bT \to\bT_{Y_a}\to\bT\]
by $\on{co-Loc}_a$. 

\sssec{}  \label{sss:functor and a}
Suppose $\bT_1$ and $\bT_2$ are two cocomplete triangulated categories equipped with actions of $A$.
Let $F:\bT_1\to \bT_2$ be a continuous triangulated functor between triangulated categories, compatible
with the $A$-actions (in the obvious sense). It is clear that it sends $(\bT_1)_{\Spec(A)-Y_a}$ to $(\bT_2)_{\Spec(A)-Y_a}$.

\medskip

In addition, since $F$ is continuous and triangulated, it preserves homotopy colimits. Thus,
if $\bt_1\in \bT_1$ satisfies $\on{Loc}_a(\bt_1)=0$,
then $\on{Loc}_a(F(\bt_1))=0$. Hence, the functor $F$ sends $(\bT_1)_{Y_a}$ to 
$(\bT_2)_{Y_a}$. 

\medskip

This formally implies that the diagram 
$$
\xymatrix{
(\bT_1)_{Y_a} \ar[r] \ar[d]_F & \bT_1 \ar[r] \ar@<.7ex>[l] \ar[d]^F & (\bT_1)_{\Spec(A)-Y_a} \ar@<.7ex>[l] \ar[d]^F \\
(\bT_2)_{Y_a} \ar[r] & \bT_2 \ar[r] \ar@<.7ex>[l]  & (\bT_1)_{\Spec(A)-Y_a} \ar@<.7ex>[l] }
$$
is commutative.

\ssec{Zariski localization}

In this subsection we will show how localization with respect to individual homogeneous
elements can be organized into localization with respect to the Zariski topology on $\Spec(A)$.
This will not involve much beyond the usual constructions in commutative algebra. 

\sssec{}

Let $a_1,a_2\in A$ be two homogeneous elements. 

\begin{lem}  \label{l:localizations compatible}
The functor $\on{Loc}_{a_2}$ preserves both 
$\bT_{Y_{a_1}}$ and $\bT_{\Spec(A)-Y_{a_1}}$
\end{lem}

\begin{proof}

Follows from \secref{sss:functor and a} applied to the tautological
embeddings $\bT_{Y_{a_1}}\hookrightarrow \bT$ and $\bT_{\Spec(A)-Y_{a_1}}\hookrightarrow \bT$
for the action of $a_2$.

\end{proof}

\begin{rem}
Note that \lemref{l:localizations compatible} did not use the fact that the actions of $a_1$ and $a_2$
commute.
\end{rem}

\sssec{}

From \lemref{l:localizations compatible} we obtain that the short exact sequences
$$\bT_{Y_{a_1}}\rightleftarrows \bT\rightleftarrows \bT_{\Spec(A)-Y_{a_1}} \text{ and }
\bT_{Y_{a_2}}\rightleftarrows \bT\rightleftarrows \bT_{\Spec(A)-Y_{a_2}}$$
are compatible in the sense of \cite{BeVo}, Sect. 1.3. Thus we obtain a 
commutative diagram in which every row and every column is a short exact sequence:
\begin{equation} \label{e:localizations commute}
\xymatrix{
\bT_{Y_{a_1}}\cap \bT_{Y_{a_2}} \ar[r] \ar[d] & \bT_{Y_{a_1}} \ar[r] \ar@<.7ex>[l] \ar[d] & \bT_{Y_{a_1}}\cap  \bT_{\Spec(A)-Y_{a_2}} 
\ar@<.7ex>[l] \ar[d]  \\
\bT_{Y_{a_2}} \ar[r] \ar[d] \ar@<.7ex>[u]  & \bT  \ar[r] \ar[d] \ar@<.7ex>[u] \ar@<.7ex>[l] & \bT_{\Spec(A)-Y_{a_2}} \ar[d] \ar@<.7ex>[u] \ar@<.7ex>[l]  \\
\bT_{\Spec(A)-Y_{a_1}}\cap \bT_{Y_{a_2}}  \ar[r]  \ar@<.7ex>[u]  &  \bT_{\Spec(A)-Y_{a_1}} \ar[r] \ar@<.7ex>[u] \ar@<.7ex>[l]
& \bT_{\Spec(A)-Y_{a_1}}\cap \bT_{\Spec(A)-Y_{a_2}}. \ar@<.7ex>[u] \ar@<.7ex>[l]}
\end{equation}

In particular, the functors $\on{Loc}_{a_1}$, $\on{Loc}_{a_2}$, $\on{co-Loc}_{a_1}$ and $\on{co-Loc}_{a_2}$ pairwise
commute. 

\sssec{}

Note that the fact that $a_1$ and $a_2$ commute implies that
$$\bT_{\Spec(A)-Y_{a_1}}\cap \bT_{\Spec(A)-Y_{a_2}}=\bT_{\Spec(A)-Y_{a_1\cdot a_2}}.$$

\medskip

Our next goal is to prove the following:

\begin{prop}  \label{p:ideals}
If $a\in A$ is a homogeneous element contained in the radical of the ideal generated by
$a_1,\dots,a_n$, then 
$$\bT_{Y_{a_1}}\cap\dots\cap \bT_{Y_{a_n}} \subset \bT_{Y_a}.$$
\end{prop}

Prior to giving the proof, we will need the following more explicit description of the category $\bT_{Y_a}$. 

\sssec{}

We start with a remark about  $A$-modules, valid for an arbitrary commutative graded ring $A$.

\medskip

Consider the (DG) category of graded $A$-modules, i.e., 
$(A\mod)^{\BG_m}$. Fix a homogeneous element $a\in A$. We identify the DG category of graded modules 
over the localization, i.e., $(A_a\mod)^{\BG_m}$, with a full subcategory of $(A\mod)^{\BG_m}$. 

\begin{lem} \label{l:cosupp}
For $M\in (A\mod)^{\BG_m}$, the following conditions are equivalent:

\smallskip

\noindent {\em (a)} $\Hom_A(A_a(j),M[i])=0,\, \forall i,j\in \BZ$. Here $A_a(j)$ refers to $A_a$ with 
grading shifted by $j$.

\smallskip

\noindent {\em (b)} For any $N\in(A_a\mod)^{\BG_m}$, $\Hom_A(N,M)=0$.

\smallskip

\noindent {\em (c)} The map 
\[\prod_{i=0}^\infty M\to\prod_{i=0}^\infty M:(m_0,m_1,\dots)\mapsto(m_0-a(m_1),m_1-a(m_2),\dots)\]
is an isomorphism. Here $\prod$ stands for the product in the category of graded $A$-modules.

\smallskip

\noindent {\em (d)} The homotopy limit 
\[holim\left(M\overset{a}\leftarrow M\overset{a}\leftarrow\dots\right)\]
vanishes.
\end{lem}

\begin{proof} Since the objects $A_a(j)$ generate $(A_a\mod)^{\BG_m}$, (a) and (b) are equivalent. Moreover, the space
$\Hom_A(A_a(j),M[i])$ identifies with the $j$-th graded component of the $i$-th cohomology of the cone of the map from
(c); therefore, (a) and (c) are equivalent. Finally, (c) and (d) are equivalent by definition.
\end{proof} 

\sssec{} Let us denote the full subcategory of $(A\mod)^{\BG_m}$, satisfying the equivalent conditions of \lemref{l:cosupp}
by \[(A\mod)^{\BG_m}_{\langle a\rangle}\subset(A\mod)^{\BG_m}.\] 
By condition (b), this subcategory is the right orthogonal of $(A_a\mod)^{\BG_m}$.
Note that $(A\mod)^{\BG_m}_{\langle a\rangle}$ is a thick subcategory that is closed under products, but
not co-products. 

\begin{rem}
One should think of $(A\mod)^{\BG_m}_{\langle a\rangle}$ as the category of $a$-adically complete $A$-modules.
\end{rem} 

\medskip

We now return to the setting of \secref{sss:localization a}.

\begin{lem} \label{l:colocalization via cosupport}
For an object $\bt\in\bT$, the following conditions are equivalent:

\smallskip

\noindent{\em (a)} $\bt\in\bT_{Y_a}$.

\noindent{\em (b)} For any $\bt'\in\bT$, the graded $A$-module $\Hom_\bT^\bullet(\bt,\bt')$ belongs to
$(A\mod)^{\BG_m}_{\langle a\rangle}$.
\end{lem}

\begin{proof} Indeed, using \lemref{l:cosupp}(c), we see that (b) is equivalent to $$\Hom_\bT^\bullet(\on{Loc}_a(\bt),\bt')=0.$$
\end{proof}

\sssec{}

We are now ready to prove \propref{p:ideals}:

\begin{proof}
By \lemref{l:colocalization via cosupport}, it suffices to show that
\[(A\mod)^{\BG_m}_{\langle a_1\rangle}\cap\dots\cap (A\mod)^{\BG_m}_{\langle a_n\rangle}\subset 
(A\mod)^{\BG_m}_{\langle a\rangle}.\]
Using \lemref{l:cosupp}(b), we see that it is enough to prove that $(A_a\mod)^{\BG_m}$ is contained in 
the full subcategory of $(A\mod)^{\BG_m}$ generated by the subcategories $(A_{a_i}\mod)^{\BG_m}$ for $i=1,\dots,n$. 
But this is obvious from the \v{C}ech resolution. (In fact, the latter subcategory identifies
with the category of $\BG_m$-equivariant modules on the scheme $\Spec(A)-\bigcap_i Y_{a_i}$.)
\end{proof}

\ssec{The definition of support}

In this subsection we will finally define the support of an object, and study how
this notion behaves under functors between triangulated categories and morphisms
between algebras. 

\sssec{}  \label{sss:subcategory for closed}

Let $Y$ be a conical (i.e., $\BG_m$-invariant) Zariski-closed subset of $\Spec(A)$.

\medskip

We define the full subcategory
$$\bT_Y:=\underset{a}\cap\, \bT_{Y_a},$$
where the intersection is taken over the set of homogeneous elements of $a\in A$ that vanish on $Y$.

\medskip

Suppose $Y_1, Y_2\subset\Spec(A)$ are two closed conical subsets. By \propref{p:ideals},   
\begin{equation} \label{e:intersection of closed}
\bT_{Y_1\cap Y_2}=\bT_{Y_1}\cap \bT_{Y_2}.
\end{equation}

\sssec{}

We give the following definitions:

\begin{defn} \label{defn:support}
Given a conical Zariski-closed subset $Y\subset\Spec(A)$ and $\bt\in\bT$,
we say that 
$$\on{supp}_A(\bt)\subset Y$$
if $\bt\in \bT_Y$.
\end{defn}

\begin{defn} 
Given $\bt\in\bT$, we define $\on{supp}_A(\bt)$ to be the minimal 
conical Zariski-closed subset $Y\subset\Spec(A)$ such that $\bt\in\bT_Y$.
\end{defn}

\begin{rem} \label{r:suppkrause}
The definition of support given in Definition \ref{defn:support} differs from the one in \cite{Kr}.
When $\bT$ is compactly generated, so that the definition of \cite{Kr} applies, 
what we call ``support" is the Zariski closure of the support from \cite{Kr}.
\end{rem}

\sssec{}

It is clear that 
$$\on{supp}_A(\bt)=\underset{a}\cap\, Y_a,$$
where the intersection is taken over the set of homogeneous elements $a$
such that $\bt\in \bT_{Y_a}$.

\begin{lem} \label{l:T_Y}
Let $Y\subset\Spec(A)$ be a conical Zariski-closed subset whose complement  
$\Spec(A)-Y$ is quasi-compact. 
(If $A$ is Noetherian, this condition is automatic.) 
Then the embedding $\bT_Y\hookrightarrow \bT$ admits a continuous right adjoint. 
\end{lem}
\begin{proof} By the assumption, there 
exists a finite collection of homogeneous elements $a_1,\dots,a_n\in A$
such that
\begin{equation} \label{e:Y cut}
Y=Y_{a_1}\cap\dots\cap Y_{a_n}.
\end{equation}

\medskip

By \propref{p:ideals}, we then have:
$$\bT_Y=\bT_{Y_{a_1}}\cap\dots\cap \bT_{Y_{a_n}}.$$
Iterating the diagram \eqref{e:localizations commute}, we see that 
the embedding
\[\bT_{Y_{a_1}}\cap\dots\cap \bT_{Y_{a_n}}\hookrightarrow\bT\]
admits a continuous right adjoint such that the composed functor
$$\bT\to \bT_{Y_{a_1}}\cap\dots\cap \bT_{Y_{a_n}}\hookrightarrow \bT$$
is isomorphic to the composition 
$$\on{co-Loc}_{a_1}\circ\dots\circ \on{co-Loc}_{a_n}.$$
\end{proof}

\sssec{} Let $Y\subset\Spec(A)$ be a conical Zariski-closed subset whose
complement is quasi-compact. From \lemref{l:T_Y}, we obtain a short exact sequence of categories
\[\bT_Y\rightleftarrows \bT\rightleftarrows \bT_{\Spec(A)-Y},\]
where $\bT_{\Spec(A)-Y}$ is the right orthogonal to $\bT_Y$. We also see that $\bT_{\Spec(A)-Y}$
is generated by categories $\bT_{\Spec(A)-Y_a}$, where $a\in A$ runs over homogeneous elements such that
$Y_a\supset Y$. (In fact, it suffices to consider $a=a_i$ for a finite collection of homogeneous 
elements $a_1,\dots,a_n\in A$ satisfying \eqref{e:Y cut}.) 

\begin{cor} \label{c:intersection of closed}
Suppose $Y_1,Y_2\subset\Spec(A)$ are conical Zariski-closed subsets whose
complements are quasi-compact. Then the category 
$\bT_{Y_1\cup Y_2}$ is generated by $\bT_{Y_1}$ and $\bT_{Y_2}$.
\end{cor}
\begin{proof}
Similar to \eqref{e:localizations commute}, we have a diagram 
\begin{equation}
\xymatrix{
\bT_{Y_1}\cap \bT_{Y_2} \ar[r] \ar[d] & \bT_{Y_1} \ar[r] \ar@<.7ex>[l] \ar[d] & \bT_{Y_1}\cap  \bT_{\Spec(A)-Y_2} 
\ar@<.7ex>[l] \ar[d]  \\
\bT_{Y_2} \ar[r] \ar[d] \ar@<.7ex>[u]  & \bT  \ar[r] \ar[d] \ar@<.7ex>[u] \ar@<.7ex>[l] & \bT_{\Spec(A)-Y_2} \ar[d] \ar@<.7ex>[u] \ar@<.7ex>[l]  \\
\bT_{\Spec(A)-Y_1}\cap \bT_{Y_{a_2}}  \ar[r]  \ar@<.7ex>[u]  &  \bT_{\Spec(A)-Y_1} \ar[r] \ar@<.7ex>[u] \ar@<.7ex>[l]
& \bT_{\Spec(A)-Y_1}\cap \bT_{\Spec(A)-Y_2} \ar@<.7ex>[u] \ar@<.7ex>[l]}
\end{equation}
with exact rows and columns. In order to prove the corollary, it suffices to check that
\[\bT_{\Spec(A)-Y_1}\cap \bT_{\Spec(A)-Y_2}=\bT_{\Spec(A)-Y_1\cup Y_2}.\]
Clearly, the right-hand side is contained in the left-hand side. On the other hand, 
\[\bT_{\Spec(A)-Y_1}\cap \bT_{\Spec(A)-Y_2}=\left(\bT_{\Spec(A)-Y_1}\right)_{\Spec(A)-Y_2}\]
is generated by the essential images 
\[\on{Loc}_{a_1\cdot a_2}(\bT)=\on{Loc}_{a_2}\circ\on{Loc}_{a_1}(\bT),\]
where $a_1,a_2\in A$ run over homogeneous elements such that $Y_1\subset Y_{a_1}$ and $Y_2\subset Y_{a_2}$.
This proves the converse inclusion.
\end{proof}

\sssec{}  \label{sss:change of category}

Let $F:\bT_1\to \bT_2$ be a continuous triangulated functor compatible with the actions of $A$. 
Let $Y\subset\Spec(A)$ be a conical Zariski-closed subset whose complement is quasi-compact. 
 It is clear from \secref{sss:functor and a} 
that $F$ induces a commutative diagram of functors:
$$
\xymatrix{
(\bT_1)_{Y} \ar[r] \ar[d]_F & \bT_1 \ar[r] \ar@<.7ex>[l] \ar[d]^F & (\bT_1)_{\Spec(A)-Y} \ar@<.7ex>[l]  \ar[d]^F\\
(\bT_2)_{Y} \ar[r] & \bT_2 \ar[r] \ar@<.7ex>[l] & (\bT_2)_{\Spec(A)-Y}. \ar@<.7ex>[l]  }
$$

\medskip

Thus, for any $\bt\in\bT_1$, we have $\on{supp}_A(\bt)\supset\on{supp}_A(F(\bt))$. If we assume that $F$ is conservative,
then $\on{supp}_A(\bt)=\on{supp}_A(F(\bt))$.

\medskip

In particular if $\bT'\subset \bT$ is a full triangulated subcategory closed under direct sums, then
$$\bT'_Y=\bT_Y\cap \bT' \text{ and } \bT'_{\Spec(A)-Y}=\bT_{\Spec(A)-Y}\cap \bT' $$
as subcategories of $\bT$.

\sssec{} \label{sss:change of algebras}

The notion of support behaves functorially under homomorphisms of algebras. Namely, let $\phi:A'\to A$
be a homomorphism of evenly graded algebras. Let $\Phi$ denote the resulting map
$\Spec(A)\to \Spec(A')$. For $\bT$ as above, the algebra $A'$ maps to the graded
center of $\bT$ by composing with $\phi$.

\medskip

We have:

\begin{lem}  \label{l:change of algebras}
For $\bt\in \bT$ and $Y'\subset \Spec(A')$ and $Y:=\Phi^{-1}(Y')$, 
$$\on{supp}_{A'}(\bt)\subset Y'\Leftrightarrow \on{supp}_A(\bt)\subset Y.$$
\end{lem}

Equivalently for $\bt\in \bT$,
$$\on{supp}_{A'}(\bt)=\overline{\Phi(\on{supp}_{A}(\bt))}.$$

\ssec{The compactly generated case}

In this subsection we will show that if $\bT$ is compactly generated, we can measure
supports of objects more explicitly. 

\sssec{}

Assume that $\bT$ is compactly generated. 

\begin{lem} \label{l:compact generated}
Let $Y\subset \Spec(A)$ be a conical Zariski-closed subset whose complement $\Spec(A)-Y$ is quasi-compact.
Then the category $\bT_Y$ is compactly generated.
\end{lem}

\begin{proof}
By induction and \eqref{e:localizations commute}, we can assume that $Y$ is cut out by one
homogeneous element $a$. It is easy to see that  the objects
$$\on{Cone}(\bt\overset{a}\longrightarrow \bt),\quad \bt\in \bT^c$$
generate $\bT_{Y_a}$. Indeed, the right orthogonal to the class of these objects coincides with $\bT_{\Spec(A)-Y_a}$.

\end{proof}

\sssec{}

One can use compact objects to rewrite the definition of support:

\begin{lem}  \label{l:supp via comp}
Let $Y$ be an arbitrary conical Zariski-closed subset of $\Spec(A)$. 

\smallskip

\noindent{\em(a)} For $\bt\in \bT$, its support is contained in $Y$ if and only
if for a set of compact generators $\bt_\alpha\in \bT$, the support of
the $A$-module
$$\Hom^\bullet_{\bT}(\bt_\alpha,\bt)$$
is contained in $Y$ for all $\alpha$ (cf. \cite[Corollary~5.3]{Kr}.)

\smallskip

\noindent{\em(b)} If $\bt$ is compact, its support is contained in $Y$ if and
only if the support of the $A$-module $\Hom^\bullet_{\bT}(\bt,\bt)$ is contained in $Y$.

\smallskip

\noindent{\em(c)} If $\bt$ is compact, and $a\in A$ is a homogeneous element that vanishes
on $\on{supp}_A(\bt)$, then there exists an integer $i$ such that $\bt\overset{a^i}\to \bt[2k\cdot i]$
vanishes. Here $2k=\deg(a)$. 
\end{lem}

\begin{proof}
Let $a$ be a homogeneous element of $A$ of degree $2k$. Suppose that $a$ vanishes on $Y$. The fact that
$\on{supp}(\bt)\subset Y_a$ is equivalent to 
the colimit
$$\bt\overset{a}\to \bt[2k]\overset{a}\to\dots$$
being zero, which can be tested by mapping the generators $\bt_\alpha[m]$, $m\in \BZ$ into this colimit. Since the $\bt_\alpha$'s 
are compact, the above Hom is isomorphic to the colimit
$$\Hom_{\bT}(\bt_\alpha,\bt[-m])\overset{a}\to 
\Hom_{\bT}(\bt_\alpha,\bt[2k-m])\overset{a}\to\dots,$$
taken in the category $\Vect^\heartsuit$.
The vanishing of the latter is equivalent to 
$$\Hom^\bullet_{\bT}(\bt_\alpha,\bt)$$
being supported over $Y$ as an $A$-module, which is the assertion of point (a) of the lemma.

\medskip

For point (b), the ``only if" direction follows from point (a). The ``if" direction
holds tautologically for any $\bt$ (with no compactness hypothesis).

\medskip

Point (c) follows from point (b): the unit element in $\Hom_{\bT}(\bt,\bt)$
is annihilated by some power of $a$. 

\end{proof}

\ssec{Support in DG categories}  \label{ss:support in DG}

From now on we will assume that $\bT$ is the homotopy category of a DG category $\bC$,
equipped with an action of an $\BE_2$-algebra $\CA$ (see \secref{sss:action of E2 on category},
where the notion of action of an $\BE_2$-algebra on a DG category is discussed). 

\medskip 

We will show that the notion of support in $\bC$ can be expressed in terms of the universal
situation, namely, for $\bC=\CA\mod$. In addition, we will study the behavior of support under
tensor products of the $\bC$'s.  

\sssec{}  \label{sss:support via E2}

Set 
$$A:=\underset{n}\bigoplus\, H^{2n}(\CA).$$
Since $\CA$ has an $\BE_2$-structure, the algebra $A$ is commutative. The action of $\CA$ on $\bC$
gives rise to a homomorphism from $A$ to the graded center of $\bT$.

\sssec{}

For a conical Zariski-closed subset $Y\subset \Spec(A)$, we let 
$$\bC_Y\subset \bC$$
be the full DG subcategory of $\bC$ defined as the preimage of $\bT_Y\subset\bT$.

\sssec{}  \label{sss:E2 on itself}

In particular, we can consider $\bC=\CA\mod$. It is clear that the resulting subcategory
$$\CA\mod_Y\subset \CA\mod$$
is a (two-sided) monoidal ideal. (In fact, any full cocomplete subcategory of $\CA\mod$ is a monoidal ideal, since
$\CA\mod$ is generated by $\CA$, which is the unit object.) 

\sssec{}

The following assertion will play a crucial role:

\begin{prop}  \label{p:support via monoidal}
Let $Y$ be such that its complement is quasi-compact. Then for any DG category $\bC$ equipped with an 
action of $\CA$, we have 
$$\bC_Y=\CA\mod_Y\underset{\CA\mod}\otimes\bC $$
as full subcategories of 
$$\bC\simeq \CA\mod\underset{\CA\mod}\otimes \bC.$$
\end{prop}

\begin{proof}

Note that if
$$\bC_1\rightleftarrows \bC_2\rightleftarrows \bC_3$$
is a short exact sequence of right modules over a DG monoidal category $\bO$,
and $\bC'$ is a left module, then 
$$\bC_1\underset{\bO}\otimes \bC'\rightleftarrows 
\bC_2\underset{\bO}\otimes \bC'\rightleftarrows \bC_3\underset{\bO}\otimes \bC'$$
is a short exact sequence of DG categories.

\medskip

This observation together with \eqref{e:localizations commute} for $\bC$ and $\CA\mod$ reduces
the proposition to the case when $Y=Y_a$ for some homogeneous element $a\in A$. In this case,
it is sufficient to show that
$$\bC_{\Spec(A)-Y_a} \text{ and } \CA\mod_{\Spec(A)-Y_a}\underset{\CA\mod}\otimes\bC $$
coincide as subcategories of $\bC$.

\medskip

First, let us show the inclusion $\supset$, i.e., we have to show that the element $a$ acts
as an isomorphism on objects from $\CA\mod_{\Spec(A)-Y_a}\underset{\CA\mod}\otimes\bC$.
This property is enough to establish on the generators, which we can take to be of the form
$\CM\otimes \bc$, where $\bc\in \bC$ and $\CM\in \CA\mod_{\Spec(A)-Y_a}$. The action of
$a$ on such an object equals
$$a_{\CM}\otimes\on{id}_{\bc},$$
and the assertion follows from the fact that $a$ is an isomorphism on $\CM$.

\medskip

In particular, we obtain a natural transformation of endofunctors 
$$\on{Loc}_{a,\bC}\to \on{Loc}_{a,\CA\mod}\otimes\on{Id}_\bC,$$
viewed as acting on 
$$\bC\simeq \CA\mod\underset{\CA\mod}\otimes\bC.$$
It suffices to show that this natural transformation is an isomorphism.
The latter follows immediately from \eqref{e:loc}.

\end{proof}

\sssec{}  \label{sss:abs product}

Suppose now we have two $\BE_2$-algebras $\CA_i$ acting on DG categories $\bC_i$, respectively ($i=1,2$).
Let $Y_i\in \Spec(A_i)$ be conical Zariski-closed subsets whose complements are quasi-compact. 

\medskip

Set $\bC:=\bC_1\otimes \bC_2$, $\CA=\CA_1\otimes \CA_2$. We then have a natural graded homomorphism
\[\phi:A_1\otimes A_2\to A,\]
where 
\[A_i:=\underset{n}\bigoplus\, H^{2n}(\CA_i)\qquad(i=1,2).\]
It induces a map
$$\Spec(A)\to \Spec(A_1)\times \Spec(A_2);$$
let $Y\subset \Spec(A)$ be the preimage of $Y_1\times Y_2\subset \Spec(A_1)\times \Spec(A_2)$. 

\medskip

As in \propref{p:support via monoidal}, one shows:

\begin{prop} \label{p:abs product}
The subcategories
$$(\bC_1)_{Y_1}\otimes (\bC_2)_{Y_2} \text{ and } \bC_Y$$
of $\bC_1\otimes \bC_2=\bC$
coincide.
\end{prop}

\sssec{} Let $\bC_i$, $\CA_i$, $A_i$ ($i=1,2$) be as in \secref{sss:abs product}. Suppose that $\bC_1$ is dualizable.
Let $F:\bC_1\to \bC_2$ be a continuous functor. Such functors are in a bijection with objects
$$\bc'\in \bC':=\bC_1^\vee\otimes \bC_2.$$

Note that $\bC_1^\vee$ is acted on by $\CA_1^{\on{op}}$ (see \secref{sss:E2 action on dual}).

\medskip

We can regard $\bC_1^\vee\otimes \bC_2$ 
as acted on by the $\BE_2$-algebra $\CA':=\CA^{\on{op}}_1\otimes \CA_2$. Let $A$ be the corresponding 
graded algebra; we have a natural morphism
$\phi:A_1\otimes A_2\to A$. (Note that the graded algebra corresponding to $\CA_1^{\on{op}}$ coincides with 
$A_1$.) Let $p_1,p_2$ be the two components of the corresponding map
\[(p_1,p_2):\Spec(A)\to\Spec(A_1)\times\Spec(A_2).\]
We have:

\begin{prop} \label{p:supp functoriality}
Let $Y_i\subset \Spec(A_i)$ (for $i=1,2$) be conical Zariski-closed subsets such that
the complement of $Y_1$ is quasi-compact. Let $\bc'$ be the object of $\bC_1^\vee\otimes \bC_2$
corresponding to $F$. Suppose that  
$$p_2(\on{supp}_A(\bc')\cap p_1^{-1}(Y_1))\subset Y_2.$$ 
Then the functor $F$ maps $(\bC_1)_{Y_1}$ to $(\bC_2)_{Y_2}$.
\end{prop}

\begin{proof} Set $Y_1'=p_1^{-1}(Y_1)\subset\Spec(A)$. It is a conical Zariski-closed subset whose complement is
quasi-compact. Consider the corresponding exact sequence of categories
\[\bC'_{Y_1'}\rightleftarrows\bC'\rightleftarrows\bC'_{\Spec(A)-Y_1'}.\] 
It is clear that the objects of $\bC'_{\Spec(A)-Y_1'}$ correspond to functors $\bC_1\to\bC_2$ that vanish on
$(\bC_1)_{Y_1}$. Therefore, we may replace $\bc'$ by its colocalization and assume that $\bc'\in\bC'_{Y_1'}$. 
We then have $p_2(\on{supp}_A(\bc'))\subset Y_2$. For such $\bc'$, it is clear that the essential image of the corresponding
functor $\bC_1\to\bC_2$ is contained in $(\bC_2)_{Y_2}$. 
\end{proof}

\ssec{Grading shift for $\BE_2$-algebras}  

In this subsection we will show how to relate the notion of support developed in the previous
subsections to the more algebro-geometric notion of support \emph{over} an algebraic stack. 

\medskip

This subsection may be skipped on the first reading, and returned to when necessary.

\sssec{}  \label{sss:good grading}

Suppose that the $\BE_2$-algebra $\CA$ carries an action of $\BG_m$ such that the corresponding $\BE_2$-algebra
$\CA^{\on{shift}}$ (see \secref{sss:shift of grading alg}) is classical. 

\medskip

In particular, $\CA^{\on{shift}}$ has a canonical $E_\infty$ (i.e., commutative algebra) structure, which restricts to the initial $\BE_2$-structure.
Hence, the same is true for $\CA$. 

\medskip

We thus obtain 
a canonical isomorphism
$$\CA^{\on{shift}}\simeq A$$
as classical commutative algebras, which is compatible with the grading after scaling the grading on the left-hand side by $2$.

\sssec{}
Consider the stack $\CS_{\CA}=\Spec(\CA^{\on{shift}}/\BG_m)$. Given a 
conical Zariski-closed subset $$Y\subset\Spec(A),$$ we regard $Y/\BG_m$ as a closed
substack in $\CS_{\CA}$. 

\medskip

Recall that by \secref{sss:shift of grading alg}, we have a 
canonical equivalence of DG categories
$$(\CA\mod)^{\BG_m}\simeq \QCoh(\CS_{\CA}).$$

This equivalence naturally extends to an equivalence of (symmetric) monoidal categories.

\begin{prop}  \label{p:support via graded}
We have 
$$\CA\mod_Y=\CA\mod\underset{\QCoh(\CS_{\CA})}\otimes \QCoh(\CS_{\CA})_{Y/\BG_m}$$
as full subcategories of
$$\CA\mod= \CA\mod\underset{(\CA\mod)^{\BG_m}}\otimes (\CA\mod)^{\BG_m}\simeq
\CA\mod\underset{\QCoh(\CS_{\CA})}\otimes \QCoh(\CS_{\CA}).$$
\end{prop}

\begin{proof}

Note first that an action of $\BG_m$ on an associative DG algebra $\CA$ defines a bi-grading on $A$. 
Suppose that $\CA$ is an $\BE_2$-algebra, and let $Y\subset\Spec(A)$ be a Zariski-closed subset
conical with respect to the cohomological grading. Then the corresponding subcategory
$\CA\mod_Y$ is $\BG_m$-invariant (see \secref{sss:de-eq subcategory}) if and only if $Y$ is conical with respect
to both gradings. 

\medskip

Note, however, that the assumptions of the proposition imply that the two gradings on $A$
coincide. Hence, for any conical $Y$, the subcategory $\CA\mod_Y$ is $\BG_m$-invariant.

\medskip

By \secref{sss:de-eq reconstr}, it suffices to show that 
$$(\CA\mod_Y)^{\BG_m} \text{ and } \QCoh(\CS_{\CA})_{Y/\BG_m}$$
coincide as subcategories of $(\CA\mod)^{\BG_m}\simeq \QCoh(\CS_{\CA})$. 

\medskip

Note that by \eqref{e:eq subcategory}, the category $(\CA\mod_Y)^{\BG_m}$ identifies 
with the full subcategory of $(\CA\mod)^{\BG_m}$ consisting of modules supported on $Y$
as plain $\CA$-modules. This makes the required assertion manifest. 
\end{proof}

\sssec{}

Let $\bC$ be a DG category acted on by $\CA$, where $\CA$ is as in \secref{sss:good grading}. In particular,
we obtain that $\bC$ is a module category over $\QCoh(\CS_\CA)$. Let $Y\subset \Spec(A)$
be a conical Zariski-closed subset such that its complement is quasi-compact. 

\medskip

Combining Propositions \ref{p:support via monoidal} and \ref{p:support via graded}, we obtain:

\begin{cor} \label{c:support by tensor}
$\bC_Y\simeq \bC\underset{\QCoh(\CS_\CA)}\otimes \QCoh(\CS_\CA)_{Y/\BG_m}$
as subcategories of $\bC$.
\end{cor}

In particular, for $\bc\in \bC$ we can express $\on{supp}_A(\bc)$ in terms of the more
familiar notion of support of an object in a category tensored over $\QCoh$ of an
algebraic stack.

\sssec{} Suppose that the algebra $A$ is Noetherian.
Let us show that in this case we can use fibers to study supports of objects.

\medskip

Let $i_s:\Spec(k')\to\Spec(A)$
be a geometric point of $\Spec(A)$. We have natural monoidal functors
\[\QCoh(\CS_\CA)\to A\mod\to\Vect_{k'},\]
where $\Vect_{k'}$ is the category of vector spaces over $k'$. This defines an action of 
the monoidal category $\QCoh(\CS_\CA)$ on $\Vect_{k'}$.

\medskip

Given $\bc\in\bC$, we define $i_s^*(\bc)$ to be the object 
\[\bc\otimes k'\in\bC\underset{\QCoh(\CS_\CA)}\otimes\Vect_{k'}.\]
By Noetherian induction, one proves the following:

\begin{lem} \label{l:fiberwise vanishing}
If $i_s^*(\bc)=0$ for all geometric points $s$ of $A$, then $\bc=0$. 
\end{lem}

\medskip

As a consequence, we obtain:

\begin{cor} \label{c:support by category fibers}
Let $Y\subset\Spec(A)$ be a conical Zariski-closed subset. Fix
$\bc\in\bC$. Then

\smallskip
\noindent {\em(a)} $\bc\in\bC_Y$ if and only if $i_s^*(\bc)=0$ for all $s\notin Y$;

\smallskip
\noindent {\em(b)} $\bc\in\bC_{\Spec(A)-Y}$ if and only if $i_s^*(\bc)=0$ for all
$s\in Y$;

\smallskip
\noindent {\em(c)} $\on{supp}_A(\bc)$ is the Zariski closure of the set 
\[\{s\in\Spec(A):i_s^*(\bc)\ne 0\}.\] 
\end{cor}
\begin{proof} 
Recall that we have an exact sequence of categories
\begin{equation}\label{e:exactseq}
\bC_{Y}\rightleftarrows\bC\rightleftarrows\bC_{\Spec(A)-Y},
\end{equation}
which identifies with
\[\bC\underset{\QCoh(\CS_\CA)}\otimes\QCoh(\CS_\CA)_{Y/\BG_m}\rightleftarrows
\bC\underset{\QCoh(\CS_\CA)}\otimes\QCoh(\CS_\CA)\rightleftarrows
\bC\underset{\QCoh(\CS_\CA)}\otimes\QCoh(\CS_\CA-Y/\BG_m).\]
The ``only if'' direction in part (a) follows because
\[\QCoh(\CS_\CA)_{Y/\BG_m}\underset{\QCoh(\CS_\CA)}\otimes\Vect_{k'}=0\]
for any point $i_s:\Spec(k')\to\Spec(A)$ not contained in $Y$. Similarly,
the ``only if'' direction in part (b) follows because 
\[\QCoh(\CS_\CA-Y/\BG_m)\underset{\QCoh(\CS_\CA)}\otimes\Vect_{k'}=0\]
for any point $i_s:\Spec(k')\to\Spec(A)$ contained in $Y$.
Now the ``if'' directions in both parts follow from the sequence \eqref{e:exactseq} and \lemref{l:fiberwise vanishing}.
Part (c) follows from part (a).
\end{proof}

\bigskip

\bigskip

\centerline{\bf Part II: The theory of singular support}

\section{Singular support of ind-coherent sheaves}  \label{s:sing}

For the rest of the paper, we will be working with DG schemes locally almost
of finite type over a ground field $k$, which is assumed to have characteristic $0$. 

\medskip

In this section we introduce the notion of singular support for objects of $\IndCoh(Z)$, where $Z$ is
a quasi-smooth DG scheme, and study the basic properties of the corresponding categories $\IndCoh_Y(Z)$,
where $Y\subset \Sing(Z)$ is a conical Zariski-closed subset.

\ssec{The definition of singular support}  \label{ss:sing}

Throughout this section, $Z$ will be a quasi-smooth DG scheme.  It will be assumed \emph{affine}
until \secref{sss:nonaffine}.

\sssec{}  \label{l:hochschild}
Consider the $\BE_2$-algebra of Hochschild cochains $\on{HC}(Z)$; see \secref{ss:two isoms for HH}
where the definition of $\on{HC}(Z)$ is recalled (and see also \secref{s:E2} for some background
material on $\BE_2$-algebras).

\medskip

As is explained in \secref{ss:two isoms for HH}, the $\BE_2$-algebra $\on{HC}(Z)$ 
identifies canonically with the $\BE_2$-algebra of Hochschild cochains $\on{HC}(Z)^\IndCoh$
of the category $\IndCoh(Z)$. 

\medskip

Let  $\on{HH}^\bullet(Z)$ denote the classical graded associative algebra
\begin{equation} \label{e:HC}
\underset{n}\bigoplus\, H^n\left(\on{HC}(Z))\right).
\end{equation}

Let $\on{HH}^{\on{even}}(Z)$ denote the even part of $\on{HH}^\bullet(Z)$, viewed as a classical 
graded associative algebra. As was mentioned in \secref{sss:support via E2}, the algebra 
$\on{HH}^{\on{even}}(Z)$ is commutative, and $\on{HH}^{\on{even}}(Z)$ maps to the graded center of $\on{Ho}(\IndCoh(Z))$.

\sssec{}

Since $Z$ was assumed quasi-smooth, $T^*(Z)$ is perfect, and we can regard $T(Z)[-1]$ as a 
Lie algebra\footnote{For reasons of tradition, while we call Lie algebras in an arbitrary symmetric monoidal
$\infty$-category $\bO$ ``Lie algebras", we refer to Lie algebras in $\Vect$ as ``DG Lie algebras".}
in $\QCoh(Z)$, see \corref{c:HH as univ env}. Note that from \corref{c:HH as univ env} we obtain a canonical map of commutative algebras
$$\Gamma(Z,\CO_{^{cl}\!Z})\to \on{HH}^0(Z),$$
and of $\Gamma(Z,\CO_{^{cl}\!Z})$-modules
$$\Gamma\left(Z,H^1(T(Z))\right)\to \on{HH}^2(Z).$$

\medskip

It induces a homomorphism of graded algebras
\begin{equation}
\label{e:Sing and Hochschild}
\Gamma\left(\Sing(Z),\CO_{\Sing(Z)}\right)=\Gamma\left(Z,\on{Sym}_{\CO_{^{cl}\!Z}} \left(H^1(T(Z))\right)\right)\to\on{HH}^{\on{even}}(Z),
\end{equation}
where we assign to $\Gamma\left(Z,H^1(T(Z))\right)$ degree $2$. 

\sssec{}


We are now ready to give the main definitions of this paper:

\begin{defn} The \emph{singular support} of $\CF\in\IndCoh(Z)$, denoted $\on{SingSupp}(\CF)$, is 
$$\on{supp}_{\Gamma(\Sing(Z),\CO_{\Sing(Z)})}(\CF)\subset \Sing(Z).$$
\end{defn}

\begin{defn}
Let $Y$ be a conical Zariski-closed subset of $\Sing(Z)$. We let
$$\IndCoh_Y(Z)\subset \IndCoh(Z)$$
denote the full subcategory spanned by objects whose singular supports are contained 
in $Y$.
\end{defn}

\sssec{}

The following assertion (borrowed from  \cite[Theorem~11.3]{Kr})
gives an explicit expression for singular support: 

\begin{lem}  \label{l:sing supp via end}
For $\CF\in\IndCoh(Z)^c:=\Coh(Z)$, its singular support is equal to the support of the graded $\Gamma(\Sing(Z),\CO_{\Sing(Z)})$-module 
$\End^\bullet_{\Coh(Z)}(\CF)$.
\end{lem}

\begin{proof} Follows immediately from \lemref{l:supp via comp}(b).
\end{proof}

In addition, we have the following result:

\begin{thm} \label{t:finiteness}
For two objects $\CF_1,\CF_2\in \Coh(Z)$, the graded vector space $\Hom^\bullet_{\Coh(Z)}(\CF_1,\CF_2)$,
regarded as a module over $\Gamma(\Sing(Z),\CO_{\Sing(Z)})$, is finitely generated. 
\end{thm}

In particular, for $\CF\in \Coh(Z)$, the $\Gamma(\Sing(Z),\CO_{\Sing(Z)})$-module $\End^\bullet_{\Coh(Z)}(\CF)$ appearing
in \lemref{l:sing supp via end} is finitely generated.

\begin{rem} \label{r:finiteness}
If $Z$ is a classical scheme, the assertion of \thmref{t:finiteness} is due to Gulliksen \cite{Gu}; also see
references in the proof of \cite[Theorem~11.3]{Kr}. For completeness, we will present
a proof of \thmref{t:finiteness} in Appendix \ref{s:proof of finiteness}.
\end{rem} 

\ssec{Basic properties}

\sssec{}

First, we note:

\begin{lem}  \label{l:tensored over qcoh}
The subcategory $\IndCoh_Y(Z)\subset \IndCoh(Z)$ is stable under the monoidal action of 
$\QCoh(Z)$ on $\IndCoh(Z)$.
\end{lem}

\begin{proof}
For any module category $\bC$ over $\QCoh(Z)$, any full cocomplete subcategory $\bC'\subset \bC$
is stable under the action, since $\QCoh(Z)$ is generated by its unit object, $\CO_Z$.
\end{proof}

\sssec{}

It is easy to see that the dualizing sheaf $\omega_Z\in \IndCoh(Z)$ 
belongs to $\IndCoh_{\{0\}}(Z)$, where $\{0\}\subset \Sing(Z)$ denotes 
the zero-section. 

\medskip

Indeed, the construction of the isomorphism of 
\corref{c:HH as univ env} shows that the action of 
$\Gamma(Z,T(Z)[-1])\to \on{HC}(Z)$ on $\omega_Z\in \IndCoh(Z)$ is trivial.

\sssec{}

Recall the fully faithful functor 
$$\Xi_Z:\QCoh(Z)\to \IndCoh(Z)$$
(see \cite[Proposition 1.5.3]{IndCoh}), which is defined since $Z$ is eventually coconnective.

\medskip

Note that since $Z$ is
quasi-smooth, and hence Gorenstein (see \corref{c:Gorenstein abs}),
$\omega_Z$ is the image of a line bundle under $\Xi_Z$. Hence, by 
\lemref{l:tensored over qcoh}, the essential image of
all of $\QCoh(Z)$ under $\Xi_Z$ is contained in $\IndCoh_{\{0\}}(Z)$. 

\sssec{}

In \secref{ss:proof of zero sect} we will prove the converse inclusion:

\begin{thm}  \label{t:zero sect}
The subcategory $\IndCoh_{\{0\}}(Z)\subset \IndCoh(Z)$
coincides with the essential image of $\QCoh(Z)$ under $\Xi_Z$.
\end{thm}

\ssec{Compact generation}

\sssec{}

Lemma~\ref{l:compact generated} (or \cite[Theorem~6.4]{Kr}) immediately implies the following claim.

\begin{cor} \label{c:with support comp gen}
Let $Z$ be a quasi-smooth affine DG scheme. For any conical Zariski-closed subset $Y\subset\Sing(Z)$, the category 
$\IndCoh_Y(Z)$ is compactly generated. 
\end{cor}

Define 
$$\Coh_Y(Z):=\IndCoh_Y(Z)\cap \Coh(Z).$$

\corref{c:with support comp gen} can be rephrased as follows:

\begin{cor}
$\IndCoh_Y(Z)\simeq \Ind(\Coh_Y(Z))$. \qed
\end{cor}

\sssec{}  \label{sss:two Ys}

Let $Y_1\subset Y_2$ be two conical Zariski-closed subsets.

\medskip

We have a tautologically defined fully faithful functor
$$\Xi_Z^{Y_1,Y_2}:\IndCoh_{Y_1}(Z)\to \IndCoh_{Y_2}(Z)$$
that sends $\Coh_{Y_1}(Z)$ identically into $\Coh_{Y_2}(Z)$. 

\medskip

Since this functor sends compact objects to compacts, it admits a continuous
right adjoint. We denote this right adjoint by $\Psi_Z^{Y_1,Y_2}$. Thus,
the functor $\Psi_Z^{Y_1,Y_2}$ realizes $\IndCoh_{Y_1}(Z)$ as a colocalization 
of $\IndCoh_{Y_2}(Z)$.

\medskip

Note that the functor $\Xi_Z^{Y_1,Y_2}$ is (tautologically) compatible with the action of
$\QCoh(Z)$.  Hence the functor $\Psi_Z^{Y_1,Y_2}$ acquires a structure of being lax-compatible. 
However, we claim: 

\begin{lem} \label{l:commutes with action}
The functor $\Psi_Z^{Y_1,Y_2}$ is compatible with the action of
$\QCoh(Z)$, i.e., the lax compatibility is strict.
\end{lem}

\begin{proof}
Same as that of \lemref{l:tensored over qcoh}.
\end{proof} 

\sssec{}

In particular, we can take $Y_2$ to be all of $\Sing(Z)$, in which case 
$\IndCoh_{Y_2}(Z)$ is all of $\IndCoh(Z)$.

\medskip

We will denote the resulting pair of adjoint functors
$$\IndCoh_Y(Z)\rightleftarrows \IndCoh(Z)$$
by $(\Xi_Z^{Y,\on{all}},\Psi_Z^{Y,\on{all}})$. 

\sssec{}

Similarly, for any $Y$ that contains the zero-section, we obtain the corresponding pair of adjoint functors 
$$\Xi_Z^Y:\QCoh(Z)\rightleftarrows \IndCoh_Y(Z):\Psi_Z^Y$$
with $\Xi_Z^Y$ being fully faithful. 

\medskip

By definition, the functor $\Xi_Z^Y$ is the ind-extension of the natural embedding
$$\QCoh(Z)^c=\QCoh(Z)^{\on{perf}}\hookrightarrow \Coh_Y(Z)\hookrightarrow \IndCoh_Y(Z),$$
and $\Psi_Z^Y$ is the ind-extension of
$$\Coh_Y(Z)\hookrightarrow \Coh(Z)\hookrightarrow \QCoh(Z).$$

\medskip

Note that according to \thmref{t:zero sect}, the latter functors are a particular case
of $(\Xi_Z^{Y_1,Y_2},\Psi_Z^{Y_1,Y_2})$ for $Y_1=\{0\}$ and $Y_2=Y$.

\ssec{The t-structure}  \label{ss:t-structure}

Let $Y$ be a conical Zariski-closed subset of $\Sing(Z)$ containing the zero-section. 

\sssec{}

We define a t-structure on $\IndCoh_Y(Z)$ by declaring that
$$\CF\in \IndCoh_Y(Z)^{\leq 0} \, \Leftrightarrow \Psi^Y_Z(\CF)\in \QCoh(Z)^{\le 0}.$$

\medskip

Note that for $Y=\Sing(Z)$, this t-structure coincides with the canonical t-structure
on $\IndCoh(Z)$ of \cite[Sect. 1.2]{IndCoh}. 

\begin{lem}
$$\CF\in \IndCoh_Y(Z)^{\leq 0} \, \Leftrightarrow \Xi_Z^{Y,\on{all}}(\CF)\in \IndCoh(Z)^{\leq 0}.$$
\end{lem}

\begin{proof}
Follows from the fact that
$$\Psi^Y_Z\simeq \Psi^Y_Z(\CF)\circ \Psi^{Y,\on{all}}_Z\circ \Xi^{Y,\on{all}}_Z\simeq \Psi_Z\circ \Xi^{Y,\on{all}}_Z.$$
\end{proof}

\begin{cor}
The functors $\Psi_Z^{Y,\on{all}}$ and $\Psi_Z^Y$ are t-exact.
\end{cor}

\begin{proof}
The previous lemma implies that $\Xi_Z^{Y,\on{all}}$ is right t-exact. Hence, 
$\Psi_Z^{Y,\on{all}}$ is left t-exact by adjunction. The fact that $\Psi_Z^{Y,\on{all}}$ 
is right t-exact follows from the fact that
$$\Psi_Z^Y\circ \Psi_Z^{Y,\on{all}}\simeq \Psi_Z.$$

\medskip

The functor $\Psi_Z^Y$ is right t-exact by definition. To show that it is left t-exact,
it is enough to show that $\Xi_Z^Y$ is right t-exact. The latter is equivalent, by
definition, to the fact that $\Psi_Z^Y\circ \Xi_Z^Y$ is right t-exact. However, the latter
functor is isomorphic to the identity.
\end{proof}

\sssec{}

We will now prove:

\begin{prop}
The functors $\Psi_Z^{Y,\on{all}}$ and $\Psi_Z^Y$ induce equivalences
$$\IndCoh(Z)^{\geq 0}\to \IndCoh_Y(Z)^{\geq 0}\to \QCoh(Z)^{\geq 0}.$$
\end{prop}

\begin{proof}
Note that the fact that the composed functor $\Psi_Z^Y\circ \Psi_Z^{Y,\on{all}}\simeq \Psi_Z$
induces an equivalence 
$$\IndCoh(Z)^{\geq 0}\to \QCoh(Z)^{\geq 0}$$
is \cite[Proposition 1.2.4]{IndCoh}. In particular, the functor 
$$\Psi_Z^{Y,\on{all}}:\IndCoh(Z)^{\geq 0}\to \IndCoh_Y(Z)^{\geq 0}$$
is conservative. 

\medskip

The left adjoint to $\Psi_Z^{Y,\on{all}}|_{\IndCoh(Z)^{\geq 0}}$ is given by
$$\CF\mapsto \tau^{\geq 0}\left(\Xi_Z^{Y,\on{all}}(\CF)\right).$$

It is enough to show that this left adjoint is fully faithful, i.e., that the
adjunction map
$$\CF\to \Psi_Z^{Y,\on{all}}\left(\tau^{\geq 0}\left(\Xi_Z^{Y,\on{all}}(\CF)\right)\right)$$
is an isomorphism. 

\medskip

Since $\Psi_Z^{Y,\on{all}}$ is t-exact, we have:
$$\Psi_Z^{Y,\on{all}}\left(\tau^{\geq 0}\left(\Xi_Z^{Y,\on{all}}(\CF)\right)\right)\simeq
\tau^{\geq 0}\left(\Psi_Z^{Y,\on{all}}\circ \Xi_Z^{Y,\on{all}}(\CF)\right).$$

However, since $\Xi_Z^{Y,\on{all}}$ is fully faithful, the latter expression is isomorphic to
$\tau^{\geq 0}(\CF)\simeq \CF$, as required.

\end{proof}

\sssec{}

Recall that a t-structure on a triangulated category $\bT$ is said to be \emph{compactly generated} if
$$\CF\in \bT^{>0}\, \Leftrightarrow\, \Hom_{\bT}(\CF_1,\CF)=0,\, \forall\, \CF_1\in \bT^{\leq 0}\cap \bT^c.$$

\begin{prop}
The t-structure on $\IndCoh_Y(Z)$ is compactly generated.
\end{prop}

\begin{proof}

Let $\CF\in \IndCoh_Y(Z)$ be an object that is right-orthogonal to 
$$\Coh_Y(Z)\cap \Coh(Z)^{\leq 0}.$$
Let us prove that $\CF\in\IndCoh_Y(Z)^{>0}$. 
Truncating, we may assume that $\CF\in \IndCoh_Y(Z)^{\leq 0}$; we need to show that $\CF=0$. 

\medskip

By assumption, $\CF$ is right-orthogonal to the essential image of $\QCoh(Z)^{\on{perf}}\cap \Coh(Z)^{\leq 0}$
under $\Xi_Z^Y$. Hence, $\Psi^Y_Z(\CF)\in \QCoh(Z)^{>0}$. However, since $\CF\in \IndCoh_Y(Z)^{\leq 0}$, we have also that
$\Psi^Y_Z(\CF)\in \QCoh(Z)^{\leq 0}$, so $\Psi^Y_Z(\CF)=0$.

\medskip

Thus, $\CF$ is right-orthogonal to the essential image of all of $\QCoh(Z)^{\on{perf}}$ under $\Xi_Z^Y$. To prove that $\CF=0$,
we need to show that 
$$\Hom_{\IndCoh_Y(Z)}(\CF_1,\CF)=0$$
for any $\CF_1\in \Coh_Y(Z)$. However, for any such $\CF_1$, there exists $\CF_2\in \QCoh(Z)^{\on{perf}}$ and a map
$\CF_2\to \CF_1$, such that
$$\on{Cone}(\CF_2\to \CF_1)\in \Coh_Y(Z)\cap \Coh(Z)^{\leq 0}.$$

This implies the required assertion by the long exact sequence. 

\end{proof}

\ssec{Localization with respect to $Z$} \label{sss:localization}

\sssec{}

Let $V\overset{i}\hookrightarrow  Z$ be a closed DG subscheme. Let $\IndCoh(Z)_V$ be the corresponding
full subcategory of $\IndCoh(Z)$ (see \cite[Sect. 4.1.2]{IndCoh}), i.e., $\IndCoh(Z)_V$ consists of those objects
that vanish when restricted to $Z-V$.

\medskip

Equivalently, we can then
produce $\IndCoh(Z)_V$ in the framework of \secref{ss:support in DG}, using the action of
$\Gamma(Z,\CO_Z)$ on the category $\IndCoh(Z)$.

\medskip

Consider the scheme
$$\Sing(Z)_V={}^{cl}(\Sing(Z)\underset{Z}\times V)\subset \Sing(Z).$$

Consider the corresponding subcategory
$$\IndCoh_{\Sing(Z)_V}(Z)\subset \IndCoh(Z).$$

The next assertion results immediately from \lemref{l:change of algebras}:

\begin{cor} \label{c:supp on subscheme}
The subcategories
$\IndCoh_{\Sing(Z)_V}(Z)$ and $\IndCoh(Z)_V$ of $\IndCoh(Z)$ coincide.
\end{cor}

\sssec{}

Let $Y\subset \Sing(Z)$ be a conical Zariski-closed subset. Set
$$Y_V:={}^{cl}(Y\underset{Z}\times V).$$

From \corref{c:supp on subscheme} and \eqref{e:intersection of closed} we obtain:

\begin{cor} \label{c:on closed} 
The subcategories
$$\IndCoh_{Y_V}(Z)  \text{ and } \IndCoh(Z)_V\cap \IndCoh_Y(Z)$$
of $\IndCoh(Z)$ coincide. 
\end{cor}

\sssec{}

Let now $U\overset{j}\hookrightarrow Z$ be an open affine. By \cite[Lemma 4.1.1]{IndCoh}, 
we have a pair of adjoint functors 
$$j^{\IndCoh,*}:\IndCoh(Z)\rightleftarrows \IndCoh(U):j_*^{\IndCoh},$$
which realize $\IndCoh(U)$ as a localization of $\IndCoh(Z)$. By \cite[Corollary 4.4.5]{IndCoh},
we have a commutative diagram with vertical arrows being equivalences:
$$
\CD
\IndCoh(Z)\underset{\QCoh(Z)}\otimes \QCoh(U)  @>{\on{Id}\otimes j_*}>>  \IndCoh(Z)\underset{\QCoh(Z)}\otimes \QCoh(Z)  \\
@VVV     @VVV    \\
\IndCoh(U)   @>{j_*^{\IndCoh}}>>   \IndCoh(Z).
\endCD
$$

\medskip

In particular, we obtain that $\IndCoh(U)$ can be interpreted as the full subcategory of $\IndCoh(Z)$
corresponding to $^{cl}U\subset {}^{cl}\!Z$ with respect to the action of
$\Gamma(Z,\CO_{^{cl}\!Z})$ on $\IndCoh(Z)$ in the sense of \secref{sss:abs category on open}.

\medskip

Let $Y\subset \Sing(Z)$ be a conical Zariski-closed subset. Set
$$Y_U:={}^{cl}(Y\underset{Z}\times U)\subset \Sing(Z)_U\simeq \Sing(U).$$

From \eqref{e:localizations commute}, we obtain:

\begin{cor} \label{c:on open} \hfill

\smallskip

\noindent{\em(a)}
The functors $j_*^{\IndCoh}$ and $j^{\IndCoh,*}$ define an equivalence between $\IndCoh_{Y_U}(U)$
and the intersection of $\IndCoh_Y(Z)$ with the essential image of $\IndCoh(U)$
under $j_*^{\IndCoh}$. 

\smallskip

\noindent{\em(b)} The functors $(j^{\IndCoh,*},j_*^{\IndCoh})$ map the categories
$$\IndCoh_Y(Z)\rightleftarrows \IndCoh_{Y_U}(U)$$
to one another, and are mutually adjoint.

\smallskip

\noindent{\em(c)} We have a commutative diagram with vertical arrows being isomorphisms:
$$
\CD
\IndCoh_Y(Z)\underset{\QCoh(Z)}\otimes \QCoh(U)  @>{\on{Id}\otimes j_*}>>  \IndCoh_Y(Z)\underset{\QCoh(Z)}\otimes \QCoh(Z)  \\
@VVV     @VVV    \\
\IndCoh_{Y_U}(U)   @>{j_*^{\IndCoh}}>>   \IndCoh_Y(Z).
\endCD
$$

\end{cor}

\begin{cor} \label{c:on open +}
Let $U_i$ be a cover of $Z$ by open affine subsets. An object $\CF\in \IndCoh(Z)$ belongs to $\IndCoh_Y(Z)$
if and only if $\CF|_{U_i}$ belongs to $\IndCoh_{Y_{U_i}}(U_i)$ for every $i$.
\end{cor}

\begin{proof}
The ``only if" direction follows immediately from \corref{c:on open}(b). The ``if" direction follows
from the \v{C}ech complex.
\end{proof}

\sssec{}

For $U$ as above, let $V$ be a complementary
closed DG subscheme. By \cite[Corollary 4.1.5]{IndCoh}, we have a short exact sequence
of categories
\begin{equation} \label{e:Cousin sos IndCoh}
\IndCoh(Z)_V\rightleftarrows \IndCoh(Z)\overset{j^{\IndCoh,*}}\rightleftarrows \IndCoh(U).
\end{equation}
It can be obtained from the short exact sequence of categories
\begin{equation} \label{e:Cousin sos}
\QCoh(Z)_V\rightleftarrows \QCoh(Z)\overset{j^*}\rightleftarrows \QCoh(U),
\end{equation}
by the operation
$$\IndCoh(Z)\underset{\QCoh(Z)}\otimes-.$$
(Here $\QCoh(Z)_V\subset \QCoh(Z)$
is the full subcategory consisting of objects set-theoretically supported on $V$.)

\medskip

Hence, we obtain:

\begin{cor}\label{c:Cousin}
There exists a short exact sequence of DG categories
$$\IndCoh_{Y_V}(Z)\rightleftarrows \IndCoh_Y(Z)\rightleftarrows \IndCoh_{Y_U}(U),$$
which can be obtained from the short exact sequence \eqref{e:Cousin sos}
by the operation $$\IndCoh_Y(Z)\underset{\QCoh(Z)}\otimes-.$$
\end{cor}

\begin{cor} 
Let $Y_1\subset Y_2$ be two conical Zariski-closed subsets. Then the functors
$$\Xi_Z^{Y_1,Y_2}:\IndCoh_{Y_1}(Z)\rightleftarrows \IndCoh_{Y_2}(Z):\Psi_Z^{Y_1,Y_2}$$
induce (mutually adjoint) functors
$$\IndCoh_{(Y_1)_V}(Z)\rightleftarrows \IndCoh_{(Y_2)_V}(Z) \text{ and }
\IndCoh_{(Y_1)_U}(U)\rightleftarrows \IndCoh_{(Y_2)_U}(U).$$
\end{cor}

\sssec{} \label{sss:nonaffine}

From \corref{c:on open +}, we obtain that the notion of singular support of an object of $\IndCoh(Z)$ makes sense
for any quasi-smooth DG scheme $Z$ (not necessarily affine). 

\medskip

Namely, we choose an affine cover $U_\alpha$, and we set
$$\on{SingSupp}(\CF)\cap \Sing(Z)_{U_\alpha}:=\on{SingSupp}(\CF|_{U_\alpha}),$$
where we identify
$$\Sing(Z)_{U_\alpha}\simeq \Sing(U_\alpha),$$ and
$\CF|_{U_\alpha}:=j_\alpha^{\IndCoh,*}(\CF)$, where $j_\alpha:U_\alpha\hookrightarrow Z$.

\medskip

\corref{c:on open +} implies that $\on{SingSupp}(\CF)$ is well-defined, an in particular, independent
of the choice of the cover. 

\medskip

Furthermore, to $Y\subset \Sing(Z)$ we can attach a full subcategory
$$\IndCoh_Y(Z)\subset \IndCoh(Z),$$
by the requirement
$$\CF\in \IndCoh_Y(Z) \, \Leftrightarrow\, \CF|_{U_\alpha}\in \IndCoh_{Y_{U_\alpha}}(U_\alpha),\, \forall \alpha.$$




\ssec{Behavior with respect to products}

\sssec{}

Recall (see \cite[Proposition 4.6.2]{IndCoh}) that if $Z_i$, $i=1,2$ are two quasi-compact DG schemes, 
the natural functor
\begin{equation} \label{e:product}
\IndCoh(Z_1)\otimes \IndCoh(Z_2)\to \IndCoh(Z_1\times Z_2)
\end{equation}
is an equivalence.

\begin{rem}
This simple statement uses the assumption that $k$ is perfect (recall that we assume $\on{char}(k)=0$).
\end{rem}

\sssec{}

Assume now that $Z_i$ are both quasi-smooth. It is easy to see that we have a natural isomorphism
$$\Sing(Z_1)\times \Sing(Z_2)\simeq \Sing(Z_1\times Z_2).$$

Let $Y_i\subset \Sing(Z_i)$ be conical Zariski-closed subsets, and consider the corresponding
subset 
$$Y_1\times Y_2\subset \Sing(Z_1\times Z_2).$$

\begin{lem} \label{l:products}
We have:
$$\IndCoh_{Y_1}(Z_1)\otimes \IndCoh_{Y_2}(Z_2)=
\IndCoh_{Y_1\times Y_2}(Z_1\times Z_2)$$
as subcategories of $\IndCoh(Z_1\times Z_2)$.
\end{lem}

\begin{proof}

This follows immediately from \propref{p:abs product}.

\end{proof}

\ssec{Compatibility with Serre duality}

\sssec{}

Recall (see, e.g., \cite[Sect. 9.5]{IndCoh}) 
that the category $\Coh(Z)$ carries a canonical
anti-involution given by Serre duality, denoted $\BD^{\on{Serre}}_Z$.

\begin{prop}  \label{p:Serre}
For any conical Zariski-closed subset $Y\subset \Sing(Z)$, 
the anti-involution $\BD^{\on{Serre}}_Z$ preserves the subcategory $\Coh_Y(Z)\subset \Coh(Z)$.
\end{prop}

One proof is given in \secref{sss:proof of Serre abs}. Another (in a sense more hands-on, but logically equivalent)
proof is given in \secref{ss:proof of Serre}.

\medskip

\begin{cor} \label{c:Serre}
For $\CF\in \Coh(Z)$, we have
$$\on{SingSupp}(\CF)=\on{SingSupp}(\BD_Z(\CF)).$$
\end{cor}

\sssec{}

We obtain that there exists a canonical identification
$$(\IndCoh_Y(Z))^\vee\simeq \IndCoh_Y(Z),$$
obtained by extending $\BD^{\on{Serre}}_Z|_{\Coh_Y(Z)}$.

\medskip

Let $Y_1\subset Y_2$ be two conical Zariski-closed subsets, and 
consider the pair of adjoint functors
$$(\Psi_Z^{Y_1,Y_2})^\vee:\IndCoh_{Y_1}(Z)\rightleftarrows \IndCoh_{Y_2}(Z):(\Xi_Z^{Y_1,Y_2})^\vee,$$
obtained from
$$\Xi_Z^{Y_1,Y_2}:\IndCoh_{Y_1}(Z)\rightleftarrows \IndCoh_{Y_2}(Z):\Psi_Z^{Y_1,Y_2}$$
by passing to the dual functors.

\begin{lem}  \label{l:Serre dual of embedding}
We have canonical isomorphisms
$$(\Psi_Z^{Y_1,Y_2})^\vee\simeq\Xi_Z^{Y_1,Y_2} \text{ and }
(\Xi_Z^{Y_1,Y_2})^\vee \simeq \Psi_Z^{Y_1,Y_2}.$$
\end{lem}

\begin{proof}
Follows from \cite{DG}, Lemma 2.3.3 using the fact that
$$\Xi_Z^{Y_1,Y_2}\circ \BD^{\on{Serre}}_Z\simeq \BD^{\on{Serre}}_Z\circ \Xi_Z^{Y_1,Y_2}.$$
\end{proof}

\begin{rem}
Note that if $Y=\{0\}$ is the zero-section, 
the resulting self duality on $\IndCoh_{\{0\}}(Z)\simeq \QCoh(Z)$ is different
from the ``naive" self-duality: the two differ by the automorphism of $\QCoh(Z)$ given by
tensoring with $\omega_Z$.
\end{rem} 

\section{Singular support on a global complete intersection and Koszul duality}   \label{s:gci}

In this section we analyze the behavior of singular support on a DG scheme $Z$ which
is a ``global complete intersection." Recall that this means that $Z$ is presented as a Cartesian square
$$
\CD
Z  @>{\iota}>> \CU  \\ 
@VVV    @VVV \\
\on{pt}  @>>>  \CV,
\endCD
$$
where $\CU$ and $\CV$ are smooth affine schemes.
In this case, we will reinterpret the notion of singular support in terms of Koszul duality. 

\medskip

Our basic tool will be the group DG scheme $\CG_{\on{pt}/\CV}:=\on{pt}\underset{\CV}\times \on{pt}$. 

\ssec{Koszul duality}  \label{ss:KD for V}

\sssec{}

Consider the groupoid 
$$\CG_{\on{pt}/\CV}:=\on{pt}\underset{\CV}\times \on{pt}$$ over $\on{pt}$, that is, a group DG scheme:

\begin{gather}
\xy
 (-15,0)*+{\on{pt}}="X";
(15,0)*+{\on{pt}.}="Y";
(0,15)*+{\CG_{\on{pt}/\CV}}="Z";
{\ar@{->}_{p_1} "Z";"X"};
{\ar@{->}^{p_2} "Z";"Y"};
\endxy
\end{gather}

Let $\Delta_{\on{pt}}$ denote the diagonal map
$$\on{pt}\to \on{pt}\underset{\CV}\times \on{pt}=\CG_{\on{pt}/\CV}.$$ 
The object 
$$(\Delta_{\on{pt}})_*^{\IndCoh}(k)\in\IndCoh(\CG_{\on{pt}/\CV})$$ is the unit in the
monoidal category $\IndCoh(\CG_{\on{pt}/\CV})$. 

\sssec{}

By \secref{ss:E2 and grpds},  the above groupoid gives rise to an $\BE_2$-algebra 
$$\CA_{\CG_{\on{pt}/\CV}}=:\on{HC}(\on{pt}/\CV),$$
whose underlying associative DG algebra identifies canonically with
$$\CMaps_{\IndCoh(\CG_{\on{pt}/\CV})}\left((\Delta_{\on{pt}})_*^{\IndCoh}(k),(\Delta_{\on{pt}})_*^{\IndCoh}(k)\right).$$

\begin{rem}
The $\BE_2$-algebra that we denote here by $\CA_{\CG_{\on{pt}/\CV}}$
is what should be properly denoted $\CA^{\IndCoh}_{\CG_{\on{pt}/\CV}}$. This is done for
the purpose of unburdening the notation. As was explained in Remark \ref{r:E2 and grpds},
it is the IndCoh version (and not the QCoh one) that we use in this paper.
\end{rem}

\medskip

Let $V$ denote the tangent space to $\CV$ at $\on{pt}$. We have:

\begin{lem} \label{l:koszul as E_1}
The associative DG algebra underlying $\on{HC}^\bullet(\on{pt}/\CV)$
identifies canonically with
$\Sym(V[-2])$.
\end{lem}
\begin{proof} This is a special case of \corref{c:HC as E_1 for groups}, since a groupoid over $\on{pt}$ is canonically
a group DG scheme. 
\end{proof}

\medskip

In particular, we see that 
$$\on{HH}^\bullet(\on{pt}/\CV):=\underset{n}\bigoplus\, H^n\left(\on{HC}(\on{pt}/\CV)\right)$$
identifies canonically with $\on{Sym}(V)$ as a classical graded algebra, where the elements of $V$ have degree $2$. 
Geometrically,
$$\Spec\left(\on{HH}^\bullet(\on{pt}/\CV)\right)\simeq V^*.$$

\sssec{}  

Note that according to \cite[Proposition 4.1.7(b)]{IndCoh}, $(\Delta_{\on{pt}})_*^{\IndCoh}(k)$ is a compact generator
of $\IndCoh(\CG_{\on{pt}/\CV})$.  

\medskip

Hence, from \secref{sss:recovering monoidal} we obtain
a natural monoidal equivalence 
\begin{equation} \label{e:koszul over pt}
\on{HC}(\on{pt}/\CV)^{\on{op}}\mod\to \IndCoh(\CG_{\on{pt}/\CV}).
\end{equation}

We will denote the inverse functor $\IndCoh(\CG_{\on{pt}/\CV})\to \on{HC}(\on{pt}/\CV)^{\on{op}}\mod$
by $\on{KD}_{\on{pt}/\CV}$, and refer to is as the Koszul duality functor. Explicitly, 
\begin{equation}  \label{e:KD for V}
\on{KD}_{\on{pt}/\CV}=\CMaps_{\IndCoh(\CG_{\on{pt}/\CV})}((\Delta_{\on{pt}})_*^{\IndCoh}(k),-).
\end{equation}

The functor $\on{KD}_{\on{pt}/\CV}$ intertwines the forgetful functor $\on{HC}(\on{pt}/\CV)\mod\to \Vect$
and 
$$\Delta_{\on{pt}}^!:\IndCoh(\CG_{\on{pt}/\CV})\to \IndCoh(\on{pt})=\Vect.$$

\sssec{}

Recall the set-up of \secref{sss:E2 on itself} with the $\BE_2$-algebra
being $\on{HC}(\on{pt}/\CV)^{\on{op}}$. 

\medskip

In particular, to an object 
$$\CF\in \IndCoh(\CG_{\on{pt}/\CV})\simeq \on{HC}(\on{pt}/\CV)^{\on{op}}\mod$$ we can
associate its support, which is a conical Zariski-closed subset of $V^*$. Conversely, to a
conical Zariski-closed subset $Y\subset V^*$ we associate a full subcategory
$$\left(\IndCoh(\CG_{\on{pt}/\CV})\right)_{Y}\subset \IndCoh(\CG_{\on{pt}/\CV}).$$

\begin{lem} \label{l:KD supp}
The support of $\CF\in\IndCoh(\CG_{\on{pt}/\CV})$ is equal to the support of the
graded $\on{Sym}(V[-2])$-module 
\[\on{Hom}^\bullet((\Delta_{\on{pt}})_*^{\IndCoh}(k),\CF)=H^\bullet\left(\on{KD}_{\on{pt}/\CV}(\CF)\right).\]
\end{lem}
\begin{proof} Since $(\Delta_{\on{pt}})_*^{\IndCoh}(k)$ is a compact generator of $\IndCoh(\CG_{\on{pt}/\CV})$,
this follows from \lemref{l:supp via comp}.
\end{proof}

\begin{rem}\label{r:KD as enh!} 
The DG scheme $\CG_{\on{pt}/\CV}$ is quasi-smooth and $V^*=\Sing(\CG_{\on{pt}/\CV})$.
It is easy to see that the support of $\CF\in \IndCoh(\CG_{\on{pt}/\CV})$ is  
nothing but $\on{SingSupp}(\CF)\subset V^*$ (see \lemref{l:sing supp via kosz} for a more general statement).
\end{rem}

\sssec{}

Note that by combining Lemma~\ref{l:koszul as E_1} and \eqref{e:koszul over pt}, we obtain:

\begin{cor} \label{c:koszul as E_1}
The categories $\IndCoh(\CG_{\on{pt}/\CV})$ and $\Sym(V[-2])\mod$ are canonically equivalent as plain 
DG categories.
\end{cor}

\begin{rem}   \label{r:non-monoidal}
Both sides in \corref{c:koszul as E_1} are naturally monoidal categories. However, 
the equivalence of \corref{c:koszul as E_1} does not come with a monoidal structure (cf. 
Remark \ref{r:not E2}). We will see in \corref{c:parallelized E2} that a choice of a \emph{parallelization}
of $\CV$ at $\on{pt}$ upgrades the above equivalence to a monoidal one.
\end{rem}

\ssec{Functoriality of Koszul duality}   \label{ss:funct of KD}

The material in this section will be needed for the proof of some key properties
of singular support in \secref{s:funct}.

\medskip

However, as it will not be used for the discussion of singular support in the rest of this section,
the reader might choose to skip it on the first pass. 

\sssec{}

Let $f=f_\CV:\CV_1\to \CV_2$ be a map between smooth classical schemes. Fix a point $\on{pt}\to\CV_1$.
Set $V_i=T_{\on{pt}}(\CV)$ and $g=(df_{\on{pt}})^*:V_2^*\to V_1^*$. 
Now consider the morphism of DG group schemes
$$f_\CG:\CG_{\on{pt}/\CV_1}\to \CG_{\on{pt}/\CV_2}.$$
(Note that $g=\Sing(f_\CG)$.) 
The following lemma is obvious.

\begin{lem} The following three conditions are equivalent:

\smallskip
\noindent
{\em (a)} $g$ is injective;

\smallskip
\noindent
{\em (b)} $f$ is smooth at $\on{pt}$;

\smallskip
\noindent
{\em (c)} $f_\CG$ is quasi-smooth.
\smallskip

The following two conditions are also equivalent:

\smallskip
\noindent
{\em (a')} $g$ is surjective;

\smallskip
\noindent
{\em (b')} $f$ is unramified (and then a regular immersion) at $\on{pt}$.

\end{lem}

\sssec{} The map $f_\CG$ induces a monoidal functor
$$(f_\CG)_*^{\IndCoh}:\IndCoh(\CG_{\on{pt}/\CV_1})\to \IndCoh(\CG_{\on{pt}/\CV_2}),$$
and a homomorphism of $\BE_2$-algebras
$$f_{\on{HC}}:\on{HC}(\on{pt}/\CV_1)\to \on{HC}(\on{pt}/\CV_2).$$

\medskip

It is easy to see that the corresponding homomorphism of graded algebras 
$$\on{HH}^\bullet(\on{pt}/\CV_1)\to \on{HH}^\bullet(\on{pt}/\CV_2)$$
corresponds to the homomorphism
\[g^*:\Gamma(V_1^*,\CO_{V_1^*})\to\Gamma(V_2^*,\CO_{V_2^*})\]
under the isomorphism
\[\on{HH}^\bullet(\on{pt}/\CV_i)\simeq\Sym(V_i)\simeq\Gamma(V_i^*,\CO_{V_i^*})\qquad (i=1,2).\]

\medskip

In terms of the equivalence of \eqref{e:koszul over pt}, the functor $(f_\CG)_*^{\IndCoh}$
corresponds to the extension of scalars functor
$$\on{HC}(\on{pt}/\CV_1)^{\on{op}}\mod\to \on{HC}(\on{pt}/\CV_2)^{\on{op}}\mod.$$

\sssec{}

Consider now the right adjoint functor
$$(f_\CG)^!:\IndCoh(\CG_{\on{pt}/\CV_2})\to \IndCoh(\CG_{\on{pt}/\CV_1}),$$
which can also be thought of as the forgetful functor
\begin{equation} \label{e:! as forget}
\on{HC}(\on{pt}/\CV_2)^{\on{op}}\mod\to \on{HC}(\on{pt}/\CV_1)^{\on{op}}\mod
\end{equation}
along the homomorphism $f_{\on{HC}}$.

\medskip

We have:

\begin{prop} \label{p:support KD} \hfill

\smallskip

\noindent{\em(a)} Suppose the support of $\CF_2\in\IndCoh(\CG_{\on{pt}/\CV_2})$ equals $Y_2\subset V_2^*$.
Then the support of \[(f_\CG)^!(\CF_2)\in\IndCoh(\CG_{\on{pt}/\CV_1})\] equals $\overline{g(Y_2)}\subset V_1^*$.

\smallskip

\noindent{\em(b)} Suppose the support of $\CF_1\in\IndCoh(\CG_{\on{pt}/\CV_1})$ equals $Y_1\subset V_1^*$.
Then the support of \[(f_\CG)^{\IndCoh}_*(\CF_1)\in\IndCoh(\CG_{\on{pt}/\CV_2})\] is 
contained in $g^{-1}(Y_1)\subset V_2^*$. 

\smallskip

\noindent{\em(b')} The support in (b) equals all of $g^{-1}(Y_1)$ if $Y_1\subset g(V_2^*)$.

\smallskip

\noindent{\em(b'')} The support in (b) equals all of $g^{-1}(Y_1)$ if $\CF_1\in\Coh(\CG_{\on{pt}/\CV_1})$.

\end{prop}

\begin{rem}
As we will see in the proof, point (b'') of the proposition is the only non-tautological one, but it is not
essential for the main results of this paper.
\end{rem}

\begin{proof} 

By \lemref{l:koszul as E_1},
the statement reduces the corresponding assertion about modules over symmetric algebras.
Let $f:W_1\to W_2$ be a map of finite-dimensional vector spaces, and consider the corresponding
homomorphism
$$\Sym(W_1[-2])\to \Sym(W_2[-2]).$$

Let $M_2$ be an object of $\Sym(W_2[-2])\mod$, and let $M_1$ be its image under the forgetful functor
$$\Sym(W_2[-2])\mod\to \Sym(W_1[-2])\mod.$$

It is clear that $$\on{supp}_{\Sym(W_1)}(H^\bullet(M_1))\subset W_1^*$$ equals the closure of the image of
$$\on{supp}_{\Sym(W_2)}(H^\bullet(M_2))\subset W_2^*$$
under the map $g:W_2^*\to W_1^*$. This proves point (a) of the proposition.

\medskip

Let now $M_1$ be an object of $\Sym(W_1[-2])\mod$ and set
$$M_2:=\Sym(W_2[-2])\underset{\Sym(W_1[-2])}\otimes M_1.$$

First, it is clear that
\begin{equation} \label{e:supp cont}
\on{supp}_{\Sym(W_2)}(H^\bullet(M_2))\subset g^{-1}\left(\on{supp}_{\Sym(W_1)}(H^\bullet(M_1))\right).
\end{equation}

\medskip

It remains to show that the above containment is an equality if either $$\on{supp}_{\Sym(W_1)}(H^\bullet(M_1))\subset g(W_2^*)$$
or if $M_1$ is perfect. 

\medskip

Both assertions are easy when $g$ is surjective. Moreover, if they hold for two composable maps
$f'$ and $f''$, then they hold for their composition. This allows us to reduce the statement to the case when $g$ is an embedding
of codimension one. Thus, let us assume that $W_2$ is the cokernel of $t:k\to W_1$.

\medskip

Denote $N_1:=H^\bullet(M_1)$, viewed as an object in $\Sym(W_1)\mod^{\heartsuit}$. Denote
$$N_2:=\Sym(W_2)\underset{\Sym(W_1)}\otimes N_1\in \Sym(W_2)\mod.$$
The object $N_2$ has cohomologies in two degrees: 
$$N'_2:=H^{0}(N_2)=\on{coker}(t:N_1\to N_1) \text{ and }
N''_2:=H^{-1}(N_2)=\on{ker}(t:N_1\to N_1).$$

\medskip

We claim that
$$\on{supp}_{\Sym(W_2)}(H^\bullet(M_2))=\on{supp}_{\Sym(W_2)}(N'_2)\cup \on{supp}_{\Sym(W_2)}(N''_2).$$
This follows from the short exact sequence in $\Sym(W_2)\mod^{\heartsuit}$
$$0\to N'_2\to H^\bullet(M_2)\to N''_2\to 0.$$

\medskip

Now, if 
$$\on{supp}_{\Sym(W_1)}(H^\bullet(M_1))=\on{supp}_{\Sym(W_1)}(N_1)\subset W_2^*,$$ 
we have $$\on{supp}_{\Sym(W_1)}(N_1)=\on{supp}_{\Sym(W_2)}(N'_2)\cup \on{supp}_{\Sym(W_2)}(N''_2),$$
which implies the desired equality in \eqref{e:supp cont} in this case. 

\medskip

If $M_1$ is perfect, it is easy to see that the module $N_1$
is finitely generated (this is in fact a particular case of \thmref{t:finiteness}). Then, by the Nakayama Lemma,
$$\on{supp}_{\Sym(W_2)}(N'_2)=\left(\on{supp}_{\Sym(W_1)}(H^\bullet(M_1))\right)\cap W_2^*,$$
which again implies the equality in \eqref{e:supp cont}. 
 
\end{proof}  

\begin{cor} \label{c:cons KD} \hfill

\smallskip

\noindent{\em(a)} $(f_\CG)^!:\IndCoh(\CG_{\on{pt}/\CV_2})\to\IndCoh(\CG_{\on{pt}/\CV_1})$ is conservative.

\smallskip

\noindent{\em(b)} Set $Y_{1,can}=g(V_2^*)\subset V_1^*$. Then the restriction of $(f_\CG)^{\IndCoh}_*$ to $\left(\IndCoh(\CG_{\on{pt}/\CV_1})\right)_{Y_{1,can}}$ is conservative.
\end{cor}

\begin{proof} Follows immediately from \propref{p:support KD}. (It is easy also to give a direct proof.)
\end{proof}

\begin{cor} \label{c:support KD} \hfill

\smallskip

\noindent{\em(a)} 
For a conical Zariski-closed subset $Y_1\subset V^*_1$, set $Y_2:=g^{-1}(Y_1)\subset V_2^*$.
Then the essential image of $\left(\IndCoh(\CG_{\on{pt}/\CV_1})\right)_{Y_1}$
under the functor $(f_\CG)^\IndCoh_*$ is contained in the subcategory
$$\left(\IndCoh(\CG_{\on{pt}/\CV_2})\right)_{Y_2}\subset \IndCoh(\CG_{\on{pt}/\CV_2})$$
and generates it.

\smallskip

\noindent{\em(b)} Suppose $f$ is a smooth morphism at $\on{pt}$, so that $g$ is injective. 
For a conical Zariski-closed subset 
$Y_2\subset V^*_2$, set $Y_1:=g(Y_2)\subset V^*_1$. Then 
the essential image of $\left(\IndCoh(\CG_{\on{pt}/\CV_2})\right)_{Y_2}$
under the functor $(f_\CG)^!$ is contained in the subcategory
$$\left(\IndCoh(\CG_{\on{pt}/\CV_1})\right)_{Y_1}\subset \IndCoh(\CG_{\on{pt}/\CV_1})$$
and generates it.

\smallskip
\end{cor}

\begin{proof} For part (a), note that we have a pair of adjoint functors
\[(f_\CG)^\IndCoh_*:\IndCoh(\CG_{\on{pt}/\CV_1})\rightleftarrows\IndCoh(\CG_{\on{pt}/\CV_2}):(f_\CG)^!.\]
By \propref{p:support KD}, they restrict to a pair of adjoint functors  
\[(f_\CG)^\IndCoh_*:\IndCoh(\CG_{\on{pt}/\CV_1})_{Y_1}\rightleftarrows\IndCoh(\CG_{\on{pt}/\CV_2})_{Y_2}:(f_\CG)^!.\]
Since $f_\CG^!$ is conservative by \corref{c:cons KD}(a), the claim follows. (Note that the categories 
involved are compactly generated.)

\medskip

For part (b), note that $f_\CG$ is quasi-smooth, therefore, we have a pair of adjoint functors
\[(f_\CG)^{\IndCoh,*}:\IndCoh(\CG_{\on{pt}/\CV_2})\rightleftarrows\IndCoh(\CG_{\on{pt}/\CV_1}):(f_\CG)^{\IndCoh}_*.\]
Moreover, $f_\CG^{\IndCoh,*}$ can be obtained from $f_\CG^!$ by twisting by a cohomologically shifted line bundle
(by \corref{c:Gorenstein}).
By \propref{p:support KD}, restriction produces a pair of adjoint functors 
\[(f_\CG)^{\IndCoh,*}:\IndCoh(\CG_{\on{pt}/\CV_2})_{Y_2}\rightleftarrows\IndCoh(\CG_{\on{pt}/\CV_1})_{Y_1}:(f_\CG)^{\IndCoh}_*.\]
Note that $Y_1\subset g(V_2^*)=Y_{1,can}$, so by \corref{c:cons KD}(b'), $(f_\CG)^{\IndCoh}_*$ 
is conservative on $\IndCoh(\CG_{\on{pt}/\CV_1})_{Y_1}$.
The claim follows.  
\end{proof}

\sssec{}

The particular usefulness of \propref{p:support KD} and Corollaries \ref{c:cons KD} and \ref{c:support KD}
for us is explained by the following observation:

\begin{lem}  \label{l:shape of point} \hfill

\smallskip

\noindent{\em(a)} Let $Z$ be a quasi-smooth DG scheme such that $^{cl}\!Z\simeq \on{pt}$. Then
$Z$ is (non-canonically) isomorphic to $\on{pt}\underset{\CV}\times \on{pt}$ for some
smooth classical scheme $\CV$.

\smallskip

\noindent{\em(b)} If $Z_i=\on{pt}\underset{V_i}\times \on{pt}$, $i=1,2$ where $V_i$ are vector spaces,
then any map $Z_1\to Z_2$ can be realized as coming from a linear map $V_1\to V_2$.

\end{lem}

\begin{proof}

We have 
$$\CO_Z\simeq \on{C}^\cdot(L),$$
where $L$ is the DG Lie algebra $T_z(Z)[-1]$, where $z$ is the unique $k$-point of $Z$. By assumption,
$L$ has only cohomology in degree $+2$. Hence,
$$H^\bullet(\CO_Z)\simeq \Sym(V[1]),$$
where $V$ is the vector space dual to $H^2(L)$. This implies that $\CO_Z$ is itself non-canonically
isomorphic to $\Sym(V[1])$. This establishes point (a). 

\medskip

Point (b) follows from the fact that the space of maps of commutative DG algebras 
$$\Sym(V_1[1])\to \Sym(V_2[1])$$
is isomorphic to
$$\Maps(V_1[1],\Sym(V_2[1])).$$
In particular, the set of homotopy classes of such maps is in bijection with $\Hom(V_1,V_2)$.

\end{proof}

\ssec{Singular support via Koszul duality}   \label{ss:singular support via smooth}

In this subsection we let $Z$ be a quasi-smooth DG scheme,  
presented as a fiber product in the category of DG schemes 
\begin{equation}  \label{e:gci for Koszul}
\CD
Z  @>{\iota}>> \CU  \\ 
@VVV    @VVV \\
\on{pt}  @>>>  \CV,
\endCD
\end{equation}
with smooth $\CU$ and $\CV$, as in \secref{sss:expl presentations}. We will also assume that $\CU$ and $\CV$ are affine.

\sssec{}   \label{sss:action on Z}

Note that we have a Cartesian diagram 
\begin{gather} \label{e:groupoids morphism} 
\xy
(-25,0)*+{Z}="X";
(25,0)*+{Z}="Y";
(0,15)*+{Z\underset{\CU}\times Z}="Z";
(-25,-20)*+{\on{pt}}="X_1";
(25,-20)*+{\on{pt}.}="Y_1";
(0,-5)*+{\on{pt}\underset{\CV}\times \on{pt}}="Z_1";
{\ar@{->}"Z";"X"};
{\ar@{->}"Z";"Y"};
{\ar@{->}"Z_1";"X_1"};
{\ar@{->}"Z_1";"Y_1"};
{\ar@{->}"Z";"Z_1"};
{\ar@{->}"X";"X_1"};
{\ar@{->}"Y";"Y_1"};
\endxy
\end{gather}

In particular, we obtain that the group DG scheme $\CG_{\on{pt}/\CV}$ canonically acts on $Z$,
preserving its map to $\CU$. 

\medskip

In other words, we have a canonical isomorphism of groupoids
$$\CG_{Z/\CU}\simeq \CG_{\on{pt}/\CV}\times Z$$
acting on $Z$, commuting with the map to $\CU$:

$$
\xy
(-30,0)*+{Z}="X";
(30,0)*+{Z}="Y";
(0,20)*+{\CG_{Z/\CU}}="Z";
(-30,-20)*+{Z}="X_1";
(30,-20)*+{Z.}="Y_1";
(0,0)*+{\CG_{\on{pt}/\CV}\times Z}="Z_1";
{\ar@{->}"Z";"X"};
{\ar@{->}"Z";"Y"};
{\ar@{->}^{\on{pr}}"Z_1";"X_1"};
{\ar@{->}_{\on{act}_{\CG_{\on{pt}/\CV},Z}}"Z_1";"Y_1"};
{\ar@{->}_{\sim}"Z";"Z_1"};
{\ar@{->}_{\on{id}}"X";"X_1"};
{\ar@{->}^{\on{id}}"Y";"Y_1"};
\endxy
$$
Here $\on{act}_{\CG_{\on{pt}/\CV},Z}:\CG_{\on{pt}/\CV}\times Z\to Z$ is the action morphism.

\sssec{}

In particular, we obtain a canonical homomorphism of monoidal categories 
\begin{equation} \label{e:tensor product acts}
\IndCoh(\CG_{\on{pt}/\CV})\otimes \QCoh(\CU)\to \IndCoh(\CG_{Z/\CU}),
\end{equation}
and hence a homomorphism of $\BE_2$-algebras

\begin{equation} \label{e:tensor product acts E2}
\CA:=\on{HC}(\on{pt}/\CV) \otimes \Gamma(\CU,\CO_{\CU})\to \on{HC}(Z/U)\to  \on{HC}(Z),
\end{equation} 
where $\on{HC}(Z/U):=\CA_{\CG_{Z/U}}$, see \secref{sss:rel hh}.

\medskip

Note that
\[A=\underset{k}\bigoplus\, H^{2k}(\CA)=\on{HH}^{\on{even}}(\on{pt}/\CV)\otimes \Gamma(\CU,\CO_{\CU})=
\on{Sym}(V)\otimes \Gamma(\CU,\CO_{\CU}).\]

\sssec{}  Thus, we find ourselves in the paradigm of \secref{sss:support via E2} with the DG category in question being
$\IndCoh(Z)$.

\medskip

In particular, to an object
$\CF\in\IndCoh(Z)$, we can assign a conical Zariski-closed subset 
\[\supp_A(\CF)\subset\Spec(A)=V^*\times \CU.\]
The following lemma shows that this recovers the singular support of $\CF$.

\begin{lem}   \label{l:sing supp via kosz}
For any $\CF$, the support $\supp_A(\CF)\subset V^*\times \CU$ is the image of
$$\on{SingSupp}(\CF)\subset\Sing(Z)$$ under the embedding 
$\Sing(Z)\hookrightarrow V^*\times Z\hookrightarrow V^*\times \CU$, 
where the first map is given by \eqref{e:expl presentation}.
\end{lem}
\begin{proof}
It suffices to verify that the diagram 
\begin{gather*}  
\xy
(-20,0)*+{\on{Sym}(V)\otimes\Gamma(\CU,\CO_\CU)}="Z";
(20,0)*+{\Gamma(\Sing(Z),\CO_{\Sing(Z)})}="W";
(0,-15)*+{\on{HH}^{\on{even}}(Z)}="T";
{\ar@{->}"Z";"W"};
{\ar@{->}"Z";"T"};
{\ar@{->}"W";"T"};
\endxy
\end{gather*}
commutes.
This is straightforward.
\end{proof}

From \propref{p:support via monoidal}, we obtain:

\begin{cor} \label{c:sing supp via ten kosz}
For a conical Zariski-closed subset $Y\subset \Sing(Z)$, we have:
$$\IndCoh_Y(Z)\simeq \IndCoh(Z)\underset{\on{HC}(\on{pt}/\CV)\mod\otimes \QCoh(\CU)} \otimes 
\left(\on{HC}(\on{pt}/\CV)\mod\otimes \QCoh(\CU)\right)_Y.$$
\end{cor}

\ssec{Koszul duality in the parallelized situation}  \label{ss:parallel}

In this subsection we will assume that the formal completion of $\CV$ at $\on{pt}$ has been 
parallelized, i.e., that it is identified with the formal completion of $V$ at $0$. 

\sssec{}

In this case we have:

\begin{lem}\label{l:parallel}
The $\BE_2$-algebra structure $\on{HC}^\bullet(\on{pt}/\CV)$ is canonically commutative
\emph{(}i.e., comes by restriction from a canonically defined $\BE_\infty$-structure\emph{)}, and, as such,
identifies with $\Sym(V[-2])$.
\end{lem}

\begin{proof}
The $\BG_m$-action on $V$ by dilations gives rise to a $\BG_m$-action on the $\BE_2$-algebra 
$\on{HC}(\on{pt}/\CV)$. The computation of the cohomology of $\on{HC}(\on{pt}/\CV)$
puts us in the framework of \secref{sss:good grading}.
\end{proof}

\begin{cor}  \label{c:parallelized E2}
A parallelization of $\CV$ upgrades the monoidal structure on the category $\IndCoh(\CG_{\on{pt}/\CV})\simeq 
\on{HC}(\on{pt}/\CV)^{\on{op}}\mod$ 
to a symmetric monoidal structure, and as such
it is canonically equivalent to $\Sym(V[-2])\mod$.
\end{cor}

\sssec{} 

Thus, we see that in the parallelized situation, we can use \corref{c:support by tensor} to study
support in the category $\IndCoh(\CG_{\on{pt}/\CV})$ as follows. 

\medskip

The category 
$$\IndCoh(\CG_{\on{pt}/\CV})\simeq  \on{HC}(\on{pt}/\CV)^{\on{op}}\mod\simeq \Sym(V[-2])\mod$$
is naturally a module category over $\QCoh(V^*/\BG_m)$. For a conical Zariski-closed subset $Y\subset V^*$
we have:

\begin{cor}
$$\IndCoh(\CG_{\on{pt}/\CV})_Y=\IndCoh(\CG_{\on{pt}/\CV})\underset{\QCoh(V^*/\BG_m)}\otimes \QCoh(V^*/\BG_m)_{Y/\BG_m}.$$
\end{cor}

\sssec{}\label{sss:global intersection parallelized}

Let $Z$ be as in \secref{ss:singular support via smooth}. By \secref{sss:good grading}, the category 
$\IndCoh(Z)$ carries an action of the monoidal category $\QCoh(\CS_\CA)$ for the
stack 
\[\CS_\CA=\Spec(\on{Sym}(V)\otimes \Gamma(\CU,\CO_\CU))/\BG_m=V^*/\BG_m\times \CU.\]
This allows us to study the singular support of objects in $\IndCoh(Z)$ using 
Corollaries~\ref{c:support by tensor} and \ref{c:support by category fibers}. In particular, we have:

\begin{cor} \label{c:sing supp via ten kosz parallelized}
For a conical Zariski-closed subset $Y\subset \Sing(Z)$, we have:
$$\IndCoh_Y(Z)\simeq \IndCoh(Z)\underset{\QCoh(V^*/\BG_m\times \CU)}\otimes 
\QCoh\left(V^*/\BG_m\times \CU\right)_{Y/\BG_m}.$$
\end{cor}

\medskip

In the rest of this section we will use the view point on singular support via Koszul duality and prove some
results stated earlier in the paper.

\ssec{Proof of \propref{p:Serre}}  \label{ss:proof of Serre}

Recall that \propref{p:Serre} says that for $\CF\in \Coh(Z)$, its singular support equals that of $\BD_Z^{\on{Serre}}(\CF)$. 

\sssec{}

With no restriction of generality, we can assume that $Z$ is affine and fits into the Cartesian square \eqref{e:gci for Koszul}.
By \corref{l:sing supp via kosz}, it suffices to show that the supports of $\CF$ and $\BD_Z^{\on{Serre}}(\CF)$ for the action of
the graded algebra $A:=\on{Sym}(V)\otimes \Gamma(\CU,\CO_{\CU})$ are equal, where the grading on $\on{Sym}(V)$
is such that $\deg(V)=2$. 

\sssec{}

By \lemref{l:supp via comp}, is suffices to show that the supports of the $A$-modules 
$\Hom^\bullet_{\Coh(Z)}(\CF,\CF)$ and $\Hom^\bullet_{\Coh(Z)}(\BD_Z^{\on{Serre}}(\CF),\BD_Z^{\on{Serre}}(\CF))$
are equal. This follows from the next assertion:

\begin{lem}  \label{l:action of vector fields}
The diagram 
$$
\CD
\on{Sym}(V)\otimes \Gamma(\CU,\CO_{\CU})   @>>>   \Hom^\bullet_{\Coh(Z)}(\CF,\CF)  \\
@VVV   @VVV    \\
\on{Sym}(V)\otimes \Gamma(\CU,\CO_{\CU})   @>>>  \Hom^\bullet_{\Coh(Z)}(\BD_Z^{\on{Serre}}(\CF),\BD_Z^{\on{Serre}}(\CF))
\endCD
$$
is commutative, where the right vertical arrow is the isomorphism given by the anti-equivalence 
$$\BD_Z^{\on{Serre}}:\Coh(Z)^{\on{op}}\to \Coh(Z),$$
and the left vertical arrow is the automorphism that acts as identity on $\Gamma(\CU,\CO_{\CU})$ and
as $-1$ on $V\subset \on{Sym}(V)$.
\end{lem}

\sssec{}

Since the action map $\CG_{\on{pt}/\CV}\times Z\to Z$ is proper, and hence commutes with Serre duality, the 
assertion of \lemref{l:action of vector fields} follows from the next one: 

\begin{lem}
The equivalence $\BD^{\on{Serre}}_{\CG_{\on{pt}/\CV}}:\Coh(\CG_{\on{pt}/\CV})^{\on{op}}\to \Coh(\CG_{\on{pt}/\CV})$
induces the automorphism of 
$$\Hom^\bullet_{\CG_{\on{pt}/\CV}}\left((\Delta_{\on{pt}})_*^{\IndCoh}(k),(\Delta_{\on{pt}})_*^{\IndCoh}(k)\right)\simeq \Sym(V[-2]),$$
given by $v\mapsto -v:V\to V$.
\end{lem} 

\ssec{Constructing objects with a given singular support}   \label{ss:koszul functors}

The material in this subsection is not strictly speaking necessary for the rest of the paper.

\medskip

Let $Z$ fit into a Cartesian square as in \eqref{e:gci for Koszul}. In this subsection we will give 
an explicit procedure for producing compact objects in $\IndCoh_Y(Z)$.

\sssec{}    \label{sss:Hecke group}    \label{sss:F&G}

Consider the action of the group DG scheme $\CG_{\on{pt}/\CV}$ on $Z$ as in 
\secref{sss:action on Z}, and recall that $\on{act}_{\CG_{\on{pt}/\CV},Z}$ denotes the corresponding 
action map
$$\CG_{\on{pt}/\CV}\times Z\to Z.$$

\medskip

Clearly, the map $\on{act}_{\CG_{\on{pt}/\CV},Z}$  is proper.
This gives rise to a pair of adjoint functors:
\begin{equation}  \label{e:*!}
(\on{act}_{\CG_{\on{pt}/\CV},Z})_*^{\IndCoh}:\IndCoh(\CG_{\on{pt}/\CV}\times Z)
\rightleftarrows \IndCoh(Z):(\on{act}_{\CG_{\on{pt}/\CV},Z})^!.
\end{equation}

We will also use the notation 
$$\sF:=(\on{act}_{\CG_{\on{pt}/\CV},Z})_*^{\IndCoh}  \text{ and } \sG:=(\on{act}_{\CG_{\on{pt}/\CV},Z})^!,$$
and identify
$$\IndCoh(\CG_{\on{pt}/\CV}\times Z)\simeq \IndCoh(\CG_{\on{pt}/\CV})\otimes \IndCoh(Z).$$

\medskip

The functor $\sG$ is conservative because it admits a retract (given by the pullback
along the unit of $\CG_{\on{pt}/\CV}$). 

\medskip

Therefore, the essential image of $\sF$ generates $\IndCoh(Z)$.

\sssec{} 

Let us regard $\IndCoh(\CG_{\on{pt}/\CV})\otimes \IndCoh(Z)$ as a module category over $\QCoh(\CU)$
via the second factor. 

\medskip

In addition, we can regard $\IndCoh(\CG_{\on{pt}/\CV})\otimes \IndCoh(Z)$ as acted on by 
$\IndCoh(\CG_{\on{pt}/\CV})$ via the first factor. 

\medskip

Combining, we obtain that the monoidal category
$$\IndCoh(\CG_{\on{pt}/\CV})\otimes \QCoh(\CU)$$ 
acts on $\IndCoh(\CG_{\on{pt}/\CV})\otimes \IndCoh(Z)$. 

\medskip

In particular, we obtain a map of $\BE_2$-algebras
\[\CA=\on{HC}(\on{pt}/\CV) \otimes \Gamma(\CU,\CO_{\CU})\to \on{HC}(\IndCoh(\CG_{\on{pt}/\CV})\otimes \IndCoh(Z)).\]
Thus, by \secref{sss:support via E2}, to an object 
$$\CF\in \IndCoh(\CG_{\on{pt}/\CV})\otimes \IndCoh(Z)$$
we can assign its support
$$\on{supp}_A(\CF)\subset \Spec(A)\simeq V^*\times \CU.$$

Conversely, a conical Zariski-closed subset $Y\subset V^*\times \CU$ yields a full subcategory
\[\left(\IndCoh(\CG_{\on{pt}/\CV})\otimes \IndCoh(Z)\right)_Y=\{\CF\in\IndCoh(\CG_{\on{pt}/\CV}\times Z),\,\,\supp_A(\CF)\subset Y\}.\]

\sssec{Warning} The above notation may seem abusive, 
because $\CG_{\on{pt}/\CV}\times Z$ is itself a quasi-smooth affine DG scheme, which
has its own notion of singular support. Note that 
\[\Sing(\CG_{\on{pt}/\CV}\times Z)\simeq \Sing(\CG_{\on{pt}/\CV})\times \Sing(Z)\simeq V^*\times \Sing(Z)
\subset V^*\times V^*\times \CU.\] Thus,
for 
$$\CF\in\IndCoh(\CG_{\on{pt}/\CV}\times Z)\simeq \IndCoh(\CG_{\on{pt}/\CV})\otimes \IndCoh(Z),$$ 
its singular support is a conical subset 
$$\on{SingSupp}(\CF)\subset V^*\times V^*\times \CU.$$ 
It follows from \lemref{l:change of algebras} that $\supp_A(\CF)$ is the closure of the 
projection $p_{13}(\on{SingSupp}(\CF))$.
To avoid confusion, we never consider this singular support for objects of $$\IndCoh(\CG_{\on{pt}/\CV})\otimes \IndCoh(Z)$$ and only
deal with the ``coarse" support contained in $V^*\times \CU$. 

\sssec{}

It is clear that the functor $\sF$ is compatible with the action of the monoidal category 
$\IndCoh(\CG_{\on{pt}/\CV})\otimes \QCoh(\CU)$ on $\IndCoh(\CG_{\on{pt}/\CV})\otimes \IndCoh(Z)$
given above, 
and the action of $\IndCoh(\CG_{\on{pt}/\CV})\otimes \QCoh(\CU)$
on $\IndCoh(Z)$ given by \eqref{e:tensor product acts}. 

\medskip

Hence, the functor $\sG$, being the
right adjoint of $\sF$, is lax-compatible with the above actions.

\begin{lem}
The lax compatibility of $\sG$ with the actions of $\IndCoh(\CG_{\on{pt}/\CV})\otimes \QCoh(\CU)$ on
$\IndCoh(Z)$ and $\IndCoh(\CG_{\on{pt}/\CV})\otimes \IndCoh(Z)$ is strict.
\end{lem}

\begin{proof}
It is easy to see that the monoidal category $\IndCoh(\CG_{\on{pt}/\CV})\otimes \QCoh(\CU)$ is rigid 
(see \cite{DG}, Sect. 6, where the notion of rigidity is discussed). Now, the assertion follows
from the fact that if $\sF:\bC_1\to \bC_2$ is a functor between module categories over a rigid
monoidal category, which admits a continuous right adjoint as a functor between plain DG categories,
then the lax compatibility structure on this right adjoint is automatically strict.
\end{proof}

From \secref{sss:change of category}, we obtain:

\begin{cor} \label{c:Hecke action abstract}
Let $Y\subset V^*\times \CU$ be a conical Zariski-closed subset. Then the functors
$\sF$ and $\sG$ restrict to an adjoint pair of functors
$$\left(\IndCoh(\CG_{\on{pt}/\CV})\otimes \IndCoh(Z)\right)_Y \rightleftarrows \IndCoh_{Y\cap \Sing(Z)}(Z).$$
Moreover, the diagram
$$
\xymatrix{
\IndCoh(\CG_{\on{pt}/\CV})\otimes \IndCoh(Z) \ar[r] \ar[d] & \IndCoh(Z) \ar@<.7ex>[l] \ar[d]^{\Psi_Z^{Y,\on{all}}} \\
\left(\IndCoh(\CG_{\on{pt}/\CV})\otimes \IndCoh(Z)\right)_Y  \ar[r] & \IndCoh_{Y\cap \Sing(Z)}(Z). \ar@<.7ex>[l] }$$
commutes as well,
where the left vertical arrow is the right adjoint to the inclusion
$$\left(\IndCoh(\CG_{\on{pt}/\CV})\otimes \IndCoh(Z)\right)_Y\hookrightarrow \IndCoh(\CG_{\on{pt}/\CV})\otimes \IndCoh(Z).$$
\end{cor}

\begin{cor}\label{c:Koszul generation} Suppose $Y$ is a conical Zariski-closed subset of $\Sing(Z)\subset V^*\times \CU$. 

\smallskip

\noindent{\em (a)} For any $\CF\in\IndCoh(Z)$, we have:
\[\CF\in \IndCoh_Y(Z)\,\Leftrightarrow\, \sG(\CF)\in\left(\IndCoh(\CG_{\on{pt}/\CV})\otimes \IndCoh(Z)\right)_{Y}.\]

\smallskip

\noindent{\em (b)} The essential image under $\sF$ of the category
$\left(\IndCoh(\CG_{\on{pt}/\CV})\otimes \IndCoh(Z)\right)_{Y}$ generates $\IndCoh_Y(Z)$.
\end{cor}
\begin{proof} Both claims formally follow from the conservativeness of $\sG$. Indeed, $\CF\in\IndCoh_Y(Z)$ if and only if
the natural morphism $\Psi_Z^{Y,\on{all}}(\CF)\to\CF$ is an isomorphism. Since $\sG$ is conservative, this happens if and only if
the morphism $\sG(\Psi_Z^{Y,\on{all}}(\CF))\to\sG(\CF)$ is an isomorphism; by \corref{c:Hecke action abstract}, this is equivalent to
$\sG(\CF)\in\left(\IndCoh(\CG_{\on{pt}/\CV})\otimes \IndCoh(Z)\right)_{Y}$. 
We have therefore proved part (a). Now note that
the restriction 
\[\sF:\left(\IndCoh(\CG_{\on{pt}/\CV})\otimes \IndCoh(Z)\right)_{Y}\to \IndCoh_{Y\cap \Sing(Z)}(Z)\] is left adjoint to a conservative
functor; this proves part (b).
\end{proof}

\sssec{} \label{sss:interp Koszul}

Note that by ~\eqref{e:koszul over pt}, we can identify the category $\IndCoh(\CG_{\on{pt}/\CV})$ with the category
$\on{HC}(\on{pt}/\CV)^{\on{op}}\mod$. Hence, we obtain an equivalence
\begin{multline}\label{e:ZG}
\IndCoh(\CG_{\on{pt}/\CV}) \otimes \IndCoh(Z)\simeq\on{HC}(\on{pt}/\CV)^{\on{op}}\mod\otimes \IndCoh(Z)\simeq \\
\simeq \on{HC}(\on{pt}/\CV)^{\on{op}}\mod(\IndCoh(Z)).
\end{multline}

Using the equivalence \eqref{e:ZG}, we can translate the functors \eqref{e:*!} into the language of $\BE_2$-algebras.
Recall homomorphisms of $\BE_2$-algebras
\[\on{HC}(\on{pt}/\CV)\to\on{HC}(Z/\CU)\to\on{HC}(Z).\]

\medskip

Thus any $\CF\in\IndCoh(Z)$ carries a natural structure of $\on{HC}(\on{pt}/\CV)^{\on{op}}$-module, i.e., 
we have a functor
$$\IndCoh(Z)\to \on{HC}(\on{pt}/\CV)^{\on{op}}\mod(\IndCoh(Z)).$$
It follows from the definitions that this functor identifies with the functor $\sG$. 

\medskip

Since $\IndCoh(Z)$ is tensored over $\on{HC}(\on{pt}/\CV)$, we obtain a functor
$$\on{HC}(\on{pt}/\CV)^{\on{op}}\mod\otimes \IndCoh(Z)\to \IndCoh(Z).$$
This is our functor $\sF$. 

\ssec{Proof of \thmref{t:zero sect}}  \label{ss:proof of zero sect}

The proof of \thmref{t:zero sect}, given in this subsection, relies on the material of \secref{ss:koszul functors}.
The reader who skipped \secref{ss:koszul functors} will find a proof of the most essential point of the argument in
\remref{r:zero sect}. Yet another proof, which uses a different idea, is given in \secref{sss:zero sect alt}.

\sssec{}

Let us recall that \thmref{t:zero sect} asserts that for an affine quasi-smooth DG scheme $Z$, the essential image of 
$$\Xi_Z:\QCoh(Z)\to \IndCoh(Z)$$
coincides with $\IndCoh_{\{0\}}(Z)$. Here by a slight abuse of notation $\{0\}$ denotes
the zero-section of $\Sing(Z)$. 
We are now ready to prove this theorem, using \corref{c:Koszul generation}. 

\sssec{}

Note that the statement is local on $Z$. Indeed, recall that $\Xi_Z$ is fully faithful and has a right adjoint 
$$\Psi_Z:\IndCoh(Z)\to\QCoh(Z).$$ 
\thmref{t:zero sect} is equivalent
to conservativeness of the restriction 
\[\Psi_Z|_{\IndCoh_{\{0\}}(Z)}:\IndCoh_{\{0\}}(Z)\to\QCoh(Z),\]
which can be verified locally. Thus, we may assume that $Z$ fits into a Cartesian diagram \eqref{e:gci for Koszul}.

\medskip

By \corref{c:Koszul generation}(b), it is enough to show that the essential image of $\Xi_Z$ contains the essential
image of the functor 
$$\sF:\left(\IndCoh(\CG_{\on{pt}/\CV})\otimes \IndCoh(Z)\right)_{\{0\}\times \CU}\to \IndCoh_{\{0\}}(Z).$$

\sssec{}

Consider the projection $p_{\CG_{\on{pt}/\CV}}:\CG_{\on{pt}/\CV}\to\on{pt}$. We claim 
that the essential image of the functor 
\[(p_{\CG_{\on{pt}/\CV}}\times \on{id}_Z)^!:\IndCoh(Z)\to\IndCoh(\CG_{\on{pt}/\CV}\times Z)\simeq \IndCoh(\CG_{\on{pt}/\CV})
\otimes \IndCoh(Z)\]
is contained in the category $\left(\IndCoh(\CG_{\on{pt}/\CV})\otimes \IndCoh(Z)\right)_{\{0\}}$
and generates it. 

\medskip

By \propref{p:abs product},
$$\left(\IndCoh(\CG_{\on{pt}/\CV})\otimes \IndCoh(Z)\right)_{\{0\}}=\IndCoh(\CG_{\on{pt}/\CV})_{\{0\}}\otimes \IndCoh(Z).$$
So, it is sufficient to see that the essential image of
\[p_{\CG_{\on{pt}/\CV}}^!:\Vect=\IndCoh(\on{pt})\to\IndCoh(\CG_{\on{pt}/\CV})\]
is contained in the category $\IndCoh(\CG_{\on{pt}/\CV})_{\{0\}}$ and generates it. However, this is 
a particular case of \corref{c:support KD}(b). 

\sssec{}

Hence, we obtain that it is sufficient to show that the essential image of the
composed functor
\begin{equation} \label{e:KD of trivial}
\IndCoh(Z)
\overset{(p_{\CG_{\on{pt}/\CV}}\times \on{id}_Z)^!}\longrightarrow
\IndCoh(\CG_{\on{pt}/\CV})\otimes \IndCoh(Z)\overset{\sF}\longrightarrow \IndCoh(Z)
\end{equation}
is contained in the essential image of $\Xi_Z$.

\medskip

We have the following assertion:

\begin{lem}  \label{l:calc comp smooth}
The composition \eqref{e:KD of trivial} is canonically isomorphic to $\iota^!\circ \iota^\IndCoh_*$, where
$\iota:Z\hookrightarrow \CU$.
\end{lem}

\begin{proof}
By the definition of the functor $\sF$, the lemma follows by base change along the Cartesian square
$$
\CD
\CG_{Z/\CU}  @>>>  Z \\
@VVV    @VV{\iota}V  \\
Z  @>{\iota}>>  \CU.
\endCD
$$
\end{proof}

\sssec{}

By \lemref{l:calc comp smooth}, it is sufficient to show that the essential image of the functor
$\iota^!$ is contained in the essential image of $\Xi_Z$. 

\medskip

However, since $\CU$ is smooth,
the monoidal action of $\QCoh(\CU)$ on $\omega_\CU\in \IndCoh(\CU)$ generates the latter category.
So, it is enough to show that $\iota^!(\omega_\CU)$ belongs to the essential image of $\Xi_Z$.
However, $\iota^!(\omega_\CU)=\omega_Z$, and the assertion follows from the fact that $Z$ is Gorenstein.

\section{A point-wise approach to singular support}   \label{s:pointwise}

\ssec{The functor of enhanced fiber}    \label{ss:at a point}

\sssec{} \label{sss:action on fiber}

Let $Z$ be an affine DG scheme with a perfect cotangent complex, and let 
$i_z:\on{pt}\hookrightarrow Z$ be a $k$-point. 

\medskip

Consider the functor
\[i^!_z:\IndCoh(Z)\to\Vect.\]

We claim that this functor can be naturally enhanced to a functor
\begin{equation} \label{e:enhanced fiber}
i^{\on{enh},!}_z:\IndCoh(Z)\to T_z(Z)[-1]\mod,
\end{equation}
where the DG Lie algebra $T_z(Z)[-1]$ is the fiber of $T(Z)[-1]$ at $z$,
and where we remind that $T(Z)[-1]$ is a Lie algebra by \corref{c:HH as univ env}.

\medskip

Indeed, let us interpret $i^!_z$ as
\begin{equation} \label{e:fiber as maps}
\CMaps_{\IndCoh}(i^\IndCoh_*(k),-).
\end{equation}
Now, it is easy to see that the canonical action of the DG Lie algebra $\Gamma(Z,T(Z)[-1])$ on $i^\IndCoh_*(k)$,
given by \corref{c:HH as univ env}, factors through 
$$k\underset{\Gamma(Z,\CO_Z)}\otimes \Gamma(Z,T(Z)[-1])\simeq T_z(Z)[-1].$$
This endows the functor in \eqref{e:fiber as maps} with an action of the DG Lie algebra $T_z(Z)[-1]$, as desired.

\medskip

We will refer to the functor \eqref{e:enhanced fiber} as that of \emph{enhanced fiber}. 

\sssec{}

Let us reinstate the assumption that $Z$ is quasi-smooth. 

\medskip

For an object $M\in T_z(Z)[-1]\mod$, we can consider
the graded vector space of its cohomologies $H^\bullet(M)$ as a module over the graded
Lie algebra $H^\bullet(T_z(Z)[-1]\mod)$. 

\medskip

In particular, $H^\bullet(M)$ is a module over
$\Sym\left(H^1(T_z(Z))\right)$, viewed as a graded commutative algebra whose generators are placed in 
degree $2$. Therefore, to $M$ we can associate the support
$$\on{supp}_{\Sym\left(H^1(T_z(Z))\right)}\left(H^\bullet(M)\right)\subset 
\Spec\left(\Sym\left(H^1(T_z(Z))\right)\right).$$

\medskip

Note that
$$\Spec\left(\Sym\left(H^1(T_z(Z))\right)\right)\simeq \Sing(Z)_{\{z\}}:={}^{cl}\!\left(\{z\}\underset{Z}\times \Sing(Z)\right).$$

\begin{lem}  \label{l:pointwise}
For $\CF\in \IndCoh(Z)$, set $M=i^{\on{enh},!}_z(\CF)\in T_z(Z)[-1]\mod$. Then:

\smallskip

\noindent{\em(a)} $\on{SingSupp}(\CF)\cap \Sing(Z)_{\{z\}}\supset 
\on{supp}_{\Sym\left(H^1(T_z(Z))\right)}\left(H^\bullet(M)\right)$.

\smallskip

\noindent{\em(b)} If $\CF\in \IndCoh(Z)_{\{z\}}$, then 
$$\on{SingSupp}(\CF)=\on{supp}_{\Sym\left(H^1(T_z(Z))\right)}\left(H^\bullet(M)\right).$$
\end{lem}

\begin{proof}

With no restriction of generality we can can replace $Z$ by an open affine that contains the point $z$.

\medskip

Consider the graded vector space $H^\bullet(i^!_z(\CF))$ as acted on by $\on{HH}^{\on{even}}(Z)$. 
Since the action of the subalgebra $\Gamma(Z,\CO_{^{cl}\!Z})\subset\on{HH}^{\on{even}}(Z)$ on this space
factors through the morphism 
\[\Gamma(Z,\CO_{^{cl}\!Z})\to k,\]
the action of $\on{HH}^{\on{even}}(Z)$ factors through the quotient
\[{}^{cl}(\on{HH}^{\on{even}}(Z)\underset{\Gamma(Z,\CO_{^{cl}\!Z})}\otimes k)\]
of $\on{HH}^{\on{even}}(Z)$. Here $k$ is considered as a $\Gamma(Z,\CO_{^{cl}\!Z})$-module
via $i_z$.

\medskip

Similarly, the resulting action of $\Gamma(\Sing(Z),\CO_{\Sing(Z)})$ on $H^\bullet(i^!_z(\CF))$ factors
through 
$$\Gamma(\Sing(Z)_{\{z\}},\CO_{\Sing(Z)_{\{z\}}})\simeq \Sym\left(H^1(T_z(Z)))\right),$$
which is equal to the action of the latter on
$$H^\bullet(i^!_z(\CF))\simeq H^\bullet(M).$$

\medskip

Now, point (a) of the lemma follows from the interpretation of $H^\bullet(i^!_z(\CF))$ as 
$$\Hom^\bullet_{\IndCoh(Z)}((i_z)_*(k),\CF),$$ since $(i_z)_*(k)$ is compact in $\IndCoh(\CF)$
(see \lemref{l:supp via comp}).

\medskip

Suppose now that $\CF\in \IndCoh(Z)_{\{z\}}$. Then $\on{SingSupp}(\CF)$ coincides with the support 
of $\CF$ computed using the action of $\Gamma(\Sing(Z),\CO_{\Sing(Z)})$ on
$\IndCoh(Z)_{\{z\}}$ (by \secref{sss:change of category}). But
$(i_z)_*(k)$ generates $\IndCoh(Z)_{\{z\}}$ (see 
\cite[Proposition~4.1.7(b)]{IndCoh}), and the required equality follows from
\lemref{l:supp via comp}. 
\end{proof} 
 
\sssec{An alternative proof of \thmref{t:zero sect}}  \label{sss:zero sect alt}

Let us sketch an alternative proof of Theorem~\ref{t:zero sect}. 

\medskip

By \corref{c:with support comp gen}, $\IndCoh_{\{0\}}(Z)$ is generated by $\Coh_{\{0\}}(Z)$.
Therefore, it suffices to check that $\Coh_{\{0\}}(Z)$ coincides with the essential image
$\Xi_Z(\QCoh(Z)^{\on{perf}})$. It suffices to show that for $\CF\in \Coh_{\{0\}}(Z)$ and 
every $k$-point $z$ of $Z$,
the object $i_z^!(\CF)\in \Vect$ is perfect.

\medskip

Consider the action of $\on{Sym}(H^1(T_z(Z)))$ on $H^\bullet(i_z^!(\CF))$. On the one hand, this
module is finitely generated by \thmref{t:finiteness}. On the other hand, by assumption
and \lemref{l:pointwise}(a), it is supported at 
\[0\in H^1(T_z(Z))^*=\Spec\left(\on{Sym}(H^1(T_z(Z)))\right).\]
Hence, it is finite-dimensional, as desired. 

\qed

\ssec{Estimates from below}
 
\lemref{l:pointwise} says that the support of $i^{\on{enh},!}_z(\CF)$ is bounded from above by the
singular support of $\CF$. In this subsection we will prove some converse estimates. 
 
\sssec{} 

The next assertion describes the singular support of an arbitrary object $\CF\in \IndCoh(Z)$
in terms of its !-fibers. 

\medskip

Note that for every geometric point $i_z:\Spec(k')\to Z$, we can consider the DG Lie algebra
$T_z(Z)[-1]\in \Vect_{k'}$ and the functor
$$i^{\on{enh},!}_z:\IndCoh(Z)\to T_z(Z)[-1]\mod.$$
We can do this by viewing $Z':=Z\underset{\Spec(k)}\times \Spec(k')$ 
as a quasi-smooth DG scheme over $k'$ and viewing $z'$ as a $k'$-rational point $i_{z'}:\Spec(k')\to Z'$.
The functor $i_z^!$ is the composition of $i_{z'}^!$ preceded by the tensoring-up functor 
$\IndCoh(Z)\to \IndCoh(Z')$.

\begin{prop}  \label{p:pointwise} Let $Y$ be a conical Zariski-closed subset $Y\subset \Sing(Z)$.
An object $\CF\in \IndCoh(Z)$ belongs to $\IndCoh_Y(Z)$ if and only if for every geometric point
$i_z:\Spec(k')\to Z$, the object
$$i^{\on{enh},!}_z(\CF)\in U(T_z(Z)[-1])\mod$$
is such that the subset
$$\on{supp}_{\Sym\left(H^1(T_z(Z))\right)}\left(H^\bullet(i^{\on{enh},!}_z(\CF))\right)\subset \Sing(Z)_{\{z\}}$$
is contained in 
$$Y_{\{z\}}:=\{z\}\underset{Z}\times Y\subset  \Sing(Z)_{\{z\}}.$$
\end{prop}

\begin{proof}

The ``only if" direction was established in \lemref{l:pointwise}. Let us prove the ``if" direction.
With no restriction of generality,
we can assume that $Y$ is cut out by one homogeneous element $a\in  \Gamma(\Sing(Z),\CO_{\Sing(Z)})$.
Set $2n=\deg(a)$.

\medskip

Recall that we have a pair of adjoint functors
$$\Xi_Z^{Y,\on{all}}:\IndCoh_Y(Z)\rightleftarrows \IndCoh(Z):\Psi_Z^{Y,\on{all}}.$$
Without loss of generality, we may replace $\CF$ with the cone of the morphism
\[\Xi_Z^{Y,\on{all}}\circ\Psi_Z^{Y,\on{all}}(\CF)\to\CF.\]
Note that this cone is the localization $\on{Loc}_a(\CF)$ introduced in 
\secref{sss:abs category on open}. Thus, we can assume that $a$ acts as an isomorphism
$$\CF\overset{a}\to \CF[2n].$$ Let us show that in this case $\CF=0$.

\medskip

Fix a point $z$ as above, and consider $H^\bullet(i^{\on{enh},!}_z(\CF))$ as a quasi-coherent sheaf
on $\Sing(Z)_{\{z\}}$. On the one hand, the action of $a$ on it is invertible; on the other,
its support is contained in $Y_{\{z\}}$, which is the zero locus of 
$a\in\Gamma(\Sing(Z)_{\{z\}},\CO_{\Sing(Z)_{\{z\}}})$. Hence, the quasi-coherent sheaf vanishes and
 $i^!_z(\CF)=0$. 

\medskip

The assertion now follows from the next lemma (which in turn follows from \cite[Proposition~4.1.7(a)]{IndCoh}). \end{proof}

\begin{lem} \label{l:pointwise vanishing}
If $\CF\in \IndCoh(Z)$ is such that $i^!_z(\CF)=0$ for all geometric points $z$, then $\CF=0$. \qed
\end{lem}

\sssec{The coherent case}

We have the following variant of \propref{p:pointwise}:
\footnote{Which also follows from \cite[Theorem~11.3 and Remark~11.4]{Kr}.}

\begin{prop}  \label{p:pointwise coh}
For $\CF\in \Coh(Z)$, and $z\in Z(k)$, the inclusion
$$\on{supp}_{\Sym\left(H^1(T_z(Z))\right)}\left(H^\bullet(i^{\on{enh},!}_z(\CF))\right)\subset \on{SingSupp}(\CF)\cap \Sing(Z)_{\{z\}}$$
is an equality.
\end{prop}

We note that \propref{p:pointwise coh} is \emph{not} essential for the main results of this paper. 

\medskip

\begin{rem}  \label{r:pointwise} We can reformulate Propositions \ref{p:pointwise} and \ref{p:pointwise coh} as follows.
Fix $\CF\in\IndCoh(Z)$, and consider the union
\[Y':=\underset{z\in Z}\bigcup\on{supp}_{\Sym\left(H^1(T_z(Z))\right)}\left(H^\bullet(i^{\on{enh},!}_z(\CF))\right)\subset \Sing(Z).\]
Here the union is over all (not necessarily closed) points of $Z$. \propref{p:pointwise} says that 
\[\on{SingSupp}(\CF)=\overline{Y'}.\] 
\propref{p:pointwise coh} says that $\on{SingSupp}(\CF)=Y'$ for $\CF\in\Coh(Z)$.
\end{rem}

\sssec{}

The rest of this subsection is devoted to the proof of \propref{p:pointwise coh}

\medskip

Let $f$ be a function on $Z$ such that $z$ belongs to the set of its zeros; let $Z'$
denote the corresponding DG subscheme of $Z$,
$$Z':=\on{pt}\underset{\BA^1}\times Z,$$ where
$Z\to \BA^1$ is given by $f$.

\medskip

Consider the closed subset
$$\Sing(Z)_{Z'}:={}^{cl}\!(\Sing(Z)\underset{Z}\times Z')\subset \Sing(Z).$$
For $\CF\in \Coh(Z)$, let $\CF'\in \Coh(Z)$ denote the object $\on{Cone}(f:\CF\to \CF)\in \Coh(Z)$.

\medskip

Taking into account \lemref{l:pointwise}(b), the statement of the proposition follows by induction 
from the next assertion:

\begin{lem}  \label{l:sing supp restr} For $\CF\in \Coh(Z)$,
$$\on{SingSupp}(\CF)\cap \Sing(Z)_{Z'}=\on{SingSupp}(\CF')$$
as subsets of $\Sing(Z)_{Z'}$.
\end{lem}

\sssec{Proof of \lemref{l:sing supp restr}}

The assertion is local, so with no restriction of generality, we can assume that $Z$ fits into a Cartesian diagram as in 
\eqref{e:gci for Koszul}. 

\medskip

Note that the map 
$$\on{act}_{\CG_{\on{pt}/\CV},Z}:\CG_{\on{pt}/\CV}\times Z\to Z$$
is quasi-smooth; therefore, its Tor-dimension is bounded. Hence, for $\CF\in \Coh(Z)$, the object 
$$\on{act}_{\CG_{\on{pt}/\CV},Z}^!(\CF)=:\sG(\CF)\in \IndCoh(\CG_{\on{pt}/\CV}\times \IndCoh(Z))$$
is compact by \cite[Lemma 7.1.2]{IndCoh}. 

\medskip

Now, the required assertion follows from the next one:

\begin{lem}
For a compact object $\CM\in \Sym(V[-2])\mod\otimes \IndCoh(Z)$, the support of
$\on{Cone}(f:\CM\to \CM)$ in $V^*\times Z'$ equals 
$$(V^*\times Z')\cap \supp_{V^*\times Z}(\CM).$$
\end{lem}

The lemma is proved by the argument given in the proof of
\propref{p:support KD}(b'').

\qed

\begin{rem}
Denote by  $i$ the closed embedding $Z'\hookrightarrow Z$. Clearly, $i$
is quasi-smooth. In particular, the map
$$\Sing(i):\Sing(Z)_{Z'}\to \Sing(Z')$$
is a closed embedding. 
Clearly, the object $\CF'$ above is canonically isomorphic to
$i^{\IndCoh}_*(i^!(\CF))[1]$.

\medskip

Thus, \lemref{l:sing supp restr} computes the singular support of $i^{\IndCoh}_*(i^!(\CF))$. 

\medskip

More generally,
for any morphism of quasi-smooth DG schemes $f:Z'\to Z$ and any $\CF\in\Coh(Z)$ (resp., $\CF'\in \Coh(Z')$),
there is a formula for $\on{SingSupp}(i^!(\CF))$ (resp., $\on{SingSupp}(i_*(\CF))$, assuming that $f$ is finite), see
Theorems \ref{t:ss of pullback} and \ref{t:under finite}, respectively. 
\end{rem}

\ssec{Enhanced fibers and Koszul duality} \label{ss:other enhanced}

Let $i_z:\on{pt}\to Z$ be a quasi-smooth DG scheme and a $k$-point. In this subsection we will
assume that $Z$ is written 
as a fiber product
\begin{equation} \label{e:Cart prod diag}
\CD
Z  @>{\iota}>> \CU  \\ 
@VVV    @VVV \\
\on{pt}  @>>>  \CV,
\endCD
\end{equation}
as in \eqref{e:gci for Koszul}.

\sssec{}

The action of $\CG_{\on{pt}/\CV}$ on $Z$ gives rise to a map of Lie algebras 
$$V[-2]\otimes \CO_Z\to T(Z)[-1]$$
in $\QCoh(Z)$, where $V$ is the tangent space to $\CV$ at $\on{pt}$.  This follows
from the functoriality of the construction in \corref{c:groupoid Lie} with respect to the groupoid. 

\medskip

In particular, we obtain a map of DG Lie algebras $V[-2]\to T_z(Z)[-1]$. 

\medskip

Composing, from $i_z^{\on{enh},!}$ we obtain a functor
\begin{equation} \label{e:first enhance}
\IndCoh(Z)\to V[-2]\mod.
\end{equation}

In this subsection we will give a different interpretation of the functor \eqref{e:first enhance}.

\sssec{} \label{sss:dg point}

Consider the morphism
\[\left(\on{id}\times(\iota\circ i_z)\right):\CG_{\on{pt}/\CV}=\on{pt}\underset{\CV}\times\on{pt}\to\on{pt}\underset{\CV}\times\CU=Z;\]
we denote it by $i_{z,\CV}$. It can be viewed as the action of the group DG scheme $\CG_{\on{pt}/\CV}$ on the point $z$. 
It is easy to see that $i_{z,\CV}$ is quasi-smooth.

\medskip

Thus, restriction defines a functor
$$i^!_{z,\CV}:\IndCoh(Z)\to\IndCoh(\CG_{\on{pt}/\CV}).$$

\medskip

Note that $i_z$ can be written as a composition
\begin{equation}\label{e:through the never}
i_z=i_{z,\CV}\circ \Delta_{\on{pt}},
\end{equation}
where \[\Delta_{\on{pt}}:\on{pt}\to\on{pt}\underset{\CV}\times\on{pt}=\CG_{\on{pt}/\CV}\] is the diagonal.
Let us observe the following:

\begin{lem} \label{l:enhanced by V} 
For $\CF\in\IndCoh(Z)$, $i^!_z(\CF)=0$ if and only if $i^!_{z,\CV}(\CF)=0$.
\end{lem}
\begin{proof}
Indeed, by \eqref{e:through the never}, we have $i^!_z=\Delta_{\on{pt}}^!\circ i^!_{z,\CV}$.
Therefore, it suffices to check that the functor $\Delta_{\on{pt}}^!$ is conservative. Equivalently, we need to prove that
the essential image of $(\Delta_{\on{pt}})_*^{\IndCoh}$ generates the category 
$\IndCoh(\CG_{\on{pt}/\CV})$.
This follows from \cite[Proposition 4.1.7(b)]{IndCoh} (or, in the case at hand, from \corref{c:support KD}(a)).
\end{proof}

\sssec{}

Combining the functor $i^!_{z,\CV}$ with the equivalence
$$\IndCoh(\CG_{\on{pt}/\CV})\simeq \Sym(V[-2])\mod\simeq V[-2]\mod,$$
we thus obtain a functor 
\begin{equation} \label{e:other enhance}
\IndCoh(Z)\to\IndCoh(\CG_{\on{pt}/\CV})\to  V[-2]\mod.
\end{equation}

\begin{prop}  \label{p:other enhanced}
The functors \eqref{e:other enhance} and \eqref{e:first enhance} are canonically isomorphic.
\end{prop}

\begin{proof}

The lemma easily reduces to the case when $Z=\CG_{\on{pt}/\CV}$, and $z$ is
given by $\Delta_{\on{pt}}$. In this case, the assertion is tautological from the definitions.

\end{proof}

\sssec{}  \label{sss:other enh}
As was mentioned in \remref{r:KD as enh!}, $\CG_{\on{pt}/\CV}$ is a quasi-smooth DG scheme and 
\[\Sing(\CG_{\on{pt}/\CV})=V^*,\] where $V$ is the tangent space to $\CV$ at $\on{pt}$.
In addition, for 
$$\CF\in \IndCoh(\on{pt}\underset{\CV}\times\on{pt})\simeq \on{HC}(\on{pt}/\CV)^{\on{op}},$$
we have an equality of Zariski-closed conical subsets:
\[\on{SingSupp}(\CF)=\on{supp}_{\Sym(V[-2])}(\CF)\subset V^*.\]

\sssec{}
Note that the diagram \eqref{e:Cart prod diag} yields an embedding 
\[\Sing(Z)_{\{z\}}\hookrightarrow V^*,\]
which can be viewed as the singular codifferential of $i_{z,\CV}$. (The fact that $\Sing(i_{z,\CV})$ is an embedding 
also follows from quasi-smoothness of $i_{z,\CV}$ by \lemref{l:codiff for q-smooth}.)

\medskip

From \propref{p:other enhanced} we obtain that for $\CF\in \IndCoh(Z)$, the support of $H^\bullet(i_z^!(\CF))$ as a module over
$\on{Sym}\left(H^1(T_z(Z))\right)$, considered as a subset of
$$\Spec\left(\on{Sym}\left(H^1(T_z(Z))\right)\right)=\Sing(Z)_{\{z\}}\subset  V^*,$$
equals the singular support of 
$$i^!_{z,\CV}(\CF)\in \IndCoh(\CG_{\on{pt}/\CV}).$$

\section{Functorial properties of the category $\IndCoh_Y(Z)$}   \label{s:funct}

So far, we have studied the category $\IndCoh_Y(Z)$ for a given quasi-smooth DG scheme $Z$. 
In this section we will establish a number of results on how these categories interact
under pullback and pushforward functors for maps between quasi-smooth DG schemes. The main results of this section
are stated in the introduction as Theorems~\ref{t:funct preview}, \ref{t:tensor up preview}, and \ref{t:proper preview},
corresponding to \propref{p:singsupp functoriality}, \corref{c:quasi-smooth tensor up}, and \thmref{t:prop cons} below. 

\ssec{Behavior under direct and inverse images}  \label{ss:functoriality}

Recall that a map $f$ between DG schemes induces the following functors between the categories
of ind-coherent sheaves: the ``ordinary'' pushforward, which is denoted by $f_*^{\IndCoh}$ (the notation
$f^*$ is reserved for the pushforward on the category of quasi-coherent sheaves), and the
``extraordinary'' pullback $f^!$. (Actually, we need $f$ to be quasi-compact for $f_*^{\IndCoh}$ to exist.)
If $f$ is eventually coconnective, the ``ordinary'' pullback, denoted by $f^{\IndCoh,*}$, makes sense as well.

\sssec{} Let $f:Z_1\to Z_2$ be a map between quasi-smooth DG schemes. Consider the functor
$$f^!:\IndCoh(Z_2)\to \IndCoh(Z_1),$$
(see \cite[Sect. 5.2.3]{IndCoh})
and, assuming that $f$ is quasi-compact, the functor
$$f^\IndCoh_*:\IndCoh(Z_1)\to \IndCoh(Z_2)$$
(see \cite[Sect. 3.1]{IndCoh}). 

\medskip

Recall also that if $f$ is proper, the above functors $(f^\IndCoh_*,f^!)$ are naturally adjoint.

\sssec{}

Recall that $f$ gives rise to the singular codifferential
\[\Sing(f):\Sing(Z_2)_{Z_1}\to\Sing(Z_1),\]
where
\[\Sing(Z_2)_{Z_1}={}^{cl}\left(\Sing(Z_2)\underset{Z_2 }\times Z_1\right).\] 

\begin{prop}\label{p:singsupp functoriality}
  Let $Y_i\subset \Sing(Z_i)$ be conical Zariski-closed subsets.
\smallskip

\noindent{\em(a)} Suppose 
\[\Sing(f)(Y_2\underset{Z_2}\times Z_1)\subset Y_1.\] 
Then
$f^!$ sends $\IndCoh_{Y_2}(Z_2)$ to $\IndCoh_{Y_1}(Z_1)$.

\smallskip

\noindent{\em(b)} Suppose that $f$ is quasi-compact, and that 
\[\Sing(f)^{-1}(Y_1)\subset Y_2\underset{Z_2}\times Z_1.\] 
Then
$f^\IndCoh_*$ sends $\IndCoh_{Y_1}(Z_1)$ to $\IndCoh_{Y_2}(Z_2)$.
\end{prop}

\begin{proof}
First of all, in both claims we may assume that $Z_1$ and $Z_2$ are affine. Indeed, claim (a) is clearly local on both $Z_1$ and $Z_2$.
On the other hand, claim (b) is clearly local on $Z_2$. By \corref{c:on open}(b), it is also local on $Z_1$: an open cover of $Z_1$
can be used to compute the direct image $f^\IndCoh_*$ using the \v{C}ech resolution. 

\smallskip

Since the functor
\[f^!:\IndCoh(Z_2)\to\IndCoh(Z_1)\] is continuous, it
corresponds to an object of
\[\IndCoh(Z_2)^\vee\otimes\IndCoh(Z_1).\]
By Serre's duality,
\[\IndCoh(Z_2)^\vee\otimes\IndCoh(Z_1)\simeq \IndCoh(Z_2)\otimes\IndCoh(Z_1)\simeq\IndCoh(Z_1\times Z_2),\]
and it is clear that $f^!$ corresponds to the object
\[\Gamma(f)^\IndCoh_*(\omega_{Z_1})\in \IndCoh(Z_1\times Z_2),\]
where $\Gamma(f):Z_1\to Z_1\times Z_2$ is the graph of $f$.
Similarly, the continuous functor 
\[f^\IndCoh_*:\IndCoh(Z_1)\to\IndCoh(Z_2)\] 
corresponds to the same object under the identification
\[\IndCoh(Z_1)^\vee\otimes\IndCoh(Z_2)\simeq \IndCoh(Z_1)\otimes\IndCoh(Z_2)\simeq\IndCoh(Z_1\times Z_2).\]

Now the assertion follows from \propref{p:supp functoriality} and \lemref{l:singsupp of graph} below. 
\end{proof}

\begin{lem} \label{l:singsupp of graph} The singular support
\[\on{SingSupp}\left(\Gamma(f)^\IndCoh_*(\omega_{Z_1})\right)\subset\Sing(Z_1\times Z_2)=\Sing(Z_1)\times\Sing(Z_2)\] 
is contained in the image of
\[\Sing(Z_2)_{Z_1}=Z_1\underset{Z_2}\times \Sing(Z_2)\] 
under the natural map of the latter to 
$\Sing(Z_1)\times \Sing(Z_2)$.
\end{lem}

\begin{proof} The statement is clearly local on both $Z_1$ and $Z_2$, so we may assume that $Z_1$ and $Z_2$
are affine without losing generality.
Since $\Gamma(f)^\IndCoh_*(\omega_{Z_1})$ is compact, by \lemref{l:supp via comp}(b)
it is enough to show that the homomorphisms
$$\Gamma(\Sing(Z_i),\CO_{\Sing(Z_i)})\to \End^\bullet_{\IndCoh(Z_1\times Z_2)}\left(\Gamma(f)^\IndCoh_*(\omega_{Z_1})\right)$$
for $i=1,2$
factor through a map 
$$\Gamma(\Sing(Z_2)_{Z_1},\CO_{\Sing(Z_2)_{Z_1}})\to \End^\bullet_{\IndCoh(Z_1\times Z_2)}\left(\Gamma(f)^\IndCoh_*(\omega_{Z_1})\right)$$
and the natural homomorphisms
$$\Gamma(\Sing(Z_i),\CO_{\Sing(Z_i)})\to \Gamma(\Sing(Z_2)_{Z_1},\CO_{\Sing(Z_2)_{Z_1}}).$$

\medskip

We have
$$\CMaps_{\IndCoh(Z_1\times Z_2)}\left(\Gamma(f)^\IndCoh_*(\omega_{Z_1}),\Gamma(f)^\IndCoh_*(\omega_{Z_1})\right)\simeq
\Gamma\left(Z_1,U_{\CO_{Z_1}}(T(Z_2)[-1]|_{Z_1})\right)$$
(established in the course of the proof of \propref{p:groupoid 1} due to the retraction $Z_1\times Z_2\to Z_1$).

\medskip

Moreover, the homomorphisms of $\BE_1$-algebras
$$\on{HC}(Z_i)\to \CMaps_{\IndCoh(Z_1\times Z_2)}\left(\Gamma(f)^\IndCoh_*(\omega_{Z_1}),\Gamma(f)^\IndCoh_*(\omega_{Z_1})\right)$$
identify with the naturally defined maps
$$\Gamma\left(Z_i,U_{\CO_{Z_i}}(T(Z_i)[-1])\right)\to \Gamma\left(Z_1,U_{\CO_{Z_1}}(T(Z_2)[-1]|_{Z_1})\right).$$

This establishes the desired assertion.
\end{proof} 

\sssec{}
Assume now that both $Z_1$ and $Z_2$ are quasi-compact. 
Recall (see \cite[Sect. 9.2.3]{IndCoh}), that under the self-duality
$$\bD^{\on{Serre}}_{Z_i}:\IndCoh(Z_i)^\vee\simeq \IndCoh(Z_i),$$
the dual of the functor $f^!$ is $f^\IndCoh_*$, and vice versa.

\medskip

Hence, from \propref{p:singsupp functoriality} and \lemref{l:Serre dual of embedding}, we obtain:

\begin{prop} \label{p:functoriality quotient}
Let $Y_i\subset \Sing(Z_i)$ be conical Zariski-closed subsets.
\smallskip

\noindent{\em(a)} Suppose 
\[\Sing(f)(Y_2\underset{Z_2}\times Z_1)\subset Y_1.\] 
Then  we have a commutative
diagram of functors:
$$
\CD
\IndCoh(Z_1)   @>{\Psi_{Z_1}^{Y_1,\on{all}}}>>  \IndCoh_{Y_1}(Z_1)  \\
@V{f^\IndCoh_*}VV    @VVV   \\
\IndCoh(Z_2)   @>{\Psi_{Z_2}^{Y_2,\on{all}}}>>  \IndCoh_{Y_2}(Z_2).
\endCD
$$
That is, the counter-clockwise composition functor factors through the colocalization $\Psi_{Z_1}^{Y_1,\on{all}}$.

\smallskip

\noindent{\em(b)} Suppose that 
\[\Sing(f)^{-1}(Y_1)\subset Y_2\underset{Z_2}\times Z_1.\] 
Then
we have have a commutative
diagram of functors:
$$
\CD
\IndCoh(Z_1)   @>{\Psi_{Z_1}^{Y_1,\on{all}}}>>  \IndCoh_{Y_1}(Z_1)  \\
@A{f^!}AA    @AAA   \\
\IndCoh(Z_2)   @>{\Psi_{Z_2}^{Y_2,\on{all}}}>>  \IndCoh_{Y_2}(Z_2).
\endCD
$$
That is, the clockwise composition functor factors through the colocalization $\Psi_{Z_2}^{Y_2,\on{all}}$.

\hfill\qed
\end{prop}

\ssec{Singular support and preservation of coherence}

In this subsection all DG schemes will be quasi-compact. 

\sssec{}

Let $Z$ be a quasi-smooth DG scheme. In turns out that the knowledge of the singular 
support of an object $\CF\in \Coh(Z)$ allows one to predict when certain functors applied to it produce a coherent 
object. Namely, we will prove the following assertion:

\begin{prop}  \label{p:preserve coherence} \hfill

\smallskip

\noindent{\em(a)} For $\CF',\CF''\in \Coh(Z)$ such that, set-theoretically,
$$\on{SingSupp}(\CF')\cap \on{SingSupp}(\CF'')\subset \{0\},$$
their internal Hom object $$\ul\Hom_{\QCoh(Z)}(\CF',\CF'')\in \QCoh(Z)$$
belongs to $\Coh(Z)$ (equivalently, is cohomologically bounded above). 

\smallskip

\noindent{\em(b)} Under the assumptions of point \emph{(a)}, the tensor product
$$\CF'\otimes \CF''\in \QCoh(Z)$$
belongs to $\Coh(Z)$ (equivalently, is cohomologically bounded below). 

\smallskip

\noindent{\em(c)} Let $f:Z_1\to Z_2$ be a morphism of quasi-smooth DG schemes. Let us
denote by
$\on{ker}(\Sing(f))\subset\Sing(Z_2)_{Z_1}$ the preimage of the zero section under 
\[\Sing(f):\Sing(Z_2)_{Z_1}\to\Sing(Z_1).\] 
For any $\CF_2\in \Coh(Z_2)$ such that, set-theoretically,
\[\left(\on{SingSupp}(\CF_2)\underset{Z_2}\times Z_1\right)\cap\on{ker}(\Sing(f))\subset \{0\},\]
we have $f^!(\CF_2)\in \Coh(Z_1)$.

\smallskip

\noindent{\em(d)} Under the assumptions of point \emph{(c)}, we have
$f^*(\CF_2)\in \Coh(Z_1)$ (equivalently, $f^*(\CF_2)\in \QCoh(Z_1)$ is bounded below).

\smallskip

\noindent{\em(e)} Under the assumptions of point \emph{(c)}, the partially defined left adjoint
$f^{\IndCoh,*}$ to 
$$f^{\IndCoh}_*:\IndCoh(Z_1)\to \IndCoh(Z_2),$$ 
is defined on $\CF_2$.  

\end{prop}

\begin{rem}
One can show, mimicking the proof of \thmref{t:ss of pullback} below, that the assertions in the above
proposition are actually ``if and only if".
\end{rem}

The rest of this subsection is devoted to the proof of the above proposition. First, we notice that all
assertions are local in the Zariski topology, so we can assume that the DG schemes involved are
affine.

\sssec{Proof of point \emph{(a)}}

Since $Z$ is affine, it suffices to show that the graded vector space 
$$\Hom^\bullet(\CF',\CF'')$$
is cohomologically bounded above. 

\medskip

By \thmref{t:finiteness}, $\Hom^\bullet(\CF',\CF'')$ is finitely generated as a module over
$\Gamma(\Sing(Z),\CO_{\Sing(Z)})$. Note that the $\Gamma(\Sing(Z),\CO_{\Sing(Z)})$-action on
$\Hom^\bullet(\CF',\CF'')$ factors through both its action on $\End^\bullet(\CF')$ and 
$\End^\bullet(\CF'')$. Hence, we obtain that
$$\on{supp}_{\Gamma(\Sing(Z),\CO_{\Sing(Z)})}(\Hom^\bullet(\CF',\CF''))\subset 
\on{SingSupp}(\CF')\cap \on{SingSupp}(\CF'')\subset \{0\}.$$

The latter implies that $\Hom^\bullet(\CF',\CF'')$ is finitely generated as a module over
$\Gamma(Z,\CO_Z)$. This implies the desired assertion.

\qed

\sssec{Proof of point \emph{(c)}}

Replacing $Z_2$ by $Z_1\times Z_2$ and $\CF_2$ by $\omega_{Z_1}\boxtimes \CF_2$, we can
assume that $f$ is a closed embedding. 

\medskip

It suffices to show that $f^!(\CF_2)$ is cohomologically bounded above. The latter is equivalent to
$$\Hom^\bullet_{\Coh(Z_2)}(f_*(\CO_{Z_1}),\CF_2)$$
living in finitely many cohomological degrees. 

\medskip

Now, \propref{p:singsupp functoriality}(b) implies that $\on{SingSupp}(f_*(\CO_{Z_1}))$ is contained
in the image of $\on{ker}(\Sing(f))$ under the projection
$$\Sing(Z_2)_{Z_1}\to Z_2.$$

Therefore, the condition on $\on{SingSupp}(\CF_2)$ implies that
$$\on{SingSupp}(\CF_2)\cap \on{SingSupp}(f_*(\CO_{Z_1}))=\{0\}_{Z_2}.$$

Hence, the required assertion follows from point (a) of the proposition.

\qed

\sssec{Proof of point \emph{(e)}}

By \corref{c:Serre}, the object $\BD_{Z_2}^{\on{Serre}}(\CF_2)\in \Coh(Z_2)$ satisfies the condition of point (c). 
We claim that the object
$$\BD^{\on{Serre}}_{Z_1}\left(f^!(\BD_{Z_2}^{\on{Serre}}(\CF_2))\right)\in \Coh(Z_1)$$
satisfies the required adjunction property. Indeed, for $\CF_1\in \IndCoh(Z_1)$, we have
\begin{multline*}
\Hom_{\IndCoh(Z_1)}(\BD^{\on{Serre}}_{Z_1}\left(f^!(\BD_{Z_2}^{\on{Serre}}(\CF_2))\right),\CF_1)\simeq
\langle f^!(\BD_{Z_2}^{\on{Serre}}(\CF_2)),\CF_1\rangle_{\IndCoh(Z_1)}\simeq \\
\simeq \langle \BD_{Z_2}^{\on{Serre}}(\CF_2),f^{\IndCoh}_*(\CF_1)\rangle \simeq
\Hom_{\IndCoh(Z_2)}(\CF_2,f^{\IndCoh}_*(\CF_1)),
\end{multline*}
where 
$$\langle-,-\rangle_{\IndCoh(Z_1)} \text{ and } \langle-,-\rangle_{\IndCoh(Z_2)}$$
denote the canonical pairings corresponding to the Serre duality equivalences
$$\bD^{\on{Serre}}_{Z_i}:\IndCoh(Z_i)^\vee\simeq \IndCoh(Z_i),$$
$i=1,2$.

\qed

\sssec{Proof of point \emph{(d)}}

Consider the object
$$f^{\IndCoh,*}(\CF_2)\in \IndCoh(Z_1),$$
whose existence is guaranteed by point (e).
In particular, for $\CF_1\in \IndCoh(\CF_1)^+$ we have a functorial isomorphism
$$\Hom_{\IndCoh(Z_1)}(f^{\IndCoh,*}(\CF_2),\CF_1)\simeq \Hom_{\IndCoh(Z_2)}(\CF_2,f^{\IndCoh}_*(\CF_1)).$$

\medskip

By construction, $f^{\IndCoh,*}(\CF_2)\in\Coh(Z_1)$ (this also follows because it is the value on a compact object 
of a partially defined left adjoint to a continuous functor).

\medskip

From the commutative diagram
$$
\CD
\IndCoh(Z_1)^+  @>{\sim}>> \QCoh(Z_1)^+ \\
@V{f^{\IndCoh,*}}VV    @VV{f^*}V   \\
\IndCoh(Z_2)^+  @>{\sim}>> \QCoh(Z_2)^+
\endCD
$$
we obtain an adjunction
$$\Hom_{\QCoh(Z_1)}(f^{\IndCoh,*}(\CF_2),\CF'_1)\simeq \Hom_{\QCoh(Z_2)}(\CF_2,f_*(\CF'_1))$$
for $\CF'_1\in \QCoh(Z_1)^+$. Now, the fact that $\QCoh(Z_i)$, $i=1,2$ is left-complete in its
t-structure implies that the above adjunction remains valid for any $\CF'_1\in \QCoh(Z_1)$. 
Hence, $f^{\IndCoh,*}(\CF_2)$, viewed as an object of
$$\Coh(Z_1)\subset \QCoh(Z_1),$$
is isomorphic to $f^*(\CF_2)$. In particular, the latter belongs to $\Coh(Z_1)$, as desired.

\qed

\sssec{Proof of point \emph{(b)}}

This follows formally from point (d) applied to the diagonal morphism $Z\to Z\times Z$ and
$$\CF'\boxtimes \CF''\in \Coh(Z\times Z).$$

\qed

\ssec{Direct image for finite morphisms}  \label{ss:finite}

In this subsection, we let $f:Z_1\to Z_2$ be a finite morphism between quasi-smooth 
DG schemes ($Z_1$ and $Z_2$ need not be quasi-compact). 
For instance, $f$ may be a closed embedding. 

\sssec{}   

Define $Y_{1,can}\subset\Sing(Z_1)$ to be the image of the singular codifferential
\[\Sing(f):\Sing(Z_2)_{Z_1}=\Sing(Z_2)\underset{Z_2}\times Z_1\to\Sing(Z_1).\]
Note that $Y_{1,can}$ is constructible, but not necessarily Zariski closed. If $f$ is quasi-smooth,
$\Sing(f)$ is a closed embedding and $Y_{1,can}$ is closed.

\sssec{}
Let $\CF$ be an object of $\IndCoh(Z_1)$, and let $Y_1\subset \Sing(Z_1)$ be its singular support. 
Let $Y_2\subset \Sing(Z_2)$ be the conical Zariski-closed subset equal to the projection of
$$\Sing(f)^{-1}(Y_1)\subset \Sing(Z_2)\underset{Z_2}\times Z_1$$
under $p:\Sing(Z_2)\underset{Z_2}\times Z_1\to \Sing(Z_2)$. It is automatically closed since $p$
is finite and therefore proper.

\medskip

Consider the object $f^\IndCoh_*(\CF)\in \IndCoh(Z_2)$. Note that by 
\propref{p:singsupp functoriality}(b), we have:
$$\on{SingSupp}(f^\IndCoh_*(\CF))\subset Y_2.$$

\begin{thm}  \label{t:under finite} \hfill

\smallskip

\noindent{\em(a)}
Suppose that $Y_1\subset Y_{1,can}$. Then $\on{SingSupp}(f^\IndCoh_*(\CF))=Y_2$.

\smallskip

\noindent{\em(b)}
Suppose that $\CF\in\Coh(Z_1)$. Then $\on{SingSupp}(f^\IndCoh_*(\CF))=Y_2$. 
\end{thm}

\begin{rem}
Point (b) of the theorem will be used in \corref{c:Drinfeld}
to give an explicit characterization of singular support of coherent
sheaves, due to Drinfeld. However, it is not essential for the main results
of this paper.
\end{rem}

\begin{proof} As in \remref{r:pointwise}, consider the union
\[Y_1':=\underset{z_1\in Z_1}\bigcup\on{supp}_{\Sym\left(H^1(T_{z_1}(Z_1))\right)}\left(H^\bullet(i^{\on{enh},!}_{z_1}(\CF))\right)
\subset \Sing(Z_1).\]
By \propref{p:pointwise}, $Y_1=\overline{Y_1'}$. Similarly, consider the union
\[Y_2':=\underset{z_2\in Z_2}\bigcup\on{supp}_{\Sym\left(H^1(T_{z_2}(Z_2))\right)}
\left(H^\bullet(i^{\on{enh},!}_{z_2}(f^\IndCoh_*(\CF)))\right)\subset \Sing(Z_2);\]
then $\on{SingSupp}(f^\IndCoh_*(\CF))=\overline{Y_2'}$. It suffices to verify that under the hypotheses of the theorem, 
$Y_2'$ is equal to the projection of $\Sing(f)^{-1}(Y_1')$.

\medskip

Let us reduce the assertion of the theorem
to the case when $^{cl}\!Z_2$ is a single point. Let $z_2\in Z_2$ be a point of $Z_2$, which we may 
assume to be a $k$-point after extending scalars. Choose
a quasi-smooth map $$i_2:Z'_2\to Z_2,$$ as in \secref{sss:dg point}, so that $Z'_2$ is a DG scheme of the form 
$\on{pt}\underset{\CV_2}\times \on{pt}$,
with $\CV_2$ smooth, and such that the unique $k$-point of $Z'_2$ goes to $z_2$. 

\medskip

Denote
$$Z'_1:=Z_1\underset{Z_2}\times Z'_2,$$ and let $i_1$ denote the corresponding map $Z'_1\to Z_1$.
The map $i_1$ is also quasi-smooth by base change. Since $Z_1$ itself is quasi-smooth, we obtain that $Z_1'$ is quasi-smooth. 
Note also that $Z_1'$ is finite; therefore, by \lemref{l:shape of point}(a),
$Z_1'$ is isomorphic to a finite disjoint union of DG schemes of the form $\on{pt}\underset{\CV_1}\times \on{pt}$.

\medskip

By \secref{sss:other enh}, we know that 
\[\on{SingSupp}(i_2^!(f^\IndCoh_*(\CF)))=Y_2'\cap\Sing(Z_2)_{Z_2'}\subset \Sing(Z_2)_{Z_2'}=\Sing(Z_2)_{\{z_2\}}\subset\Sing(Z_2'),\]
and 
\[\on{SingSupp}(i_1^!(\CF))=Y_1'\cap\Sing(Z_1)_{Z_1'}\subset\Sing(Z_1)_{Z_1'}\subset\Sing(Z_1').\]
Base change allows us to replace $Z_1$, $Z_2$, and $\CF$ by
$Z_1'$, $Z_2'$, and $i_1^!(\CF)$, respectively. Note that $i_1^!(\CF)$ satisfies the hypotheses of the theorem
(this relies on $i_1$ being eventually coconnective, so that $i_1^!$ preserves coherence). 

\medskip

Thus, we assume that $^{cl}\!Z_2$ is a single point. It suffices to check the claim with $Z_1$ replaced by each of 
its connected components, so we
may assume that $^{cl}\!Z_1$  is a single point as well. Now the claim follows from \lemref{l:shape of point}(b)
and \propref{p:support KD}(b' and b'').
\end{proof}

\sssec{}
From \thmref{t:under finite}(b), we can derive an explicit characterization of singular support for objects of 
$\Coh(Z)\subset \IndCoh(Z)$.

\medskip

Let $(z,\xi)$ be a point of $\Sing(Z)$, where $z\in Z(k)$ and $0\neq \xi \in H^{-1}(T^*_z(Z))$.
We would like to determine when this point belongs to $\on{SingSupp}(\CF)$ for a given $\CF\in \Coh(Z)$.

\medskip

Let $Z$ be written as 
\[
\CD
Z  @>{\iota}>> \CU  \\ 
@VVV    @VVV \\
\on{pt}  @>>>  \CV,
\endCD
\]
with smooth $\CU$ and $\CV$, as in \secref{sss:expl presentations}. 

\medskip

Using the embedding $\Sing(Z)\hookrightarrow V^*\times Z$, we can view 
$\xi$ as a cotangent vector to $\CV$ at $\on{pt}$. Choose a function $\CV\to \BA^1$ 
that sends $\on{pt}\mapsto 0$, and whose differential equals $\xi$. Let $Z'$ be the
fiber product
$$
\CD
Z'  @>>>  \CU \\
@VVV   @VVV  \\
\on{pt}  @>>>  \BA^1.
\endCD
$$

Let $f$ denote the closed embedding $Z\hookrightarrow Z'$. 

\medskip

We have the following characterization of singular support, suggested to us by V.~Drinfeld:

\begin{cor} \label{c:Drinfeld}
The element $(z,\xi)$ belongs to $\on{SingSupp}(\CF)$ if and only if
$f_*(\CF)\in \Coh(Z')$ is \emph{not perfect} on a Zariski neighborhood of $z$.
\end{cor}

\begin{proof}

Note that $\Sing(Z')_{\{z\}}=\on{Span}(\xi)$. 

\medskip

Let is first prove the ``only if" direction. As was mentioned above, \propref{p:singsupp functoriality}(b)
implies that if $(z,\xi)\notin \on{SingSupp}(\CF)$, then $(z,\xi)\notin \on{SingSupp}(f^\IndCoh_*(\CF))$.
Hence, on a Zariski neighborhood of $z$, we have
$$\on{SingSupp}(f^\IndCoh_*(\CF))\subset \{0\}.$$
Therefore, by \thmref{t:zero sect}, $f^\IndCoh_*(\CF)$ belongs to the essential image of the functor
$$\Xi_{Z'}:\QCoh(Z')\to \IndCoh(Z').$$ Now, the assertion follows from the following general lemma
(\cite[Lemma 1.5.8]{IndCoh}):

\begin{lem}
For an eventually coconnective DG scheme $Z$, the intersection $$\Coh(Z)\cap \Xi_Z(\QCoh(Z))\subset \IndCoh(Z)$$
equals $\Xi_Z(\QCoh(Z)^{\on{perf}})$.
\end{lem}

\begin{proof}
The assertion is local, so we can assume that $Z$ is quasi-compact.  
Since the functor $\Xi_Z$ is fully faithful and continuous, if $\Xi_Z(\CF)$ is compact in $\IndCoh(Z)$,
then $\CF$ is compact in $\QCoh(Z)$, i.e., $\CF\in \QCoh(Z)^{\on{perf}}$.
\end{proof}

\medskip

For the ``if" direction, assume that $(z,\xi)\in \on{SingSupp}(\CF)$. By \thmref{t:under finite}, we obtain that 
$(z,\xi)$ belongs to $\on{SingSupp}(f_*^\IndCoh(\CF))$, considered as an object
of $\Coh(Z')$. Hence, $f_*^\IndCoh(\CF)$ is not perfect on any Zariski neighborhood of $z$ by 
the easy direction in \thmref{t:zero sect}.

\end{proof}

\begin{rem}
We note that the assertion of \corref{c:Drinfeld} makes sense also when $\xi=0$. It is easy to adapt
the proof to show that it is valid in this case as well.
\end{rem}

\sssec{}

Let $Y_1\subset \Sing(Z_1)$ be conical Zariski-closed subset, and assume that $Y_1$
is contained in the image of $\Sing(Z_1)_{Z_2}$ under $\Sing(f)$. (Recall that the image 
is constructible, but not necessarily closed.)

\medskip

From \thmref{t:under finite}(a), we obtain the following corollary.

\begin{cor} \label{c:under finite}
The restricted functor 
\[f^{\IndCoh}_*|_{\IndCoh_{Y_1}(Z_1)}:\IndCoh_{Y_1}(Z_1)\to \IndCoh(Z_2)\]
is conservative.
\end{cor}

\ssec{Conservativeness for finite quasi-smooth maps} \label{ss:finite qs maps}

\sssec{}

Let us remain in the setting of \secref{ss:finite}, and 
let us assume in addition that $f$ is quasi-smooth. For instance, $f$ could be
a quasi-smooth closed embedding, so that $Z_1$ is a ``locally
complete intersection in $Z_2$." 

\medskip

In this case, $\Sing(f)$ is a closed embedding. As before, let $Y_{1,can}\subset\Sing(Z_1)$ be the image of
$\Sing(f)$, which is a conical Zariski-closed subset. 

\sssec{} Recall now that $f$ is eventually coconnective, so by \corref{c:f^*}, the functor $f^\IndCoh_*$ admits a left adjoint, $f^{\IndCoh,*}$.
Moreover, $f$ is Gorenstein by \corref{c:Gorenstein}, so the functors 
$f^{\IndCoh,*}$ and $f^!$ can be obtained from one another by tensoring
by a cohomologically shifted line bundle (see \cite[Proposition 7.3.8]{IndCoh}). 

\medskip

By \propref{p:singsupp functoriality}(a), we have two pairs of adjoint functors
$$f^\IndCoh_*:\IndCoh_{Y_{1,can}}(Z_1)\rightleftarrows \IndCoh(Z_2):f^!$$
and
\[f^{\IndCoh,*}:\IndCoh(Z_2)\rightleftarrows\IndCoh_{Y_{1,can}}(Z_1):f^\IndCoh_*.\]
 
\begin{prop} \label{p:closed embed cons}
Suppose $f:Z_1\to Z_2$ is a finite quasi-smooth morphism between quasi-smooth DG schemes. Then 
the essential image of $\IndCoh(Z_2)$ under the functor $f^!$ generates $\IndCoh_{Y_{1,can}}(Z_1)$.
\end{prop}

\begin{proof}
Since the functors $f^!$ and $f^{\IndCoh,*}$ differ by tensoring by a cohomologically shifted line bundle, 
the statement is equivalent to the claim that the restriction $f^\IndCoh_*|_{\IndCoh_{Y_{1,can}}(Z_1)}$
is conservative. This is a particular case of \corref{c:under finite}.
\end{proof}

\begin{rem} \label{r:zero sect}
\propref{p:closed embed cons} is a generalization of \thmref{t:zero sect}. Indeed, if we assume that 
$Z$ is a global complete intersection, it admits a quasi-smooth closed embedding $\iota:Z\to\CU$, where $\CU$ is smooth.
The key step in the proof of \thmref{t:zero sect} (see \secref{ss:proof of zero sect}) is to show that 
the essential image $\iota^!(\IndCoh(\CU))$ generates $\IndCoh_{\{0\}}(Z)$. But this is exactly the 
assertion of \propref{p:closed embed cons} applied to $\iota$. 
\end{rem}

\sssec{}  \label{sss:q-smooth with supports}

Let now $Y_2$ be an arbitrary conical Zariski-closed subset of $\Sing(Z_2)$. Let 
$Y_1\subset \Sing(Z_1)$
be the image of $Y_2\underset{Z_2}\times Z_1$ under the singular codifferential
\[\Sing(f):\Sing(Z_2)_{Z_1}\to \Sing(Z_1).\]

By \propref{p:singsupp functoriality} we have two pairs of adjoint functors:
$$f^\IndCoh_*:\IndCoh_{Y_1}(Z_1)\rightleftarrows \IndCoh_{Y_2}(Z_2):f^!$$
and
\[f^{\IndCoh,*}:\IndCoh_{Y_2}(Z_2)\rightleftarrows\IndCoh_{Y_1}(Z_1):f^\IndCoh_*.\]

Then from \corref{c:under finite} we obtain:

\begin{cor} \label{c:closed embed cons}  Under the above circumstances: 

\smallskip

\noindent{\em(a)} The functor $f^\IndCoh_*:\IndCoh_{Y_1}(Z_1)\to \IndCoh_{Y_2}(Z_2)$
is conservative.

\smallskip

\noindent{\em(b)} The essential image of $\IndCoh_{Y_2}(Z_2)$ under $f^!$ (or under $f^{\IndCoh,*}$) generates
$\IndCoh_{Y_1}(Z_1)$.
\end{cor}

\sssec{A digression}

Let $f:W_1\to W_2$ be a locally eventually coconnective morphism and that $W_2$ is quasi-compact.
Let $\QCoh(W_2)$ act on $\IndCoh(W_1)$ via the homomorphism 
of monoidal categories $f^*:\QCoh(W_2)\to \QCoh(W_1)$.
The functors $f^!$ and $f^{\IndCoh,*}$ are $\QCoh(W_2)$-linear, and therefore induce two functors
\[\QCoh(W_1)\underset{\QCoh(W_2)}\otimes \IndCoh(W_2)\rightrightarrows \IndCoh(W_1).\]

\begin{lem}  \label{l:tensor up ff}
Let $f:W_1\to W_2$ be a locally eventually coconnective morphism with $W_2$ quasi-compact. Then the
functors
$$\QCoh(W_1)\underset{\QCoh(W_2)}\otimes \IndCoh(W_2)\rightrightarrows \IndCoh(W_1)$$
induced by the functors $f^!$ and $f^{\IndCoh,*}$ are fully faithful.
\end{lem}

\begin{proof}
This is \cite[Propositions 4.4.2 and 7.5.9]{IndCoh}. 
\end{proof}

\sssec{} Let $f:Z_1\to Z_2$ be as above, with $Z_1$ (and hence $Z_2$) quasi-compact. 
As in \secref{sss:q-smooth with supports}, we let $Y_2$ be a conical Zariski-closed subset of $\Sing(Z_2)$, and 
let $Y_1\subset \Sing(Z_1)$ be the image of $Y_2\underset{Z_2}\times Z_1$ under $\Sing(f)$.

\medskip

Since the restriction
$$f^!:\IndCoh_{Y_2}(Z_2)\to \IndCoh_{Y_1}(Z_1)$$
is $\QCoh(Z_2)$-linear, 
it induces a functor
\begin{equation} \label{e:closed embedding base change supports}
\QCoh(Z_1)\underset{\QCoh(Z_2)}\otimes \IndCoh_{Y_2}(Z_2)\to \IndCoh_{Y_1}(Z_1).
\end{equation}

\begin{cor}  \label{c:closed embedding base change}
The functor \eqref{e:closed embedding base change supports} is an equivalence. 
\end{cor}

\begin{proof}
The functor is fully faithful by \lemref{l:tensor up ff}, and its essential image generates
$\IndCoh_{Y_1}(Z_1)$ by \corref{c:closed embed cons}(b).
\end{proof}

\ssec{Behavior under smooth morphisms}

\sssec{}

Let $f:Z_1\to Z_2$ be a smooth map between DG schemes, and assume that $Z_2$ is quasi-compact. 
Recall (see \cite[Proposition 4.5.3]{IndCoh}) that the functor 
$$f^{\IndCoh,*}:\IndCoh(Z_2)\to \IndCoh(Z_1)$$ gives rise to an equivalence of categories
$$\IndCoh(Z_2)\underset{\QCoh(Z_2)}\otimes \QCoh(Z_1)\to \IndCoh(Z_1).$$

Similarly,the functor $f^!:\IndCoh(Z_2)\to \IndCoh(Z_1)$ gives rise to an equivalence
\begin{equation} \label{e:smooth tensor-up}
\IndCoh(Z_2)\underset{\QCoh(Z_2)}\otimes \QCoh(Z_1)\to \IndCoh(Z_1)
\end{equation}
see \cite[Corollary 7.5.7]{IndCoh}. 

\sssec{}

Assume now that $Z_2$ (and, hence, $Z_1$) is quasi-smooth. Recall from \lemref{l:codiff for smooth}
that in this case, the singular codifferential
\[
\Sing(f):\Sing(Z_2)_{Z_1}:={}^{cl}(\Sing(Z_2)\underset{Z_2}\times Z_1)\simeq \Sing(Z_2)\underset{Z_2}\times Z_1
\to \Sing(Z_1)\]
is an isomorphism. 

\medskip

Fix a conical Zariski-closed subset $Y_2\subset \Sing(Z_2)$, and let
$Y_1\subset \Sing(Z_1)$
be the image
\[\Sing(f)\left(Y_2\underset{Z_2}\times Z_1\right).\]

\medskip

We have:
\begin{prop} \label{p:smooth pullback}
Under the equivalence of \eqref{e:smooth tensor-up}, we have
$$\IndCoh_{Y_2}(Z_2)\underset{\QCoh(Z_2)}\otimes \QCoh(Z_1)=\IndCoh_{Y_1}(Z_1)$$
as subcategories of $\IndCoh(Z_1)$.
\end{prop}

\begin{proof} Since $\QCoh(Z_2)$ is rigid 
and $\IndCoh_{Y_2}(Z_2)$ is dualizable,
the formation of
$$\IndCoh_{Y_2}(Z_2)\underset{\QCoh(Z_2)}\otimes -$$
commutes with limits (see \cite[Corollary 4.3.2 and 6.4.2]{DG}). 
Hence, the assertion is local on $Z_1$. 
Similarly, it is easy to see that the assertion is local on $Z_2$. 

\medskip

Hence, by \corref{c:shape of smooth}, we can assume that $f$ fits 
into a commutative diagram 
$$
\CD
Z_1  @>>>  \CU_1 \\
@V{f}VV   @VV{f_\CU}V  \\
Z_2  @>>>  \CU_2 \\
@VVV   @VVV   \\
\on{pt} @>>>  \CV,
\endCD
$$
where $\CU_1,\CU_2$, and $\CV$ are smooth affine schemes, $f_\CU$ is a smooth morphism, and
all squares are Cartesian.

\medskip

We can view 
$$\IndCoh(Z_2)\underset{\QCoh(Z_2)}\otimes \QCoh(Z_1)\simeq \IndCoh(Z_1)$$
as a category tensored over
$$\QCoh(\CU_1)\otimes \on{HC}(\on{pt}/\CV)^{\on{op}}\mod,$$
and both subcategories in the proposition correspond to the condition that the support
be contained in $Y_1\subset \CU_1\times V^*$ (see \propref{p:abs product}).
\end{proof}

\sssec{}

Let $Y_2$ and $Y_1$ be as above. From \propref{p:smooth pullback}, we obtain:

\begin{cor}    \label{c:pullback under smooth}
We have the following commutative diagrams:
$$\xymatrix{
\IndCoh_{Y_1}(Z_1)\ar[rr]<2pt>^{\hskip0.7cm \Xi^{Y_1,\on{all}}_{Z_1}} \ar[d] && 
\IndCoh(Z_1)  \ar[ll]<2pt>^{\hskip0.7cm \Psi^{Y_1,\on{all}}_{Z_1}} \ar[d]^{f^\IndCoh_*} \\
\IndCoh_{Y_2}(Z_2)\ar[rr]<2pt>^{\hskip0.7cm \Xi^{Y_2,\on{all}}_{Z_1}} && 
\IndCoh(Z_2)  \ar[ll]<2pt>^{\hskip0.7cm \Psi^{Y_2,\on{all}}_{Z_1}} \\}$$
and 
$$\xymatrix{
\IndCoh_{Y_1}(Z_1)\ar[rr]<2pt>^{\hskip0.7cm \Xi^{Y_1,\on{all}}_{Z_1}} && 
\IndCoh(Z_1)  \ar[ll]<2pt>^{\hskip0.7cm \Psi^{Y_1,\on{all}}_{Z_1}}  \\
\IndCoh_{Y_2}(Z_2)\ar[rr]<2pt>^{\hskip0.7cm \Xi^{Y_2,\on{all}}_{Z_1}} \ar[u] && 
\IndCoh(Z_2).  \ar[ll]<2pt>^{\hskip0.7cm \Psi^{Y_2,\on{all}}_{Z_1}} \ar[u]^{f^!} \\}$$
\end{cor}

\sssec{}
As another corollary of \propref{p:smooth pullback}, we obtain the following. Let $f:Z_1\to Z_2$
be as above, and assume moreover that it is surjective, i.e., $f$ is a smooth cover. 

\medskip

Let $Z^\bullet$ denote the \v{C}ech
nerve of $f$. Fix $Y_2\subset \Sing(Z_2)$, and for each $i$, let $Y^i\subset \Sing(Z^i)$ be the 
corresponding subset of $\Sing(Z^i)$. 

\medskip

We can form the cosimplicial category $\IndCoh_{Y^\bullet}(Z^\bullet)$ using either the !-pullback
or $(\IndCoh,*)$-pullback functors. In each case the resulting cosimplicial category is augmented
by $\IndCoh_{Y_2}(Z_2)$.

\begin{prop} \label{p:descent of condition}
Under the above circumstances the augmentation functor
$$\IndCoh_{Y_2}(Z_2)\to \on{Tot}\left(\IndCoh_{Y^\bullet}(Z^\bullet)\right)$$
is an equivalence.
\end{prop}

\begin{proof}
This follows from the fact that
$$\QCoh(Z_2)\to \on{Tot}\left(\QCoh(Z^\bullet)\right)$$
is an equivalence, combined with the fact that the operation 
$$\IndCoh_{Y_1}(Z_2)\underset{\QCoh(Z_2)}\otimes -$$
commutes with limits. 
\end{proof}

\begin{cor}  \label{c:sing supp via smooth cover}
For $\CF\in \IndCoh(Z_2)$, we have
$$\on{SingSupp}(\CF)\subset Y_2 \Leftrightarrow \on{SingSupp}(f^!(\CF))\subset Y_2\underset{Z_2}\times Z_1,$$
and also
$$\on{SingSupp}(\CF)\subset Y_2 \Leftrightarrow \on{SingSupp}(f^{\IndCoh,*}(\CF))\subset Y_2\underset{Z_2}\times Z_1.$$
\end{cor}

\begin{rem} From \thmref{t:ss of pullback} one can derive a more precise statement: if $f:Z_1\to Z_2$ is a smooth map
and $\CF\in\IndCoh(Z_2)$, then 
\[\on{SingSupp}(f^!(\CF))= \on{SingSupp}(f^{\IndCoh,*}(\CF))=\on{SingSupp}(\CF)\underset{Z_2}\times Z_1.\]
\end{rem} 

\ssec{Quasi-smooth morphisms, revisited}

In this subsection we will establish a generalization of \corref{c:closed embedding base change} for arbitrary 
quasi-smooth maps. We will treat the case of the !-pullback, while the $(\IndCoh,*)$-pullback is similar.

\sssec{}

Let $f:Z_1\to Z_2$ be a quasi-smooth morphisms between quasi-smooth DG schemes, and assume that $Z_2$
is quasi-compact. For a conical Zariski-closed $Y_2\subset \Sing(Z_2)$ let 
$$Y_1=\Sing(f)(Y_2\underset{Z_2}\times Z_1)\subset \Sing(Z_1),$$
where we regard $Y_2\underset{Z_2}\times Z_1$ as a subset of $\Sing(Z_2)_{Z_1}$.

\medskip

By \propref{p:singsupp functoriality}(a), we have a well-defined functor
$$f^!:\IndCoh_{Y_2}(Z_2)\to \IndCoh_{Y_1}(Z_1).$$
It extends by $\QCoh(Z_2)$-linearity to a functor
\begin{equation} \label{e:closed base change supports}
\QCoh(Z_1)\underset{\QCoh(Z_2)}\otimes \IndCoh_{Y_2}(Z_2)\to \IndCoh_{Y_1}(Z_1).
\end{equation}

\begin{cor}  \label{c:quasi-smooth tensor up}
The functor \eqref{e:closed base change supports} is an equivalence.
\end{cor}

\begin{proof}

As in the proof of \propref{p:smooth pullback}, the statement is local on both $Z_1$ and $Z_2$. By \lemref{l:factoring quasi-smooth}, 
locally, the map $f$ can be decomposed as a composition of a quasi-smooth closed embedding followed by a smooth
map. Now the assertion follows by combining \propref{p:smooth pullback} and \corref{c:closed embedding base change}.

\end{proof}

\sssec{} We can now generalize the results of \secref{ss:finite qs maps} as follows.

\begin{prop} \label{p:affine cons}
Suppose that $Z_1$ and $Z_2$ are quasi-smooth DG schemes and let $f:Z_1\to Z_2$ be an affine quasi-smooth morphism.
For a conical Zariski-closed $Y_2\subset \Sing(Z_2)$ let 
$$Y_1=\Sing(f)(Y_2\underset{Z_2}\times Z_1)\subset \Sing(Z_1),$$
where we regard $Y_2\underset{Z_2}\times Z_1$ as a subset of $\Sing(Z_2)_{Z_1}$.
Then \corref{c:closed embed cons} holds:

\smallskip

\noindent{\em(a)} The functor $f^\IndCoh_*:\IndCoh_{Y_1}(Z_1)\to \IndCoh_{Y_2}(Z_2)$
is conservative.

\smallskip

\noindent{\em(b)} The essential image of $\IndCoh_{Y_2}(Z_2)$ under $f^!$ (or under $f^{\IndCoh,*}$) generates
$\IndCoh_{Y_1}(Z_1)$.
\end{prop}

\begin{proof} As in the proof of \corref{c:closed embed cons}, the claim is local on $Z_2$, so we may assume that 
$Z_2$ (and, therefore, $Z_1$) is quasi-compact. Assertions (a) and (b) are equivalent, so it suffices to verify (b). But
it follows from \corref{c:quasi-smooth tensor up}, because the essential image of 
$\QCoh(Z_2)$ under 
$f^*$ generates $\QCoh(Z_1)$.
\end{proof}

\ssec{Inverse image}

Let us now consider the behavior of singular support under the operation of inverse image.
Let $f:Z_1\to Z_2$ be a morphism between quasi-smooth DG schemes ($Z_1$ and $Z_2$ need not be quasi-compact). 

\sssec{}
Let $\CF$ be an object of $\IndCoh(Z_2)$, and let $Y_2\subset \Sing(Z_2)$ be its singular support. 
Let $Y_1\subset \Sing(Z_1)$ be the Zariski closure of $\Sing(f)(p^{-1}(Y_2))$, where
\[p:\Sing(Z_2)_{Z_1}=\Sing(Z_2)\underset{Z_2}\times Z_1\to\Sing(Z_2)\]
is the projection. 

\medskip

Consider the object $f^!(\CF)\in \IndCoh(Z_1)$. Note that by 
\propref{p:singsupp functoriality}(a), we have:
$$\on{SingSupp}(f^!(\CF))\subset Y_1.$$

\begin{thm}  \label{t:ss of pullback}
Suppose that either $\CF\in\Coh(Z_2)$ or $f$ is a topologically open morphism (e.g. flat). Then \[\on{SingSupp}(f^!(\CF))=Y_1.\]
\end{thm}

\begin{rem}
The assertion of \thmref{t:ss of pullback} is not necessary for the main results of this paper.
\end{rem}

\begin{proof} As in \remref{r:pointwise}, consider the union
\[Y_2':=\underset{z_2\in Z_2}\bigcup
\on{supp}_{\Sym\left(H^1(T_{z_2}(Z_2))\right)}\left(H^\bullet(i^{\on{enh},!}_{z_2}(\CF))\right)\subset \Sing(Z_2).\]
By \propref{p:pointwise}, $Y_2=\overline{Y_2'}$; if $\CF\in\Coh(Z_2)$, then $Y_2=Y_2'$ by \propref{p:pointwise coh}. 
Similarly, consider the union
\[Y_1':=\underset{z_1\in Z_1}\bigcup\on{supp}_{\Sym\left(H^1(T_{z_1}(Z_1))\right)}
\left(H^\bullet(i^{\on{enh},!}_{z_1}(f^!(\CF)))\right)\subset \Sing(Z_1);\]
then $\on{SingSupp}(f^!(\CF))=\overline{Y_1'}$.

\medskip

Let us show that \[Y_1'=\Sing(f)(p^{-1}(Y_2')),\]
which would imply the assertion of the theorem. 

\medskip

Let $z_1\in Z_1$ be a point of $Z_1$, which we may assume to be a $k$-point after extending scalars. 
Set $z_2=f(z_2)\in Z_2$. Choose
a quasi-smooth map $$i_2:Z'_2\to Z_2$$ as in \secref{sss:dg point}, so that $Z'_2$ is a DG scheme of the form 
$\on{pt}\underset{\CV_2}\times \on{pt}$,
with $\CV_2$ smooth, and such that the unique $k$-point of $Z'_2$ goes to $z_2$.  

\medskip

Set $$Z'_1:=Z_1\underset{Z_2}\times Z'_2.$$ 
Since DG scheme $Z_1$ and the morphism $Z'_1\to Z_1$ are quasi-smooth, $Z_1'$ is quasi-smooth. Also, $z_1\in Z_1'$. We can therefore
choose a quasi-smooth map $$Z''_1\to Z_1'$$
as in \secref{sss:dg point}, so that $Z''_1$ is a DG scheme of the form $\on{pt}\underset{\CV_1}\times \on{pt}$,
with $\CV_1$ smooth, and such that the unique $k$-point of $Z''_1$ goes to $z_1$. Let $i_1$ be the composition $Z_1''\to Z_1'\to Z_1$.
Being a composition of quasi-smooth maps, $i_1$ is quasi-smooth.

\medskip

By \secref{sss:other enh}, we know that 
\[\on{SingSupp}(i_1^!(f^!(\CF)))=Y_1'\cap\Sing(Z_1)_{Z_1''}\subset \Sing(Z_1)_{Z_1''}=\Sing(Z_1)_{\{z_1\}}\subset\Sing(Z_1''),\]
and 
\[\on{SingSupp}(i_2^!(\CF))=Y_2'\cap\Sing(Z_2)_{Z_2'}\subset \Sing(Z_2)_{Z_2'}=\Sing(Z_2)_{\{z_2\}}\subset\Sing(Z_2').\]
Now we can replace $Z_1$, $Z_2$, and $\CF$ with 
$Z_1''$, $Z_2'$, and $i_2^!(\CF)$, respectively. Note that if $\CF$ is coherent, then so is $i_2^!(\CF)$.

\medskip

Thus, we assume that $^{cl}\!Z_1$ and $^{cl}\!Z_2$ are both isomorphic to $\on{pt}$. 
The claim now follows from \lemref{l:shape of point}(b) and 
\propref{p:support KD}(a).
\end{proof}

\ssec{Conservativeness for proper maps}

\sssec{} Suppose now that $f:Z_1\to Z_2$ is a proper morphism between quasi-smooth DG schemes.
Let $Y_1\subset \Sing(Z_1)$ be a conical Zariski-closed subset, and let $Y_2\subset \Sing(Z_2)$ be the image of $(\Sing(f))^{-1}(Y_1)$ 
under the projection
$$\Sing(Z_2)_{Z_1}\to\Sing(Z_2).$$
The subset $Y_2$ is automatically closed since the above map is proper.

\medskip

By \propref{p:singsupp functoriality}(b), the functor $f^\IndCoh_*$ sends $\IndCoh_{Y_1}(Z_1)$
to $\IndCoh_{Y_2}(Z_2)$. Our goal is to prove the following result.

\begin{thm} \label{t:prop cons}
Under the above circumstances, the essential image of $\IndCoh_{Y_1}(Z_1)$ under $f^\IndCoh_*$
generates $\IndCoh_{Y_2}(Z_2)$. 
\end{thm}

\medskip

We will derive \thmref{t:prop cons} from the following more general statement:

\begin{prop} \label{p:prop cons}
Let $f:Z_1\to Z_2$ be a (not necessarily proper) morphism of quasi-smooth DG schemes. Let $Y_1\subset\Sing(Z_1)$
and $Y_2\subset\Sing(Z_2)$ be conical Zariski-closed subsets. Suppose that $Y_2$ is contained in the image
of $\Sing(f)^{-1}(Y_1)$ under the projection \[\Sing(Z_2)_{Z_1}\to\Sing(Z_2).\] 
Suppose $\CF\in\IndCoh(Z_2)$ is such that $\Psi_{Z_2}^{Y_2,\on{all}}(\CF)\ne 0$.
Then \[\Psi_{Z_1}^{Y_1,\on{all}}\circ f^!(\CF)\ne 0.\]  
\end{prop}

\begin{proof}[Proof of \thmref{t:prop cons}]
The statement of the theorem is equivalent to the claim that the functor right adjoint to
\[f^\IndCoh_*:\IndCoh_{Y_1}(Z_1)\to \IndCoh_{Y_2}(Z_2)\]
is conservative. 

\medskip

The right adjoint 
in question is equal to the composition
\begin{equation} \label{e:right adjoint to f with supp}
\IndCoh_{Y_2}(Z_2)\overset{\Xi_{Z_2}^{Y_2,\on{all}}}\longrightarrow \IndCoh(Z_2)
\overset{f^!}\longrightarrow \IndCoh(Z_1)\overset{\Psi_{Z_1}^{Y_1,\on{all}}}\longrightarrow 
\IndCoh_{Y_1}(Z_1).
\end{equation}

\medskip

Suppose $\CF\in\IndCoh_{Y_2}(Z_2)$ is annihilated by the composition \eqref{e:right adjoint to f with supp}. But by 
\propref{p:prop cons}, the vanishing
\[\Psi_{Z_1}^{Y_1,\on{all}}\circ f^!\circ \Xi_{Z_2}^{Y_2,\on{all}}(\CF)=0\]
implies
\[0=\Psi_{Z_2}^{Y_2,\on{all}}\circ\Xi_{Z_2}^{Y_2,\on{all}}(\CF)=\CF,\]
as required.
\end{proof}

The rest of this subsection is devoted to the proof of \propref{p:prop cons}.

\sssec{Step 1}

We are going to reduce the statement of the proposition 
to the case when $Z_2$ is of the form 
$\on{pt}\underset{\CV_2}\times\on{pt}$ for a smooth scheme $\CV_2$ and a point $\on{pt}\hookrightarrow\CV_2$.

\medskip
  
Indeed, by \lemref{l:pointwise vanishing}, there exists a geometric point $z_2$ of $Z_2$ such that
such that $i_{z_2}^!(\Psi_{Z_2}^{Y_2,\on{all}}(\CF))\ne 0$. Extending the ground field,
we may assume that $i_{z_2}:\on{pt}\hookrightarrow Z_2$ is a rational point.

\medskip

As in \secref{sss:dg point}, we now extend the morphism $i_{z_2}$ to a quasi-smooth morphism of quasi-smooth 
DG schemes
\[i_2'=i_{z_2,\CV_2}:Z_2'\to Z_2,\]
where 
\[Z_2'=\on{pt}\underset{\CV_2}\times\on{pt}\]
for a smooth scheme $\CV_2$ with a marked point $\on{pt}\hookrightarrow\CV_2$.
Then 
\[(i_2')^!(\Psi_{Z_2}^{Y_2,\on{all}}(\CF))\ne 0\]
by \lemref{l:enhanced by V}. 

\medskip

Recall that the singular codifferential $\Sing(i_2')$ is an embedding 
\[\Sing(Z_2)_{\{z_2\}}\hookrightarrow V_2^*=\Sing(Z_2'),\]
where $V_2=T_{\on{pt}}(\CV_2)$.
Set \[Y_2'=\Sing(i_2')\left(Y_2\cap\Sing(Z_2)_{\{z_2\}}\right)\subset V_2^*.\]

\medskip 

Set now $Z_1'=Z_1\underset{Z_2}\times Z_2'$. The morphism $i_1':Z_1'\to Z_1$ is quasi-smooth, so by 
\lemref{l:codiff for q-smooth},
\[\Sing(i_1'):\Sing(Z_1)_{Z_1'}=\Sing(Z_1)\underset{Z_1}\times Z_1'\to\Sing(Z_1')\]
is a closed embedding. Set 
\[Y_1'=\Sing(i_1')\left(Y_1\underset{Z_1}\times Z_1'\right)\subset \Sing(Z_1').\]

\medskip

From \propref{p:functoriality quotient}(b), we obtain a commutative diagram of functors
\[
\CD
\IndCoh_{Y_2}(Z_2) @<{\Psi_{Z_2}^{Y_2,\on{all}}}<< \IndCoh(Z_2) @>{f^!}>> \IndCoh(Z_1) 
@>{\Psi_{Z_1}^{Y_1,\on{all}}}>> \IndCoh_{Y_1}(Z_1)\\
@V{(i_2')^!}VV @V{(i_2')^!}VV  @V{{i_1'}^!}VV @V{{i_1'}^!}VV\\
\IndCoh_{Y'_2}(Z'_2) @<{\Psi_{Z'_2}^{Y'_2,\on{all}}}<< \IndCoh(Z_2') @>{(f')^!}>> \IndCoh(Z_1')
@>{\Psi_{Z'_1}^{Y'_1,\on{all}}}>> \IndCoh_{Y'_1}(Z'_1),
\endCD
\]
where $f':Z'_1\to Z'_2$ is the natural morphism. Hence, it suffices to show that
$$\Psi_{Z'_1}^{Y'_1,\on{all}}\circ (f')^!\circ (i'_2)^!(\CF)\neq 0.$$

\medskip

Note that $f'$ satisfies the conditions of \propref{p:prop cons} with respect $Y'_2\subset \Sing(Z'_2)$
and $Y'_1\subset \Sing(Z'_1)$. 

\medskip

Thus, we obtain that the statement of proposition is reduced to the case when $Z_2$ is replaced by $Z'_2$, 
$Y_2$ with $Y_2'$, $\CF$ by $(i_2')^!(\CF)$, $Z_1$ by $Z_1'$,
and $Y_1$ by $Y_1'$. 

\medskip

In other words, we can assume that $Z_2=\on{pt}\underset{\CV_2}\times\on{pt}$, as desired.

\sssec{Step 2}

We are now going to reduce the assertion of the proposition to the case when 
$Z_1$ is also of the form $\on{pt}\underset{\CV_1}\times\on{pt}$. 

\medskip

To do so, let us fix a parallelization of the 
formal neighborhood of $\on{pt}$ in $\CV_2$. 
As explained in \secref{sss:global intersection parallelized}, this equips $\IndCoh(Z_2)$ 
with an action of the monoidal category $\QCoh(V_2^*/\BG_m)$.  
By \lemref{l:fiberwise vanishing}, there exists a geometric point $y_2\in V_2^*$ 
such that the fiber $i_{y_2}^*(\Psi_{Z_2}^{Y_2,\on{all}}(\CF))\ne 0$. 
By \corref{c:support by category fibers}(a), we see that $$y_2\in Y_2\subset V_2^*.$$
Extending the ground field,
we may assume that 
\[i_{y_2}:\on{pt}\hookrightarrow V_2^*=\Sing(Z_2)\] is a $k$-rational point.

\medskip

Since $Y_2$ is contained in the image of $\Sing(f)^{-1}(Y_1)$ under the projection 
$$\Sing(Z_2)_{Z_1}\to\Sing(Z_2),$$ there is a $k$-point $z_1\in Z_1$ such that $\Sing(f)$
sends \[(y_2,z_1)\in V_2^*\times{}^{cl}\!(Z_1)={}^{cl}\!\left(\Sing(Z_2)\underset{Z_2}\times Z_1\right)=\Sing(Z_2)_{Z_1}\]
into $Y_1\subset\Sing(Z_1)$.

\medskip

Now extend the morphism $i_{z_1}:\on{pt}\to Z_1$ to a quasi-smooth morphism of quasi-smooth 
DG schemes
\[\tilde i_1=i_{z_1,\CV_1}:\tilde Z_1\to Z_1,\]
where 
\[\tilde Z_1=\on{pt}\underset{\CV_1}\times\on{pt}\]
for a smooth scheme $\CV_1$ with a marked point $\on{pt}\hookrightarrow\CV_1$.
Set $\tilde f=f\circ\tilde i_1$.

\medskip

The singular codifferential $\Sing(\tilde i_1)$ is an embedding 
\[\Sing(Z_1)_{\{z_1\}}\hookrightarrow V_1^*=\Sing(\tilde Z_1).\]
Let $\tilde Y_1$ be the image
\[\Sing(\tilde i_1)(Y_1\cap\Sing(Z_1)_{\{z_1\}})\subset V_1^*.\]
The singular codifferential $\Sing(\tilde f)$ is a linear map $V_2^*\to V_1^*$.
Set   
\[\tilde Y_2:=Y_2\cap\Sing(\tilde f)^{-1}(\tilde Y_1)\subset V_2^*.\]

\medskip

By construction, $y_2\in\tilde Y_2$. By \corref{c:support by category fibers}(b),
\[i_{y_2}^*(\Psi_{Z_2}^{\tilde Y_2,\on{all}}(\CF))\simeq i_{y_2}^*(\CF)\simeq 
i_{y_2}^*(\Psi_{Z_2}^{Y_2,\on{all}}(\CF))\ne 0,\]
and hence
$$\Psi_{Z_2}^{\tilde Y_2,\on{all}}(\CF)\ne 0.$$

\medskip

Thus, it suffices to prove the assertion of the proposition after replacing 
$Y_2$ by $\tilde Y_2$, $Z_1$ by $\tilde Z_1$,
and $Y_1$ by $\tilde Y_1$, while keeping $\CF$ and $Z_2$ the same. 

\sssec{Step 3}

Thus, we can assume that $Z_i\simeq \on{pt}\underset{\CV_i}\times\on{pt}$
for $i=1,2$. In this case, the assertion of the proposition 
follows from \lemref{l:shape of point}(b) and \corref{c:support KD}(a).

\qed[\propref{p:prop cons}]

\section{Singular support on stacks}   \label{s:stacks}

In this section we develop the notion of singular support for objects of $\IndCoh(\CZ)$, where $\CZ$
is a quasi-smooth Artin stack. This will not be difficult, given the good functorial properties of $\IndCoh(-)$
on DG schemes under smooth maps. 

\medskip

Essentially, all this section amounts to is showing that for Artin stacks things work just as well as for DG schemes.
For this reason, this section, as well as \secref{s:gci stacks} may be skipped on the first pass.

\medskip

Recall that all schemes and stacks are assumed derived by default. To simplify the terminology, from now on we discard
the words ``differential graded" for stacks. Thus, ``Artin stack" stands for ``DG Artin stack."

\ssec{Quasi-smoothness for stacks}

\sssec{}  \label{sss:q-smooth via cotangent stacks}

Let $\CZ$ be an Artin stack (see \cite{Stacks}, Sect. 4). We say that it is quasi-smooth if for every affine DG scheme $Z$
equipped with a smooth map $Z\to \CZ$, the DG scheme $Z$ is quasi-smooth. 

\medskip

Equivalently, $\CZ$ is quasi-smooth if for some (equivalently, every) smooth atlas 
$f:Z\to \CZ$, the DG scheme $Z$ is quasi-smooth.

\medskip

Recall that for a $k$-Artin stack $\CZ$ its cotangent complex $T^*(\CZ)$ is an object of $\QCoh(\CZ)^{\leq k}$. We have:

\begin{lem} 
A $k$-Artin stack $\CZ$ is quasi-smooth if and only if $T^*(\CZ)$ is perfect of Tor-amplitude $[-1,k]$.
\end{lem}

\begin{proof} Let $f:Z\to\CZ$ be a smooth atlas. Then $T^*(\CZ)$ is perfect of Tor-amplitude $[-1,k]$ if and only if
$f^*(T^*(\CZ))$ has this property. Besides, we have an exact triangle
$$f^*(T^*(\CZ))\to T^*(Z)\to T^*(Z/\CZ),$$
where $T^*(Z/\CZ)$ is perfect of Tor-amplitude $[0,k-1]$. Thus, $f^*(T^*(\CZ))$ is perfect of Tor-amplitude $[-1,k]$ if and only if
$T^*(Z)$ is perfect of Tor-amplitude $[-1,k]$. The latter condition is equivalent to $Z$ being quasi-smooth.
\end{proof}

\medskip

We say that a map $f:\CZ_1\to \CZ_2$ between Artin stacks is quasi-smooth if $T^*(\CZ_1/\CZ_2)$ is perfect
of Tor-amplitude bounded from below by $-1$. 

\sssec{}

Recall the property of local eventual coconnectivity for a morphism between DG schemes
(see \secref{sss:loc ev coconn}). Clearly, this
property is local in the smooth topology on the source and on the target. Hence, it makes sense 
for morphisms between Artin stacks. 

\begin{lem} \label{l:q-smooth eventually coconn stacks} 
A quasi-smooth morphism $f:\CZ_1\to\CZ_2$ of Artin stacks is locally eventually coconnective. In 
particular, a quasi-smooth Artin stack is locally eventually coconnective.
\end{lem}
\begin{proof}
Follows from \corref{c:q-smooth eventualy coconn}.
\end{proof}

\sssec{}

If $\CZ$ is a quasi-smooth Artin stack, we introduce a classical Artin stack $\Sing(\CZ)$, equipped with
an affine (and, in particular, schematic) map to $^{cl}\CZ$. We call $\Sing(\CZ)$ the stack of singularities of $\CZ$. It is
constructed as follows: 

\medskip

For every affine DG scheme $Z$ with a smooth map
to $\CZ$ we set
$$\Sing(\CZ)\underset{\CZ}\times Z:=\Sing(Z),$$
and this assignment satisfies the descent conditions because of \lemref{l:codiff for smooth}. Equivalently,
$\Sing(\CZ)$ can be defined as 
\begin{equation}\label{e:stacky sing as spec}
^{cl}\!\left(\Spec_{\CZ}\left(\on{Sym}_{\CO_{\CZ}}(T({\CZ})[1])\right)\right).
\end{equation}

\sssec{The singular codifferential for stacks}

Let $f:\CZ_1\to \CZ_2$ be a map between quasi-smooth Artin stacks. We claim that we have a naturally
defined singular codifferential map 

\begin{equation} \label{e:codiff stacks}
\Sing(f):\Sing(\CZ_2)_{\CZ_1}:=\Sing(\CZ_2)\underset{\CZ_2}\times \CZ_1\to \Sing(\CZ_1).
\end{equation}

It can be obtained from the differential of $f$ using \eqref{e:stacky sing as spec}. 

\medskip

As in the case of DG schemes, it is easy to see that a map $f$ is quasi-smooth (resp., smooth)
if and only if the map $\Sing(f)$ is a closed embedding (resp., isomorphism).

\ssec{Definition of the category with supports for stacks}

\sssec{}

Recall (\cite[Sect. 11.2]{IndCoh}) that we have a well-defined category $\IndCoh(\CZ)$,
and that it can be recovered as
\begin{equation} \label{e:IndCoh via smooth}
\IndCoh(\CZ)\simeq \underset{Z\in \affdgSch_{/\CZ,\on{smooth}}}{lim} \IndCoh(Z),
\end{equation}
where $\affdgSch_{/\CZ,\on{smooth}}$ denotes the non-full subcategory of $(\affdgSch_{\on{aft}})_{/\CZ}$,
where the objects are restricted to pairs $(Z\in \affdgSch_{\on{aft}},f:Z\to \CZ)$ where $f$ is smooth, 
and where $1$-morphisms are restricted to maps $g:Z_1\to Z_2$ that are smooth as well.

\medskip

In the formation of the above limit we can use either the !-pullback functors or the $(\IndCoh,*)$-pullback functors,
as the two differ by the twist by a cohomologically shifted line bundle (this is due to the smoothness assumption on the morphisms).

\sssec{}

We let $\Coh(\CZ)\subset \IndCoh(\CZ)$ be the full (but not cocomplete) subcategory defined as
$$\Coh(\CZ)\simeq \underset{Z\in \affdgSch_{/\CZ,\on{smooth}}}{lim} \Coh(Z).$$

Note that we can think of $\Coh(\CZ)$ also as a full subcategory of $\QCoh(\CZ)$, where the latter,
according to \cite[Proposition 5.1.2]{QCoh}, is isomorphic to
$$\underset{Z\in \affdgSch_{/\CZ,\on{smooth}}}{lim} \QCoh(Z).$$

\sssec{}

Let $Y$ be a conical Zariski-closed subset of $\Sing(\CZ)$. We define the full subcategory
$$\IndCoh_Y(\CZ)\subset \IndCoh(\CZ)$$
as
\begin{equation} \label{e:IndCoh w/ support via smooth}
\IndCoh_Y(\CZ)\simeq \underset{Z\in \affdgSch_{/\CZ,\on{smooth}}}{lim} \IndCoh_{Y\underset{\CZ}\times Z}(Z),
\end{equation}
where we view $Y\underset{\CZ}\times Z$ as a closed subset of
$$\Sing(\CZ)\underset{\CZ}\times Z\simeq \Sing(Z).$$

From Lemmas \ref{l:tensored over qcoh} and \ref{l:commutes with action}, we obtain:

\begin{cor} \label{c:tensor stacks}
The action of $\QCoh(\CZ)$ on $\IndCoh(\CZ)$ preserves
$\IndCoh_Y(\CZ)$. 
\end{cor}  

\sssec{}

From \corref{c:pullback under smooth}, we obtain:

\begin{cor}  \label{c:colocalization stacks}
There exists a pair of adjoint functors
$$\Xi_\CZ^{Y,\on{all}}:\IndCoh_Y(\CZ)\rightleftarrows \IndCoh(Z):\Psi_\CZ^{Y,\on{all}},$$
with $\Xi_\CZ^{Y,\on{all}}$ being fully faithful. 
Moreover, for a smooth map $f:Z\to \CZ$, we have commutative diagrams
$$
\CD
\IndCoh_{Y\underset{\CZ}\times Z}(Z)  @>{\Xi_Z^{Y\underset{\CZ}\times Z,\on{all}}}>> \IndCoh(Z) \\
@AAA   @AA{f^!}A   \\
\IndCoh_Y(\CZ)   @>{\Xi_\CZ^{Y,\on{all}}}>> \IndCoh(\CZ)
\endCD
$$
and
$$
\CD
\IndCoh_{Y\underset{\CZ}\times Z}(Z)  @<{\Psi_Z^{Y\underset{\CZ}\times Z,\on{all}}}<< \IndCoh(Z) \\
@AAA   @AA{f^!}A   \\
\IndCoh_Y(\CZ)   @<{\Psi_\CZ^{Y,\on{all}}}<< \IndCoh(\CZ)
\endCD
$$
\end{cor}

\sssec{}

Recall from \cite[Sect. 11.7.3]{IndCoh} that for an eventually coconnective Artin stack,
we have a fully faithful functor
$$\Xi_\CZ:\QCoh(\CZ)\to \IndCoh(\CZ).$$

From \thmref{t:zero sect} we obtain:

\begin{cor}   \label{c:zero sing supp stacks}
If $Y$ is the zero-section, the subcategory
$$\IndCoh_{\{0\}}(\CZ)\subset \IndCoh(\CZ)$$
coincides with the essential image of $\QCoh(\CZ)$ under the functor
$$\Xi_\CZ:\QCoh(\CZ)\to \IndCoh(\CZ).$$
\end{cor}

\sssec{}

Let $\CV\hookrightarrow \CZ$ be a closed substack (not necessarily quasi-smooth), and 
let $j:\CU\hookrightarrow \CZ$ be the complementary open. 

\begin{cor}\label{c:on open stacks}
Let $Y\subset \Sing(\CZ)$ be a closed conical subset. Set 
\[Y_\CV={}^{cl}\!\left(Y\underset\CZ\times\CV\right)\subset\Sing(\CZ).\] 

\smallskip

\noindent{\em (a)} The subcategory
\[\IndCoh_Y(\CZ)\cap \IndCoh(\CZ_\CV)\subset \IndCoh(\CZ)\]
is equal to $\IndCoh_{Y_\CV}(\CZ)$.

\smallskip

\noindent{\em (b)} We have a short exact sequence of categories
$$\IndCoh_{Y_\CV}(\CZ)\rightleftarrows \IndCoh_Y(\CZ)\rightleftarrows \IndCoh_{Y\underset{\CZ}\times \CU}(\CU).$$

\end{cor}
\begin{proof} 
The two claims follow from Corollaries~\ref{c:supp on subscheme} and \ref{c:Cousin}
\end{proof}

\sssec{}

We have no reason to expect that the category $\IndCoh_Y(\CZ)$ is compactly generated for 
an arbitrary $\CZ$.

\medskip

Assume now that $\CZ$ is a QCA algebraic stack in the sense of \cite{DrGa0} Definition 1.1.8\footnote{QCA means quasi-compact,
and the automorphism group of any geometric point is affine.} (in particular, $\CZ$ is a $1$-Artin stack). 

\medskip

It is shown in {\it loc.cit.}, Theorem 3.3.4, that in 
this case the category $\IndCoh(\CZ)$ is compactly 
generated by $\Coh(\CZ)$. In particular, $\IndCoh(\CZ)$ is dualizable. 

\medskip

By \corref{c:colocalization stacks}, the category $\IndCoh_Y(\CZ)$ is a retract of $\IndCoh(\CZ)$.
Hence, by \cite[Lemma 4.3.3]{DrGa0}, we obtain:

\begin{cor}  \label{c:dualizable}
Under the above circumstances, the category $\IndCoh_Y(\CZ)$ is dualizable.
\end{cor}

\begin{rem}
We do not know whether under the assumptions of \corref{c:dualizable}, the category 
$\IndCoh_Y(\CZ)$ is compactly generated. In fact, we do not know this even for 
$Y=\{0\}$, i.e., we do not know whether $\QCoh(\CZ)$ is compactly generated.
We will describe two cases when this holds: one is proved in Appendix \ref{a:B} 
(when $\CZ=Z$ is a quasi-compact DG scheme) and the other in \secref{ss:expl presentation stacks}.
\end{rem}

\ssec{Smooth descent}

\sssec{Smooth maps of stacks}  \label{sss:smooth maps of stacks}

It follows from \lemref{l:codiff for q-smooth} that if $\CZ_1\to \CZ_2$ is smooth,
the singular codifferential
$$\Sing(f):\Sing(\CZ_2)_{\CZ_1}=\Sing(\CZ_2)\underset{\CZ_2}\times \CZ_1\to \Sing(\CZ_1)$$
is an isomorphism. 

\medskip

For a conical closed subset $Y_2\subset \Sing(\CZ_2)$, set
$$Y_1:=\Sing(f)\left(Y_2\underset{\CZ_2}\times \CZ_1\right)\subset\Sing(\CZ_1).$$

\begin{lem} \label{l:pullback under smooth stacks}
Let $f:\CZ_1\to \CZ_2$ be a smooth map between quasi-smooth Artin stacks. Then
we have the following commutative diagram
$$\xymatrix{
\IndCoh_{Y_1}(\CZ_1)\ar[rr]<2pt>^{\hskip0.7cm \Xi^{Y_1,\on{all}}_{\CZ_1}} && 
\IndCoh(\CZ_1)  \ar[ll]<2pt>^{\hskip0.7cm \Psi^{Y_1,\on{all}}_{\CZ_1}}  \\
\IndCoh_{Y_2}(\CZ_2)\ar[rr]<2pt>^{\hskip0.7cm \Xi^{Y_2,\on{all}}_{\CZ_1}} \ar[u] && 
\IndCoh(\CZ_2).  \ar[ll]<2pt>^{\hskip0.7cm \Psi^{Y_2,\on{all}}_{\CZ_1}} \ar[u]^{f^!} \\}$$
If in addition $f$ is quasi-compact and schematic, we also have a commutative diagram
$$\xymatrix{
\IndCoh_{Y_1}(\CZ_1)\ar[rr]<2pt>^{\hskip0.7cm \Xi^{Y_1,\on{all}}_{\CZ_1}} \ar[d] && 
\IndCoh(\CZ_1)  \ar[ll]<2pt>^{\hskip0.7cm \Psi^{Y_1,\on{all}}_{\CZ_1}} \ar[d]^{f^\IndCoh_*} \\
\IndCoh_{Y_2}(\CZ_2)\ar[rr]<2pt>^{\hskip0.7cm \Xi^{Y_2,\on{all}}_{\CZ_1}} && 
\IndCoh(\CZ_2)  \ar[ll]<2pt>^{\hskip0.7cm \Psi^{Y_2,\on{all}}_{\CZ_1}}. \\}$$
\end{lem}
\begin{proof} 
Both assertions follow formally from \corref{c:pullback under smooth}.
\end{proof}

\sssec{}

Let $f:\CZ_1\to \CZ_2$ be a smooth map between quasi-smooth Artin stacks. Let
$\CZ_1^\bullet$ denote its \v{C}ech nerve. Consider the co-simplicial DG category
$\IndCoh(\CZ_1^\bullet)$ formed by using either !- or $(\IndCoh,*)$-pullback functors,
augmented by $\IndCoh(\CZ_2)$.

\medskip

Let $Y_2\subset \Sing(\CZ_2)$ be a conical Zariski-closed subset. Set $Y_1^\bullet\subset \CZ_1^\bullet$
to be equal to
$$\CZ_1^\bullet\underset{\CZ_2}\times \Sing(\CZ_2).$$

According to \lemref{l:pullback under smooth stacks}, we have a well-defined full cosimplicial 
subcategory
$$\IndCoh_{Y_1^\bullet}(\CZ_1^\bullet)\subset \IndCoh(\CZ_1^\bullet),$$
augmented by $\IndCoh_{Y_2}(\CZ_2)$.

\begin{prop}  \label{p:Cech for stacks}
Suppose that $f$ is surjective on $k$-points. Then the augmentation functor
$$\IndCoh_{Y_2}(\CZ_2)\to \on{Tot}\left(\IndCoh_{Y_1^\bullet}(\CZ_1^\bullet)\right)$$
is an equivalence.
\end{prop}

\begin{proof} The statement formally follows from \propref{p:descent of condition}: smooth descent
for schemes implies smooth descent for stacks. To make the argument precise, we consider an auxiliary
category of ``inputs for $\IndCoh$". Namely, the objects of the category are pairs $(Z,Y)$, where $Z$ is a 
quasi-smooth affine DG scheme, and $Y\subset \Sing(Z)$ is a conical Zariski-closed subset. 
Morphisms $(Z_1,Y_1)$ to $(Z_2,Y_2)$ are smooth maps $Z_1\to Z_2$
whose singular codifferential induces an isomorphism $Y_2\underset{Z_2}\times Z_1\to Y_1$.
Denote this category by $\affdgSch_{\on{aft,q-smooth}+\on{supp}}$. 

\medskip

Let $\CZ'$ be a quasi-smooth Artin stack, and let $Y'\subset \Sing(\CZ')$ be a conical Zariski-closed subset.
This pair defines a presheaf $(\CZ',Y')$ on $\affdgSch_{\on{aft,q-smooth}+\on{supp}}$. Namely, for any 
$(Z,Y)\in\affdgSch_{\on{aft,q-smooth}+\on{supp}}$, the groupoid
$$\Maps((Z,Y),(\CZ',Y'))$$ is the full subgroupoid in $\Maps(Z,\CZ')$ consisting of maps $Z\to\CZ'$ that are smooth and 
whose singular
codifferential induces an isomorphism
$$Z\underset{\CZ'}\times Y'\to Y.$$

\medskip

The assignment $(Z,Y)\mapsto \IndCoh_Y(Z)$ is a functor
$$\IndCoh_{\on{supp}}:(\affdgSch_{\on{aft,q-smooth}+\on{supp}})^{\on{op}}\to \StinftyCat_{\on{cont}}.$$
It follows from \cite[Proposition~11.2.2]{IndCoh} that the category $\IndCoh_{Y'}(\CZ')$ identifies with the value on
$$(\CZ',Y')\in \on{PreShv}(\affdgSch_{\on{aft,q-smooth}+\on{supp}})$$
of the right Kan extension of $\IndCoh_{\on{supp}}$ along the Yoneda embedding
$$(\affdgSch_{\on{aft,q-smooth}+\on{supp}})^{\on{op}}\hookrightarrow (\on{PreShv}(\affdgSch_{\on{aft,q-smooth}+\on{supp}}))^{\on{op}}.$$

Now, \propref{p:descent of condition} says that the functor $\IndCoh_{\on{supp}}$ satisfies descent
with respect to surjective maps. This implies the assertion of the lemma by \cite[6.2.3.5]{Lu0}.

\end{proof}

\sssec{}  \label{sss:descent princ}

\propref{p:Cech for stacks} allows us to reduce statements concerning morphisms of Artin stacks $f:\CZ'\to \CZ$ 
to the case when $\CZ$ is a DG scheme. Such proofs proceed by induction along the hierarchy
$$\dgSch_{\on{aft}}\to \on{Alg.Spaces}\to \on{Stk}_{\on{1-Artin}}\to \on{Stk}_{\on{2-Artin}}\to\dots.$$

Namely, we choose an atlas $Z\to \CZ$ with $Z$ being a DG scheme that is locally almost of finite type.
Now if $\CZ$ is a $k$-Artin stack,
then the terms of the \v{C}ech nerve $Z^\bullet$ are $(k-1)$-Artin stacks.

\ssec{Functorial properties}

Let $\CZ_1$ and $\CZ_2$ be two quasi-smooth Artin stacks, and let $$f:\CZ_1\to \CZ_2$$ be a map.

\sssec{Functoriality under pullbacks}

Let $Y_i\subset \Sing(\CZ_i)$ be conical Zariski-closed subsets. 

\begin{lem} \label{l:pullback stacks}
Assume that the image of $Y_2\underset{\CZ_2}\times \CZ_1$ under the singular codifferential
\eqref{e:codiff stacks} is contained in $Y_1$. Then the functor $f^!$
sends $\IndCoh_{Y_2}(\CZ_2)$ to $\IndCoh_{Y_1}(\CZ_1)$.
\end{lem}

\begin{proof}
By \secref{sss:descent princ} we reduce the statement to the case when $\CZ_2$ is a DG scheme. In the latter
case, the statement from \propref{p:singsupp functoriality}(a).
\end{proof}

Similarly, we have:

\begin{lem} \label{l:pullback quotient stacks}
Assume that the preimage of $Y_1$ under the singular codifferential
\eqref{e:codiff stacks} is contained in 
$Y_2\underset{\CZ_2}\times \CZ_1$. Then the functor
$$\IndCoh(\CZ_2)\overset{f^!}\longrightarrow \IndCoh(\CZ_1)\overset{\Psi_{\CZ_1}^{Y_1,\on{all}}}\longrightarrow 
\IndCoh_{Y_1}(\CZ_1)$$
factors through the colocalization 
$$\IndCoh(\CZ_2)\overset{\Psi_{\CZ_2}^{Y_2,\on{all}}}\longrightarrow \IndCoh_{Y_2}(\CZ_2).$$
\end{lem}

\begin{proof}
Again, by \secref{sss:descent princ} we reduce the statement to the case when $\CZ_2$ is a DG scheme. In the latter
case, the statement from  \propref{p:functoriality quotient}(b).
\end{proof}

\sssec{Functoriality under pushforwards}

Let now $$f:\CZ_1\to \CZ_2$$ be schematic and quasi-compact. Recall (see \cite[Sect. 10.6]{IndCoh}) that in this case,
we have a well-defined functor
$$f_*^\IndCoh:\IndCoh(\CZ_1)\to \IndCoh(\CZ_2),$$
which satisfies a base-change property with respect to !-pullbacks for maps $\CZ'_2\to \CZ_2$.

\medskip

\begin{lem}  \label{l:direct image stacks}
Let $f:\CZ_1\to \CZ_2$ be schematic and quasi-compact. Assume that the preimage of $Y_1$ under the singular codifferential
\eqref{e:codiff stacks} is contained in 
$Y_2\underset{\CZ_2}\times \CZ_1$. Then the functor $f_*^\IndCoh$ sends
$\IndCoh_{Y_1}(\CZ_1)$ to $\IndCoh_{Y_2}(\CZ_2)$.
\end{lem}

\begin{proof}
Follows from \propref{p:singsupp functoriality}(b) by base change.
\end{proof}

Similarly, we have:

\begin{lem}  \label{l:direct image quotient stacks}
Let $f:\CZ_1\to \CZ_2$ be schematic and quasi-compact. Assume that the image of $Y_2\underset{\CZ_2}\times \CZ_1$ under the singular codifferential
\eqref{e:codiff stacks} is contained in $Y_1$. Then the functor 
$$\IndCoh(\CZ_1)\overset{f^\IndCoh_*}\longrightarrow \IndCoh(\CZ_2)\overset{\Psi_{\CZ_2}^{Y_2,\on{all}}}\longrightarrow 
\IndCoh_{Y_2}(\CZ_2)$$
factors through the colocalization
$$\IndCoh(\CZ_1)\overset{\Psi_{\CZ_1}^{Y_1,\on{all}}}\longrightarrow \IndCoh_{Y_1}(\CZ_1).$$
\end{lem}
\begin{proof}
Follows from \propref{p:functoriality quotient}(a) by base change.
\end{proof}

\sssec{Preservation of coherence}

We have the following generalization of \propref{p:preserve coherence}:

\begin{cor}  \hfill

\smallskip

\noindent{\em(a)} Let $\CF',\CF''\in \Coh(\CZ)$ be such that, set-theoretically, 
$$\on{SingSupp}(\CF')\cap \on{SingSupp}(\CF'')=\{0\}.$$
Then both
$$\CF'\otimes \CF'' \text{ and } \ul{\Hom}(\CF',\CF'')$$
belong to $\Coh(\CZ)$.

\smallskip

\noindent{\em(b)} Let $f:\CZ_1\to \CZ_2$ be a morphism and $\CF_2\in \Coh(\CZ_2)$ 
such that, set-theoretically, 
$$\left(\on{SingSupp}(\CF_2)\underset{\CZ_2}\times \CZ_1\right)\cap
\on{ker}\Bigl(\Sing(f):\Sing(\CZ_2)_{\CZ_1}\to \Sing(\CZ_1)\Bigr)\subset \{0\}\underset{\CZ_2}\times \CZ_1.$$
Then $f^!(\CF_2)\in \IndCoh(\CZ_1)$ belongs to $\Coh(\CZ_1)\subset \IndCoh(\CZ_1)$, and 
$f^*(\CF_2)\in \QCoh(\CZ_1)$ belongs to $\Coh(\CZ_1)\subset \QCoh(\CZ_1)$.

\end{cor}

\sssec{Conservativeness for finite maps}

Let now $f:\CZ_1\to \CZ_2$ be a finite (and, in particular, affine) map of quasi-smooth Artin stacks. Let
$Y_1\subset \Sing(\CZ_1)$ be a conical Zariski-closed subset contained in the image
of $\Sing(f)$. 

\medskip

From \corref{c:under finite} we obtain:

\begin{cor}
The functor $f^\IndCoh_*|_{\IndCoh_{Y_1}(\CZ_1)}:\IndCoh_{Y_1}(\CZ_1)\to \IndCoh(\CZ_2)$
is conservative. 
\end{cor}

\sssec{Conservativeness for quasi-smooth affine maps}
Let us prove an extension of \propref{p:affine cons}.
Suppose $\CZ_1$ and $\CZ_2$ are quasi-smooth Artin stacks, and $f:\CZ_1\to \CZ_2$ is a quasi-smooth 
affine map; in particular, $f$ is schematic and quasi-compact. 
As in the case of schemes,
the singular codifferential 
\[\Sing(f):\Sing(\CZ_2)_{\CZ_1}\to \Sing(\CZ_1)\]
is a closed embedding; this follows from \lemref{l:codiff for q-smooth} by base change.



\begin{prop} \label{p:affine cons stacks}
Let $Y_2\subset \Sing(\CZ_2)$ be a conical closed subset. 
Set \[Y_1=\Sing(f)\left(Y_2\underset{\CZ_2}\times\CZ_1\right)\subset\Sing(Z_1).\] 

\smallskip

\noindent{\em(a)} The essential image of $\IndCoh_{Y_2}(\CZ_2)$
under the functor $f^!$ generates $\IndCoh_{Y_1}(\CZ_1)$.

\smallskip

\noindent{\em(b)} The restriction of the functor $f_*^\IndCoh$
to $\IndCoh_{Y_1}(\CZ_1)$ is conservative.
\end{prop}

\begin{proof}
By \cite[Proposition 10.7.7]{IndCoh}, we have a pair of adjoint functors
\[f^\IndCoh_*:\IndCoh(\CZ_1)\rightleftarrows \IndCoh(\CZ_2):f^!.\]
From Lemmas~\ref{l:pullback stacks} and \ref{l:direct image stacks}, we see that they restrict to a pair
of functors
$$f^\IndCoh_*:\IndCoh_{Y_1}(\CZ_1)\rightleftarrows \IndCoh_{Y_2}(\CZ_2):f^!.$$
Moreover, since $f$ is locally eventually coconnective and Gorenstein, we have another pair of adjoint functors
\[f^{\IndCoh,*}:\IndCoh_{Y_2}(\CZ_2)\rightleftarrows\IndCoh_{Y_1}(\CZ_1):f^\IndCoh_*,\]
where $f^!$ differs from $f^{\IndCoh,*}$ by tensoring with the relative dualizing 
sheaf (see \cite[Proposition 7.3.8]{IndCoh}). Therefore, the two claims of the proposition are equivalent. 

\medskip

By \secref{sss:descent princ}, claim (b) is local in smooth topology on 
$\CZ_2$; hence we may assume that $\CZ_2$ is a DG scheme. 
This reduces the proposition to \propref{p:affine cons}.
\end{proof} 

\sssec{Quasi-smooth maps of stacks}

Let $f:\CZ_1\to \CZ_2$ be a quasi-smooth map of Artin stacks. Assume now that $\CZ_2$ is quasi-compact and
has an affine diagonal. In particular $\CZ_2$ is QCA, and by \cite[Corollary 4.3.8]{DrGa0}, the
category $\QCoh(\CZ_2)$ is rigid as a monoidal category. 

\medskip

Let $Y_2\subset \Sing(\CZ_2)$ and let
$$Y_1:=\Sing(f)\left(Y_2\underset{\CZ_2}\times \CZ_1\right)\subset \Sing(\CZ_1).$$

\begin{prop} \label{p:quasi-smooth pullback ten prod stacks}
Under the above circumstances, the functor
$$\IndCoh_{Y_2}(\CZ_2)\underset{\QCoh(\CZ_2)}\otimes \QCoh(\CZ_1)\to \IndCoh_{Y_1}(\CZ_1),$$
induced by the $\QCoh(\CZ_2)$-linear functor
$$f^!:\IndCoh(\CZ_2)\to \IndCoh(\CZ_1),$$
is an equivalence.
\end{prop}

\begin{proof}
By definition we have
$$\IndCoh_{Y_1}(\CZ_1)\simeq 
\underset{Z_1\in \affdgSch_{/\CZ,\on{smooth}}}{lim} \IndCoh_{Y_1\underset{\CZ_1}\times Z_1}(Z_1).$$
In addition, by \cite[Proposition 5.1.2(b)]{QCoh}
$$\QCoh(\CZ_1)\simeq \underset{Z_1\in \affdgSch_{/\CZ,\on{smooth}}}{lim} \QCoh(Z_1).$$

Since the category $\IndCoh_{Y_2}(\CZ_2)$ is dualizable, and $\QCoh(\CZ_2)$ is rigid, by 
\cite[Corollaries 4.3.2 and 6.4.2]{DG}, the formation of
$$\IndCoh_{Y_2}(\CZ_2)\underset{\QCoh(\CZ_2)}\otimes -$$
commutes with limits. This reduces the assertion of the proposition to the case when $\CZ_1=Z_1$ is an affine
DG scheme.

\medskip

Choose a smooth atlas $Z_2\to \CZ_2$, and let $Z_2^\bullet$ be its \v{C}ech nerve. Note that the assumption on $\CZ_2$
implies that the terms of $Z_2^\bullet$ are DG schemes (and not Artin stacks). 

\medskip

By \propref{p:Cech for stacks}, we obtain that
$\IndCoh_{Y_2}(\CZ_2)$ is the totalization of $\IndCoh_{Y_2^\bullet}(Z_2^\bullet)$,
where
$$Y_2^\bullet:=Y_2\underset{\CZ_2}\times Z_2^\bullet.$$

Since $\QCoh(Z_1)$ is dualizable and $\QCoh(\CZ_2)$ is rigid, we obtain that
$$\IndCoh_{Y_2}(\CZ_2)\underset{\QCoh(\CZ_2)}\otimes \QCoh(Z_1)$$
maps isomorphically to the totalization of 
\begin{equation} \label{e:totalization}
\IndCoh_{Y^\bullet_2}(Z^\bullet_2)\underset{\QCoh(\CZ_2)}\otimes \QCoh(Z_1).
\end{equation}

However,
$$\IndCoh_{Y^\bullet_2}(Z^\bullet_2)\underset{\QCoh(\CZ_2)}\otimes \QCoh(Z_1)
\simeq \IndCoh_{Y^\bullet_2}(Z^\bullet_2)\underset{\QCoh(Z^\bullet_2)}\otimes 
\QCoh(Z^\bullet_2)\underset{\QCoh(\CZ_2)}\otimes \QCoh(Z_1).$$

Now, we claim that the natural functor
$$\QCoh(Z^\bullet_2)\underset{\QCoh(\CZ_2)}\otimes \QCoh(Z_1)\to \QCoh(Z^\bullet_2\underset{\CZ_2}\times Z_1)$$
is an equivalence. This follows from \lemref{l:ten prod over stack} below.

\medskip

Thus, we obtain that the cosimplicial category \eqref{e:totalization}
identifies with
$$
\IndCoh_{Y_2\underset{\CZ_2}\times Z^\bullet_2}(Z^\bullet_2)
\underset{\QCoh(Z^\bullet_2)}\otimes \QCoh(Z^\bullet_2\underset{\CZ_2}\times Z_1)
\simeq \IndCoh_{Y_2\underset{\CZ_2}\times (Z^\bullet_2\underset{\CZ_2}\times Z_1)}(Z^\bullet_2\underset{\CZ_2}\times Z_1),$$
where the last isomorphism takes place due to \propref{p:smooth pullback}. 

\medskip

Now, $Z^\bullet_2\underset{\CZ_2}\times Z_1$ is the \v{C}ech nerve of the smooth cover
$Z_2\underset{\CZ_2}\times Z_1\to Z_1$, and by \propref{p:descent of condition},
the totalization of 
$$\IndCoh_{Y_2\underset{\CZ_2}\times (Z^\bullet_2\underset{\CZ_2}\times Z_1)}(Z^\bullet_2\underset{\CZ_2}\times Z_1)$$
is isomorphic to $\IndCoh_{Y_1}(Z_1)$, as required.

\end{proof}

\begin{lem} \label{l:ten prod over stack}
Let $\CZ$ be a quasi-compact stack with an affine diagonal. Then for any two prestacks
$\CZ_1$ and $\CZ_2$ mapping to $\CZ$, the naturally defined functor
$$\QCoh(\CZ_1)\underset{\QCoh(\CZ)}\otimes \QCoh(\CZ_2)\to \QCoh(\CZ_1\underset{\CZ}\times \CZ_2)$$
is an equivalence, provided that one of the categories $\QCoh(\CZ_1)$ or $\QCoh(\CZ_2)$ is dualizable.
\end{lem}

\begin{proof}

This follows by combining \cite[Proposition 3.3.3]{QCoh} and \cite[Corollary 4.3.8]{DrGa0}.

\end{proof}

\sssec{}

Let $\CZ$ be again a quasi-compact stack with an affine diagonal. Let $\CV\subset \CZ$ and $Y\subset \Sing(\CZ)$ be
as in \corref{c:on open stacks}.  

\medskip

In a way analogous to the proof of \propref{p:quasi-smooth pullback ten prod stacks} one shows:

\begin{prop}    \label{p:on open stacks}
Under the above circumstances, the short exact sequence of categories 
$$\IndCoh_{Y_\CV}(\CZ)\rightleftarrows \IndCoh_Y(\CZ)\rightleftarrows \IndCoh_{Y\underset{\CZ}\times \CU}(\CU)$$
is obtained from 
$$\QCoh(\CZ)_\CV\rightleftarrows \QCoh(\CZ)\overset{j^*}\rightleftarrows \QCoh(\CU)$$
by tensoring with $\IndCoh_Y(\CZ)$ over $\QCoh(\CZ)$.
\end{prop}

\sssec{Conservativeness for proper maps of stacks}
Suppose now that $f:\CZ_1\to\CZ_2$ is a schematic proper morphism between quasi-smooth Artin stacks.
Let $Y_1\subset \Sing(\CZ_1)$ be a conical closed subset, and let 
$Y_2$ be the image of 
\[(\Sing(f))^{-1}(Y_1)\subset \Sing(\CZ_2)_{\CZ_1}\] under the projection
$$\Sing(\CZ_2)_{\CZ_1}\to\Sing(\CZ_2).$$
Since the projection is proper, $Y_2\subset\Sing(\CZ_2)$ is a closed subset.

\medskip

By \lemref{l:direct image stacks}, the functor $f^\IndCoh_*$ sends $\IndCoh_{Y_1}(\CZ_1)$
to $\IndCoh_{Y_2}(\CZ_2)$. We have the following generalization of Theorem~\ref{t:prop cons}.

\begin{prop} \label{p:prop cons stacks}
Under the above circumstances, the essential image of $\IndCoh_{Y_1}(\CZ_1)$ under $f^\IndCoh_*$
generates $\IndCoh_{Y_2}(\CZ_2)$. 
\end{prop}

\begin{proof}
It is enough to verify that the claim is local on $\CZ_2$ in the smooth topology; one can then use \thmref{t:prop cons}. Indeed, 
the proposition is equivalent to the claim that the functor right adjoint to
\[f^\IndCoh_*:\IndCoh_{Y_1}(Z_1)\to \IndCoh_{Y_2}(Z_2)\]
is conservative. The right adjoint in question is the composition
\[
\IndCoh_{Y_2}(\CZ_2)\overset{\Xi_{\CZ_2}^{Y_2,\on{all}}}\longrightarrow \IndCoh(\CZ_2)
\overset{f^!}\longrightarrow \IndCoh(\CZ_1)\overset{\Psi_{\CZ_1}^{Y_1,\on{all}}}\longrightarrow 
\IndCoh_{Y_1}(\CZ_1).\]
The locality of this assertion follows from \secref{sss:descent princ}.
\end{proof}

\section{Global complete intersection stacks}  \label{s:gci stacks}

In this section, we adapt the approach of \secref{s:gci} to stacks. Our main objective
is to show that for a quasi-compact algebraic stack $\CZ$, globally given as a ``complete
intersection," and $Y\subset \Sing(\CZ)$, the corresponding category $\IndCoh_Y(\CZ)$
is compactly generated. The precise meaning of the words ``global complete intersection stack''
is explained in \secref{ss:expl presentation stacks}.

\medskip

As was mentioned earlier, this section may be skipped on the first pass.

\medskip

In this section all Artin stacks
will be quasi-compact with an affine diagonal (in particular, they all are QCA algebraic stacks in the sense of \cite{DrGa0}). 

\ssec{Relative Koszul duality}  \label{ss:rel Koszul}

\sssec{}   \label{sss:groupoid over stacks}

Let us consider a relative version of the setting of \secref{s:gci}. Let $\CX$ be a smooth stack, $\CV\to \CX$ a smooth schematic map, and let $\CX\to \CV$ be a section. 

\medskip

Consider the fiber product
\[\CG_{\CX/\CV}=\CX\underset{\CV}\times\CX.\]

As in \secref{sss:F&G}, it is naturally a group DG scheme over $\CX$. 

\sssec{}

The group structure on $\CG_{\CX/\CV}$ over $\CX$ turns $\IndCoh(\CG_{\CX/\CV})$ into a monoidal category over the symmetric monoidal category $\QCoh(\CX)$; the operation on $\IndCoh(\CG_{\CX/\CV})$ is the convolution. The unit object of $\IndCoh(\CG_{\CX/\CV})$
is \[(\Delta_{\CX})_*^\IndCoh(\omega_\CX)\in\IndCoh(\CG_{\CX/\CV}),\]
where $\omega_\CX$ is the dualizing complex on $\CX$.
Its endomorphisms naturally form an $\BE_2$-algebra in the symmetric monoidal category $\QCoh(\CX)$
(see \secref{ss:E2 and monoidal}). Denote this $\BE_2$-algebra by $\on{HC}(\CX/\CV)$. 

\medskip

Moreover, $(\Delta_{\CX})_*^\IndCoh(\omega_\CX)$ generates $\IndCoh(\CG_{\CX/\CV})$ over $\QCoh(\CX)$. 
Therefore, taking maps from $(\Delta_{\CX})_*^\IndCoh(\omega_\CX)$ defines
an equivalence of monoidal categories:
\[\on{KD}_{\CX/\CV}:\IndCoh(\CG_{\CX/\CV})\to \on{HC}(\CX/\CV)^{\on{op}}\mod.\]
This equivalence is the relative version of the Koszul duality \eqref{e:KD for V}.

\begin{lem}  \label{l:E2 mod c g}
The monoidal category $\on{HC}(\CX/\CV)^{\on{op}}\mod$ is rigid and compactly generated.
\end{lem}

\begin{proof}
First, note that $\QCoh(\CX)$ is compactly generated, that is, that $\CX$ is a perfect stack. Indeed, $\CX$ is smooth, so
it suffices to show that $\IndCoh(\CX)$ is compactly generated. The latter statement holds because 
$\CX$ is a QCA stack,
(see \cite[Theorem 0.4.5]{DrGa0}). 

\medskip

To show that $\on{HC}(\CX/\CV)^{\on{op}}\mod$ is rigid and compactly generated, we must show that it admits 
a family of compact dualizable generators
(see \cite[Lemma 5.1.1 and Proposition 5.2.3]{DG}). It is not hard to see 
that $\on{HC}(\CX/\CV)^{\on{op}}$-modules induced from perfect objects
of $\QCoh(\CX)$ form such a family. Let us prove the corresponding general statement:

\medskip

Let $\bO$ be a symmetric monoidal category that is compactly generated and rigid as a monoidal category. 
Then for any $\BE_2$-algebra $\CA$ in $\bO$, the monoidal category $\CA\mod(\bO)$ is rigid and compactly generated.

\medskip
Indeed, 
for $\bo\in \bO^c$, the object $\CA\otimes \bo$ is
compact in $\CA\mod(\bO)$. Clearly, such compact objects generate $\CA\mod(\bO)$.

\medskip

Since $\bO$ is rigid, its compact objects are dualizable. Hence, $\CA\otimes \bo$ is also
dualizable: its dual is $\CA\otimes \bo^\vee$. Thus, $\CA\mod(\bO)$ admits a family of compact dualizable generators, as required. 
\end{proof}

\sssec{}

As in \lemref{l:koszul as E_1}, we obtain that the $\BE_1$-algebra underlying $\on{HC}(\CX/\CV)$
is canonically isomorphic to $\Sym_{\CO_\CX}(V[-2])$, where $V$ is the pullback along $\CX\to \CV$
of the relative tangent sheaf to $\CV\to \CX$. 

\medskip

In particular, we have a canonical identification 
\begin{equation} \label{e:non-par mod E2 stack}
\on{HC}(\CX/\CV)^{\on{op}}\mod\simeq \Sym_{\CO_\CX}(V[-2])\mod
\end{equation}
as module categories over $\QCoh(\CX)$. 

\begin{rem}
A remark parallel to Remark \ref{r:non-monoidal} applies in the present situation.
\end{rem}

\sssec{} \label{sss:sing of Hecke stack}

Clearly, $\Sing(\CG_{\CX/\CV})\simeq V^*$, where $V^*$ denotes the total space of the corresponding
vector bundle over $\CX$.

\medskip

Let $Y\subset V^*$ be a conical Zariski-closed subset. Let us denote by
$$\on{HC}(\CX/\CV)^{\on{op}}\mod_Y\subset  \on{HC}(\CX/\CV)^{\on{op}}\mod\simeq \Sym_{\CO_\CX}(V[-2])\mod$$
the full subcategory of objects supported on $Y$. If $\CX$ is an affine scheme, it can be defined via the formalism of
\secref{sss:support via E2}; in general, we define it using an affine atlas $X\to\CX$. 

\begin{cor}\label{c:Koszul at point stacks} 
The functor $\on{KD}_{\CX/\CV}$ provides an equivalence between $\IndCoh_Y(\CG_{\CX/\CV})$ and
$\on{HC}(\CX/\CV)^{\on{op}}\mod_Y$.
\end{cor}
\begin{proof} The claim is local in the smooth topology on $\CX$. Therefore, we may assume that $\CX$ is affine,
and the assertion follows from the definition.
\end{proof}

\ssec{Explicit presentation of a quasi-smooth stack}  \label{ss:expl presentation stacks}

\sssec{} 

Let $\CZ$ be an Artin stack, and assume that we have a commutative diagram
\begin{gather}  \label{e:map of q-smooth w/ embeddings stacks}
\xy
(-10,20)*+{\CZ}="X";
(10,20)*+{\CU}="Y";
(-10,0)*+{\CX}="Z";
(10,0)*+{\CV}="W";
(0,-15)*+{\CX}="T";
{\ar@{->}^{\iota}"X";"Y"};
{\ar@{->}"X";"Z"};
{\ar@{->}"Y";"W"};
{\ar@{->}"Z";"W"};
{\ar@{->}_{\on{id}}"Z";"T"};
{\ar@{->}"W";"T"};
\endxy
\end{gather}
where the upper square is Cartesian, the lower portion of the diagram is as in \secref{ss:rel Koszul}, 
and $\CU$ is a smooth stack. Recall that $\CX$ is assumed to be smooth, and that $\CV$ is smooth and schematic
over $\CX$.

\medskip

In this situation, we say that $\CZ$ is presented as a global complete intersection stack. 
It is easy to see that such $\CZ$ is quasi-smooth. 

\sssec{}

We have a commutative diagram:

\begin{gather} \label{e:groupoids morphism stacks} 
\xy
(-25,0)*+{\CZ}="X";
(25,0)*+{\CZ}="Y";
(0,15)*+{\CG_{\CZ/\CU}}="Z";
(-25,-20)*+{\CX}="X_1";
(25,-20)*+{\CX,}="Y_1";
(0,-5)*+{\CG_{\CX/\CV}}="Z_1";
{\ar@{->}"Z";"X"};
{\ar@{->}"Z";"Y"};
{\ar@{->}"Z_1";"X_1"};
{\ar@{->}"Z_1";"Y_1"};
{\ar@{->}"Z";"Z_1"};
{\ar@{->}"X";"X_1"};
{\ar@{->}"Y";"Y_1"};
\endxy
\end{gather}
in which both parallelograms are Cartesian.

\medskip

In particular, as in \secref{sss:action on Z}, we obtain that 
the relative group DG scheme $\CG_{\CX/\CV}$ canonically acts on $\CZ$.

\sssec{}

We obtain homomorphisms of monoidal categories
\[\on{HC}(\CX/\CV)\mod\underset{\QCoh(\CX)}\otimes\QCoh(\CU)\simeq 
\IndCoh(\CG_{\CX/\CV})\underset{\QCoh(\CX)}\otimes\QCoh(\CU)\to\IndCoh(\CG_{\CZ/\CU}).\]
This allows us to view $\IndCoh(\CZ)$ as a category tensored over
\[\on{HC}(\CX/\CV)\mod\underset{\QCoh(\CX)}\otimes\QCoh(\CU).\]

\sssec{}

Let $Y\subset V^*\underset{\CX}\times\CU$ be a conical Zariski-closed subset. We can attach to it 
a full subcategory
\begin{equation} \label{e:abstract Hecke support}
\left(\on{HC}(\CX/\CV)\mod\underset{\QCoh(\CX)}\otimes \IndCoh(\CU)\right)_Y
\subset \on{HC}(\CX/\CV)\mod\underset{\QCoh(\CX)}\otimes \IndCoh(\CU)
\end{equation}
by interpreting
$$\on{HC}(\CX/\CV)\mod\underset{\QCoh(\CX)}\otimes \IndCoh(\CU)\simeq
\IndCoh(\CG_{\CX/\CV})\underset{\QCoh(\CX)}\otimes \QCoh(\CU)\simeq
\IndCoh(\CG_{\CX/\CV}\underset{\CX}\times \CU),$$
where the latter equivalence follows from \lemref{l:ten prod over stack} and \propref{p:quasi-smooth pullback ten prod stacks}
by
\begin{multline*}
\IndCoh(\CG_{\CX/\CV})\underset{\QCoh(\CX)}\otimes \QCoh(\CU)=
\IndCoh(\CG_{\CX/\CV})\underset{\QCoh(\CG_{\CX/\CV})}\otimes \QCoh(\CG_{\CX/\CV})
\underset{\QCoh(\CX)}\otimes \QCoh(\CU)\simeq \\
\simeq \IndCoh(\CG_{\CX/\CV})\underset{\QCoh(\CG_{\CX/\CV})} \otimes \QCoh(\CG_{\CX/\CV}\underset{\CX}\times \CU)
\simeq \IndCoh(\CG_{\CX/\CV}\underset{\CX}\times \CU).
\end{multline*}

Finally, we note that 
$$\Sing(\CG_{\CX/\CV}\underset{\CX}\times \CU)\simeq V^*\underset{\CX}\times\CU,$$
and we let the subcategory \eqref{e:abstract Hecke support} correspond to
$$\IndCoh_Y(\CG_{\CX/\CV}\underset{\CX}\times \CU)\subset \IndCoh(\CG_{\CX/\CV}\underset{\CX}\times \CU).$$

\sssec{}

We have a canonical closed embedding
\begin{equation} \label{e:embed Sing stacks}
\Sing(\CZ)\hookrightarrow V^*\underset{\CX}\times \CZ.
\end{equation}

The following assertion is parallel to \corref{c:sing supp via ten kosz}:

\begin{lem} \label{l:category w support as tensor}
For a conical Zariski-closed subset $Y\subset \Sing(\CZ)$, 
$$\IndCoh_Y(\CZ)\simeq 
\IndCoh(\CZ)\underset{\on{HC}(\CX/\CV)\mod\underset{\QCoh(\CX)}\otimes\QCoh(\CU)}\otimes
\left(\on{HC}(\CX/\CV)\mod\underset{\QCoh(\CX)}\otimes\QCoh(\CU)\right)_Y.$$
\end{lem}

\begin{proof}

First, \lemref{l:ten prod over stack} reduces the assertion to the case when $\CU$ is an affine DG scheme.

\medskip

Note that for any morphism $f:\CX_1\to \CX_2$ of prestacks and an associative algebra $\CA_2\in \QCoh(\CX_2)$,
the natural functor
$$\CA_2\mod\underset{\QCoh(\CX_2)}\otimes \QCoh(\CX_1)\to \CA_1\mod$$
is an equivalence (here $\CA_1:=f^*(\CA_2)$). This follows from \cite[Proposition 4.8.1]{DG}.

\medskip

This observation, combined with \lemref{l:ten prod over stack}, reduces the assertion to the case when $\CX$
is an affine DG scheme. In the latter case, the assertion follows from \corref{c:sing supp via ten kosz}.

\end{proof}
 
\begin{cor} \label{c:category comp gen stacks}
For any conical Zariski-closed subset $Y\subset \Sing(\CZ)$, 
the category $\IndCoh_Y(\CZ)$ is compactly generated.
\end{cor}

\begin{proof}

Since the monoidal category $\on{HC}(\CX/\CV)\mod\underset{\QCoh(\CX)}\otimes\QCoh(\CU)$ is rigid,
and $\IndCoh(\CZ)$ is compactly generated, it suffices to show that 
$$\left(\on{HC}(\CX/\CV)\mod\underset{\QCoh(\CX)}\otimes\QCoh(\CU)\right)_Y$$
is compactly generated. By \eqref{e:non-par mod E2 stack}, the latter is equivalent to 
$$\left(\Sym_{\CO_\CX}(V[-2])\mod \underset{\QCoh(\CX)}\otimes\QCoh(\CU)\right)_Y$$
being compactly generated, which in turn would follow from the compact generation of
$$\left(\Sym_{\CO_\CX}(V[-2])\mod \underset{\QCoh(\CX)}\otimes\QCoh(\CU)\right)^{\BG_m}_Y,$$
where $\BG_m$ acts on $V$ by dilations. 

\medskip

Let us now compare the categories $\Sym_{\CO_\CX}(V[-2])\mod$ and $\Sym_{\CO_\CX}(V)\mod$. While the two categories
are not equivalent, they differ only by a shift of grading; this implies that their categories of $\BG_m$-equivariant objects 
are naturally equivalent (we discuss the framework of the grading shift in \secref{sss:shift of grading alg}). The same applies 
to the categories 
\begin{align*}
\Sym_{\CO_\CX}(V[-2])&\mod \underset{\QCoh(\CX)}\otimes\QCoh(\CU)\\ \intertext{and}
\Sym_{\CO_\CX}(V)&\mod \underset{\QCoh(\CX)}\otimes\QCoh(\CU).
\end{align*} 

\medskip

Hence,
\begin{multline*}
\left(\Sym_{\CO_\CX}(V[-2])\mod \underset{\QCoh(\CX)}\otimes\QCoh(\CU)\right)^{\BG_m}_Y \simeq \\
\simeq \left(\Sym_{\CO_\CX}(V)\mod \underset{\QCoh(\CX)}\otimes\QCoh(\CU)\right)^{\BG_m}_Y
\simeq \QCoh((V^*/\BG_m)\underset{\CX}\times \CU)_{Y/\BG_m},
\end{multline*}
and the latter is easily seen to be compactly generated by 
$$\Coh((V^*/\BG_m)\underset{\CX}\times \CU)_{Y/\BG_m}=\QCoh((V^*/\BG_m)\underset{\CX}\times \CU)^{\on{perf}}\cap 
\QCoh((V^*/\BG_m)\underset{\CX}\times \CU)_{Y/\BG_m},$$
since the stack $(V^*/\BG_m)\underset{\CX}\times \CU$ is smooth. 

\end{proof}

\begin{cor}
Under the circumstances of \corref{c:category comp gen stacks} we have:
$$\IndCoh_Y(\CZ)\simeq \Ind(\Coh_Y(\CZ)),$$
where $\Coh_Y(\CZ):=\IndCoh_Y(\CZ)\cap \Coh(Z)$.
\end{cor}

\begin{proof}
By \corref{c:category comp gen stacks}, it suffices to show that
$$(\IndCoh_Y(\CZ))^c=\IndCoh_Y(\CZ)\cap \Coh(Z),$$
as subcategories of $\IndCoh_Y(\CZ)$. 

\medskip

However, this follows from the fact that the functor $\IndCoh_Y(\CZ)\hookrightarrow \IndCoh(\CZ)$
admits a continuous right adjoint and hence sends compacts to compacts, is fully faithful, and 
$$\Coh(\CZ)=\IndCoh(\CZ)^c$$ (the latter is \cite[Proposition 3.4.2(b)]{DrGa0}).
\end{proof}

\ssec{Parallelized situation} \label{ss:parallel stacks}

Assume now that in diagram \eqref{e:map of q-smooth w/ embeddings stacks}, the map $\CV\to \CX$
has been parallelized. That is, assume that $\CV$ is a vector bundle $V$ over $\CX$, and the section $\CX\to \CV$
is the zero-section.

\sssec{} \label{sss:universal} The diagram \eqref{e:map of q-smooth w/ embeddings stacks} can then be simplified, at least assuming
that the rank of the vector bundle $V$ is constant on $\CX$ (for instance, this is 
true if $\CX$ is connected). Indeed, suppose that $\on{rk}(V)=n$. By definition, the vector 
bundle $V$ on $\CX$ defines a morphism from $\CX$ into the classifying stack $\on{pt}/\on{GL}(n)$.
Clearly,
\[V=\CX\underset{\on{pt}/\on{GL}(n)}\times(\BA^n/\on{GL}(n)).\] 
Consider the composition \[\CU\to V\to(\BA^n/\on{GL}(n));\]
we then have
\[\CZ=\CU\underset{\BA^n/\on{GL}(n)}\times(\on{pt}/\on{GL}(n)),\]
where we embed $\on{pt}$ into $\BA^n$ as the origin. 
In other words, we may assume that $\CX=\on{pt}/\on{GL}(n)$ and $V=\BA^n/\on{GL}(n)$ in 
\eqref{e:map of q-smooth w/ embeddings stacks}.

\medskip

In more explicit terms, $\CU$ is equipped with a rank $n$ vector bundle and a section, and $\CZ$ is the zero 
locus of this section. That is, we may assume that $\CX=\CU$ in \eqref{e:map of q-smooth w/ embeddings stacks}, 
and that the composition $\CU\to V\to\CX$ is the identity map.

\sssec{}  \label{sss:tensored over conormal stacks parallel}

As in \lemref{l:parallel}, we obtain that $\on{HC}(\CX/\CV)$ is canonically isomorphic to 
$\Sym_{\CO_\CX}(V[-2])$ as an $\BE_2$-algebra. 

\medskip

In particular, we obtain that for $\CZ$ in \eqref{e:map of q-smooth w/ embeddings stacks}, the category
$\IndCoh(\CZ)$ is tensored over the monoidal category
\[\QCoh(V^*/\BG_m)\underset{\QCoh(\CX)}\otimes\QCoh(\CU)\simeq \QCoh((V^*/\BG_m)\underset{\CX}\times \CU).\]

Moreover, we have the following version of \lemref{l:category w support as tensor}: 

\begin{cor} \label{c:category as tensor product stacks}
For a conical Zariski-closed subset $Y\subset\Sing(\CZ)$, we have
\begin{equation} \label{e:category as tensor product stacks}
\IndCoh_Y(\CZ)=\IndCoh(\CZ)\underset{\QCoh((V^*/\BG_m)\underset{\CX}\times \CU)}\otimes 
\QCoh((V^*/\BG_m)\underset{\CX}\times \CU)_{Y/\BG_m}
\end{equation}
as full subcategories of $\IndCoh(\CZ)$.
\end{cor}

\begin{proof}
Follows from the fact that
$$\Vect\underset{\QCoh(\on{pt}/\BG_m)}\otimes \QCoh((V^*/\BG_m)\underset{\CX}\times \CU)\simeq 
\QCoh(V^*\underset{\CX}\times \CU)$$
as monoidal categories, and 
$$\Vect\underset{\QCoh(\on{pt}/\BG_m)}\otimes \QCoh((V^*/\BG_m)\underset{\CX}\times \CU)_{Y/\BG_m}\simeq 
\QCoh(V^*\underset{\CX}\times \CU)_Y$$
as modules over them.
\end{proof}

\ssec{Generating the category defined by singular support on a stack}

\sssec{}  \label{sss:Koszul generation stacks}

As in \secref{sss:F&G}, we have a tautologically defined functor 
$$\sG:\IndCoh(\CZ)\to \on{HC}(\CX/\CV)\mod\underset{\QCoh(\CX)}\otimes\IndCoh(\CZ) ,$$
and its left adjoint
$$\sF:\on{HC}(\CX/\CV)\mod\underset{\QCoh(\CX)}\otimes\IndCoh(\CZ) \to \IndCoh(\CZ).$$

\medskip

These functors are obtained as pullback and pushforward, respectively, for the action map
\[\on{act}_{\CG_{\CX/\CV},\CZ}:\CG_{\CX/\CV}\underset{\CX}\times\CZ\to\CZ.\]

\medskip

We have the following versions of Corollaries~\ref{c:Hecke action abstract} and \ref{c:Koszul generation}.

\begin{prop} \label{p:Hecke action abstract stacks}
For any conical Zariski-closed subset $Y\subset \CV^*\underset{\CX}\times \CU$, the functors
$\sF$ and $\sG$ restrict to a pair of adjoint functors
\[\sF:\left(\on{HC}(\CX/\CV)\mod\underset{\QCoh(\CX)}\otimes \IndCoh(\CZ)\right)_Y\rightleftarrows\IndCoh_{Y\cap\Sing(\CZ)}(\CZ):\sG.\]
Moreover, the diagram
$$
\xymatrix{
\on{HC}(\CX/\CV)\mod\underset{\QCoh(\CX)}\otimes \IndCoh(\CZ) \ar[r] \ar[d] & \IndCoh(\CZ) \ar@<.7ex>[l] \ar[d]^{\Psi_\CZ^{Y,\on{all}}} \\
\left(\on{HC}(\CX/\CV)\mod\underset{\QCoh(\CX)}\otimes \IndCoh(\CZ)\right)_Y  \ar[r] & \IndCoh_{Y\cap \Sing(\CZ)}(\CZ). \ar@<.7ex>[l] }$$
commutes.
Here the left vertical arrow is the right adjoint to the inclusion
$$\left(\on{HC}(\CX/\CV)\mod\underset{\QCoh(\CX)}\otimes \IndCoh(\CZ)\right)_Y
\hookrightarrow \on{HC}(\CX/\CV)\mod\underset{\QCoh(\CX)}\otimes \IndCoh(\CZ).$$
\end{prop}

\begin{proof}

Reduces to the case of DG schemes as in the proof of \lemref{l:category w support as tensor}.

\end{proof}

\begin{cor}\label{c:Koszul generation stacks} Suppose $Y$ is a conical Zariski-closed subset of $\Sing(\CZ)\subset\CV^*\underset\CX\times\CU$. 

\smallskip

\noindent{\em (a)} For any $\CF\in\IndCoh(\CZ)$, we have:
\[\CF\in \IndCoh_Y(\CZ)\,\Leftrightarrow\, \sG(\CF)\in\left(\on{HC}(\CX/\CV)\mod\underset{\QCoh(\CX)}\otimes \IndCoh(\CZ)\right)_{Y}.\]

\smallskip

\noindent{\em (b)} The essential image under $\sF$ of the category
$\left(\on{HC}(\CX/\CV)\mod\underset{\QCoh(\CX)}\otimes \IndCoh(\CZ)\right)_{Y}$ generates $\IndCoh_Y(\CZ)$.
\end{cor}
\begin{proof} The corollary formally follows from the conservativeness of $\sG$, similarly to the proof of \corref{c:Koszul generation}.
Namely, $\CF\in\IndCoh_Y(\CZ)$ if and only if
the natural morphism $\Psi_\CZ^{Y,\on{all}}(\CF)\to\CF$ is an isomorphism. Since $\sG$ is conservative, this happens if and only if
the morphism $\sG(\Psi_\CZ^{Y,\on{all}}(\CF))\to\sG(\CF)$ is an isomorphism; by \propref{p:Hecke action abstract stacks}, 
this is equivalent to
\[\sG(\CF)\in\left(\on{HC}(\CX/\CV)\mod\underset{\QCoh(\CX)}\otimes \IndCoh(\CZ)\right)_{Y}.\] 
We have therefore proved part (a). Also, by \propref{p:Hecke action abstract stacks},
the restriction 
\[\sF:\left(\on{HC}(\CX/\CV)\mod\underset{\QCoh(\CX)}\otimes \IndCoh(\CZ)\right)_Y\to\IndCoh_{Y\cap\Sing(\CZ)}(\CZ)\] 
is left adjoint to a conservative functor; this proves part (b).
\end{proof}

\bigskip

\bigskip

\centerline{\bf Part III: The geometric Langlands conjecture.}

\section{The stack $\LocSys_G$: recollections}   \label{s:LocSys}

In this section $G$ is an arbitrary affine algebraic group. Given a DG scheme $X$, 
we will recall the construction of the stack $\LocSys_G(X)$ of $G$-local systems on $X$.
We will compute its tangent and cotangent complexes. When $X$ is a smooth and complete
curve, we will show that $\LocSys_G(X)$ is quasi-smooth and calculate the corresponding
classical stack $\Sing(\LocSys_G(X))$. We will also show that $\LocSys_G(X)$ can in fact 
be written as a ``global complete intersection" as in \secref{s:gci stacks}.

\medskip

This section may be skipped if the reader is willing to take the existence of the stack $\LocSys_G(X)$ 
and its basic properties on faith.

\ssec{Definition of $\LocSys_G$} \hfill

\medskip

\noindent As the stack $\LocSys_G(X)$ of local systems is in general 
an object of derived algebraic geometry, some extra care is required. 
In this subsection we give the relevant definitions. Since the discussion 
will be purely technical, the reader can skip this subsection and return to it
when necessary.

\medskip

For the duration of this subsection we remove the a priori assumption that all prestacks
are locally almost of finite type. 

\sssec{}

Let $X$ be an arbitrary DG scheme almost of finite type. We define the prestacks $\Bun_G(X)$
and $\LocSys_G(X)$ using the general framework of Appendix \ref{s:Weil}. 

\medskip

Namely, for
$S\in \affdgSch$, we set
$$\Maps(S,\Bun_G(X)):=\Maps(S\times X,\on{pt}/G)$$ and  
$$\Maps(S,\LocSys_G(X)):=\Maps(S\times X_\dr,\on{pt}/G),$$
respectively. Here $X_\dr$ denotes the de Rham prestack of $X$, see \cite[Sect. 1.1.1]{Crys}. 

\medskip

The natural projection $X\to X_\dr$ defines the forgetful map 
\begin{equation} \label{e:proj to Bun}
\LocSys_G(X)\to \Bun_G(X).
\end{equation}

\sssec{}

The following is tautological:

\begin{lem}  \label{l:Bun G cl}
Suppose that $X$ is classical. Then the classical prestack $^{cl}\!\Bun_G(X)$ is the \emph{usual}
prestack of $G$-bundles on $X$, i.e., for $S\in \affSch$,
$$\Maps(S,\Bun_G(X))=\Maps_{^{cl}\!\on{PreStk}}(S\times X,\on{pt}/G).$$
\end{lem}

\begin{rem}
Note that when $X$ is \emph{not} classical, the classical prestack $^{cl}\!\Bun_G(X)$ cannot be recovered
from classical algebraic geometry. For instance, take $X$ to be a DG scheme with $^{cl}\!X=\on{pt}$. Then the
$\infty$-groupoid $\Maps(\on{pt},\Bun_G(X))$ is that of $G$-bundles on $X$, which can can have 
non-zero homotopy groups up to degree $n+1$ if $X$ is $n$-coconnective.
\end{rem}

\sssec{}

Similarly, we have:
\begin{lem}  \label{l:Loc G cl}
Let $X$ be arbitrary (but still almost of finite type). Then the classical prestack $^{cl}\!\LocSys_G(X)$ is the \emph{usual}
prestack of $G$-local systems on $X$, i.e., for $S\in \affSch$,
$$\Maps(S,\LocSys_G(X))=\Maps_{^{cl}\!\on{PreStk}}(S\times X_\dr,\on{pt}/G).$$
\end{lem}

\begin{proof}
Follows from the fact that $X_\dr\in \on{PreStk}$ is classical, see \cite[Proposition 1.3.3(b)]{Crys}.
\end{proof}

\begin{rem}
For a more familiar description of $^{cl}\!\LocSys_G(X)$ via the Tannakian formalism, see 
\secref{sss:Tannaka 3}.
\end{rem}

\ssec{A Tannakian description}  

One can describe the $\infty$-groupoids 
$$\Maps(S\times X,\on{pt}/G) \text{ and } \Maps(S\times X_\dr,\on{pt}/G)$$
in more intuitive terms using Tannakian duality.

\medskip

The material in this subsection will not be used elsewhere in the paper. 

\sssec{}  \label{sss:Tannaka 1}

Let $\CX$ be any prestack (which we will take to be either $S\times X$ or $S\times X_\dr$). Then
by \cite[Theorem 3.4.2]{Lu3}, the $\infty$-groupoid $\Maps(\CX,\on{pt}/G)$ identifies
with the full subcategory of symmetric monoidal functors
$$\Rep(G):=\QCoh(\on{pt}/G)\overset{\Phi}\to \QCoh(\CX)$$ 
that satisfy:
\begin{itemize}

\item $\Phi$ is continuous;

\item $\Phi$ is right t-exact, i.e., sends $\Rep(G)^{\leq 0}$ to $\QCoh(\CX)^{\leq 0}$
(see \cite[Sect. 1.2.3]{QCoh} for the definition of the t-structure on $\QCoh$ over an
arbitrary prestack);

\item $\Phi$ sends \emph{flat} objects to \emph{flat} objects (an object of $\QCoh$
on a prestack is said to be flat if its pullback to an arbitrary affine DG scheme is
flat.\footnote{If $A$ is a connective ring, an $A$-module $M$ is flat if $M$ is connective,
and $H^0(A)\underset{A}\otimes M$ is acyclic off degree $0$ and is flat over $H^0(A)$.})

\end{itemize}

\begin{rem}  \label{r:Tannaka 1}
The above description of $\Maps(\CX,-)$ is valid for $\on{pt}/G$ replaced by any \emph{geometric stack}
(see \cite[Definition 3.4.1]{Lu3}). Note that any quasi-compact algebraic stack with
an affine diagonal is geometric.
\end{rem}

\sssec{} \label{sss:Tannaka 2}

Since $\Rep(G)$ identifies with the derived category of its heart, 
the above category of symmetric monoidal functors can be identified with that of symmetric monoidal functors
$$\Rep(G)^{\heartsuit,c}=\Rep(G)^\heartsuit\cap \Rep(G)^c \overset{\Phi}\to \QCoh(\CX)^{\leq 0},$$
such that:

\begin{itemize}

\item $\Phi$ takes short exact sequences in $\Rep(G)^{\heartsuit,c}$ to distinguished triangles in $\QCoh(\CX)^{\leq 0}$.

\end{itemize}

\begin{rem}  \label{r:Tannaka 1.5}
Note that every object of $\Rep(G)^{\heartsuit,c}$ is dualizable. Therefore, the same holds for its image under 
such $\Phi$. From this, it is easy to see that the essential image of $\Rep(G)^{\heartsuit,c}$ 
automatically belongs to the full subcategory $\QCoh(\CX)^{\on{loc.free}}$
of $\QCoh(\CX)$ spanned by locally free sheaves of finite rank (i.e., those objects whose pullback to any affine DG 
scheme $S$ is a direct summand of $\CO_S^{\oplus n}$ for some integer $n$).
\end{rem}


\sssec{} \label{sss:Tannaka 3}

Assume for a moment that $X$ is classical. We obtain that \secref{sss:Tannaka 2} recovers the usual description of
the classical prestack $^{cl}\!\Bun_G(X)$. Namely, for $S\in \affSch$, the groupoid $\Maps(S,\Bun_G(X))$
is that of exact symmetric monoidal functors
$$\Phi:\Rep(G)^{\heartsuit,c}\to \QCoh(S\times X)^\heartsuit.$$
Note that by Remark \ref{r:Tannaka 1.5}, the essential image of such $\Phi$ consists of vector bundles on $S\times X$.

\medskip

Similarly, for an arbitrary $X$, we obtain the usual description of 
the classical prestack $^{cl}\!\LocSys_G(X)$. Namely, for $S\in \affSch$, the groupoid $\Maps(S,\LocSys_G(X))$
is that of exact symmetric monoidal functors
$$\Phi:\Rep(G)^{\heartsuit,c} \to \QCoh(S\times X_\dr)^\heartsuit.$$
By Remark \ref{r:Tannaka 1.5}, the essential image of such $\Phi$ belongs to the subcategory of $\QCoh(S\times X_\dr)$
that consists of $S$-families of local systems on $X$.

\begin{rem}  \label{r:Tannaka 3}
To obtain the above Tannakian description of $^{cl}\!\Bun_G(X)$ (for $X$ classical) and
$^{cl}\!\LocSys_G(X)$ (for $X$ arbitrary) one needs something weaker than the full strength of
\cite[Theorem 3.4.2]{Lu3}. Namely, one can make do with its classical version, given by 
\cite[Theorem 5.11]{Lu4}. 
\end{rem}

\ssec{Basic properties of $\Bun_G$ and $\LocSys_G$}

In this subsection we will calculate the (pro)-cotangent spaces of $\Bun_G(X)$ and $\LocSys_G(X)$.
In addition, we will show that if $X$ is classical and proper, the forgetful map $\LocSys_G(X)\to \Bun_G(X)$
is schematic and affine.

\medskip

Throughout this subsection we will assume that $X$ is eventually coconnective.  

\sssec{}   \label{sss:cotangent to LocSys}

Let $X$ be $l$-coconnected. It follows from 
\propref{p:maps deform} that $\Bun_G(X)$ admits a $(-l-1)$-connective deformation
theory. 

\medskip

Similarly, since $X_\dr$ is classical, we obtain that $\LocSys_G(X)$ admits a $(-1)$-connective 
deformation theory.

\sssec{}

The pro-cotangent spaces of $\Bun_G(X)$ and $\LocSys_G(X)$ can be described as follows:

\medskip

The stack $\on{pt}/G$ admits co-representable $(-1)$-connective deformation theory, and its
cotangent complex identifies with $\fg^*_{\CP^{univ}}[-1]$, where $\CP^{univ}$ is the universal
$G$-bundle on $\on{pt}/G$, and $\fg^*_{\CP^{univ}}$ is the vector bundle on $\on{pt}/G$ associated with the
coadjoint representation. 

\medskip

Let $\CP$ (resp., $(\CP,\nabla)$) be an $S$-point of $\Bun_G(X)$ (resp., $\LocSys_G(X)$). Then
by \propref{p:maps deform}(b), the pro-cotangent space to $\Bun_G(X)$ (resp., $\LocSys_G(X)$) at the 
above point, viewed as a functor
$$\QCoh(S)^{\leq 0}\to \inftygroup,$$
identifies with
\begin{equation} \label{e:cotangent Bun}
\CM\mapsto \Gamma(S\times X,\CM\otimes \fg_\CP)[1]
\end{equation}
and 
\begin{equation} \label{e:cotangent Loc}
\CM\mapsto \Gamma(S\times X_\dr,\CM\otimes \fg_\CP)[1],
\end{equation}
respectively, where $\fg_\CP$ denotes the bundle associated with the adjoint representation. 

\medskip

The relative pro-cotangent space to the map \eqref{e:proj to Bun} at the above point is the functor
\begin{equation} \label{e:rel cotangent LocSys}
\CM\mapsto \on{ker}\Bigl(\Gamma(S\times X_\dr,\CM\otimes \fg_\CP)[1]\to \Gamma(S\times X,\CM\otimes \fg_\CP)[1]\Bigr).
\end{equation}

\medskip

All of the above functors commute with colimits. 

\sssec{}

Note that by \corref{c:maps deform}, if $X$ is proper, the pro-cotangent spaces to $\Bun_G(X)$ and $\LocSys_G(X)$
are co-representable by objects of $\QCoh(S)^-$. 

\medskip

In other words, $\Bun_G(X)$ (resp., $\LocSys_G(X)$)
admits a \emph{co-representable} $(-l-1)$-connective (resp., $(-1)$-connective) deformation theory.

\medskip 

For $S\to \Bun_G(X)$ (resp., $S\to \LocSys_G(X)$), 
we will denote the resulting cotangent spaces, viewed as objects of $\QCoh(S)$, by  

\begin{equation} \label{e:cotangent spaces pr}
T^*(\Bun_G(X))|_S,\,\,\,T^*(\LocSys_G(X))|_S, \text{ and }  
T^*(\LocSys_G(X)/\Bun_G(X))|_S,
\end{equation}
respectively. 

\sssec{}

We claim:

\begin{lem}  \label{l:schematicity}
Assume that $X$ is classical. Then the relative pro-cotangent spaces of $\LocSys_G(X)\to \Bun_G(X)$
are connective.
\end{lem}

\begin{proof}

We need to show that the functor \eqref{e:rel cotangent LocSys}, viewed as a functor
$\QCoh(S)\to \Vect$, is left t-exact. 

\medskip

By \cite[Proposition 3.4.3]{Crys}, the object $\Gamma(S\times X_\dr,\CM\otimes \fg_\CP)\in \Vect$
can be calculated as the totalization of the co-simplicial object of $\Vect$ whose $n$-simplices are
$$\Gamma\left(S\times (X^n)^\wedge_X,\CM\otimes \fg_\CP|_{S\times (X^n)^\wedge_X}\right),$$
where $(X^n)^\wedge_X$ is the DG indscheme equal to the formal completion of $X^i$ along
the main diagonal.

\medskip

In particular, the projection onto the $0$-simplices is the canonical map
$$\Gamma(S\times X_\dr,\CM\otimes \fg_\CP)\to \Gamma(S\times X,\CM\otimes \fg_\CP).$$

\medskip

Hence, it suffices to show that for every $n$ and $\CM\in \QCoh(S)^{\geq 0}$, we have
\begin{equation} \label{e:above 0}
\Gamma\left(S\times (X^n)^\wedge_X,\CM\otimes \fg_\CP|_{S\times (X^n)^\wedge_X}\right)\in \Vect^{\geq 0}.
\end{equation}

The key observation is that the assumption that $X$ be classical implies that 
the DG indscheme $(X^n)^\wedge_X$ is \emph{classical}, see \cite[Proposition 6.8.2]{IndSch}.
That is, it can be written as a colimit of classical schemes $Z_\alpha$.  

\medskip

Hence,
$$\Gamma\left(S\times (X^n)^\wedge_X,\CM\otimes \fg_\CP|_{S\times (X^n)^\wedge_X}\right)\simeq
\underset{\alpha}{lim}\, \Gamma\left(S\times Z_\alpha,\CM\otimes \fg_\CP|_{S\times Z_\alpha}\right).$$

Now \eqref{e:above 0} follows from the fact that, for each $\alpha$, 
$$\Gamma\left(S\times Z_\alpha,\CM\otimes \fg_\CP|_{S\times Z_\alpha}\right) \in \Vect^{\geq 0}.$$

\end{proof}

\begin{cor}
For a proper (classical) scheme $X$ and any $S\to \LocSys_G(X)$, the object 
$T^*(\LocSys_G(X)/\Bun_G(X))|_S$ belongs to $\QCoh(S)^{\leq 0}$.
\end{cor}

\sssec{}

Finally, we claim:

\begin{prop} \label{p:schematicity}  \hfill

\smallskip

\noindent{\em(a)} If $X$ is a classical scheme, 
the map \eqref{e:proj to Bun} is ind-schematic and in fact ind-affine. 

\smallskip

\noindent{\em(a')} If $X$ is a proper classical scheme,  
the map \eqref{e:proj to Bun} is schematic and affine. 

\smallskip

\noindent{\em(b)} The prestacks $\Bun_G(X)$ and $\LocSys_G(X)$ are locally almost of finite type.

\end{prop}

\begin{proof}

The fact that $\Bun_G(X)$ is locally almost of finite type is a particular case of \corref{c:maps laft}.

\medskip

Assume now that $X$ is classical. Let $S$ be an affine DG scheme almost of finite type equipped with a map to $\Bun_G(X)$. 
We will prove that
\begin{equation} \label{e:conn on bundle}
S\underset{\Bun_G(X)}\times \LocSys_G(X)
\end{equation}
is an ind-affine DG indscheme locally almost of finite type, and is in fact an affine DG scheme if $X$ proper.
This implies (a) and (a'). It also implies (b): to show that $\LocSys_G(X)$ is locally almost
of finite type we can replace the initial $X$ by $^{cl}\!X$. 

\medskip

The fact that
$$^{cl}(S\underset{\Bun_G(X)}\times \LocSys_G(X))$$
is an ind-affine indscheme (resp., affine scheme for $X$ proper) and that it is locally of finite type follows from 
Lemmas \ref{l:Bun G cl} and \ref{l:Loc G cl}, since the corresponding assertions in classical algebraic 
geometry are well known. \footnote{In the case when $X$ is smooth, this is obvious, using the description
of local systems as bundles with a connection. For a general $X$, it is enough to consider the case of $G=GL_n$. Then the ``locally of 
finite type" assertion is a general property of the category of D-modules.  The (ind)-representability can be proved 
using the infinitesimal groupoid as in \lemref{l:schematicity}.}

\medskip

To prove that \eqref{e:conn on bundle} is a DG indscheme/DG scheme we use \thmref{t:repr}. 
Indeed, the required condition on the pro-cotangent spaces follows from \lemref{l:schematicity}.

\medskip

The fact that \eqref{e:conn on bundle} is locally almost of finite type follows from \lemref{l:laft crit}.

\end{proof}

\ssec{The case of curves}

From now on let us assume that $X$ is a smooth, complete and connected curve. In what follows we will omit $X$
from the notation in $\Bun_G(X)$ and $\LocSys_G(X)$, unless an ambiguity is likely to occur.

\medskip

We will show that $\Bun_G$ is a smooth classical stack. We will also show that $\LocSys_G$ is quasi-smooth,
and compute the corresponding classical stack $\on{Sing}(\LocSys_G)$. 

\sssec{}

First, we note that $\Bun_G$ is an algebraic stack (a.k.a. $1$-Artin stack in the terminology of \cite{Stacks}).
Indeed, the usual proof that $^{cl}\!\Bun_G$ is a classical algebraic stack (see, e.g., \cite{Av})
applies in the context of derived algebraic geometry to show the corresponding property of $\Bun_G$.
\footnote{Another way to see this is to choose sufficiently deep level structure
(over every fixed quasi-compact open substack in $\Bun_G$) and apply \thmref{t:repr}.}

\medskip

We claim: 

\begin{lem}
The stack $\Bun_G$ is smooth (and, in particular, classical).
\end{lem}

\begin{proof}

By \lemref{l:smooth sch}, this follows from
the fact that $\Bun_G$ is an Artin stack locally almost of finite type and 
from the description of the cotangent spaces to $\Bun_G$ given by \eqref{e:cotangent Bun}. 

\end{proof}

\sssec{}

From the fact that $\Bun_G$ is an algebraic stack and the fact that the map \eqref{e:proj to Bun} is schematic
(see \propref{p:schematicity}), we obtain that $\LocSys_G$ is also an algebraic stack. Since $\Bun_G$ has an affine 
diagonal, we obtain that the same is true for $\LocSys_G$, since the map \eqref{e:proj to Bun} is separated (in fact, it is affine).

\medskip

Moreover, it is easy to see that the image of $\LocSys_G$ in $\Bun_G$ is contained in a quasi-compact open substack
of $\Bun_G$. \footnote{Here is a sketch of the proof. It is enough to consider the case of $G=GL_n$. Now, if a rank $n$-bundle $\CE$
splits as a direct sums $\CE=\CE_1\oplus \CE_2$, then a connection on $\CE$ gives rise to connections on $\CE_i$.
However, it follows from the Riemann-Roch Theorem 
that every rank $n$-bundle outside a certain quasi-compact open substack 
of $\Bun_n$ admits a direct sum decomposition as above with either $\on{deg}(\CE_1)\neq 0$ or 
$\on{deg}(\CE_2)\neq 0$.}
This implies that $\LocSys_G$ itself is quasi-compact. Thus, $\LocSys_G$ is a QCA stack in the terminology
of \cite{DrGa0}. 

\medskip

Again, from \propref{p:schematicity} we obtain that $\LocSys_G$ is locally almost of finite type. 

\medskip

But, of course, $\LocSys_G$ is not smooth. 

\sssec{}

We now claim:

\begin{prop} \label{p:LocSys q-smooth}
The stack $\LocSys_G$ is quasi-smooth.
\end{prop}

\begin{proof}
This follows immediately from the description of the cotangent spaces given by \eqref{e:cotangent Loc}.
Namely, for any $S$-point of $\LocSys_G$, the object $T^*(\LocSys_G)|_S\in \QCoh(S)$ is given, by Verdier duality
along $X$, by 
\begin{equation}  \label{e:exp for cotangent}
\Gamma(S\times X_\dr,\fg^*_\CP)[1],
\end{equation}
which lives in cohomological degrees $\geq -1$, as required. 
\end{proof}

Note that by the same token we obtain a description of the tangent complex of $\LocSys_G$: for an $S$-point of
$\LocSys_G$, the object $T(\LocSys_G)|_S\in \QCoh(S)$ is given by
\begin{equation}  \label{e:exp for tangent}
\Gamma(S\times X_\dr,\fg_\CP)[1].
\end{equation} 

\sssec{The stack $\Sing(\LocSys_G)$} The above description of the tangent complex of $\LocSys_G$ implies the following description of
the \emph{classical} stack  $\Sing(\LocSys_G)$:

\begin{cor} \label{c:Sing of LocSys}
The stack $\Sing(\LocSys_G)$ admits the following description: given a classical affine scheme
$S\in \affSch$, 
$$\on{Maps}(S,\Sing(\LocSys_G))=(\CP,\nabla,A),$$
where $(\CP,\nabla)\in \on{Maps}(S,\LocSys_G)$, and $A$ is an element of
$$H^0\left(\Gamma(S\times X_\dr,\fg^*_\CP)\right).$$
\end{cor}

\begin{proof}

By definition, $\Sing(\LocSys_G)$ is the classical stack underlying 
$$\Spec_{\LocSys_G}\left(\Sym_{\CO_{\LocSys_G}}\left(T(\LocSys_G[1])\right)\right).$$
The assertion of the corollary follows from \eqref{e:exp for tangent}, since for $(\CP,\nabla)$
as in the corollary, by the Serre duality,
$$\Hom_{\Coh(S)}\left(T(\LocSys_G[1])_S,\CO_S\right)\simeq H^0\left(\Gamma(S\times X_\dr,\fg^*_\CP)\right).$$

\end{proof}

\sssec{} We denote the stack $\Sing(\LocSys_G)$ by $\Arth_G$. 

\begin{rem}
As was explained in the introduction, if $G$ is reductive,
the Arthur parameters for the automorphic side are supposed to 
correspond to points $(\CP,\nabla,A)\in\Arth_G$ where $A$ is nilpotent; see \secref{s:Langlands} for details.
\end{rem}

\ssec{Accessing $\LocSys_G$ via an affine cover}  \label{ss:LocSys via cover}

In this subsection we will show how to make sense in the DG world of such operations
as adding a $1$-form to a connection, and taking the polar part of a meromorphic 
connection. 

\sssec{}

Let $U\subset X$ be a non-empty open affine subset. Consider the prestack
$$\LocSys_G(X;U):=\LocSys_G(U)\underset{\Bun_G(U)}\times \Bun_G(X).$$

\medskip

It is clear that if $X=U_1\cup U_2$, we have
\begin{equation} \label{e:LocSys via open}
\LocSys_G:=\LocSys_G(X)\simeq \LocSys_G(X;U_1)\underset{\LocSys_G(X;U_{1,2})}\times \LocSys_G(X;U_2),
\end{equation}
where $U_{1,2}=U_1\cap U_2$.

\medskip

The description of $\LocSys_G$ via \eqref{e:LocSys via open} will be handy for establishing certain of its properties.

\sssec{}

The main observation is: 

\begin{prop} \label{p:over open classical}
The prestack $\LocSys_G(X;U)$ is classical.
\end{prop}

\begin{proof}

Since $\Bun_G$ is a smooth classical algebraic stack, it suffices to show that for any smooth classical affine 
scheme $S$ and any map $S\to\Bun_G$,
the fiber product
\begin{equation} \label{e:Cart prod over Bun_G}
S\underset{\Bun_G}\times \LocSys_G(X;U)
\end{equation}
is classical. 

\medskip

First, we claim that the DG indscheme \eqref{e:Cart prod over Bun_G} is formally smooth.
For that it suffices to show that $\LocSys_G(X;U)$ is  formally smooth \emph{over} $\Bun_G$. 
By \cite[Proposition 8.2.2]{IndSch}, this is equivalent to showing that for any $S'\in \affdgSch_{\on{aft}}$
with a map to $\LocSys_G(X;U)$ and any $\CM\in \QCoh(S')^{<0}$, we have
$$\pi_0\left(T^*(\LocSys_G(X;U)/\Bun_G)|_{S'}(\CM)\right)=0,$$
where $T^*(\LocSys_G(X;U)/\Bun_G)|_S$ is viewed as a functor
$$\QCoh(S')^{\leq 0}\to \inftygroup.$$

\medskip

By \eqref{e:rel cotangent LocSys}, we have:
\begin{multline} \label{e:cotangent to all}
T^*(\LocSys_G(X;U)/\Bun_G)|_{S'}(\CM)\simeq \\
\simeq \on{Cone}\Bigl(\Gamma(S'\times U_\dr,\CM\otimes \fg_\CP)\to \Gamma(S'\times U,\CM\otimes \fg_\CP)\Bigr)
\simeq \Gamma(U,\CM\otimes \fg_\CP\otimes \omega_X),
\end{multline}
and the required vanishing follows from the fact that $U$ is affine. \footnote{The fact that \eqref{e:Cart prod over Bun_G}
is formally smooth already implies that its classical via \cite[Theorem 9.1.6]{IndSch}.  Below we give a more explicit
proof which avoids (the non-trivial) \cite[Theorem 9.1.6]{IndSch} and instead uses (the more elementary) 
\cite[Proposition 9.1.4]{IndSch}.}

\medskip

Consider now the classical indscheme
$$^{cl}(S\underset{\Bun_G}\times \LocSys_G(X;U)).$$

I.e., this is the classical moduli problem corresponding to endowing a given $G$-bundle
(given by a point of $S$) with a connection defined over $U\subset X$. It is easy to
see that, locally on $S$, we have an isomorphism
$$^{cl}(S\underset{\Bun_G}\times \LocSys_G(X;U))\simeq S\times \BA^\infty,$$
where
$$\BA^\infty=\underset{n\in \BZ^{\geq 0}}{colim}\, \BA^n.$$

In particular, $^{cl}(S\underset{\Bun_G}\times \LocSys_G(X;U))$, \emph{when viewed as a DG indscheme},
is formally smooth. 

\medskip

The assertion of the proposition follows now from \cite[Proposition 9.1.4]{IndSch}. 
\footnote{The assertion if \cite[Proposition 9.1.4]{IndSch} says that if $\CZ$ is a formally
smooth DG indscheme, such that the underlying classical indscheme $^{cl}\CZ$ is formally
smooth \emph{when viewed as a derived indscheme}, then the canonical map
$^{cl}\CZ\to \CZ$ is an isomorphism, up to sheafification.}

\end{proof}

\sssec{}

As a first application of \propref{p:over open classical}, we will prove the following.
Consider the following group DG scheme over $\Bun_G$, which we denote by $\on{Hitch}_G$:

\medskip

For an $S$-point $\CP$ of $\Bun_G$, the $\infty$-groupoid of its lifts to an $S$-point of $\on{Hitch}_G$
is by definition
$$\tau^{\leq 0}\Bigl(\Gamma(S\times X,\fg_\CP\otimes \omega_X)\Bigr).$$

By Serre duality, we have
$$S\underset{\Bun_G}\times \on{Hitch}_G=
\Spec\Bigl(\Sym_{\CO_S}\left(\Gamma(S\times X,\fg^*_\CP)[1]\right)\Bigr).$$
Note that $\on{Hitch}_G$ is naturally a DG vector bundle\footnote{By a DG vector bundle over a prestack $\CZ$ we mean
a prestack of the form $\Spec(\Sym_{\CO_\CZ}(\CF))$ for $\CF\in \QCoh(\CZ)^{\leq 0}$.}
(and therefore a DG group scheme) over $\Bun_G$.

\begin{cor}
There exists a canonical action of $\on{Hitch}_G$ on $\LocSys(G)$ over $\Bun_G$; the action is simply transitive in the 
sense that the induced map
\[\on{Hitch}_G\underset{\Bun_G}\times\LocSys(G)\to\LocSys(G)\underset{\Bun_G}\times\LocSys(G)\]
is an isomorphism.
\end{cor}

\begin{rem}
This corollary is a triviality for the underlying classical stacks: any two connections 
on a given bundle over a curve differ by a $1$-form. However, it is less obvious
at the derived level, since the procedure of adding a $1$-form to a connection is 
difficult to make sense of in the $\infty$-categorical setting.
\end{rem}

\begin{proof}

As in the case of $\LocSys_G(X;U)$, we can define a relative indscheme, $\on{Hitch}(X;U)$,
over $\Bun_G$, whose
$S$-points are pairs $(\CP,\alpha)$, where $\CP$ is an $S$-point of $\Bun_G$ and 
$\alpha$ is a point of
$$\Gamma(S\times U,\fg_\CP\otimes \omega_X),$$
considered as an $\infty$-groupoid.  As in the case of $\LocSys_G(X;U)$, we show that $\on{Hitch}(X;U)$
is classical. Similarly,
$$\on{Hitch}(X;U)\underset{\Bun_G}\times \LocSys_G(X;U)$$
is classical. 

\medskip

Since we are dealing with classical objects, it is easy to see that $\on{Hitch}(X;U)$ acts simply transitively on
$\LocSys_G(X;U)$ over $\Bun_G$. Moreover, these actions are compatible under restrictions for $U\hookrightarrow U'$.

\medskip

Covering $X=U_1\cup U_2$, we have
$$\on{Hitch}_G\simeq \on{Hitch}(X;U_1)\underset{\on{Hitch}(X;U_{1,2})}\times \on{Hitch}(X;U_2),$$
as prestacks. Now, the required assertion follows from \eqref{e:LocSys via open}.

\end{proof}

\sssec{}

Let now $x$ be a $k$-point of $X$ outside of $U$. Consider the following relative DG indscheme over
$\Bun_G(X)$, denoted $\on{Polar}(G,x)$, whose $S$-points are pairs $(\CP,A_{\on{Polar}})$, where 
$\CP$ is an $S$-point of $\Bun_G(X)$, and $A_{\on{Polar}}$ is a point of
$$\Gamma(S\times X,\fg_\CP\otimes \omega_X(\infty\cdot x)/\omega_X),$$
considered as an $\infty$-groupoid via $\Vect^{\leq 0}\to \inftygroup$. 

\medskip

One easily shows that $\on{Polar}(G,x)$ is formally smooth
and classical as a prestack: locally in the fppf topology on $\Bun_G$, the stack $\on{Polar}(G,x)$ looks
like the product of $\Bun_G$ and $\BA^\infty$. 

\sssec{}

Note that since $\LocSys_G(X;U)$ and $\on{Polar}(G,x)$ are both classical prestacks, the usual
operation of taking the polar part of the connection defines a map
$$\LocSys_G(X;U)\to \on{Polar}(G,x).$$

\medskip

\begin{prop} \label{p:via polar part}
Set $U':=U\cup \{x\}$, and suppose that it is still affine.
Then there exists a canonical isomorphism
$$\LocSys_G(X;U')\simeq \LocSys_G(X;U)\underset{\on{Polar}(G,x)}\times \Bun_G(X),$$
where $\Bun_G(X)\to \on{Polar}(G,x)$ is the zero-section.
\end{prop}

\begin{proof}

Note that since $U'$ is affine and hence $\LocSys_G(X;U')$ is classical, there exists a canonically
defined map
$$\LocSys_G(X;U')\to \LocSys_G(X;U)\underset{\on{Polar}(G,x)}\times \Bun_G(X),$$
which is an isomorphism at the classical level. To show that this map is an isomorphism,
it is enough to show that it induces an isomorphism at the level of cotangent spaces
at $S$-points for a classical affine scheme $S\in \affSch$. The latter, in turn, follows from the computation of the
cotangent spaces in \secref{sss:cotangent to LocSys}.

\end{proof}

\sssec{}

Let now $U=X-x$. We claim:

\begin{cor}  \label{c:all poles}
There exists a canonical isomorphism
$$\LocSys_G(X)\simeq \LocSys_G(X;U)\underset{\on{Polar}(G,x)}\times \Bun_G(X).$$
\end{cor}

\begin{proof}
Let $U'$ be another open affine of $X$ that contains the point $x$. 
Applying \propref{p:via polar part}, we obtain:
$$\LocSys_G(X;U')\simeq \LocSys_G(X;U\cap U')\underset{\on{Polar}(G,x)}\times \Bun_G(X),$$
so 
\begin{multline*} 
\LocSys_G\simeq \LocSys_G(X;U)\underset{\LocSys_G(X;U\cap U')}\times \LocSys_G(X;U')\simeq \\
\simeq \LocSys_G(X;U)\underset{\LocSys_G(X;U\cap U')}\times \left(\LocSys_G(X;U\cap U')\underset{\on{Polar}(G,x)}\times 
\Bun_G(X)\right) \simeq \\
\simeq \LocSys_G(X;U)\underset{\on{Polar}(G,x)}\times \Bun_G(X),
\end{multline*}
as required.

\medskip

It is equally easy to see that the constructed map does not depend on the choice of $U'$: for $U''\subset U'$ the corresponding
diagram commutes. 

\end{proof}

\ssec{Presentation of $\LocSys_G$ as a fiber product} 

In this subsection we will show how to define the notion of connection that has a pole of order $\leq 1$
at a given point, and how to represent of $\LocSys_G$ as a fiber product of smooth stacks. 

\sssec{} Fix a point $x\in X$, and let
$$\on{Polar}^{\leq 1}(G,x)\subset \on{Polar}(G,x)$$
be the closed substack corresponding to 
$$\fg_\CP\otimes \omega_X(x)/\omega_X\subset \fg_\CP\otimes \omega_X(\infty\cdot x)/\omega_X.$$
That is, this substack corresponds to pairs $(\CP,A_{\on{Polar}})$ where $A_{\on{Polar}}$ has at most a simple pole.

\medskip

It is easy to see that we have a canonical identification (the residue map)
$$\on{Polar}^{\leq 1}(G,x)\simeq \fg/G\underset{\on{pt}/G}\times \Bun_G,$$
where $\Bun_G\to \on{pt}/G$ is the canonical map corresponding to the restriction of a $G$-bundle to $x\in X$.

\sssec{}

We define the stack $\LocSys_G^{\on{R.S.}}$ of local systems with (at most) a simple pole at $x$ by
$$\LocSys_G^{\on{R.S.}}:=\LocSys_G(X;X-x)\underset{\on{Polar}(G,x)}\times \on{Polar}^{\leq 1}(G,x).$$

\medskip

By \corref{c:all poles}, we have a canonical map 
$$\iota:\LocSys_G\hookrightarrow \LocSys_G^{\on{R.S.}}$$
and a canonical map
$$\on{res}:\LocSys_G^{\on{R.S.}}\to \fg/G\underset{\on{pt}/G}\times \Bun_G$$
that fit into a Cartesian square
\begin{equation} \label{e:residue diagram}
\CD
\LocSys_G @>>> \LocSys_G^{\on{R.S.}}  \\
@VVV    @VV{\on{res}}V    \\
\Bun_G   @>>>  \fg/G\underset{\on{pt}/G}\times \Bun_G,
\endCD
\end{equation}
where the bottom horizontal arrows comes from the zero-section map $\on{pt}/G\to \fg/G$.

\sssec{}  \label{sss:cotangent RS}

From \eqref{e:cotangent to all}, we obtain the following description of the relative cotangent spaces
of $\LocSys_G^{\on{R.S.}}$ over $\Bun_G$:

\medskip

For an $S$-point $(\CP,\nabla,A)$, the cotangent space
$T^*(\LocSys_G^{\on{R.S.}}/\Bun_G)|_S$, viewed as a functor
$$\QCoh(S)^{\leq 0}\to \inftygroup$$ is given by
$$\CM\mapsto \Gamma(S\times X,\CM\otimes \fg_\CP\otimes \omega_X(x)).$$

\medskip

In particular,
$$T^*(\LocSys_G^{\on{R.S.}}/\Bun_G)|_S\in \QCoh(S)^{\leq 0}.$$
In fact, by the Serre duality,
\begin{equation} \label{e:cotangent rel to RS}
T^*(\LocSys_G^{\on{R.S.}}/\Bun_G)|_S\simeq \Gamma(S\times X,\fg^*_\CP(-x))[1].
\end{equation}

\sssec{}

By \thmref{t:repr}, we obtain that the map
$$\LocSys_G^{\on{R.S.}}\to \Bun_G$$
is schematic. (The same argument applies to connections with poles of any fixed order instead of simple poles.)
This map is also easily seen to be separated (and, in fact, affine). This implies that $\LocSys_G^{\on{R.S.}}$
has an affine diagonal. 

\medskip

Note that, unlike $\LocSys_G$, the stack $\LocSys_G^{\on{R.S.}}$ is not quasi-compact (unless $G$ is unipotent).
However, for our applications the stack $\LocSys_G^{\on{R.S.}}$ may be replaced by a Zariski neighborhood of the image
$\iota(\LocSys_G)\subset\LocSys_G^{\on{R.S.}}$; we can choose such a neighborhood to be quasi-compact, and therefore QCA. 

\sssec{}

From \eqref{e:cotangent rel to RS}, we obtain that the map 
$$\LocSys_G^{\on{R.S.}}\to \Bun_G$$
is quasi-smooth. 

\medskip

Since $\Bun_G$ is smooth, we obtain that the stack $\LocSys_G^{\on{R.S.}}$ is quasi-smooth. We now claim:

\begin{prop} \label{p:RS smooth} \hfill

\smallskip

\noindent{\em(a)} The stack $\LocSys_G^{\on{R.S.}}$ is smooth in a Zariski neighborhood of the image of the
closed embedding
$\iota:\LocSys_G\hookrightarrow \LocSys_G^{\on{R.S.}}$. 

\smallskip

\noindent{\em(b)} If $G$ is unipotent \footnote{S.~Raskin has observed that the assertion and its proof remain valid
under the weaker assumption that the identity connected component of $G$ is solvable.}, 
then $\LocSys_G^{\on{R.S.}}$ is smooth.
\end{prop}

\begin{proof}

Let $\CZ$ be an Artin stack with a perfect cotangent complex. (For instance, this is the case if $\CZ$ is quasi-smooth.) 
It is easy to see that smoothness of $\CZ$ can be verified at $k$-points. Namely,
a point $z:\Spec(k)\to \CZ$
belongs to the smooth locus of $\CZ$ if and only if $T_z^*(\CZ)\in \Vect^{\geq 0}$. 

\medskip

A $k$-point of $\LocSys_G^{\on{R.S.}}$ is a pair $z=(\CP,\nabla)$, where
$\CP$ is a $G$-bundle on $X$, and $\nabla$ is a connection on $\CP$ with a simple
pole at $x$.

\medskip

We have the following description of 
$T^*_z(\LocSys_G^{\on{R.S.}})$, parallel to that of \eqref{e:exp for cotangent}:
\begin{equation} \label{e:log}
T^*_z(\LocSys_G^{\on{R.S.}})\simeq 
\on{Cone}\left(\nabla:\Gamma(X,\fg^*_{\CP}(-x))\to \Gamma(X,\fg^*_{\CP}\otimes \omega_X)\right).
\end{equation}

\medskip

Therefore, a point $z$ belongs to the smooth locus if and only if the map of (classical)
vector spaces 
$$\nabla:H^0\left(\Gamma(X,\fg^*_{\CP}(-x))\right)\to H^0\left(\Gamma(X,\fg^*_{\CP}\otimes \omega_X)\right)$$
is injective. 

\medskip

In other words, smooth points correspond to pairs $(\CP,\nabla)$ such that  
$\fg^*_\CP$ has no non-zero horizontal sections that vanish at $x$. Recall that the connection on $\fg^*_\CP$
has a simple pole at $x$; the condition automatically holds if
none of the eigenvalues of the coadjoint action of the residue $\on{res}(\nabla)\in\fg/G$ is a negative integer.

\medskip

In particular, if $(\CP,\nabla)$ is a point of $\iota(\LocSys_G)$, then $\on{res}(\nabla)=0$ and the condition trivially holds; 
this proves part (a). On the other hand, if $G$ is unipotent, the coadjoint action of $\fg$ is nilpotent, and the
condition is satisfied as well; this proves part (b).
\end{proof}

\sssec{}
By \eqref{e:embed Sing stacks}, from \propref{p:RS smooth}(a), we obtain a canonical closed embedding
\begin{equation}  \label{e:embedding Sing LocSys}
\Arth_G\hookrightarrow \fg^*/G\underset{\on{pt}/G}\times \LocSys_G.
\end{equation}

Recall that by \corref{c:Sing of LocSys}, $\Arth_G$ is isomorphic to 
the moduli stack (in the classical sense) of triples
$(\CP,\nabla,A)$, where $(\CP,\nabla)\in\LocSys_G$ and $A\in H^0(\Gamma(X_\dr,\fg^*_\CP))$. It is easy to see that
\eqref{e:embedding Sing LocSys} is given by
\[(\CP,\nabla,A)\mapsto(A(x),(\CP,\nabla)).\]

\section{The global nilpotent cone and formulation of the conjecture}   \label{s:Langlands}

As before, let $X$ be a connected smooth projective curve. 
From now on we assume that the algebraic group $G$ is reductive. Let $\cG$ be its Langlands dual.

\medskip

In this section we will formulate the Geometric Langlands conjecture, whose automorphic (a.k.a. geometric)
side involves the category $\Dmod(\Bun_G)$, and the Galois (a.k.a. spectral) side, an appropriate
modification of the category $\QCoh(\LocSysG)$. 

\ssec{The global nilpotent cone}  \label{sss:glob nilp cone}

\sssec{}

Recall that \propref{c:Sing of LocSys} provides an isomorphism between $\Sing(\LocSysG)$ and the (classical)
moduli stack $\ArthG$, which parametrizes triples $(\CP,\nabla,A)$.
Here $\CP$ is a $\cG$-bundle on $X$, $\nabla$ is a connection on $\CP$, and $A$ is a horizontal section of $\cg^*_\CP$.

\medskip

We define a Zariski-closed subset 
\begin{equation} \label{e:glob nilp cone}
\on{Nilp}_{glob}\subset \ArthG
\end{equation}
to correspond to triples $(\CP,\nabla,A)$ with nilpotent $A$.

\medskip

That is, we require that for every local trivialization of $\CP$, the element $A$ viewed (locally) as a map
$S\times X\to \cg^*$ hit the locus of nilpotent elements $\on{Nilp}(\cg^*)\subset \cg^*$.
The latter is defined as the image of the
locus of nilpotent elements $\on{Nilp}(\cg)\subset \cg$ under some (or any) $\cG$-invariant identification $\cg\simeq \cg^*$.

\sssec{}

Let $\fc(\cg)$ denote the characteristic variety of $\cg$, i.e., 
$$\fc(\cg):=\Spec(\Sym(\cg)^\cG),$$
and let $\varpi$ denote the Chevalley map
$$\varpi:\cg^*=\Spec(\Sym(\cg))\to \Spec(\Sym(\cg)^\cG)=\fc(\cg).$$ 

\medskip

For $(\CP,\nabla,A)\in \Maps(S,\ArthG)$ we thus obtain a map 
$$\varpi(A):S\times X\to \fc(\cg).$$

The nilpotence condition can be phrased as the requirement that $\varpi(A)$ should factor through $$\{0\}\subset \fc(\cg).$$

\sssec{}

We can also express the nilpotence condition locally:

\begin{lem} \label{l:nilpotent at point}
For an $S$-point $(\CP,\nabla,A)$ of $\ArthG$, 
the element $A$ is nilpotent if and only if for some (and then any) 
point $x\in X$, the value $A|_{S\times \{x\}}$ of $A$ at $x$
is nilpotent as a section of $\cg^*_{\CP_x}:=\cg^*_{\CP}|_{S\times \{x\}}$.
\end{lem}

\begin{proof}
The fact that $A$ is horizontal implies that the map $\varpi(A)$ is infinitesimally constant along $X$
(i.e., factors through a map $S\times X_\dr\to \fc(\cg)$), and therefore is 
constant (since $X$ is connected). This implies the assertion of the lemma.
\end{proof}

\medskip

Recall that by \eqref{e:embedding Sing LocSys}, we have a canonical closed embedding
$$\ArthG\hookrightarrow \cg^*/\cG\underset{\on{pt}/\cG}\times \LocSysG.$$

Thus, \lemref{l:nilpotent at point} can be reformulated as the equality between 
$$\on{Nilp}_{glob}\subset \ArthG$$ 
and the preimage of the closed subset
$$\on{Nilp}(\cg^*)/\cG\underset{\on{pt}/\cG}\times \LocSysG\subset \cg^*/\cG\underset{\on{pt}/\cG}\times \LocSysG$$
under the above map.

\sssec{The spectral side of the Geometric Langlands conjecture}

Our main object of study is the category 
$$\IndCoh_{\on{Nilp}_{glob}}(\LocSysG).$$

By definition, this is a full subcategory of $\IndCoh(\LocSysG)$, which contains the essential image of
$\QCoh(\LocSysG)$ under the functor
$$\Xi_{\LocSysG}:\QCoh(\LocSysG)\to \IndCoh(\LocSysG).$$

\medskip

We propose the category $\IndCoh_{\on{Nilp}_{glob}}(\LocSysG)$ as the category appearing on the spectral
side of the Geometric Langlands conjecture.

\sssec{}   \label{sss:life over Lie algebra}

By \corref{c:category comp gen stacks} and \propref{p:RS smooth}, the category $\IndCoh_{\on{Nilp}_{glob}}(\LocSysG)$
is compactly generated. 

\medskip

By \secref{sss:tensored over conormal stacks parallel}, \propref{p:RS smooth} allows us to 
view $\IndCoh(\LocSysG)$ as tensored over the monoidal category $\QCoh(\cg^*/(\cG\times \BG_m))$. 
We emphasize that the latter structure depends on the choice of a point $x\in X$. 

\medskip

By \lemref{l:nilpotent at point} and \corref{c:category as tensor product stacks}
we have:
$$\IndCoh_{\on{Nilp}_{glob}}(\LocSysG)\simeq \IndCoh(\LocSysG)\underset{\QCoh(\cg^*/(\cG\times \BG_m))}
\otimes \QCoh(\on{Nilp}(\cg^*)/(\cG\times \BG_m)).$$

\ssec{Formulation of the Geometric Langlands conjecture}

\sssec{}

We propose the following form of the Geometric Langlands conjecture:

\begin{conj} \label{conj:main}
There exists an equivalence of DG categories
$$\Dmod(\Bun_G)\simeq \IndCoh_{\on{Nilp}_{glob}}(\LocSysG).$$
\end{conj}

Since the DG categories appearing on both sides of \conjref{conj:main} are compactly generated,
it can be tautologically rephrased as follows:
\begin{conj}
There exists an equivalence of non-cocomplete DG categories
$$\Dmod(\Bun_G)^c\simeq \Coh_{\on{Nilp}_{glob}}(\LocSysG).$$
\end{conj}

\sssec{}  \label{sss:tempered}

In what follows we will refer to the essential image in $\Dmod(\Bun_G)$ of
$$\Xi_{\LocSysG}\left(\QCoh(\LocSysG)\right)\subset \IndCoh_{\on{Nilp}_{glob}}(\LocSysG)$$
under the above conjectural equivalence as the ``tempered part" of $\Dmod(\Bun_G)$,
and denote it by $\Dmod_{\on{temp}}(\Bun_G)$.

\sssec{}

Of course, one needs to specify a lot more data to fix the equivalence of \conjref{conj:main} uniquely.
This will be done over the course of several papers following this one. In the present paper we
will discuss the following aspects: 

\medskip

\noindent(i) The case when $G$ is a torus; 

\medskip

\noindent(ii) Compatibility with the Geometric Satake Equivalence (see \secref{s:Satake});

\medskip

\noindent(iii) Compatibility with the Eisenstein series (see \secref{s:eis}).

\sssec{The case of a torus}

Let $G$ be a torus $T$. This case offers nothing new. The subset $\on{Nilp}_{glob}$ is the zero-section of $\Sing(\LocSys_{\cT})$, so
by \corref{c:zero sing supp stacks},
$$\IndCoh_{\on{Nilp}_{glob}}(\LocSys_{\cT})=\Xi_{\LocSysG}\left(\QCoh(\LocSysG)\right),$$
as subcategories of $\IndCoh(\LocSys_{\cT})$.

\medskip

In this case, the equivalence
$$\QCoh(\LocSys_{\cT})\simeq \Dmod(\Bun_T)$$
is a particular case of the Fourier transform for D-modules on an abelian variety (see \cite{L1,L2} and \cite{Ro1,Ro2}),
appropriately adjusted to the DG setting. 

\begin{rem}
In more detail, for $G=\BG_m$, a choice of $x\in X$ identifies
$$\Bun_{\BG_m}\simeq \on{Pic}\times \on{pt}/\BG_m \times \BZ \text{ and }
\LocSys_{\BG_m}\simeq \wt{\on{Pic}}\times (\on{pt}\underset{\BA^1}\times \on{pt})\times \on{pt}/\BG_m,$$
where $\on{Pic}$ is the Picard \emph{scheme} and $\wt{\on{Pic}}$ is its universal additive extension. Then the
classical Fourier-Mukai-Laumon transform identifies 
$$\Dmod(\on{Pic})\simeq \QCoh(\wt{\on{Pic}}),$$
and we have explicit equivalences of categories 
$$\Dmod(\on{pt}/\BG_m)\simeq \QCoh(\on{pt}\underset{\BA^1}\times \on{pt}) \text{ and } \Dmod(\BZ)\simeq \QCoh(\on{pt}/\BG_m).$$

Note that the compact generator of $\Dmod(\on{pt}/\BG_m)$ is the direct image with compact supports
of $k\in \Vect$ under the map $\on{pt}\to \on{pt}/\BG_m$.
It corresponds to the structure sheaf of $\QCoh(\on{pt}\underset{\BA^1}\times \on{pt})$. 

\medskip

Note also that under this 
equivalence, the constant
sheaf $k_{\on{pt}/\BG_m}\in \Dmod(\on{pt}/\BG_m)$ is not compact, and it corresponds to the sky-scraper on 
$\on{pt}\underset{\BA^1}\times \on{pt}$, which is an object of $\Coh(\on{pt}\underset{\BA^1}\times \on{pt})$, but not of
$\QCoh(\on{pt}\underset{\BA^1}\times \on{pt})^{\on{perf}}$, i.e., it is \emph{not} compact in 
$\QCoh(\on{pt}\underset{\BA^1}\times \on{pt})$. 

\end{rem}

\section{Compatibility with Geometric Satake Equivalence}   \label{s:Satake}

One of the key properties of the Geometric Langlands equivalence is its behavior with respect to the 
Hecke functors on both sides of the correspondence. In this section we will study how this is compatible
with the proposed candidate for the spectral side: the category $\IndCoh_{\on{Nilp}_{glob}}(\LocSysG)$.

\ssec{Main results of this section}

\sssec{The Geometric Satake Equivalence}

As before, $X$ is a smooth connected projective curve, $G$ is a reductive group and $\cG$ is its Langlands dual. 
Fix
a point $x\in X$. The category of Hecke functors at $x$ (``the spherical Hecke category at $x$'') is the category
of $G(\wh\CO_x)$-equivariant D-modules on the affine Grassmannian $\Gr_{G,x}$, which we denote by 
\[\Sph(G,x):=\Dmod(\Gr_{G,x})^{G(\wh\CO_x)}.\]
We regard it as a monoidal category with respect to the convolution product.

\medskip

In \corref{c:Satake}, we will construct a monoidal equivalence 
\[\on{Sat}:\IndCoh_{\on{Nilp}(\cg^*)/\cG}(\on{Hecke}(\cG)_{spec})\simeq\Sph(G,x)\]
between $\Sph(G,x)$ and a certain category constructed from the group $\cG$. Explicitly,
$\IndCoh_{\on{Nilp}(\cg^*)/\cG}(\on{Hecke}(\cG)_{spec})$
is the category of ind-coherent sheaves on the stack
\[\on{Hecke}(\cG)_{spec}:=\on{pt}/\cG\underset{\cg/\cG}\times \on{pt}/\cG\]
whose singular support is contained in 
\[\on{Nilp}(\cg^*)/\cG\subset\cg^*/\cG\simeq\Sing(\on{Hecke}(\cG)_{spec}),\]
where $\on{Nilp}(\cg^*)\subset\cg^*$ is the nilpotent cone.

\medskip

We refer to $\on{Sat}$ as the Geometric Satake Equivalence. It is naturally related to the other versions of the Satake equivalence, constructed
in \cite{MV} and \cite{BF} (which are given below as \eqref{e:MV} and \thmref{t:Bezr}, respectively).

\sssec{The Geometric Langlands conjecture and the Geometric Satake Equivalence}

We show that the ``modified'' Geometric Langlands conjecture (\conjref{conj:main}) 
agrees with the Geometric Satake Equivalence $\on{Sat}$: there is a natural action of the monoidal
category \[\IndCoh_{\on{Nilp}(\cg^*)/\cG}(\on{pt}/\cG\underset{\cg/\cG}\times \on{pt}/\cG)\] 
on the category $\IndCoh_{\on{Nilp}_{glob}}(\LocSysG)$, see \corref{p:Hecke action spec}. Under the equivalence of \conjref{conj:main},
this action should correspond to the action of the monoidal category $\on{Sph}(G,x)$ on $\Dmod(\Bun_G)$ by the Hecke functors; this is \conjref{conj:comp with Satake}.

\sssec{Tempered D-modules}

Recall that \conjref{conj:main} implies the existence of a certain full subcategory 
$\Dmod_{\on{temp}}(\Bun_G)\subset \Dmod(\Bun_G)$, the ``tempered part'' of $\Dmod(\Bun_G)$,
defined as the essential image in $\Dmod(\Bun_G)$ of $$\Xi_{\LocSysG}\left(\QCoh(\LocSysG)\right)\subset 
\IndCoh_{\on{Nilp}_{glob}}(\LocSysG).$$ Assuming the compatibility of \conjref{conj:comp with Satake}, 
we will be able to describe the subcategory $\Dmod_{\on{temp}}(\Bun_G)$
in purely ``geometric'' terms, using the Hecke functors at a fixed point
$x\in X$. The description is independent of the Langlands conjecture; we denote the resulting full subcategory
by \[\Dmod^x_{\on{temp}}(\Bun_G)\subset\Dmod(\Bun_G).\] 
However, it is not clear that the subcategory $\Dmod^x_{\on{temp}}(\Bun_G)$ 
is independent of the choice of the point $x$. This is the content
of \conjref{conj:tempered independence}, which follows from Conjectures~\ref{conj:main} and \ref{conj:comp with Satake}.

\ssec{Preliminaries on the spherical Hecke category}

\sssec{}

Recall that the spherical Hecke category at a point $x\in X$ is defined as
\[\Sph(G,x):=\Dmod(\Gr_{G,x})^{G(\wh\CO_x)},\] 
regarded as a monoidal category with respect to the convolution product. 
We claim that as a DG category, $\Sph(G,x)$ is compactly generated. 

\medskip

Indeed, we can represent $\Gr_{G,x}$ as a 
union of $G(\wh\CO_x)$-invariant finite-dimensional closed subschemes $Z_\alpha$. We have
$$\Dmod(\Gr_{G,x})^{G(\wh\CO_x)}\simeq \underset{\alpha}{\underset{\longrightarrow}{colim}}\, \Dmod(Z_\alpha)^{G(\wh\CO_x)},$$
where for $\alpha_1\geq \alpha_2$, the functor 
$$\Dmod(Z_{\alpha_1})^{G(\wh\CO_x)}\to \Dmod(Z_{\alpha_2})^{G(\wh\CO_x)}$$
is given by direct image along the corresponding closed embedding. In particular, for every $\alpha$, the functor
$$\Dmod(Z_\alpha)^{G(\wh\CO_x)}\to \Dmod(\Gr_{G,x})^{G(\wh\CO_x)}$$
sends compacts to compacts. By \cite[Lemma 1.3.3]{DG}, this reduces the assertion to showing that each 
$\Dmod(Z_\alpha)^{G(\wh\CO_x)}$ is compactly generated. 

\medskip

Let $G_\alpha$ be a finite-dimensional quotient of $G(\wh\CO_x)$ through which it acts on $Z_\alpha$.
With no restriction of generality, we can assume that $\on{ker}(G(\wh\CO_x)\to G_\alpha)$ is pro-unipotent.
Hence, the forgetful functor
$$\Dmod(Z_\alpha)^{G_\alpha}\to \Dmod(Z_\alpha)^{G(\wh\CO_x)}$$
is an equivalence.

\medskip

Now, $Z_\alpha/G_\alpha$ is a QCA algebraic stack, and the compact generation of $\Dmod(Z_\alpha)^{G_\alpha}$ 
follows from \cite[Theorem 0.2.2]{DrGa0}. (Since $Z_\alpha/G_\alpha$ is a global quotient, the compact generation
follows more easily from the results of \cite{BFN}.)

\medskip

Note that the monoidal operation on $\Sph(G,x)$ preserves the subcategory of compact objects (this follows
from the properness of $\Gr_{G,x}$). Hence, $\Sph(G,x)^c$ acquires a structure of non-cocomplete monoidal
DG category. 

\sssec{}
Consider the heart $\left(\Dmod(\Gr_{G,x})^{G(\wh\CO_x)}\right){}\!^\heartsuit$
of the natural t-structure on the category $\Dmod(\Gr_{G,x})^{G(\wh\CO_x)}$. I.e., 
$\left(\Dmod(\Gr_{G,x})^{G(\wh\CO_x)}\right){}\!^\heartsuit$ is the abelian category of $G(\wh\CO_x)$-equivariant
D-modules on $\Gr_{G,x}$. 
The geometric Satake isomorphism of \cite{MV} gives an equivalence between
$\left(\Dmod(\Gr_{G,x})^{G(\wh\CO_x)}\right){}\!^\heartsuit$ and the abelian category of representations of $\cG$,
which we denote by $\Rep(\cG)^\heartsuit$.

\medskip

We let $\Sph(G,x)^{naive}$ be the derived category (considered as a DG category) of the abelian category
$$\left(\Dmod(\Gr_{G,x})^{G(\wh\CO_x)}\right){}\!^\heartsuit.$$

We have a canonical (but not fully faithful) monoidal functor
\begin{equation} \label{e:naive to all}
\Sph(G,x)^{naive}\to \Sph(G,x)
\end{equation}
(see \cite[Theorem 1.3.2.2]{Lu1}). 

\medskip

The Satake equivalence of \cite{MV} induces a canonical t-exact equivalence of monoidal categories
\begin{equation} \label{e:MV}
\on{Sat}^{naive}:\Rep(\cG)\simeq \Sph(G,x)^{naive}.
\end{equation}
We will refer to it as the ``naive" version of the Geometric Satake Equivalence. 

\medskip

In order to describe $\Sph(G,x)$, it is convenient
to first introduce and describe its ``renormalized" version. Here ``renormalization'' refers to the
process of changing (in this case, enlarging) the class of compact objects of the category.

\sssec{}

Let $\Sph(G,x)^{loc.c}$ denote the full subcategory of $\Sph(G,x)$ consisting of those
objects of $\Sph(G,x)\simeq \Dmod(\Gr_{G,x})^{G(\wh\CO_x)}$ that become compact
after applying the forgetful functor 
$$\Dmod(\Gr_{G,x})^{G(\wh\CO_x)}\to \Dmod(\Gr_{G,x}).$$
(The superscript ``loc.c" stands for ``locally compact.")

\medskip

The category $\Sph(G,x)^{loc.c}\subset \Sph(G,x)$ is stable under the monoidal operation, and hence 
acquires a structure of (non-cocomplete) monoidal DG category. We define
$$\Sph(G,x)^{\on{ren}}:=\Ind(\Sph(G,x)^{loc.c}),$$
which thus acquires a structure of monoidal DG category.

\sssec{}  \label{sss:Psi for Sph}

We have a canonically defined monoidal functor
$$\Psi^{\Sph}:\Sph(G,x)^{\on{ren}}\to \Sph(G,x),$$
obtained by ind-extending the tautological embedding $\Sph(G,x)^{loc.c}\hookrightarrow \Sph(G,x)$.

\medskip

Since $\Sph(G,x)^{loc.c}$ is closed under the truncations with respect to the (usual) t-structure on
$\Sph(G,x)$, we obtain that $\Sph(G,x)^{\on{ren}}$ acquires a unique t-structure, compatible
with colimits\footnote{A t-structure on a cocomplete DG category is called compatible with colimits
if the subcategory of coconnective objects is closed under filtered colimits.},
for which the functor $\Psi^{\Sph}$ is t-exact.

\medskip

The functor $\Psi^{\Sph}$ admits a left adjoint, denoted $\Xi^{\Sph}$, obtained by ind-extending the tautological
embedding $\Sph(G,x)^{c}\hookrightarrow \Sph(G,x)^{loc.c}$. By construction, the functor $\Xi^{\Sph}$ is fully faithful. 
So, $\Psi^{\Sph}$ makes $\Sph(G,x)$ into a colocalization of $\Sph(G,x)^{\on{ren}}$. 

\medskip

Finally, note that $\Xi^{\Sph}$ has a natural structure of a monoidal functor. Indeed, it is clearly co-lax monoidal, being the left
adjoint of a monoidal functor. Since the category $\Sph(G,x)$ is compactly generated, it suffices to check that $\Xi^{\Sph}$
is (strictly) monoidal after restriction to the category of compact objects $\Sph(G,x)^c$. However, this restriction identifies with the
embedding of monoidal categories $\Sph(G,x)^c\hookrightarrow \Sph(G,x)^{loc.c}$.

\sssec{}

We claim that the tautological functor $\Sph(G,x)^{naive}\to\Sph(G,x)$ of \eqref{e:naive to all} canonically factors as
$$\Sph(G,x)^{naive}\to \Sph(G,x)^{\on{ren}}
\overset{\Psi^{\Sph}}\longrightarrow \Sph(G,x).$$

Indeed, we construct the functor 
\begin{equation} \label{e:naive to ren}
\Sph(G,x)^{naive}\to \Sph(G,x)^{\on{ren}}
\end{equation}
as the ind-extension of a functor $(\Sph(G,x)^{naive})^c\to \Sph(G,x)^{loc.c}$.
The latter is obtained by noticing that the essential image of $(\Sph(G,x)^{naive})^c$ under the
functor \eqref{e:naive to all} is contained in $\Sph(G,x)^{loc.c}$.

\medskip

By construction, the functor \eqref{e:naive to ren} sends compact objects to compact ones. By contrast, the functor \eqref{e:naive to all} 
does not have this property. 

\medskip

We will denote by $\on{Sat}^{naive,\on{ren}}$ the resulting functor
$$\Rep(\cG)\to \Sph(G,x)^{\on{ren}}.$$

\ssec{The Hecke category on the spectral side}

\sssec{}

Consider now the stack
$$\on{Hecke}(\cG)_{spec}:=\on{pt}/\cG\underset{\cg/\cG}\times \on{pt}/\cG,$$
where both maps $\on{pt}\to \cg$ correspond to $0\in \cg$. In the notation of
\secref{sss:groupoid over stacks}, 
\[\on{Hecke}(\cG)_{spec}=\CG_{(\on{pt}/\cG)/(\cg/\cG)}.\]

\medskip

The stack $\on{Hecke}(\cG)_{spec}$ is naturally a groupoid acting on $\on{pt}/\cG$. 
This groupoid structure equips
$$\IndCoh(\on{Hecke}(\cG)_{spec})$$
with a structure of monoidal category via convolution. 

\medskip

We can also consider the 
subcategory 
$$\Coh(\on{Hecke}(\cG)_{spec})\subset \IndCoh(\on{Hecke}(\cG)_{spec}),$$ which is stable under the monoidal
operation, and thus acquires a structure of (non-cocomplete) monoidal category, whose ind-completion
identifies with $\IndCoh(\on{Hecke}(\cG)_{spec})$. 

\sssec{}

The following description of $\Sph(G,x)^{loc.c}$ is given by \cite[Theorem 5]{BF}\footnote{The statement in {\it loc.cit.}
is the combination of \thmref{t:Bezr} as stated below and \propref{p:Koszul dual for Hecke}.}:

\begin{thm} \label{t:Bezr}
There is a canonical equivalence of (non-cocomplete) monoidal categories
$$\Coh(\on{Hecke}(\cG)_{spec})\simeq \Sph(G,x)^{loc.c}.$$
\end{thm}

This equivalence tautologically extends to an equivalence between the ind-completions of these categories,
giving the following ``renormalized'' Geometric Satake Equivalence. 

\begin{cor}  \label{c:Bezr}
There exists a canonical equivalence of monoidal categories 
$$\on{Sat}^{\on{ren}}:\IndCoh(\on{Hecke}(\cG)_{spec})\simeq \Sph(G,x)^{\on{ren}}.$$
\end{cor}

Under this equivalence, the functor $\on{Sat}^{naive,\on{ren}}:\Rep(\cG)\to \Sph(G,x)^{\on{ren}}$
corresponds to the canonical functor 
$$\Rep(\cG)\to \IndCoh(\on{Hecke}(\cG)_{spec}),$$
given by the direct image along the diagonal map
$$\Delta_{\on{pt}/\cG}:\on{pt}/\cG\to \on{pt}/\cG\underset{\cg/\cG}\times \on{pt}/\cG=\on{Hecke}(\cG)_{spec}.$$

\begin{rem}  \label{r:E3}
The category $\Sph(G,x)$, as well as $\Sph(G,x)^{\on{ren}}$, has a richer structure, namely, that of
\emph{factorizable} monoidal category, when we allow the point $x$ to move along $X$.
One can see this structure on the category 
$\IndCoh(\on{Hecke}(\cG)_{spec})$ as well, and one can show that the equivalence of
\eqref{c:Bezr} can be naturally upgraded to an equivalence of factorizable monoidal
categories. \footnote{The latter statement is known as ``derived Satake"; it was conjectured by
V.~Drinfeld and proved by J.~Lurie and the second author (unpublished) by interpreting 
$\IndCoh(\on{Hecke}(\cG)_{spec})$ as the $\BE_3$-center of $\Rep(\cG)$, viewed as an
$\BE_2$-category.} 

\end{rem}

\ssec{A Koszul dual description}

\sssec{}

Consider the commutative DG algebra $\on{Sym}(\cg[-2])\mod$, which is acted on canonically by $\cG$.
Consider the category 
$$(\on{Sym}(\cg[-2])\mod)^\cG$$
as a monoidal category via the usual tensor product operation of modules over a commutative algebra.

\medskip

We claim:

\begin{prop} \label{p:Koszul dual for Hecke}
There exists a canonical equivalence of monoidal
categories
$$\on{KD}_{\on{Hecke}(\cG)_{spec}}:\IndCoh(\on{Hecke}(\cG)_{spec})\simeq (\on{Sym}(\cg[-2])\mod)^\cG.$$
\end{prop}

\begin{proof}
Consider $V=\cg/\cG$ as a vector bundle over $\CX=\on{pt}/\cG$. Then $\on{Sym}_{\CO_\CX}(V[-2])$ is a commutative
(i.e., $\BE_\infty$) algebra in
$\QCoh(\CX)$, and we have a natural equivalence
\[(\on{Sym}(\cg[-2])\mod)^\cG\simeq\Sym_{\CO_\CX}(V[-2])\mod.\]
Note that we are in the setting of \secref{ss:parallel stacks}; therefore, we have an isomorphism of $\BE_2$-algebras
in $\QCoh(\CX)$:
\[\Sym_{\CO_\CX}(V[-2])\simeq\on{HC}(\CX/V).\]
Now the claim follows from Koszul duality of \corref{c:Koszul at point stacks}.
\end{proof}

\begin{rem} \label{r:not E3}
Being a symmetric monoidal category, $(\on{Sym}(\cg[-2])\mod)^\cG$ also has a structure
of factorizable monoidal category over $X$. However, the equivalence of \propref{p:Koszul dual for Hecke}
is only between mere monoidal categories: it is not compatible with the factorizable structure.
In fact, one can show that $\IndCoh(\on{Hecke}(\cG)_{spec})$ does not admit an $\BE_2$-structure
which is compatible with the factorizable structure, even if $G$ is a torus. \footnote{This is more convenient
to see on the geometric side. The key fact is that the transgression
map $H^\bullet(\on{pt}/\BG_m)\otimes H_\bullet(X)\to 
H^\bullet(\Bun_{\BG_m})$ does not commute with the maps $H^\bullet(\Bun_{\BG_m})\to H^\bullet(\Bun_{\BG_m})$
given by translation by points of $X$. Here $H^\bullet(\on{pt}/\BG_m)$ appears as the endomorphism
algebra of the unit object of $\Sph(\BG_m,x)$.} 
\end{rem}

\sssec{}

Combining \corref{c:Bezr} and \propref{p:Koszul dual for Hecke}, we obtain:

\begin{cor}  \label{c:Satake Koszul renorm}
There exists a canonical equivalence of monoidal categories
$$\on{Sat}^{\on{ren}}\circ (\on{KD}_{\on{Hecke}(\cG)_{spec}})^{-1}:(\on{Sym}(\cg[-2])\mod)^\cG\simeq \Sph(G,x)^{\on{ren}}.$$
\end{cor}

Moreover, from \secref{sss:Psi for Sph}, we obtain:

\begin{cor} \label{c:Satake Koszul as coloc}
There exists a canonically defined monoidal functor
$$(\on{Sym}(\cg[-2])\mod)^\cG\to \Sph(G,x),$$
which is, moreover, a colocalization.
\end{cor}

\sssec{}  \label{sss:Bezr expl}

In fact, the equivalence of \corref{c:Satake Koszul renorm} can be made more explicit:

\medskip

Let $\delta_\one$ be the unit object of 
$$\Sph(G,x)^{loc.c}\subset \Sph(G,x),$$ 
given by the delta-function at $\one\in \Gr_{G,x}$. 
\thmref{t:Bezr} implies that there exists a canonical isomorphism of $\BE_2$-algebras
\begin{equation} \label{e:calc Bezr algebras}
\CMaps_{\Sph(G,x)^{loc.c}}(\delta_\one,\delta_\one)\simeq \on{Sym}(\cg[-2])^\cG,
\end{equation}
and that for any $\CM\in \Sph(G,x)^{loc.c}$, we have an isomorphism of $\on{Sym}(\cg[-2])^\cG$-modules
\begin{equation} \label{e:calc Bezr}
\left(\on{KD}_{\on{Hecke}(\cG)_{spec}}\left((\on{Sat}^{\on{ren}})^{-1}(\CM)\right)\right)^\cG\simeq
\CMaps_{\Sph(G,x)^{loc.c}}(\delta_\one,\CM).
\end{equation}

\sssec{}

Thus, $\on{KD}_{\on{Hecke}(\cG)_{spec}}\left((\on{Sat}^{\on{ren}})^{-1}(\CM)\right)$ is a $\cG$-equivariant
module over $\on{Sym}(\cg[-2])$, and \eqref{e:calc Bezr} recovers $\cG$-invariants in this module. 

\medskip

One can reconstruct the entire module $\on{KD}_{\on{Hecke}(\cG)_{spec}}\left((\on{Sat}^{\on{ren}})^{-1}(\CM)\right)$ by
considering convolutions of $\CM$ with objects of the form $\on{Sat}^{naive,\on{ren}}(\rho)$ for $\rho\in \Rep(\cG)^c$:

$$\left(\on{KD}_{\on{Hecke}(\cG)_{spec}}\left((\on{Sat}^{\on{ren}})^{-1}(\CM)\right)\otimes \rho\right)^\cG\simeq
\CMaps_{\Sph(G,x)^{loc.c}}(\delta_\one,\CM\star \on{Sat}^{naive,\on{ren}}(\rho)).$$

\sssec{}  \label{sss:equiv cohomology}

Note that since $\one$ is a closed $G(\wh\CO_x)$-invariant point of $\Gr_{G,x}$,
\begin{multline}  \label{e:equiv cohomology}
\CMaps_{\Sph(G,x)^{loc.c}}(\delta_\one,\delta_\one)\simeq 
\CMaps_{\Sph(G,x)}(\delta_\one,\delta_\one)\simeq \\
\simeq \CMaps_{\Dmod(\on{pt})^{G(\wh\CO_x)}}(\delta,\delta)\simeq
\CMaps_{\Dmod(\on{pt})^G}(\delta,\delta),
\end{multline}
where $\delta$ denotes the generator $k\in \Vect\simeq \on{pt}$, which is naturally equivariant
with respect to any group. 

\medskip

We note that the last isomorphism in \eqref{e:equiv cohomology} is due to the fact that 
$\on{ker}(G(\wh\CO_x)\to G)$ is pro-unipotent.

\medskip

By definition, the algebra $\CMaps_{\Dmod(\on{pt})^G}(\delta,\delta)$ is the equivariant cohomology
of $G$, which we denote by $H_\dr(\on{pt}/G)$, and we have a canonical isomorphism
$$H_\dr(\on{pt}/G)\simeq \on{Sym}(\fh^*[-2])^W\simeq \on{Sym}(\ch[-2])^W\simeq \on{Sym}(\cg[-2])^\cG.$$

Now, it follows from the construction of the isomorphism of \thmref{t:Bezr}, that the resulting isomorphism
$$\on{Sym}(\cg[-2])^\cG\simeq H_\dr(\on{pt}/G)\simeq \CMaps_{\Dmod(\on{pt})^G}(\delta,\delta)\simeq
\CMaps_{\Sph(G,x)^{loc.c}}(\delta_\one,\delta_\one)$$
equals one given by \eqref{e:calc Bezr algebras}.

\ssec{A description of $\Sph(G,x)$}

\sssec{}  \label{sss:Koszul and sing support for Satake}

Recall that the stack 
$$\on{Hecke}(\cG)_{spec}=\on{pt}/\cG\underset{\cg/\cG}\times \on{pt}/\cG=\CG_{(\on{pt}/\cG)/(\cg/\CG)}$$
is quasi-smooth, and
$$\Sing(\on{Hecke}(\cG)_{spec})\simeq \cg^*/\cG,$$
see \secref{sss:sing of Hecke stack}.

\medskip

Moreover, by \corref{c:Koszul at point stacks}, 
the equivalence of \propref{p:Koszul dual for Hecke}
calculates the singular support of objects of $\on{Hecke}(\cG)_{spec}$:

\medskip

\noindent For
$\CF\in \on{Hecke}(\cG)_{spec}$, we have
\begin{equation} \label{e:sing supp and supp}
\on{SingSupp}(\CF)=\on{supp}(\on{KD}_{\on{Hecke}(\cG)_{spec}}(\CF)),
\end{equation}
as subsets of $\cg^*/\cG$. 

\sssec{}

We will prove:

\begin{thm}  \label{t:Satake}
Under the equivalence 
$$\on{Sat}^{\on{ren}}:\IndCoh(\on{Hecke}(\cG)_{spec})\simeq \Sph(G,x)^{\on{ren}}$$
of \corref{c:Bezr}, the colocalization
$$\Xi^{\Sph}:\Sph(G,x)\rightleftarrows \Sph(G,x)^{\on{ren}}:\Psi^{\Sph}$$
identifies with
$$\IndCoh_{\on{Nilp}(\cg^*)/\cG}(\on{Hecke}(\cG)_{spec})
\rightleftarrows \IndCoh(\on{Hecke}(\cG)_{spec}).$$
\end{thm}

We remind that the functors $(\Xi^{\Sph}, \Psi^{\Sph})$ appearing in \thmref{t:Satake} are those from 
\secref{sss:Psi for Sph}.

\medskip

In terms of \corref{c:Satake Koszul as coloc}, the assertion of \thmref{t:Satake} can be reformulated
as follows:

\begin{cor}
The colocalization 
$$\Sph(G,x)\rightleftarrows (\on{Sym}(\cg[-2])\mod)^\cG$$
of \corref{c:Satake Koszul as coloc} identifies with
$$\left((\on{Sym}(\cg[-2])\mod)^\cG\right)_{\on{Nilp}(\cg^*)/\cG}\rightleftarrows (\on{Sym}(\cg[-2])\mod)^\cG.$$
\end{cor}

\thmref{t:Satake}, in particular, implies:
\begin{cor} \label{c:Satake}
There exists a canonical equivalence of monoidal categories
$$\on{Sat}:\IndCoh_{\on{Nilp}(\cg^*)/\cG}(\on{Hecke}(\cG)_{spec})\simeq \Sph(G,x),$$
and of non-cocomplete monoidal categories\footnote{The fact that $\Coh_{\on{Nilp}(\cg^*)/\cG}(\on{Hecke}(\cG)_{spec})$ is
preserved under the monoidal operation follows, e.g., from \propref{p:Koszul dual for Hecke}.}
$$\Coh_{\on{Nilp}(\cg^*)/\cG}(\on{Hecke}(\cG)_{spec})\simeq \Sph(G,x)^c.$$
\end{cor}

\ssec{Proof of \thmref{t:Satake}}

By \eqref{e:sing supp and supp}, we need to show that the essential image of $\Sph(G,x)$ under the equivalence
$$\on{KD}_{\on{Hecke}(\cG)_{spec}}\circ (\on{Sat}^{\on{ren}})^{-1}:\Sph(G,x)^{\on{ren}}\simeq (\on{Sym}(\cg[-2])\mod)^\cG.$$
coincides with the subcategory
$$(\on{Sym}(\cg[-2])\mod_{\on{Nilp}(\cg^*)})^\cG\subset (\on{Sym}(\cg[-2])\mod)^\cG.$$

\sssec{}

Let 
$$\Dmod(\Gr_{G,x})^{G(\wh\CO_x)\on{-mon}}\subset \Dmod(\Gr_{G,x})$$
be the full subcategory generated by the essential image of the forgetful functor
$$\Sph(G,x)^{loc.c}\to \Dmod(\Gr_{G,x}).$$ 

\medskip

Since $\Sph(G,x)^{loc.c}$ is closed under the truncation functors, we obtain that the category
$\Dmod(\Gr_{G,x})^{G(\wh\CO_x)\on{-mon}}$ is compactly generated by the essential image of
$$(\Sph(G,x)^{loc.c})^\heartsuit\subset \Sph(G,x)^{loc.c}.$$

\medskip

Since the generators of $\Dmod(\Gr_{G,x})^{G(\wh\CO_x)\on{-mon}}$ are holonomic, 
the forgetful functor
$$\oblv_{G(\wh\CO_x)}:\Sph(G,x)\to \Dmod(\Gr_{G,x})^{G(\wh\CO_x)\on{-mon}}$$
admits a \emph{left} adjoint, given by !-averaging with respect to $G(\wh\CO_x)$.  We denote this
functor by $\on{Av}_{G(\wh\CO_x),!}$.

\sssec{}

Since the functor $\oblv_{G(\wh\CO_x)}$ is conservative, the essential image of $\on{Av}_{G(\wh\CO_x),!}$ 
generates $\Sph(G,x)$. Moreover, being a left adjoint of a continuous functor, $\on{Av}_{G(\wh\CO_x),!}$ 
sends compact objects to compact ones.

\medskip

Thus, we obtain that $\Sph(G,x)$ is compactly generated by the objects 
$$\on{Av}_{G(\wh\CO_x),!}\left(\oblv_{G(\wh\CO_x)}(\CM)\right)$$
for $\CM\in (\Sph(G,x)^{loc.c})^\heartsuit$.

\sssec{}

Note also that for 
$$\CM_1\in \Dmod(\Gr_{G,x})^{G(\wh\CO_x)\on{-mon}} \text{ and }
\CM_2\in \Sph(G,x),$$
we have:
\begin{equation} \label{e:averaging and convolution}
\on{Av}_{G(\wh\CO_x),!}(\CM_1)\star \CM_2\simeq
\on{Av}_{G(\wh\CO_x),!}(\CM_1\star \CM_2).
\end{equation}
In particular, if $\CM_1\in \Dmod(\Gr_{G,x})^c$, and
$\CM_2\in \Sph(G,x)^{loc.c}$, then 
$$\CM_1\star \CM_2\in \Dmod(\Gr_{G,x})^c,$$
and therefore in this case
$$\on{Av}_{G(\wh\CO_x),!}(\CM_1)\star \CM_2\in \Sph(G,x)^c.$$ 

\sssec{}

Let $\wt\delta$ be the object of $\Dmod(\on{pt})^G$
equal to $\on{Av}_{G,!}(k)$, where $\on{Av}_{G,!}$ is the \emph{left} adjoint
to the forgetful functor
\begin{equation} \label{e:forget equiv}
\Dmod(\on{pt})^G\to \Dmod(\on{pt})=\Vect.
\end{equation}

The following is well-known:

\begin{lem}
\begin{equation} \label{e:delta tilde}
\wt\delta\simeq \delta\underset{\on{Sym}(\cg[-2])^\cG}\otimes \fl,
\end{equation}
where $\delta$ is as in \secref{sss:equiv cohomology}, and $\fl$ is a graded 
line (placed in the cohomological degree $-\dim(G)$), acted on trivially by $\on{Sym}(\cg[-2])^\cG$. 
\end{lem}

\begin{proof}
Let $\fa$ be the object of $\Vect$ such that $A:=\on{Sym}(\cg[-2])^\cG\simeq \Sym(\fa)$. It is well-known 
(see, e.g., \cite[Example 6.5.5]{DrGa0}) that $\Dmod(\on{pt})^G$, equipped with the forgetful functor
\eqref{e:forget equiv}, identifies with the category $B\mod$, where $B=\Sym(\fa^*[-1])\mod$. 
In particular, the object $\wt\delta$ corresponds to $B$ itself, and $\delta$ corresponds to 
the augmentation $B\to k$.

\medskip

This makes the assertion of the lemma manifest, where $\fl$ is the graded line such that
$$B\simeq B^*\otimes \fl,$$
where $B^*$ is the linear dual of $B$ regarded as an object of $B\mod$.

\end{proof}

\sssec{}

Let $\wt\delta_\one$ denote the corresponding object of $\Sph(G,x)$ obtained via
$$\Dmod(\on{pt}/G)\simeq \Dmod(\on{pt})^G\simeq  \Dmod(\on{pt})^{G(\wh\CO_x)}\overset{\one}\hookrightarrow
\Dmod(\Gr_{G,x})^{G(\wh\CO_x)}=\Sph(G,x).$$
via the inclusion of the point $\one\in \Gr_{G,x}$. 

\medskip

By construction,
$$\wt\delta_\one\simeq \on{Av}_{G(\wh\CO_x),!}(\delta_\one),$$
so from \eqref{e:averaging and convolution}, we obtain that the category $\Sph(G,x)$
is compactly generated by objects of the form
$$\wt\delta_\one\star \CM$$
for $\CM\in (\Sph(G,x)^{loc.c})^\heartsuit$. Such $\CM$ are of the form $\on{Sat}^{naive,\on{ren}}(\rho)$ for
$\rho\in (\Rep(\cG)^c)^\heartsuit$, by the construction of $\on{Sat}^{naive,\on{ren}}(\rho)$.

\sssec{}

By \eqref{e:delta tilde} and Sects. \ref{sss:Bezr expl} and \ref{sss:equiv cohomology}, for 
$\rho\in \Rep(\cG)^c$ we have
$$\on{KD}_{\on{Hecke}(\cG)_{spec}}\circ (\on{Sat}^{\on{ren}})^{-1}\left(\wt\delta_\one\star \on{Sat}^{naive,\on{ren}}(\rho)\right)
\simeq \left(\on{Sym}(\cg[-2])\underset{\on{Sym}(\cg[-2])^\cG}\otimes \fl\right)\otimes \rho,$$
regarded as an object of $(\on{Sym}(\cg[-2])\mod)^\cG$. 

\medskip

So, the essential image of $\Sph(G,x)$ under $\on{KD}_{\on{Hecke}(\cG)_{spec}}\circ (\on{Sat}^{\on{ren}})^{-1}$
is compactly generated by objects of form
$$\left(\on{Sym}(\cg[-2])\underset{\on{Sym}(\cg[-2])^\cG}\otimes k\right)\otimes \rho,\quad \rho\in \Rep(\cG)^c.$$

However, as
$$\Sym(\cg)\underset{(\on{Sym}(\cg))^\cG}\otimes k\simeq \CO_{\on{Nilp}(\cg^*)},$$
it is clear that the subcategory generated by such objects is
exactly 
$$(\on{Sym}(\cg[-2])\mod_{\on{Nilp}(\cg*)})^\cG.$$

\qed[\thmref{t:Satake}]

\ssec{The action on $\Dmod(\Bun_G)$ and $\IndCoh_{\on{Nilp}_{glob}}(\LocSysG)$}

Recall that the monoidal category $\Sph(G,x)$ canonically acts on $\Dmod(\Bun_G)$. In this subsection we 
will study the corresponding action on the spectral side. 

\sssec{}

First, we claim that the monoidal category $\IndCoh(\on{Hecke}(\cG)_{spec})$ canonically acts on 
$\IndCoh(\LocSysG)$. This follows from the fact that we have a commutative diagram
in which both parallelograms are Cartesian:
\begin{gather}  \label{e:Hecke action spectral}
\xy
(-25,0)*+{\LocSysG}="X";
(25,0)*+{\LocSysG}="Y";
(0,15)*+{\LocSysG\underset{\LocSysG^{\on{R.S.}}}\times \LocSysG}="Z";
(-25,-40)*+{\on{pt}/\cG}="X_1";
(25,-40)*+{\on{pt}/\cG,}="Y_1";
(0,-25)*+{\on{pt}/\cG\underset{\cg/\cG}\times \on{pt}/\cG}="Z_1";
{\ar@{->}"Z";"X"};
{\ar@{->}"Z";"Y"};
{\ar@{->}"Z_1";"X_1"};
{\ar@{->}"Z_1";"Y_1"};
{\ar@{->}"Z";"Z_1"};
{\ar@{->}"X";"X_1"};
{\ar@{->}"Y";"Y_1"};
\endxy
\end{gather}
indeed, this a special case of diagram \eqref{e:groupoids morphism stacks}.

\sssec{}

The next proposition shows that the several different ways to define an action of 
the monoidal category $\IndCoh_{\on{Nilp}(\cg^*)/\cG}(\on{Hecke}(\cG)_{spec})$ on
$\IndCoh_{\on{Nilp}_{glob}}(\LocSysG)$ give the same result. 

\begin{prop}  \label{p:Hecke action spec} \hfill

\smallskip

\noindent{\em(a)}
For any conical Zariski-closed subset $Y\subset\ArthG$, the action of the monoidal category $\IndCoh(\on{Hecke}(\cG)_{spec})$
sends $\IndCoh_{Y}(\LocSysG)$ to $\IndCoh_{Y}(\LocSysG)$. Moreover, the diagram
$$
\CD
\IndCoh(\on{Hecke}(\cG)_{spec})\otimes \IndCoh_{Y}(\LocSysG)  @>{\text{action}}>> \IndCoh_{Y}(\LocSysG)   \\
@A{\on{Id}\otimes \Psi^{Y,\on{all}}_{\LocSysG}}AA    @AA{\Psi^{Y,\on{all}}_{\LocSysG}}A  \\
\IndCoh(\on{Hecke}(\cG)_{spec})\otimes \IndCoh(\LocSysG) @>{\text{action}}>>  \IndCoh(\LocSysG)
\endCD
$$
commutes as well \emph{(}i.e., the functor $\Psi^{Y,\on{all}}$, which is a priori lax compatible with the
action of $\IndCoh(\on{Hecke}(\cG)_{spec})$, is strictly compatible\emph{)}.

\smallskip

\noindent{\em(b)} The action of $\IndCoh_{\on{Nilp}(\cg^*)/\cG}(\on{Hecke}(\cG)_{spec})$
sends $\IndCoh(\LocSysG)$ to the subcategory $\IndCoh_{\on{Nilp}_{glob}}(\LocSysG)$. 

\smallskip

\noindent{\em(c)} The composed functor
\begin{multline*}
\IndCoh(\on{Hecke}(\cG)_{spec})\otimes \IndCoh(\LocSysG)\overset{\text{action}}\longrightarrow \\
\IndCoh(\LocSysG)\overset{\text{colocalization}}\longrightarrow 
\IndCoh_{\on{Nilp}_{glob}}(\LocSysG)
\end{multline*}
factors through the colocalization 
\begin{multline*}
\IndCoh(\on{Hecke}(\cG)_{spec})\otimes \IndCoh(\LocSysG)\to \\
\to \IndCoh_{\on{Nilp}(\cg^*)/\cG}(\on{Hecke}(\cG)_{spec})\otimes \IndCoh_{\on{Nilp}_{glob}}(\LocSysG).
\end{multline*}
\end{prop}

\begin{proof}
To prove point (a), it is enough to do so on the generators of $\IndCoh(\on{Hecke}(\cG)_{spec})$,
i.e., on the essential image of 
$$(\Delta_{\on{pt}/\cG})_*^\IndCoh:\IndCoh(\on{pt}/\cG)\to \IndCoh(\on{pt}/\cG\underset{\cg/\cG}\times \on{pt}/\cG).$$
However, for $\CF\in \IndCoh(\on{pt}/\cG)$, the action of $(\Delta_{\on{pt}/\cG})_*^\IndCoh(\CF)$ on 
$\IndCoh(\LocSysG)$ is given by tensor product with the pullback of 
$$\CF\in \IndCoh(\on{pt}/\cG)\simeq \QCoh(\LocSysG)$$ under the map
$$\LocSysG\to \on{pt}/\cG,$$
corresponding to the point $x$. Hence, the assertion follows from \corref{c:tensor stacks}.

\medskip

Points (b) and (c) are a particular case of \corref{p:Hecke action abstract stacks}.
\end{proof} 

As a corollary, we obtain:

\begin{cor}  \label{c:Hecke action spec}
There exists a canonically defined action of 
$\IndCoh_{\on{Nilp}(\cg^*)/\cG}(\on{Hecke}(\cG)_{spec})$ on $\IndCoh_{\on{Nilp}_{glob}}(\LocSysG)$,
which is compatible with the $\IndCoh(\on{Hecke}(\cG)_{spec})$-action on $\IndCoh(\LocSysG)$
via any of the functors
$$\Xi_{\LocSysG}^{\on{Nilp}_{glob}}:\IndCoh_{\on{Nilp}_{glob}}(\LocSysG)\rightleftarrows \IndCoh(\LocSysG):\Psi_{\LocSysG}^{\on{Nilp}_{glob}}$$
and
$$\Xi^{\on{Nilp}(\cg^*)/\cG,\on{all}}_{\on{Hecke}(\cG)_{spec}}:
\IndCoh_{\on{Nilp}(\cg^*)/\cG}(\on{Hecke}(\cG)_{spec})\rightleftarrows \IndCoh(\on{Hecke}(\cG)_{spec}):
\Psi^{\on{Nilp}(\cg^*)/\cG,\on{all}}_{\on{Hecke}(\cG)_{spec}}.$$
\end{cor}

\sssec{}

The compatibility of \conjref{conj:main} with the Geometric Satake Equivalence reads:

\begin{conj} \label{conj:comp with Satake}
The action of $\Sph(G,x)$ on $\Dmod(\Bun_G)$ corresponds via
$$\on{Sat}:\IndCoh_{\on{Nilp}(\cg^*)/\cG}(\on{Hecke}(\cG)_{spec})\simeq \Sph(G,x)$$
to the action of $\IndCoh_{\on{Nilp}(\cg^*)/\cG}(\on{Hecke}(\cG)_{spec})$ on 
$\IndCoh_{\on{Nilp}_{glob}}(\LocSysG)$.
\end{conj}

\begin{rem}
The above conjecture is not yet the full compatibility of the Geometric Satake Equivalence with the Geometric
Langlands equivalence. The full version amounts to formulating \conjref{conj:comp with Satake} in a way
that takes into account the factorizable structure of $$\IndCoh_{\on{Nilp}(\cg^*)/\cG}(\on{Hecke}(\cG)_{spec})\simeq \Sph(G,x)$$
as $x$ moves along $X$.
\end{rem}

\ssec{Singular support via the Hecke action}

\sssec{}

If \conjref{conj:main} holds, an object $\CM\in\Dmod(\Bun_G)$ can be assigned its singular support, 
which by definition is equal to
the singular support of the corresponding object of $\IndCoh_{\on{Nilp}_{glob}}(\LocSysG)$. The singular support
is a conical Zariski-closed subset $\on{SingSupp}(\CM)\subset\on{Nilp}_{glob}$. 

\medskip

It turns out that \conjref{conj:comp with Satake} implies
certain relation between $\on{SingSupp}(\CM)$ and the action of the Hecke category $\Sph(G,x)$ on $\CM$.
Let us explain this in more detail.

\sssec{}

The equivalence
$$\on{KD}_{\on{Hecke}(\cG)_{spec}}:\IndCoh(\on{Hecke}(\cG)_{spec})\simeq (\on{Sym}(\cg[-2])\mod)^\cG$$
of \propref{p:Koszul dual for Hecke} and \propref{p:Hecke action spec}(a) makes 
the categories $\IndCoh(\LocSysG)$ and $\IndCoh_{\on{Nilp}_{glob}}(\LocSysG)$ into categories
tensored over $\QCoh(\cg^*/(\cG\times \BG_m))$. 

\medskip

By construction, this is the same structure as that given by the embedding 
$$\iota:\LocSysG\hookrightarrow \LocSysG^{\on{R.S.}}$$
in terms of \secref{sss:tensored over conormal stacks parallel}. 

\medskip

Thus, we can determine the singular support of objects of $\IndCoh(\LocSysG)$ via the above
action of $\IndCoh(\on{Hecke}(\cG)_{spec})$.

\sssec{} 
On the other hand, 
\corref{c:Satake Koszul as coloc} defines on $\Dmod(\Bun_G)$ a structure of category
tensored over $$\QCoh(\cg^*/(\cG\times \BG_m)).$$
 Hence, we can attach to an object $\CM\in\Dmod(\Bun_G)$ its support
 \[\on{supp}^{x}(\CM)\subset\cg^*/(\cG\times\BG_m),\] which is a Zariski-closed subset.
 (The superscript $x$ indicates that this support depends on the choice of the point $x\in X$.)

\medskip

Conjectures~\ref{conj:main} and \ref{conj:comp with Satake} would imply that $\on{supp}^{x}(\CM)$
is the Zariski closure of the image of
\[\on{SingSupp}(\CM)\subset\ArthG=\Sing(\LocSysG)\]
under the map
\[\ArthG\to\cg^*/(\cG\times\BG_m):(\CP,\nabla,A)\mapsto A(x).\]
(This easily follows from \lemref{l:change of algebras}.)
Here we use the explicit description of $\ArthG$ given in \corref{c:Sing of LocSys}.

\sssec{} In particular, consider the full subcategory
\[\Dmod^x_{\on{temp}}(\Bun_G):=\{\CM\in\Dmod(\Bun_G):\on{supp}^x(\CM)=\{0\}\}.\]

Equivalently, we can define it as the tensor product
\begin{equation}\label{e:temp as tensor}
\Dmod^x_{\on{temp}}(\Bun_G)=\Dmod(\Bun_G)\underset{\QCoh(\cg^*/(\cG\times \BG_m))}\otimes\QCoh(\cg^*/(\CG\times\BG_m))_{\{0\}}.
\end{equation}
Conjectures~\ref{conj:main} and \ref{conj:comp with Satake} imply that under the equivalence
\[\Dmod(\Bun_G)\simeq\IndCoh_{\on{Nilp}_{glob}}(\LocSysG),\]
the category $\Dmod^x_{\on{temp}}(\Bun_G)$ corresponds to the subcategory
\[\IndCoh_{\{0\}}(\LocSysG)\subset\IndCoh_{\on{Nilp}_{glob}}(\LocSysG),\]
which is the same as the essential image of $\QCoh(\LocSysG)$ under the functor
$$\Psi_{\LocSysG}:\QCoh(\LocSysG)\to \IndCoh(\LocSysG).$$
Thus, $\Dmod^x_{\on{temp}}(\Bun_G)$ should
be equal to the subcategory
$$\Dmod_{\on{temp}}(\Bun_G)\subset \Dmod(\Bun_G)$$
of \secref{sss:tempered}.

\medskip

In particular, we obtain:
\begin{conj}\label{conj:tempered independence}
The subcategory $\Dmod^x_{\on{temp}}(\Bun_G)\subset \Dmod(\Bun_G)$
is independent of the choice of the point $x\in X$. 
\end{conj}

\sssec{}

Let us provide a more explicit description of the subcategory $\Dmod^x_{\on{temp}}(\Bun_G)$. 

\medskip

Recall that $\Dmod(\Bun_G)$ is compactly generated (see \cite{DrGa1}). Now \eqref{e:temp as tensor}
implies that $\Dmod^x_{\on{temp}}(\Bun_G)$ is compactly generated by 
\[\Dmod^x_{\on{temp}}(\Bun_G)^c=\Dmod^x_{\on{temp}}(\Bun_G)\cap\Dmod(\Bun_G)^c.\] 
For this reason, it suffices to describe the compact objects of $\Dmod_{\on{temp}}(\Bun_G)$.

\medskip

By \thmref{t:Bezr}, there exists a canonical map in $\Sph(G,x)$:
$$\alpha:\Sat^{naive}(\cg)[-2]\to \Sat^{naive}(k)=\delta_{\one},$$
\medskip
where $k\in \Rep(\cG)$ is the trivial representation.

\medskip

Moreover, for any $n\in \BN$ we can consider its ``symmetric power"
$$\alpha_n:\Sat(\Sym^n(\cg))[-2n]\to \delta_{\one}.$$

From \lemref{l:supp via comp}(c), we obtain:

\begin{cor}  \label{c:temp}
For $\CM\in \Dmod(\Bun_G)^c$, the following conditions are equivalent:

\smallskip

\noindent {\em (a)} $\CM\in\Dmod^x_{\on{temp}}(\Bun_G)$;

\noindent {\em (b)} 
The induced map
$$(\alpha_n\star{\on{id}_\CM}):\Sat^{naive}(\Sym^n(\cg))\star \CM[-2n]\to \CM$$
vanishes for some integer $n\ge 0$;

\noindent {\em (c)} The map $\alpha_n\star{\on{id}_\CM}$
vanishes for all sufficiently large integers $n$. 
\end{cor}

\begin{rem}
One can show \emph{unconditionally} that when we view $\Dmod(\Bun_G)$ as tensored over $\QCoh(\cg^*/(\cG\times \BG_m))$, i.e., 
as of category over the stack $\cg^*/(\cG\times \BG_m)$, it is supported over
$$\on{Nilp}(\cg^*)/(\cG\times \BG_m)\subset \cg^*/(\cG\times \BG_m).$$

\medskip

This fact is equivalent to the following. In the situation of \corref{c:temp} consider the map of $\cG$-representations 
$$\on{Sym}^n(\cg)^{\cG}\otimes k\to \on{Sym}^n(\cg),$$
where $\on{Sym}^n(\cg)^{\cG}$ is regarded as a mere vector space and $k$ is the trivial representation. 
Composing, for $\CM\in  \Dmod(\Bun_G)$,  we obtain a map 
\begin{equation} \label{e:action of invaiants}
\on{Sym}^n(\cg)^{\cG}\otimes \CM[-2n]\to \CM.
\end{equation}
We claim that this map vanishes for any $\CM\in \Dmod(\Bun_G)^c$ whenever $n$ is sufficiently large. 

\medskip

To prove this, we note that by \secref{sss:equiv cohomology}, the map \eqref{e:action of invaiants} comes
from the composition
$$\on{Sym}^n(\cg)^{\cG}\hookrightarrow H^{2n}_\dr(\on{pt}/G)\to H^{2n}_\dr(\Bun_G),$$
where $H_\dr(\on{pt}/G)\to H_\dr(\Bun_G)$ is the homomorphism given by the map $\Bun_G\to \on{pt}/G$,
given, corresponding to the restriction along $x\to X$. 

\medskip

Now, one can show that for any (connected) algebraic stack $\CY$ and any $\CM\in \Dmod(\CY)^c$, the map
$$H_\dr(\CY)\otimes \CM\to \CM$$
vanishes on a sufficiently high power of the augmentation ideal of $H_\dr(\CY)$. 

\end{rem}

\sssec{}

Using \lemref{l:supp via comp}(a) one can give the following characterization of the entire subcategory
$\Dmod^x_{\on{temp}}(\Bun_G)\subset \Dmod(\Bun_G)$ (and not just its compact objects):

\begin{lem}
An object $\CM\in  \Dmod(\Bun_G)$ belongs to $\Dmod^x_{\on{temp}}(\Bun_G)$ of and only for
a set of compact generators $\CM_\alpha\in \Dmod(\Bun_G)$ and any $\CM_\alpha\to \CM$,
for all sufficiently large $n$ the composition
$$\CM_\alpha\to \CM\to \Sat^{naive}(\Sym^n(\cg^*))\star \CM[2n]$$
vanishes.
\end{lem}

From here we obtain:

\begin{cor}
The constant sheaf $k_{\Bun_G}\in \Dmod(\Bun_G)$ is \emph{not} tempered.
\end{cor}

\begin{proof}
Follows from the fact that for any $n$, the map 
$$k_{\Bun_G}[-2n]\to \Sat^{naive}(\Sym^n(\cg^*))\star k_{\Bun_G}
\simeq H\left(\Gr_{G,x},\Sat^{naive}(\Sym^n(\cg^*))\right)\otimes k_{\Bun_G}$$
is the inclusion of a direct summand.
\end{proof}

\section{Compatibility with Eisenstein series}  \label{s:eis}

A crucial ingredient in the formulation of the Geometric Langlands equivalence is the interaction of $G$
with its Levi subgroups. Such interaction is given by the functors of Eisenstein series on both
sides of the correspondence. In this section we will study how these functors act on our
category $\IndCoh_{\on{Nilp}^G_{glob}}(\LocSysG)$. 

\ssec{The Eisenstein series functor on the geometric side}

\sssec{}\label{sss:Eis!}

Let $P$ be a parabolic subgroup of $G$ with the  Levi quotient $M$. Let us recall the definition of
the Eisenstein series functor
$$\Eis^P_!:\Dmod(\Bun_M)\to \Dmod(\Bun_G)$$
(see \cite{DrGa2}).

\medskip

By definition, 
$$\Eis^P_!=(\sfp^P)_!\circ (\sfq^P)^*,$$
where $\sfp^P$ and $\sfq^P$ are the maps in the diagram
\begin{gather}  \label{basic diag}
\xy
(-15,0)*+{\Bun_G}="X";
(15,0)*+{\Bun_M.}="Y";
(0,15)*+{\Bun_P}="Z";
{\ar@{->}_{\sfp^P} "Z";"X"};
{\ar@{->}^{\sfq^P} "Z";"Y"};
\endxy
\end{gather}

We note that the functor $(\sfq^P)^*$ is defined because the morphism $\sfq^P:\Bun_P\to \Bun_M$ is smooth,
and that the functor $(\sfp^P)_!$, left adjoint to $(\sfp^P)^!$, is defined on the essential image of $(\sfq^P)^*$,
as is shown in \cite[Proposition~1.2]{DrGa2}.

\medskip

Note that the functor $\Eis^P_!$ sends compact objects in $\Dmod(\Bun_M)^c$ to compact objects
in $\Dmod(\Bun_G)^c$, since it admits a continuous right adjoint 
\[\on{CT}_*^P=(\sfq^P)_*\circ (\sfp^P)^!.\]

\sssec{}\label{sss:cusp}

Let $\Dmod_{\on{Eis}}(\Bun_G)$ be the full subcategory of $\Dmod(\Bun_G)$ generated by the
essential images of the functors $\Eis^P_!$ for all proper parabolic subgroups $P$. Let
$\Dmod_{\on{cusp}}(\Bun_G)$ denote the full subcategory of $\Dmod(\Bun_G)$ equal to the right
orthogonal of $\Dmod_{\on{Eis}}(\Bun_G)$. 

\medskip

Since the functors $\Eis^P_!$ preserve compactness, the category $\Dmod_{\on{Eis}}(\Bun_G)$ is compactly generated. 
Therefore, $\Dmod_{\on{cusp}}(\Bun_G)$ is
a localization of $\Dmod(\Bun_G)$ with respect to $\Dmod_{\on{Eis}}(\Bun_G)$, so we obtain a short exact sequence of DG
categories
\[\Dmod_{\on{Eis}}(\Bun_G)\rightleftarrows\Dmod(\Bun_G)\rightleftarrows\Dmod_{\on{cusp}}(\Bun_G).\]

\ssec{Eisenstein series on the spectral side}

\sssec{}

Fix a parabolic subgroup $P\subset G$, and consider the corresponding parabolic subgroup $\cP\subset \cG$, 
whose Levi quotient $\cM$ identifies with the Langlands dual of $M$. Consider the diagram:

\begin{gather}  \label{basic diag spec}
\xy
(-15,0)*+{\LocSysG}="X";
(15,0)*+{\LocSys_\cM.}="Y";
(0,15)*+{\LocSys_\cP}="Z";
{\ar@{->}_{\sfp^P_{spec}} "Z";"X"};
{\ar@{->}^{\sfq^P_{spec}} "Z";"Y"};
\endxy
\end{gather}

We define the functor
$$\Eis^{P}_{spec}:\IndCoh(\LocSys_{\cM})\to \IndCoh(\LocSys_\cG),$$
to be 
$$(\sfp^P_{spec})^{\IndCoh}_*\circ (\sfq^P_{spec})^!.$$

First, we note:

\begin{lem} \label{l:q and p} \hfill

\smallskip

\noindent {\em(a)} The map $\sfq^P_{spec}$ is quasi-smooth.

\smallskip

\noindent {\em(b)} The map $\sfp^P_{spec}$ is schematic and proper.
\end{lem}
\begin{proof}
For part (a), we claim that for any surjective homomorphism of algebraic groups $$\cG_1\to \cG_2,$$
the corresponding map $\LocSys_{\cG_1}\to \LocSys_{\cG_2}$ is quasi-smooth. This relative version of
\propref{p:LocSys q-smooth} can be proved in the same way as \propref{p:LocSys q-smooth}.
Part (b) is straightforward (and well known).
\end{proof}

\begin{cor}  \label{l:spectral Eis coarse}
The functor $\Eis^{P}_{spec}$ sends $\Coh(\LocSys_{\cM})$ to $\Coh(\LocSys_{\cG})$.
\end{cor}

\begin{proof}
First, we claim that the functor $(\sfq^P_{spec})^!$ sends
$\Coh(\LocSys_{\cM})$ to $\Coh(\LocSys_{\cP})$. This follows from 
\cite[Lemma 7.1.2]{IndCoh} and \lemref{l:q and p}(a).

\medskip

Now, $(\sfp^P_{spec})^{\IndCoh}_*$ sends $\Coh(\LocSys_{\cP})$ to $\Coh(\LocSys_{\cG})$
by \cite[Lemma 3.3.5]{IndCoh} and Lemma~\ref{l:q and p}(b). 
\end{proof}

\sssec{}
Let us now analyze the singular codifferential of the morphisms $\sfp^P_{spec}$ and $\sfq^P_{spec}$. 
To avoid confusion, let us introduce superscripts and write
\[\on{Nilp}^G_{glob}\subset\ArthG\qquad\text{and}\qquad\on{Nilp}^M_{glob}\subset\Arth_\cM\] to 
distinguish between the global nilpotent cones for $\cM$ and $\cG$.

\medskip

By \lemref{l:q and p}(a), the singular codifferential
$$\Sing(\sfq^P_{spec}):\Arth_\cM\underset{\LocSys_{\cM}}\times \LocSys_\cP\to 
\Arth_\cP$$
is a closed embedding. 
Consider the subset 
$$\on{Nilp}^M_{glob}\underset{\LocSys_\cM}\times \LocSys_\cP\subset 
\Arth_\cM\underset{\LocSys_\cM}\times \LocSys_\cP$$
and let
\[\on{Nilp}^P_{glob}:=\Sing(\sfq^P_{spec})\left(\on{Nilp}^M_{glob}\underset{\LocSys_\cM}\times \LocSys_\cP\right)\subset \Arth_\cP\]
be its image. Here is an explicit description:

\begin{lem}\label{l:Nilpcone for P}
Let us identify $\Arth_\cP$ with the moduli stack 
(in the classical sense) of triples $(\CP^\cP,\nabla,A^\cP\in H^0(\Gamma(X_\dr,\cp^*_{\CP^P}))$
using \corref{c:Sing of LocSys}. Then 
\[\on{Nilp}^P_{glob}=\{(\CP^\cP,\nabla,A^\cP),\,\,A^\cP\text{ is a nilpotent section of }\cm^*_{\CP^\cP}\subset \cp^*_{\CP^\cP}.\}.\]
\end{lem}
\begin{proof} Indeed, $\Arth_\cM\underset{\LocSys_{\cM}}\times \LocSys_\cP$
is identified with the classical moduli stack of triples $(\cP^\cP,\nabla,A^\cM)$, where
$(\CP^{\cP},\nabla)\in\LocSys_\cP$ and $A^\cM\in H^0(\Gamma(X_\dr,\cm^*_{\CP^\cP}))$. Here we use the natural action
of $\cP$ on $\cm$ (and $\cm^*$).

\medskip
Under this identification,
\[\Sing(\sfq^P_{spec})(\CP^\cP,\nabla,A^\cM)=(\CP^\cP,\nabla,A^\cP),\]
where $A^\cP$ is the image of $A^\cM$ under the natural embedding $\cm^*\hookrightarrow \cp^*$. The claim follows. 
\end{proof}

\begin{prop} \label{p:spectral Eis fine}
The functor $\Eis^{P}_{spec}$ sends the subcategory 
$$\IndCoh_{\on{Nilp}^M_{glob}}(\LocSys_{\cM})\subset \IndCoh(\LocSys_{\cM})$$
to the subcategory
$$\IndCoh_{\on{Nilp}^G_{glob}}(\LocSysG)\subset \IndCoh(\LocSysG)$$
\end{prop}

The proposition provides a commutative diagram of functors
$$
\CD
\IndCoh_{\on{Nilp}^M_{glob}}(\LocSys_{\cM})  @>>>  \IndCoh(\LocSys_{\cM})  \\
@VVV    @VV{\Eis^{P}_{spec}}V  \\
\IndCoh_{\on{Nilp}^G_{glob}}(\LocSysG)  @>>>  \IndCoh(\LocSysG),
\endCD$$
where the horizontal arrows are the tautological embeddings. In particular,
the resulting functor 
$$\Eis^{P}_{spec}:\IndCoh_{\on{Nilp}^M_{glob}}(\LocSys_{\cM}) \to 
\IndCoh_{\on{Nilp}^G_{glob}}(\LocSysG)$$
also sends compact objects to compact objects, that is, it restricts to a functor 
$$\Coh_{\on{Nilp}^M_{glob}}(\LocSys_{\cM})\to \Coh_{\on{Nilp}^G_{glob}}(\LocSysG).$$

\sssec{Proof of \propref{p:spectral Eis fine}}
By \lemref{l:pullback stacks}, we see that
\[(\sfq^P_{spec})^!\left(\IndCoh_{\on{Nilp}^M_{glob}}(\LocSys_\cM)\right)\subset\IndCoh_{\on{Nilp}^P_{glob}}(\LocSys_\cP).\]
Therefore, it is enough to check that
\[(\sfp^P_{spec})_*^\IndCoh\left(\IndCoh_{\on{Nilp}^P_{glob}}(\LocSys_\cP)\right)\subset\IndCoh_{\on{Nilp}^G_{glob}}(\LocSys_\cG).\]
By \lemref{l:direct image stacks}, it suffices to show that
the preimage of $\on{Nilp}^P_{glob}$ under the singular codifferential 
$$\Sing(\sfp^P_{spec}):\ArthG\underset{\LocSys_{\cG}}\times \LocSys_\cP\to \Arth_\cP$$
is contained in 
$$\on{Nilp}^G_{glob}\underset{\LocSys_{\cG}}\times \LocSys_\cP\subset
\ArthG\underset{\LocSys_{\cG}}\times \LocSys_\cP.$$ 

\medskip

Using \corref{c:Sing of LocSys}, we can identify $\ArthG\underset{\LocSys_{\cG}}\times \LocSys_\cP$
with the classical moduli stack of triples $(\CP^\cP,\nabla,A^\cG)$, where $(\CP^\cP,\nabla)\in\LocSys_\cP$ and
$A^\cG\in H^0(\Gamma(X_\dr,\cg^*_{\CP^\cP}))$. Under this identification, $\Sing(\sfp^P_{spec})$ sends such a triple to
the triple
$$(\CP^\cP,\nabla,A^\cP)\in H^0(\Gamma(X_\dr,\cp^*_{\CP^P})),$$
where $\A^\cP$ is obtained from $A^\cG$ via the natural projection $\cg^*\to \cp^*$. 

\medskip

It remains to notice that if $a\in\cg^*$ is such that its projection to $\cp^*$ is a 
nilpotent element of $\cm^*\subset\cp^*$, then $a$ itself is nilpotent.
\qed

\sssec{Compatibility between Geometric Langlands Correspondence and Eisenstein series}

The following is one of the key requirements on the equivalence of \conjref{conj:main}:

\begin{conj}  \label{conj:compat}
For every parabolic $P$ the following diagram of functors
$$
\CD
\Dmod(\Bun_G)      @>>>  \IndCoh_{\on{Nilp}^G_{glob}}(\LocSysG)  \\
@A{\Eis^P_!}AA    @AA{\Eis^{P}_{spec}}A  \\
\Dmod(\Bun_M)     @>>>  \IndCoh_{\on{Nilp}^M_{glob}}(\LocSys_\cM) 
\endCD
$$
commutes, up to an auto-equivalence of $\IndCoh_{\on{Nilp}^M_{glob}}(\LocSys_\cM)$
given by tensoring by a line bundle. 
\end{conj}

\ssec{The main result}

\sssec{}

Let $\LocSysG^{\on{red}}$ denote the closed substack of $\LocSysG$ equal to the union of the
images of the maps $\sfp^P_{spec}$ for all proper parabolics $P$, considered, say, with the reduced
structure. Let $\LocSysG^{\on{irred}}$
be the complementary open; we denote by $\jmath$ the open embedding 
$$\LocSysG^{\on{irred}}\hookrightarrow\LocSysG.$$

\medskip

By \corref{c:on open stacks} we obtain a diagram of short exact sequences of DG categories
\begin{equation} \label{e:localization on LocSys}
\CD
\IndCoh_{(\on{Nilp}^G_{glob})_{\LocSysG^{\on{red}}}}(\LocSysG)
@>>> \IndCoh(\LocSysG)_{\LocSysG^{\on{red}}}  \\
@VVV  @VVV  \\
\IndCoh_{\on{Nilp}^G_{glob}}(\LocSysG)  @>>> \IndCoh(\LocSysG) \\
@V{\jmath^{\IndCoh,*}}VV  @VV{\jmath^{\IndCoh,*}}V \\
\IndCoh_{\on{Nilp}^G_{glob}}(\LocSysG^{\on{irred}})   @>>>  \IndCoh(\LocSysG^{\on{irred}}) 
\endCD
\end{equation}
obtained from $\IndCoh_{\on{Nilp}^G_{glob}}(\LocSysG)\hookrightarrow \IndCoh(\LocSysG)$ by tensoring over $\QCoh(\LocSysG)$ with
the short exact sequence
$$\QCoh(\LocSysG)_{\LocSysG^{\on{red}}}  \rightleftarrows  \QCoh(\LocSysG)  \rightleftarrows  \QCoh(\LocSysG^{\on{irred}}).$$ 

In particular, the vertical arrows in the diagram \eqref{e:localization on LocSys} admit right adjoints,
and the horizontal arrows are fully faithful embeddings. Moreover, all the categories involved are compactly generated;
in particular,
$$\IndCoh(\LocSysG)_{\LocSysG^{\on{red}}}  \text{ and } \IndCoh_{(\on{Nilp}^G_{glob})_{\LocSysG^{\on{red}}}}(\LocSysG)$$
are compactly generated by 
$$\Coh(\LocSysG)_{\LocSysG^{\on{red}}} \text{ and } \Coh(\LocSysG)_{\LocSysG^{\on{red}}} \cap \Coh_{\on{Nilp}^G_{glob}}(\LocSysG),$$
respectively.

\sssec{}

We have:

\begin{prop}  \label{p:on irred}
The inclusion
$$\QCoh(\LocSysG^{\on{irred}})\hookrightarrow \IndCoh_{\on{Nilp}^G_{glob}}(\LocSysG^{\on{irred}})$$
is an equality. 
\end{prop}

\begin{proof}

This follows from \corref{c:zero sing supp stacks} using the following observation:

\begin{lem}
The preimage of $\LocSysG^{\on{irred}}$ in $\on{Nilp}^G_{glob}$ consists of the zero-section.
\end{lem}

\begin{proof}
Indeed, an irreducible $\cG$-local system admits no non-trivial horizontal
nilpotent sections of the associated bundle of Lie algebras.
\end{proof}

\end{proof}

\sssec{}

We are now ready to state the main result of this paper:

\begin{thm} \label{t:generation}
The subcategory 
$$\IndCoh_{(\on{Nilp}^G_{glob})_{\LocSysG^{\on{red}}}}(\LocSysG) \subset \IndCoh_{\on{Nilp}^G_{glob}}(\LocSysG)$$
is generated by the essential images of the functors  
$$\Eis^{P}_{spec}:\IndCoh_{\on{Nilp}^M_{glob}}(\LocSys_{\cM}) \to 
\IndCoh_{\on{Nilp}^G_{glob}}(\LocSysG)$$
for all proper parabolics $P$.
\end{thm}

\sssec{}

Note that from
\thmref{t:generation}, combined with \propref{p:on irred} and \eqref{e:localization on LocSys}, we obtain:

\begin{cor}
The subcategory $\QCoh(\LocSys_\cG)$ and the essential images of 
$$\Eis^{P}_{spec}|_{\IndCoh_{\on{Nilp}^M_{glob}}(\LocSys_{\cM})}$$ 
for all \emph{proper} parabolics $P$, generate the category $\IndCoh_{\on{Nilp}^G_{glob}}(\LocSysG)$.
\end{cor}

Now, the transitivity property of Eisenstein series and induction on the semi-simple rank imply:

\begin{cor} \label{c:generate by QCoh}
The subcategory $\QCoh(\LocSys_\cG)$, together with the essential images of the subcategories
$\QCoh(\LocSys_\cM)\subset \IndCoh_{\on{Nilp}^M_{glob}}(\LocSys_\cM)$ under the functors $\Eis^{P}_{spec}$ 
for all \emph{proper} parabolics $P$, generate $\IndCoh_{\on{Nilp}^G_{glob}}(\LocSysG)$.
\end{cor}

Still equivalently, we have:
\begin{cor}
The essential images of $\QCoh(\LocSys_\cM)\subset \IndCoh_{\on{Nilp}^M_{glob}}(\LocSys_\cM)$ 
under the functors $\Eis^{P}_{spec}$ 
for \emph{all} parabolic subgroups $P$  (including the case $P=G$) generate $\IndCoh_{\on{Nilp}^G_{glob}}(\LocSysG)$.
\end{cor}

\sssec{}

Let us explain the significance of this theorem from the point of view of \conjref{conj:main}. 
We are going to 
show that $\IndCoh_{\on{Nilp}^G_{glob}}(\LocSysG)$ is the \emph{smallest} subcategory of $\IndCoh(\LocSysG)$
that contains $\QCoh(\LocSysG)$, which can be equivalent to $\Dmod(\Bun_G)$, if we 
assume compatibility with the Eisenstein series as in \conjref{conj:compat}.

\medskip

More precisely, let us assume that there exists an equivalence between $\Dmod(\Bun_G)$ and \emph{some} subcategory
$$\QCoh(\LocSysG)  \subset \IndCoh_?(\LocSysG) \subset \IndCoh(\LocSysG),$$
for which the diagrams 
\begin{equation} \label{e:compat ?}
\CD
\Dmod(\Bun_G)      @<{\sim}<< \IndCoh_?(\LocSysG)  \\
@A{\Eis^P_!}AA    @AA{\Eis^{P}_{spec}}A   \\
\Dmod(\Bun_M)     @<<<  \QCoh(\LocSys_\cM) 
\endCD
\end{equation}
commute for all proper parabolics $P$, up to tensoring by a line bundle as
in \conjref{conj:compat}.

\medskip

We claim that in this case, $\IndCoh_?(\LocSysG)$ necessarily contains $\IndCoh_{\on{Nilp}^G_{glob}}(\LocSysG)$.
Indeed, this follows from \corref{c:generate by QCoh}.

\sssec{}

Recall (\secref{sss:cusp}) that $\Dmod_{\on{cusp}}(\Bun_G)\subset \Dmod(\Bun_G)$ is the full subcategory
of D-modules that are right orthogonal to the essential images of the Eisenstein series functors $\Eis^P_!$
for all proper parabolic subgroups $P\subset G$.
Let us note the following corollary of \conjref{conj:compat} and \thmref{t:generation}:

\begin{cor} \label{c:cuspidal}
The equivalence of \conjref{conj:main} gives rise to an equivalence
$$\Dmod_{\on{cusp}}(\Bun_G)\simeq \QCoh(\LocSysG^{\on{irred}}).$$
\end{cor}


\ssec{Proof of \thmref{t:generation}}

\sssec{}

We will prove a more precise result. For a given parabolic $P$, let $\LocSysG^{\on{red}_P}$ denote the closed
substack of $\LocSysG$ equal to the image of the map $\sfp^P_{spec}$, considered, say, with the reduced structure. 
Let
$$\IndCoh_{(\on{Nilp}^G_{glob})_{\LocSysG^{\on{red}_P}}}(\LocSysG)$$
denote the corresponding full subcategory of $\IndCoh_{\on{Nilp}^G_{glob}}(\LocSysG)$.

\medskip

Clearly, the functor 
$$\Eis^{P}_{spec}:\IndCoh_{\on{Nilp}^M_{glob}}(\LocSys_{\cM}) \to 
\IndCoh_{\on{Nilp}^G_{glob}}(\LocSysG)$$
factors through $\IndCoh_{(\on{Nilp}^G_{glob})_{\LocSysG^{\on{red}_P}}}(\LocSysG)$. 

\medskip

We will prove:
\begin{thm}  \label{t:generation for P}
The essential image of the functor
$$\Eis^{P}_{spec}:\IndCoh_{\on{Nilp}^M_{glob}}(\LocSys_{\cM}) \to 
\IndCoh_{(\on{Nilp}^G_{glob})_{\LocSysG^{\on{red}_P}}}(\LocSysG)$$
generates the target category.
\end{thm}

\thmref{t:generation for P} implies \thmref{t:generation} by \corref{c:intersection of closed}.

\sssec{}

\thmref{t:generation for P} follows from the combination of the following two statements:

\begin{prop}  \label{p:generation M P}
The essential image of the functor 
$$(\sfq^P_{spec})^!:\IndCoh_{\on{Nilp}^M_{glob}}(\LocSys_\cM)\to \IndCoh_{\on{Nilp}^P_{glob}}(\LocSys_\cP)$$
generates the target category.
\end{prop}

\begin{prop}  \label{p:generation P G}
The essential image of the functor 
$$(\sfp^P_{spec})^{\IndCoh}_*:\IndCoh_{\on{Nilp}^P_{glob}}(\LocSys_\cP)\to 
\IndCoh_{(\on{Nilp}^G_{glob})_{\LocSysG^{\on{red}_P}}}(\LocSys_\cG)$$
generates the target category.
\end{prop}

\sssec{Proof of \propref{p:generation M P}}
Recall that the map $\sfq^P_{spec}$ is quasi-smooth, so that its singular codifferential
$$\Sing(\sfq^P_{spec}):\Arth_{\cM}\underset{\LocSys_{\cM}}\times \LocSys_\cP\to 
\Arth_\cP$$
is a closed embedding.  Moreover, $\on{Nilp}^P_{glob}\subset \Arth_\cP$ is equal to the image
of the closed subset
\[\on{Nilp}^M_{glob}\underset{\LocSys_\cM}\times \LocSys_\cP\subset\Arth_\cM\underset{\LocSys_\cM}\times \LocSys_\cP\]
under $\Sing(\sfq^P_{spec})$. Therefore, \propref{p:quasi-smooth pullback ten prod stacks} implies that $(\sfq^P_{spec})^!$ 
induces an equivalence
\[
\QCoh(\LocSys_\cP)\underset{\QCoh(\LocSys_\cM)}\otimes \IndCoh_{\on{Nilp}_{glob}^M}(\LocSys_\cM)\to 
\IndCoh_{\on{Nilp}_{glob}^P}(\LocSys_\cP).
\]

\medskip 

It remains to show that the essential image of the usual pullback functor
$$(\sfq^P_{spec})^*:\QCoh(\LocSys_{\cM})\to\QCoh(\LocSys_{\cP})$$
generates the target category. 
\medskip

Since $(\sfq^P_{spec})^*$ is the left adjoint to $(\sfq^P_{spec})_*$, we need to show that 
the pushforward functor
$$(\sfq^P_{spec})_*:\QCoh(\LocSys_{\cP})\to \QCoh(\LocSys_{\cM})$$ is
conservative. But this is true because the map $\sfq^P_{spec}$  
can be presented as a quotient of an schematic affine map by an action of a unipotent group-scheme
(i.e., $\sfq^P_{spec}$ is cohomologically affine). 

\qed[\propref{p:generation M P}]

\sssec{Proof of \propref{p:generation P G}}

We will deduce the proposition from \propref{p:prop cons stacks}.

\medskip

We need to show that the map 
$$(\Sing(\sfp^P_{spec}))^{-1}(\on{Nilp}^P_{glob})\to \on{Nilp}^G_{glob}$$
is surjective at the level of $k$-points. 

\medskip

Concretely, this means the following: let $(\CP^\cG,\nabla)$ be a $\cG$-bundle on $X$, equipped with a connection
$\nabla$, which admits a horizontal reduction to the parabolic $\cP$. Let $A$ be a horizontal section of $\cg_{\CP^\cG}$,
which is nilpotent (here we have chosen some $\cG$-invariant identification of $\cg$ with $\cg^*$). We need to show 
that there exists a horizontal reduction of $\CP^\cG$ to $\cP$ such that $A$ belongs to $\cp_{\CP^\cG}$. 

\medskip

Let $\on{Sect}^\nabla(X,\CP^\cG/\cP)$ be the (classical) scheme of all horizontal reductions of $\CP^\cG$ to $\cP$. By assumption,
this scheme is non-empty, and it is also proper, since it embeds as a closed subscheme into $\CP^\cG_x/\cP$ for any/some 
$x\in X$. 

\medskip

The algebraic group $\on{Aut}(\CP^\cG,\nabla)$ acts naturally on $\on{Sect}^\nabla(X,\CP^\cG/\cP)$. Note that the Lie
algebra of $\on{Aut}(\CP^\cG,\nabla)$ identifies with $H^0(\Gamma(X_\dr,\cg_{\CP^\cG}))$. 
By assumption, the element $$A\in H^0(\Gamma(X_\dr,\cg_{\CP^\cG}))$$ is nilpotent (as a linear operator on the algebra of 
functions on $\on{Aut}(\CP^\cG,\nabla)$). Hence, it comes from a homomorphism $\BG_a\to \on{Aut}(\CP^\cG,\nabla)$. 

\medskip

By properness, the resulting action
of $\BG_a$ on $\on{Sect}^\nabla(X,\CP^\cG/\cP)$ has a fixed
point, which is the desired reduction.\footnote{This argument was inspired by the proof that every nilpotent element
in a Lie algebra (over a not necessarily algebraically closed field)
is contained in a minimal parabolic that we learned from J.~Lurie. It can also be used to reprove \cite[Lemma 6]{Gi}
which is a key ingredient of the proof in {\it loc.cit.} that the global nilpotent cone in $T^*(\Bun_G)$ is Lagrangian.}

\qed[\propref{p:generation P G}]

\appendix

\section{Action of groups on categories}  \label{s:eq}

Operations explained in this appendix have been used several times in the main body of the paper. They are
applicable to any affine algebraic group $G$ over a ground field of characteristic $0$. 

\ssec{Equivariantization and de-equivariantization}  \label{ss:eq}

\sssec{}

Let $BG^\bullet$ be the standard simplicial model of the classifying space of $G$. Quasicoherent sheaves on 
$BG^\bullet$ form a cosimplical monoidal category, which we denote by $\QCoh(BG^\bullet)$. By definition, a DG 
category acted on by $G$ is a cosimplicial category $\bC^\bullet$ tensored over $\QCoh(BG^\bullet)$ which is 
co-Cartesian in the sense that for every face map $[k]\to [l]$, the functor
\[\QCoh(BG^l)\underset{\QCoh(BG^k)}\otimes\bC^k\to\bC^l\]
is an equivalence.

\medskip

We will regard this as an additional structure over a plain DG
category $\bC:=\bC^0$. We denote the 2-category of DG categories acted on by $G$ (regarded
as an $(\infty,1)$-category) by $G\mmod$. 

\medskip

For $\bC$ as above, we let $\bC^G$ denote the category $\on{Tot}(\bC^\bullet)$.

\sssec{}  \label{sss:de-eq subcategory}

It is easy to see that for $\bC$ acted on by $G$ and a full subcategory $\bC'\subset \bC$, there is at most
one way to define a $G$-action on $\bC'$ in a way compatible with the embedding to $\bC$; this
condition is enough to check at the level of the underlying triangulated categories and for $1$-simplices.
If this is the case, we will say that $\bC'$ is invariant under the action of $G$.

\medskip

It is easy to see that in this case 
$(\bC')^G$ is a full subcategory of $\bC^G$ that fits into the pullback square
\begin{equation} \label{e:eq subcategory}
\CD
\bC'{}^G  @>>> \bC^G  \\
@VVV    @VVV  \\
\bC'  @>>>  \bC.
\endCD
\end{equation}

\sssec{}

Let $\bC$ be acted on by $G$. By construction, $\bC^G$ is 
a module category over
$$\on{Tot}(\QCoh(BG^\bullet))\simeq \Rep(G).$$ 

\medskip

Thus, we obtain a functor 
\begin{equation} \label{e:equiv}
\bC\mapsto \bC^G:G\mmod\to \Rep(G)\mmod,
\end{equation}
where $\Rep(G)\mmod$ is the 2-category of module categories over $\Rep(G)$. 

\medskip

The above functor
admits a left adjoint given by
\begin{equation} \label{e:de-equiv}
\wt\bC\mapsto \on{de-Eq}^G(\wt\bC):=\Vect\underset{\Rep(G)}\otimes \wt\bC,
\end{equation}
where $\Vect$ is naturally regarded as a DG category endowed with the trivial $G$-action and the trivial structure
of a $\Rep(G)$-module with the natural compatibility structure between the two.

\medskip

Note that we also have the naturally defined functors between plain DG categories
$\bC^G\to \bC$ or, equivalently, $\wt\bC\to \on{de-Eq}^G(\wt\bC)$. 

\sssec{}  \label{sss:de-eq reconstr}

We have the following assertion (\cite[Theorem 2.2.2]{GA}):

\begin{thm}  \label{t:de-eq}
The two functors \eqref{e:equiv} and \eqref{e:de-equiv} are mutually inverse. 
\end{thm}

\sssec{}

Several comments are in order:  

\medskip

The fact that for $\bC\in G\mmod$, the adjunction map
$$\on{de-Eq}^G(\bC^G)\to \bC$$
is an equivalence is easy. It follows from the fact that the functor $\on{de-Eq}^G$ commutes
with both \emph{colimits} and \emph{limits}, which in turn follows from the fact that the monoidal category $\Rep(G)$
is rigid (see \cite[Corollaries 4.3.2 and 6.4.2]{DG}). 

\medskip

For $\wt\bC\in \Rep(G)\mmod$, the fact that the adjunction map
$$\wt\bC\to (\on{de-Eq}^G(\wt\bC))^G$$
is an isomorphism is also easy to see when $\wt\bC$ is dualizable. 

\medskip

The above two observations are the only two cases of \thmref{t:de-eq} that have been used in the main
body of the text. 

\medskip

The difficult direction in \thmref{t:de-eq} implies that if $\bC$ is dualizable
(as an abstract DG category), then so is $\bC^G$. \footnote{We do not know 
whether the fact that $\bC$ is compactly generated implies the corresponding fact for $\bC^G$.}

\ssec{Shift of grading}  \label{ss:shift of grading}

\sssec{} \label{sss:shift of grading vect}

Consider the symmetric monoidal category $\Rep(\BG_m)\simeq \QCoh(\on{pt}/\BG_m)$.
We may view its objects as bigraded vector spaces equipped with differential of bidegree
$(1,0)$. Here in the grading $(i,k)$,
the first index refers to the cohomological grading, and the second index to
the grading coming from the $\BG_m$-action.

\medskip

The symmetric monoidal category $\Rep(\BG_m)$ carries a canonical automorphism which we will refer to as the ``grading shift"
\[M\mapsto M^{\on{shift}},\quad\text{where}\quad M^{\on{shift}}_{(i,k)}=M_{i+2k,k}.\]

\medskip

In particular, the 2-category $\Rep(\BG_m)\mmod$ of DG categories tensored over $\Rep(\BG_m)$ 
carries a canonical auto-equivalence, which commutes with the forgetful functor to the 2-category 
$\StinftyCat_{\on{cont}}$ of
plain DG categories. We denote it by 
\begin{equation} \label{e:shift on categories}
\wt\bC\rightsquigarrow \wt\bC^{\on{shift}}.
\end{equation}

\sssec{}  \label{sss:shift of grading alg}

For example, suppose $A$ is a $\BZ$-graded associative DG algebra, and set $\bC=(A\mod)^{\BG_m}$. That is,
$\bC$ is the DG category of \emph{graded} $A$-modules. 

\medskip

In this case, 
$$((A\mod)^{\BG_m})^{\on{shift}}\simeq (A^{\on{shift}}\mod)^{\BG_m}.$$
The categories $(A^{\on{shift}}\mod)^{\BG_m}$ and $(A\mod)^{\BG_m}$ are equivalent as DG categories
(but not as categories tensored over $\Rep(\BG_m)$) with the equivalence given by
\[M\mapsto M^{\on{shift}}\qquad \text{ for }M\in(A\mod)^{\BG_m}.\]

\sssec{}

By \thmref{t:de-eq}, the shift of grading auto-equivalence of $\Rep(\BG_m)\mmod$
induces an auto-equivalence of the 2-category $\BG_m\mmod$, which we denote by 
\begin{equation} \label{e:shift on categories mod}
\bC\rightsquigarrow \bC^{\on{shift}}.
\end{equation}

\medskip

Note, however, that the auto-equivalence \eqref{e:shift on categories} of $\BG_m\mmod$ \emph{does not}
commute with the forgetful functor to $\StinftyCat_{\on{cont}}$. 

\medskip

For example, for a graded associative DG algebra $A$ as above, we have:
$$(A\mod)^{\on{shift}}\simeq A^{\on{shift}}\mod.$$

\section{Spaces of maps and deformation theory}  \label{s:Weil}

In this appendix we drop the assumption that our DG schemes/Artin stacks/prestacks be
locally almost of finite type. 

\ssec{Spaces of maps}  \label{ss:spaces of maps}

Let $\CZ\in\on{PreStk}$ be an arbitrary prestack (see \cite[Sect. 1.1.1]{Stacks}), thought of as the target,
and let $\CX$ be another object of $\on{PreStk}$, thought of as the source. 

\sssec{}

We define a new prestack $\bMaps(\CX,\CZ)\in \on{PreStk}$, by
$$\Maps(S,\bMaps(\CX,\CZ)):=\Maps(S\times \CX,\CZ)$$
for $S\in \affdgSch$. 

\begin{rem}
Note that the above procedure is a particular case of restriction of scalars \`a la Weil: we can start
with a map $\CX_1\to \CX_2$ in $\on{PreStk}$ and $\CZ_1\in \on{PreStk}_{/\CX_1}$, and define
$$\CZ_2=\Res^{\CX_1}_{\CX_2}(\CZ_1)\in \on{PreStk}_{/\CX_2}$$ by
$$\Maps(S,\CZ_2):=\Maps_{\on{PreStk}_{/\CX_1}}(S\underset{\CX_2}\times \CX_1,\CZ_1).$$
In our case $\CX_1=\CX$ and $\CX_2=\on{pt}$.
\end{rem}

\sssec{}

For example, we define:
$$\Bun_G(\CX):=\bMaps(\CX,\on{pt}/G) \text{ and } \LocSys_G(\CX):=\bMaps(\CX_\dr,\on{pt}/G),$$
where $\CX_\dr$ is the de Rham prestack of $\CX$ (see \cite[Sect. 1.1.1]{Crys}). 

\ssec{Deformation theory}

Let us recall some basic definitions from deformation theory. We refer the reader to \cite[Sect. 2.1]{Lu2}
or \cite[Sect. 4]{IndSch} for a more detailed treatment. 

\sssec{}  \label{sss:split SqZ}

Recall that for $S\in \affdgSch$ we have a canonically defined functor
$$\QCoh(S)^{\leq 0}\to \affdgSch_{S/}$$
that assigns to $\CF\in \QCoh(S)^{\leq 0}$ the corresponding split square-zero extension $S_\CF$ of $S$, i.e., 
$$\CF\mapsto S_\CF:=\Spec\left(\Gamma(S,\CO_S)\oplus \Gamma(S,\CF)\right).$$

\sssec{}

Let $\CZ$ be an object of $\on{PreStk}$, and let $z$ be a point of $\Maps(S,\CZ)$. 
Consider the following functor
$\QCoh(S)^{\leq 0}\to \inftygroup$
\begin{equation} \label{e:cotangent space}
\CF\mapsto \Maps(S_\CF,\CZ)\underset{\Maps(S,\CZ)}\times \{z\}.
\end{equation}

\begin{defn} \label{d:cond A}  Let $k$ be a non-negative integer. 

\smallskip

\noindent{\em(a)}
We will say that $\CZ$ admits $(-k)$-connective pro-cotangent spaces if for any $(S,x)$ the functor \eqref{e:cotangent space} is
pro-representable by an object of $\on{Pro}(\QCoh(S)^{\leq k})$. 

\smallskip

\noindent{\em(b)}
We will say that $\CZ$ admits $(-k)$-connective cotangent spaces if for any $(S,x)$ the functor \eqref{e:cotangent space} is
co-representable by an object of $\QCoh(S)^{\leq k}$. 
\end{defn}

For $\CZ$ as in Definition \ref{d:cond A}(a) (resp., (b)), we will denote the resulting object of $\on{Pro}(\QCoh(S)^{\leq k})$
(resp., $\QCoh(S)^{\leq k}$) by $T^*_z(\CZ)$ and refer to it as the ``cotangent space to $\CZ$ at the point $z$."

\medskip
That is, if we regard
$$\CMaps_{\QCoh(S)}(T^*_z(\CZ),\CF)$$
as an object of $\Vect$, and consider its truncation 
$$\Maps_{\QCoh(S)}(T^*_z(\CZ),\CF):=\tau^{\leq 0}\left(\CMaps_{\QCoh(S)}(T^*_z(\CZ),\CF)\right)$$ 
as an $\infty$-groupoid via $\Vect^{\leq 0}\to \inftygroup$,
the result is canonically isomorphic to \eqref{e:cotangent space}. 

\sssec{}

Let $\alpha:S_1\to S$ be a map in $\affdgSch$. Consider the corresponding functor 
$$\on{Pro}(\alpha^*):\on{Pro}(\QCoh(S))\to \on{Pro}(\QCoh(S_1)).$$

By definition, for an object $\Phi\in \on{Pro}(\QCoh(S))$, viewed as a functor $\QCoh(S)\to \Vect$,
the object $\on{Pro}(\alpha^*)(\Phi)\in \on{Pro}(\QCoh(S_1))$, viewed as a functor $\QCoh(S_1)\to \Vect$,
is given by the left Kan extension of $\Phi$ along $\alpha^*:\QCoh(S)\to \QCoh(S_1)$.

\medskip

Let $\CZ$ be an object of $\on{PreStk}$ that admits
$(-k)$-connective pro-cotangent spaces. Then for $z:S\to \CZ$ and $z_1:=z\circ \alpha$ we obtain a map
\begin{equation} \label{e:cond B}
T^*_{z_1}(\CZ)\to \on{Pro}(\alpha^*)(T^*_z(\CZ))
\end{equation}
in $\on{Pro}(\QCoh(S_1))$. 

\medskip

If $\CZ$ admits $(-k)$-connective cotangent spaces, then $$T^*_z(\CZ)\in \QCoh(S)^{\leq k} \text{ and } T^*_{z_1}(\CZ)\in \QCoh(S_1)^{\leq k},$$
and the map \eqref{e:cond B} is a map
$$T^*_{z_1}(\CZ)\to \alpha^*(T^*_{z}(\CZ))\in \QCoh(S_1).$$

\medskip

\begin{defn}  \label{d:cond B}  \hfill

\smallskip

\noindent{\em(a)}
We will say that $\CZ$ admits a $(-k)$-connective pro-cotangent complex if it admits $(-k)$-connective
pro-cotangent spaces and the map \eqref{e:cond B}
is an isomorphism for any $z$ and $\alpha$. 

\smallskip

\noindent{\em(b)} 
We will say that $\CZ$ admits a $(-k)$-connective cotangent complex if it admits $(-k)$-connective
cotangent spaces and the map \eqref{e:cond B} is an isomorphism for any $z$ and $\alpha$. 
\end{defn}

If $\CZ$ admits a $(-k)$-connective cotangent complex, it gives rise
to a well-defined object in $\QCoh(\CZ)^{\leq -k}$ that we will denote by $T^*(\CZ)$. For a given $z:S\to \CZ$
we will also use the notation
$$T^*(Z)|_S:=T^*_z(Z).$$

\sssec{}

Let now $\CS$ be an arbitrary object of $\on{PreStk}$ and let $z:\CS\to \CZ$. Consider the functor
$$\CF\mapsto \CS_{\CF}:\QCoh(\CS)^{\leq 0}\to \on{PreStk}_{\CS/},$$
defined by
$$\Maps(U,\CS_{\CF})=\{s:U\to \CS,\,\, f:U\to U_{s^*(\CF)}\}$$
for $U\in \affdgSch$ (here the notation $U_{s^*(\CF)}$ is as in \secref{sss:split SqZ}). 

\medskip

The following is tautological from the definitions:

\begin{lem}   \label{l:cotangent on prestack}
Assume that $\CZ$ admits a $(-k)$-connective cotangent complex. For any 
$\CS\in \on{PreStk}$ and a map $z:\CS\to \CZ$, consider the functor $\QCoh(\CS)^{\leq 0}\to \inftygroup$ given by
$$\Maps(\CS_\CF,\CZ)\underset{\Maps(\CS,\CZ)}\times \{z\}.$$
This functor is canonically isomorphic to
$$\tau^{\leq 0}\left(\CMaps_{\QCoh(\CS)}(T_z^*(\CZ),\CF)\right).$$
\end{lem}

\sssec{}

Let $S'$ be a square-zero extension of $S\in \affdgSch$, not necessarily split. Such $S'$ corresponds to an object $\CI\in \QCoh(S)^{\leq 0}$
(the ideal of $S\subset S'$) and a map 
$$T^*(S)\to \CI[1],$$
where $T^*(S)$ is the cotangent complex of $S$ (see \cite[Sect. 4.5]{IndSch}). 

\medskip

Let $\CZ$ admit $(-k)$-connective pro-cotangent spaces. Let $z:S\to \CZ$ be a map. As in \cite[Lemma 4.6.6]{IndSch}, 
we obtain a map
\begin{equation}  \label{e:cond C}
\Maps(S',\CZ)\underset{\Maps(S,\CZ)}\times \{z\}\to 
\tau^{\leq 0}\left(\CMaps_{\QCoh(S)}\left(\on{Cone}(T^*_z(\CZ)\overset{(dz)^*}\longrightarrow T^*(S)),\CI[1]\right)\right),
\end{equation}
where $(dz)^*:T^*_z(\CZ)\to T^*(S)$ is the dual of the differential, see \cite[Sect. 4.4.5]{IndSch}.

\begin{defn}  \label{d:cond C}
We will say that $\CZ$ is infinitesimally cohesive if the map \eqref{e:cond C} is an isomorphism for
any $S$, $z$ and $S'$.
\end{defn}

\begin{rem}
The meaning of the above definition is that (pro)-cotangent spaces control not only maps out of split
square-zero extensions, but from all square-zero extensions. Iterating, we obtain that infinitesimal 
cohesiveness of $\CZ$ implies that its (pro)-cotangent spaces effectively control extensions of a given 
map $S\to \CZ$ to maps $S'\to \CZ$ for any nil-immersion $S\hookrightarrow S'$.
\end{rem}

\sssec{}

Finally, we define:

\begin{defn}  \label{d:deform} 
We will say that $\CZ\in \on{PreStk}$ admits a $(-k)$-connective deformation theory
\emph{(}resp., co-representable $(-k)$-connective deformation theory\emph{)} if:

\begin{itemize}

\item $\CZ$ is convergent (see \cite{Stacks}, Sect. 1.2.1);

\item $\CZ$ admits a $(-k)$-connective pro-cotangent complex \emph{(}resp., $(-k)$-connective cotangent complex\emph{)};

\item $\CZ$ is infinitesimally cohesive.

\end{itemize}

\end{defn}

\sssec{}

We record the following result for use in the main body of the paper. 

\begin{thm} \label{t:repr}  
An object $\CZ\in\on{PreStk}$ is a DG indscheme \emph{(}resp., DG scheme\emph{)} if and only if the following conditions hold:

\begin{itemize}

\item The classical prestack $^{cl}\CZ$ is a classical indscheme \emph{(}resp., scheme\emph{)}.

\item $\CZ$ admits a $0$-connective \emph{(}resp., co-representable $0$-connective\emph{)} deformation theory. 

\end{itemize}

\end{thm}

The above theorem for DG indschemes is \cite[Theorem 5.1.1]{IndSch}. The case of DG schemes
can be proved similarly (see \cite[Theorem 3.1.2]{Lu2}).

\medskip

The same proof also shows that if $^{cl}\CZ$ is an ind-affine indscheme (i.e., is a  filtered colimit of \emph{affine} classical schemes
under closed embeddings), or an affine scheme, then $\CZ$ is an ind-affine DG indscheme, or an affine DG scheme, respectively. 

\ssec{Deformation theory of spaces of maps}

\sssec{}

Let $\CX$ and $\CZ$ be as in \secref{ss:spaces of maps}.  
We are going to show:

\begin{prop}  \label{p:maps deform}

Assume that $\CZ$ admits a co-representable $(-k)$-connective deformation theory. Assume that $\CX$
is $l$-coconnective for some $l\in \BZ^{\geq 0}$.  Then:

\smallskip

\noindent{\em(a)} The prestack $\bMaps(\CX,\CZ)$ admits a $(-k-l)$-connective deformation theory. 

\smallskip

\noindent{\em(b)} For $\wt{z}\in \Maps(S,\bMaps(\CX,\CZ))$ corresponding to $z\in \Maps(S\times \CX,\CZ)$,
the pro-cotangent space $T^*_{\wt{z}}(\bMaps(\CX,\CZ))$, viewed as a functor $\QCoh(S)^{\leq 0}\to \Vect^{\leq 0}$, 
identifies with
\begin{equation} \label{e:cotangent to maps}
\CF \mapsto \tau^{\leq 0}\left(\CMaps_{\QCoh(S\times \CX)}(T^*_z(\CZ),\CF\boxtimes \CO_\CX)\right).
\end{equation}
\end{prop}

\begin{proof}

The fact that the pro-cotangent spaces of $\bMaps(\CX,\CZ)$ are given by the functor
\eqref{e:cotangent to maps} follows from \lemref{l:cotangent on prestack}. 

\medskip

Let us show that these pro-cotangent spaces are $(-k-l)$-connective. That is, we need to show
that the functor 
$$\CF \mapsto \CMaps_{\QCoh(S\times \CX)}(T^*_z(\CZ),\CF\boxtimes \CO_\CX),\quad \QCoh(S)\to \Vect$$
sends $\QCoh(S)^{\geq 0}$ to $\Vect^{\geq -k-l}$. 

\medskip

Since $\CX$ is $l$-coconnective, we can write $\CX$ as
$\underset{U\in {}^{\leq l}\!\affdgSch_{/\CX}}{colim}\, U$, 
and hence
$$\CMaps_{\QCoh(S\times \CX)}(T^*_z(\CZ),\CF\boxtimes \CO_\CX)\simeq
\underset{U\in ({}^{\leq l}\!\affdgSch_{/\CX})^{\on{op}}}{lim}\, 
\CMaps_{\QCoh(S\times U)}(T^*_z(\CZ)|_{S\times U},\CF\boxtimes \CO_U).$$

However, $\CF\in \QCoh(S)^{\geq 0}$ implies $\CF\boxtimes \CO_U\in \QCoh(S\times U)^{\geq -l}$,
and the assertion follows. 

\medskip

The fact that 
\eqref{e:cond B} and \eqref{e:cond C} are isomorphisms follows from the definitions. Finally,
the fact that $\bMaps(\CX,\CZ)$ is convergent follows tautologically from the fact that $\CZ$
is convergent. 

\end{proof}

\sssec{}  \label{sss:maps deform}

Let now $X$ be an $l$-coconnective DG scheme, proper over $\Spec(k)$. Let $\CZ$ be again a prestack
that admits a co-representable $(-k)$-connective deformation theory. 

\begin{cor}  \label{c:maps deform}
The prestack $\bMaps(X,\CZ)$ admits a 
co-representable $(-k-l)$-connective deformation theory and the prestack admits 
$\bMaps(X_\dr,\CZ)$ admits a co-representable $(-k)$-connective deformation theory.
\end{cor}

\begin{proof}

Let $\CX$ be either $X$ or $X_\dr$.

\medskip

We need to show the existence of a left adjoint functor to
$$\CF\mapsto \CF\boxtimes \CO_X:\QCoh(S)\to \QCoh(S\times \CX).$$
Since $\QCoh(S\times \CX)\simeq \QCoh(S)\otimes \QCoh(\CX)$ (see \cite[Proposition 1.4.4]{QCoh}),
it is sufficient to consider the case $S=\on{pt}$. 

\medskip

Since the category $\QCoh(\CX)$ is compactly generated, it suffices to 
construct the left adjoint on compact objects of $\QCoh(\CX)$. This amounts 
to showing that for $\CE\in \QCoh(\CX)^c$, the object
\begin{equation} \label{e:would be left}
\CMaps_{\QCoh(\CX)}(\CE,\CO_\CX)\in \Vect
\end{equation}
is compact (i.e., has finitely many non-zero cohomologies, all of which are finite-dimensional);
then the left adjoint in question sends $k\in \Vect$ to the dual of \eqref{e:would be left}.

\medskip

The compactness of \eqref{e:would be left} follows easily from the assumption on $\CX$.

\end{proof}

\ssec{The ``locally almost of finite type" condition}

\sssec{}

Recall the notion of prestack locally almost of finite type, see \cite[Sect. 1.3.9]{Stacks}). We have the following 
assertion:

\begin{lem}  \label{l:laft crit}
Let $\CZ\in \on{PreStk}$ be a prestack that admits a $(-k)$-connective deformation theory for some $k$.
Then $\CZ$ is locally almost of finite type if and only if:

\begin{itemize}

\item The underlying classical prestack $^{cl}\CZ$ is locally of finite type \emph{(}see \cite[Sect. 1.3.2]{Stacks} for what this means\emph{)}.

\item Given any classical scheme of finite type $S\in \affSch_{\on{ft}}$, a morphism $z:S\to \CZ$, and
an integer $k\geq 0$, the cotangent space $T^*_z(\CZ)$, viewed as a functor
$$\QCoh(S)^{\geq -k,\leq 0}\to \inftygroup,$$ commutes with filtered colimits. 

\end{itemize}

\end{lem}

The proof is essentially the same as that of \cite[Proposition 5.3.2]{IndSch}. 

\sssec{}

Let $\CZ$ be a prestack that is locally almost of finite type. 
Let $X$ be an eventually coconnective quasi-compact DG scheme almost of finite type.

\begin{cor}  \label{c:maps laft}
The prestack $\bMaps(X,\CZ)$ is locally almost of finite type, provided that the classical
prestack $^{cl}\CZ$ satisfies Zariski descent.
\end{cor}

\begin{proof}

The description of the cotangent spaces to $\bMaps(X,\CZ)$ given by \propref{p:maps deform}(b)
implies that the second condition of \lemref{l:laft crit} is satisfied. So, it remains to show that the classical
prestack $^{cl}\bMaps(X,\CZ)$
is locally of finite type. By definition, this means that the functor
\begin{equation} \label{e:maps on algs}
S\mapsto \Maps(S,\bMaps(X,\CZ))
\end{equation}
on the category $(\affSch)^{\on{op}}$ should commute with filtered colimits.

\medskip

Any quasi-compact and quasi-separated scheme can be expressed as a \emph{finite} colimit of affine schemes in the category
of Zariski sheaves
(i.e., the full subcategory of $\on{PreStk}$ consisting of objects that satisfy Zariski descent).
Since finite limits commute with filtered colimits,  the assumption on $\CZ$ reduces the assertion to the case 
when $X$ is affine. 

\medskip

Let $n$ be such that $X\in {}^{\leq n}\!\affdgSch$. The functor 
$$S\mapsto S\times X:\affSch\to {}^{\leq n}\!\affdgSch$$
commutes with filtered limits, and the required property of \eqref{e:maps on algs}
follows from the fact that $\CZ$ is locally almost of finite type (namely, that for
every $n$, it takes filtered colimits in $({}^{\leq n}\!\affdgSch)^{\on{op}}$ to colimits in
$\inftygroup$).

\end{proof}
 
\section{The Thomason-Trobaugh theorem ``with supports"}  \label{a:B}

In this appendix we will prove the following result: 

\medskip

\noindent{\bf Theorem C.}
{\it Let $Z$ be a quasi-compact DG scheme, and $Y\subset \Sing(Z)$ a conical Zariski-closed
subset. Then the category $\IndCoh_Y(Z)$ is compactly generated.}

\medskip

The proof is an easy adaptation of the argument of \cite{TT} for the compact generation of $\QCoh(Z)$. 
Let us sketch the proof following \cite{NTr}. 

\begin{proof}

Recall that the objects of $\Coh_Y(Z)$ are compact in $\IndCoh_Y(Z)$.
Hence, it suffices to check that $\Coh_Y(Z)$ generates $\IndCoh_Y(Z)$.

\medskip

The proof proceeds by induction on the number of affine open subsets covering $Z$. The base case is when $Z$ itself is affine.
In this case, the assertion follows from \corref{c:with support comp gen}. 

\medskip

Suppose now $Z$ is arbitrary, and let $W_i$ be an affine cover of $Z$. Fix $\CF\in\IndCoh_Y(Z)$, $\CF\ne 0$. Set
$$U_i=\underset{j\neq i}\cup\, W_j,$$ and choose $i$ so that $\CF|_{U_i}\ne 0$.
We now drop the index $i$ and write simply $U$ and $W$ for $U_i$ and $W_i$. 

\medskip

By the induction hypothesis, we can assume that $\IndCoh_{Y\underset{Z}\times U}(U)$ is compactly generated.
Therefore, there exists a compact object 
\[\wt\CG_U\in\Coh_{Y\underset{Z}\times U}(U)\]
together with a non-zero map $\wt\iota_U:\wt\CG_U\to\CF|_U$. Set $\CG_U:=\wt\CG_U\oplus\wt\CG_U[1]$
and let $\iota_U:\CG_U\to \CF|_U$ be equal to 
$$\wt\CG_U\oplus\wt\CG_U[1] \to \wt\CG_U\overset{\wt\iota_U}\to \CF|_U.$$

We will extend the pair $(\CG_U,\iota_U)$ to $\CG\in \Coh_Y(Z)$ and a map $\iota:\CG\to \CF$.

\medskip

Indeed, by \cite[Theorem~3.1]{NTr}, there exists an object 
\[\CG_W\in \Coh_{Y\underset{Z}\times W}(W)\]
together with an isomorphism 
\[(\CG_W)|_{U\cap W}\simeq(\CG_U)|_{U\cap W}\]
and a map $\iota_W:\CG_W\to\CF|_W$ whose restriction to $U\cap W$ equals $\iota_U|_{U\cap W}$. 

\medskip

Define 
\[\CG:=\on{Cone}\Bigl(j_{U,*}(\CG_U)\oplus j_{W,*}(\CG_W)\to j_{U\cap W,*}((\CG_W)|_{U\cap W})\Bigr)[-1].\]
Here $j_U$, $j_W$, and $j_{U\cap W}$ are the natural embeddings $U\hookrightarrow Z$, $W\hookrightarrow Z$, and 
$U\cap W\hookrightarrow Z$, respectively.  We have
$$\CG|_U\simeq \CG_U\text{ and } \CG|_W\simeq \CG_W,$$
so $\CG\in\Coh_Y(Z)$. The morphisms $\iota_U,\iota_W$ induce a map
$\iota:\CG\to\CF$, whose restrictions to $U$ and $W$ identify with $\iota_U$ and $\iota_W$, respectively. 

\end{proof}

\section{Finite generation of Exts} \label{s:proof of finiteness}

In this appendix, we will prove \thmref{t:finiteness}. Let us recall its formulation:

\medskip

\noindent{\bf Theorem D.}
{\it Let $Z$ be a quasi-smooth affine DG scheme $Z$. Given
$\CF_1,\CF_2\in \Coh(Z)$, consider the graded vector space $\Hom^\bullet_{\Coh(Z)}(\CF_1,\CF_2)$ 
as a module over the graded algebra $\Gamma(\Sing(Z),\CO_{\Sing(Z)})$. We claim that  
the module is finitely generated.}

\medskip

If $Z$ is classical, this is due to Gulliksen \cite{Gu}, and the extension to DG schemes 
is straightforward.

\medskip

\noindent{D.1.}  It is easy to see that the statement is Zariski-local on $Z$. So, we can assume that $Z$ is as in 
\eqref{e:Cart prod diag}. Let 
$$\on{pt}=\CV_0\subset \CV_1\subset...\subset \CV_n=\CV$$ be a flag of smooth closed subschemes whose
dimensions increase by one. With no restriction of generality, we can assume that $\CV_{i-1}$
is cut out by one function inside $\CV_i$.

\medskip

Set
$$Z_i:=\CV_i\underset{\CV}\times \CU.$$
All these DG schemes are quasi-smooth, and $Z_{i-1}$ is cut out inside of $Z_i$ by one function. 
We have $Z_0=Z$ and $Z_n=\CU$. Let $g_i$ denote the closed embedding $Z\to Z_i$. 

\medskip

\noindent{D.2.}  We will argue by descending induction on $i$, assuming that
$$\Hom^\bullet_{\Coh(Z_i)}((g_i)_*(\CF_1),(g_i)_*(\CF_2))$$
is finitely generated as a module over $\Gamma(\Sing(Z_i),\CO_{\Sing(Z_i)})$.

\medskip

The base of induction is $i=n$. In this case $Z_n=\CU$ is smooth, and the assertion is obvious. 

\medskip

To carry out the induction step we can thus assume that we have a quasi-smooth closed
embedding 
$$g:Z\hookrightarrow Z',$$
that fits into a Cartesian diagram
$$
\CD
Z  @>{g}>>  Z'  \\
@VVV   @VVV  \\
\{0\}  @>>>  \BA^1.
\endCD
$$

By induction, we can assume that
$$\Hom^\bullet_{\Coh(Z')}(g_*(\CF),g_*(\CF))$$
is finitely generated as a module over $\Gamma(\Sing(Z'),\CO_{\Sing(Z')})$.

\medskip

\noindent{D.3.} Note that the generator of $T_{\{0\}}(\BA^1)$ gives rise to an element $\eta\in \on{HH}^2(Z)$.

\medskip

Since the grading on the $\Gamma(\Sing(Z),\CO_{\Sing(Z)})$-module $\Hom^\bullet_{\Coh(Z)}(\CF_1,\CF_2)$ is 
bounded below, it is sufficient to show that 
$$\on{coker}\left(\eta:\Hom^\bullet_{\Coh(Z)}(\CF_1[2],\CF_2)\to \Hom^\bullet_{\Coh(Z)}(\CF_1,\CF_2)\right)$$
is finitely generated. 

\medskip

However, from the long exact sequence we obtain that the above cokernel is a submodule in 
\begin{equation} \label{e:Hom into cone}
\Hom^\bullet_{\Coh(Z)}\left(\on{Cone}(\CF_1\overset{\eta}\longrightarrow \CF_1[2]),\CF_2[1]\right).
\end{equation}

Since the algebra $\Gamma(\Sing(Z),\CO_{\Sing(Z)})$ is Noetherian, it is enough to show that the graded
$\Gamma(\Sing(Z),\CO_{\Sing(Z)})$-module given by \eqref{e:Hom into cone}, is finitely generated. 

\medskip

\noindent{D.4.} Note that for any $\CF\in \QCoh(Z)$, the corresponding object
$$\on{Cone}(\CF\overset{\eta}\longrightarrow \CF[2])$$
identifies with $g^*\circ g_*(\CF)$. 

\medskip

Hence, the module \eqref{e:Hom into cone} identifies with
$$\Hom^\bullet_{\Coh(Z)}(g^*\circ g_*(\CF_1),\CF_2),$$
up to a shift of grading.

\medskip

We have an isomorphism of vector spaces:
\begin{equation} \label{e:HH on Z and Z'}
\Hom^\bullet_{\Coh(Z)}(g^*\circ g_*(\CF_1),\CF_2)\simeq \Hom^\bullet_{\Coh(Z')}(g_*(\CF_1),g_*(\CF_2)).
\end{equation}

\medskip

Note that the right-hand side in \eqref{e:HH on Z and Z'} is acted on by $\Gamma(\Sing(Z'),\CO_{\Sing(Z')})$, and this action
factors through the surjection
$$\Gamma(\Sing(Z'),\CO_{\Sing(Z')})\twoheadrightarrow \Gamma(\Sing(Z')_Z,\CO_{\Sing(Z')_Z}).$$

In particular, this gives an action of $\Gamma(\Sing(Z),\CO_{\Sing(Z)})$ on right-hand side in \eqref{e:HH on Z and Z'} 
via the closed embedding 
$$\Sing(g):\Sing(Z')_{Z}\hookrightarrow \Sing(Z).$$

\medskip

Now, it follows from the construction that the above action of  $\Gamma(\Sing(Z),\CO_{\Sing(Z)})$ on
the right-hand side in \eqref{e:HH on Z and Z'} is compatible with the canonical action on the
left-hand side. 

\medskip

Now, as was mentioned above, by the induction hypothesis, $\Hom^\bullet_{\Coh(Z')}(g_*(\CF_1),g_*(\CF_2))$ is finitely generated
over $\Gamma(\Sing(Z'),\CO_{\Sing(Z')})$, which implies that $\Hom^\bullet_{\Coh(Z)}(g^*\circ g_*(\CF_1),\CF_2)$
is finitely generated over $\Gamma(\Sing(Z),\CO_{\Sing(Z)})$, as required.
\qed

\section{Recollections on $\BE_2$-algebras} \label{s:E2}

In this appendix we recall some basic facts regarding $\BE_2$-algebras and their
actions on categories.

\ssec{$\BE_2$-algebras}  

The main reference to the theory of $\BE_2$-algebras is \cite[Sect. 5.1]{Lu1}. Here we will summarize
some basic facts. All monoidal categories, $\BE_1$-algebras and $\BE_2$-algebras will be assumed unital.

\medskip

We will use the terms ``$\BE_1$-algebra" and ``associative DG algebra'' interchangeably; and similarly for 
the terms ``monoidal functor" and ``homomorphism of monoidal DG categories."  

\sssec{}

For the purposes of the paper, we set
$$\BE_2\Alg:=\BE_1\Alg(\BE_1\Alg).$$

I.e., an $\BE_2$-algebra is an associative algebra object in the category $\BE_1\Alg$ of associative
algebras. 

\medskip

The fact that this definition is equivalent to the definition of $\BE_2$-algebras as modules over the little discs
operad is the Dunn Additivity Theorem, see \cite[Theorem 5.1.2.2]{Lu1}.

\sssec{}   

For an object $\CA\in \BE_2\Alg$, we will refer to the associative algebra structure on the underlying 
object of $\BE_1\Alg$ as \emph{the interior multiplication}.  

\medskip

We will refer to the associative algebra structure
on the image of $\CA$ under the forgetful functor 
$$\BE_2\Alg=\BE_1\Alg(\BE_1\Alg)\to \BE_1\Alg(\Vect)=\BE_1\Alg$$
as \emph{the exterior multiplication}. 

\medskip

Unless specified otherwise, for $\CA\in \BE_2\Alg$, by the \emph{underlying $\BE_1$-algebra} we will mean
the result of the application of the forgetful functor that remembers \emph{the interior multiplication}. 

\medskip

In particular, we will denote by $\CA\mod$ the category of left $\CA$-modules, where $\CA$ is regarded
as an $\BE_1$-algebra with the interior multiplication.

\sssec{}    \label{sss:E2 opp}

We will denote by $\CA^{\on{int-op}}$ (resp., $\CA^{\on{ext-op}}$) the $\BE_2$-algebra obtained by
reversing the first (resp., second) multiplication. Note, however, because of \cite[Theorem 5.1.2.2]{Lu1}, 
the choice of orientation on $S^1$ (the real circle) gives rise to a canonical isomorphism
$$\CA^{\on{int-op}}\simeq \CA^{\on{ext-op}}$$
of $\BE_2$-algebras. 

\medskip

Having made the choice of orientation once and for all, we will sometimes use the notation $\CA^{\on{op}}$
for the resulting $\BE_2$-algebra, with one of the multiplications reversed.  

\ssec{$\BE_2$-algebras and monoidal categories}  \label{ss:E2 and monoidal}

\sssec{}

We recall the following construction, explained to us by J.~Lurie. 

\medskip

On the one hand, we consider the symmetric monoidal $\infty$-category $\BE_1\Alg$ of $\BE_1$-algebras. 
On the other hand, we consider the symmetric monoidal $\infty$-category $$(\StinftyCat_{\on{cont}})_{\Vect\!/},$$ i.e., the category
of pairs $(\bc,\bC)$, where $\bC\in \StinftyCat_{\on{cont}}$ and $\bc\in \bC$. 

\medskip

We have a fully faithful symmetric monoidal functor
\begin{equation} \label{e:from alg to cat}
\BE_1\Alg\to (\StinftyCat_{\on{cont}})_{\Vect\!/},\quad A\mapsto (A,A^{\on{op}}\mod).
\end{equation}

The essential image of this functor consists of those $(\bc,\bC)$ for which $\bc$ is a compact
generator of $\bC$. 

\medskip

The functor \eqref{e:from alg to cat} has a right adjoint, given by
\begin{equation} \label{e:from cat to alg}
(\bc,\bC)\mapsto \CMaps_\bC(\bc,\bc).
\end{equation}

Tautologically, the functor \eqref{e:from cat to alg} is \emph{right-lax} symmetric monoidal. In particular,
it naturally upgrades to a functor between the categories of algebras (for any operad). 

\sssec{}

Since the functor \eqref{e:from alg to cat} is symmetric monoidal, we obtain that for $\CA\in \BE_2\Alg$,
the category $\CA^{\on{int-op}}\mod$ acquires a monoidal structure, for which $\CA$ is the unit. 

\medskip

We denote this operation by
$$\CM_1,\CM_2\mapsto \CM_1\underset{\CA}\otimes \CM_2,$$

\medskip

Informally, if we view $\CA$ as an associative DG algebra with respect to the interior multiplication,
the tensor product is induction with respect to the homomorphism
$$\CA\otimes \CA\to \CA$$
given by the exterior multiplication.

\medskip

In fact, since the functor \eqref{e:from alg to cat} is fully faithful, the structure of $\BE_2$-algebra
on a given $\CA\in \BE_1\Alg$ is equivalent to a structure on $\CA^{\on{int-op}}\mod$ of monoidal
category with the unit being $\CA$. 

\sssec{}

For a given symmetric monoidal category $\bO$, the functor \eqref{e:from cat to alg}
sends it to the $\BE_2$-algebra whose underlying $\BE_1$-algebra is $\CMaps_\bO({\bf 1}_\bO,{\bf 1}_\bO)$.
We denote this $\BE_2$-algebra by $\CA_\bO$. 

\medskip

It is easy to see that we have a canonical isomorphism of $\BE_2$-algebras:
\begin{equation} \label{e:opp mon category}
\CA_{\bO^{\on{op-mon}}}\simeq \CA_\bO^{\on{ext-op}},
\end{equation}
where $\bO^{\on{op-mon}}$ denotes the monoidal category obtained from $\bO$ by reversing the monoidal structure. 

\sssec{}  \label{sss:recovering monoidal}

By adjunction, a datum of a homomorphism of $\BE_2$-algebras
$$\CA\to \CA_\bO$$
is equivalent to that of a continuous monoidal functor
$$\CA^{\on{int-op}}\mod\to \bO.$$

\medskip

In particular, the unit of the adjunction is the tautological isomorphism
$$\CA\to \CMaps_{\CA^{\on{int-op}}}(\CA,\CA).$$

The co-unit of the adjunction is a canonically defined monoidal functor
$$\CMaps_\bO({\bf 1}_\bO,{\bf 1}_\bO)^{\on{int-op}}\mod \to \bO.$$

The latter functor is an equivalence if and only if ${\bf 1}_\bO$ is a compact generator of $\bO$.

\ssec{Hochschild cochains of a DG category}

\sssec{}

For $\bC\in \StinftyCat_{\on{cont}}$, we consider the monoidal category
$$\on{Funct}_{\on{cont}}(\bC,\bC).$$

\medskip

We let $\on{HC}(\bC)$ denote the resulting $\BE_2$-algebra.

\medskip

By definition, the $\BE_1$-algebra underlying $\on{HC}(\bC)$ is
$$\CMaps_{\on{Funct}_{\on{cont}}(\bC,\bC)}(\on{Id}_\bC,\on{Id}_\bC).$$

\sssec{}    \label{sss:action of E2 on category}

Let $\CA$ be an $\BE_2$-algebra.  By adjunction
the following pieces of data are equivalent:

\begin{itemize}

\item an action of the monoidal category $\CA\mod$ on $\bC$;

\item a homomorphism of monoidal categories $\CA\mod\to \on{Funct}_{\on{cont}}(\bC,\bC)$;

\item a map of $\BE_2$-algebras $\CA^{\on{int-op}}\to \on{HC}(\bC)$.

\end{itemize}

\medskip

We will refer to such data as a (left) \emph{action of $\CA$ on $\bC$}. 

\medskip

By a right action of $\CA$ on $\bC$
we will mean a left action of $\CA^{\on{int-op}}$, or, equivalently, a homomorphism of $\BE_2$-algebras 
$\CA\to \on{HC}(\bC)$.

\sssec{}  \label{sss:E2 action on dual}

Assume now that $\bC$ is dualizable. In this case we have a canonical equivalence:
$$\on{Funct}_{\on{cont}}(\bC^\vee,\bC^\vee)\simeq \left(\on{Funct}_{\on{cont}}(\bC,\bC)\right)^{\on{op-mon}}.$$

Hence, from \eqref{e:opp mon category}, we obtain a canonical isomorphism of $\BE_2$-algebras
$$\on{HC}(\bC^\vee)\simeq \on{HC}(\bC)^{\on{ext-op}}.$$

\medskip

In particular, if an $\BE_2$-algebra $\CA$ acts on $\bC$, then $\CA^{\on{ext-op}}$ naturally acts on $\bC^\vee$. 

\section{Hochschild cochains of a DG scheme and groupoids}    \label{s:Hoch}

In this appendix we collect some facts and constructions regarding the $\BE_2$-algebra
of Hochschild cochains of a DG scheme $Z$.  We will also review a variant of this construction,
when we produce an $\BE_2$-algebra out of a groupoid acting on $Z$. 

\ssec{Hochschild cochains of a DG scheme}   \label{ss:two isoms for HH}
 
Let $Z$ be a quasi-compact DG scheme. We put 
$$\on{HC}(Z):=\on{HC}(\QCoh(Z)).$$

The main result of this subsection is that we can equivalently realize $\on{HC}(Z)$ as the
$\BE_2$-algebra of Hochschild cochains of the category $\IndCoh(Z)$. 

\begin{rem}
The assertion regarding the comparison of the two versions of $\on{HC}(Z)$
is \emph{not} essential for the contents of this paper. The object that we are really
interested in is the $\BE_2$-algebra of Hochschild cochains of the category $\IndCoh(Z)$,
denoted below by $\on{HC}^{\IndCoh}(Z)$.

\medskip

So, the reader can skip this subsection, substituting $\on{HC}^{\IndCoh}(Z)$ for
$\on{HC}(Z)$ in any occurrence of the latter in the main body of the paper, in 
particular, in \secref{s:sing}.
\end{rem}

\sssec{}

Recall that the category $\QCoh(Z)$ is canonically self-dual
\begin{equation} \label{e:self-duality on QCoh}
\bD_Z^{\on{naive}}:\QCoh(Z)^\vee\simeq \QCoh(Z),
\end{equation}
where the corresponding functor on the compact objects
$$(\QCoh(Z)^c)^{\on{op}}\to \QCoh(Z)^c$$
is the ``naive" duality functor $\BD_Z^{naive}(-)=\underline{\on{Hom}}_{\QCoh(Z)}(-,\CO_Z)$ on $\QCoh(Z)^c=\QCoh(Z)^{\on{perf}}$.

\medskip

In particular, from \secref{sss:E2 action on dual}, we obtain a canonical identification
\begin{equation} \label{e:inv on HH}
\on{HC}(Z)\simeq \on{HC}(Z)^{\on{ext-op}}.
\end{equation}

\sssec{}

Consider now the category $\IndCoh(Z)$. Denote 
$$\on{HC}^{\IndCoh}(Z):=\on{HC}(\IndCoh(Z)).$$ 

\medskip

Recall (see \cite[Sect. 9.2.1]{IndCoh}) that the category $\IndCoh(Z)$ is also canonically self-dual
\begin{equation} \label{e:self-duality on IndCoh}
\bD_Z^{\on{Serre}}:\IndCoh(Z)^\vee\simeq \IndCoh(Z),
\end{equation}
where the corresponding functor on the compact objects
$$(\IndCoh(Z)^c)^{\on{op}}\to \IndCoh(Z)^c$$
is the Serre duality functor $\BD_Z^{\on{Serre}}$ on $\IndCoh(Z)^c=\Coh(Z)$.

\medskip

Hence, we obtain a canonical isomorphism:
\begin{equation} \label{e:inv on HH Ind}
\on{HC}^{\IndCoh}(Z)\simeq \on{HC}^{\IndCoh}(Z)^{\on{ext-op}}.
\end{equation}

\sssec{}  \label{sss:two isoms for HH}

Recall now that there exists a canonically defined functor
$$\Psi_Z:\IndCoh(Z)\to \QCoh(Z),$$
obtained by ind-extending the tautological embedding $\Coh(Z)\hookrightarrow \QCoh(Z)$,
see \cite[Sect. 1.1.5]{IndCoh}.

\medskip

We let 
$$\Upsilon_Z:\QCoh(Z)\to \IndCoh(Z)$$ 
be its dual with respect to the identifications \eqref{e:self-duality on QCoh} and \eqref{e:self-duality on IndCoh}.

\medskip

Recall also (see \cite[Proposition 9.3.3]{IndCoh}) that 
\[\Upsilon_Z(\CF)=\CF\otimes\omega_Z,\quad \CF\in\QCoh(Z),\]
where $\omega_Z\in\IndCoh(Z)$ is the dualizing complex and ``$\otimes$" stands for the action of the
monoidal category $\QCoh(Z)$ on $\IndCoh(Z)$, see \cite[Sect. 1.4]{IndCoh}.

\begin{prop}  \label{p:two isoms for HH}
There exist uniquely defined isomorphisms
$$\Psi_{\on{HC}}:\on{HC}(Z)\to \on{HC}^{\IndCoh}(Z) \text{ and } \Upsilon_{\on{HC}}:\on{HC}^{\IndCoh}(Z)\to \on{HC}(Z),$$
the former compatible with the functor $\Psi_Z$, and the latter 
compatible with the functor $\Upsilon_Z$. Moreover, the following diagram commutes:
\begin{equation} \label{e:two isoms for HH}
\CD
\on{HC}(Z)^{\on{ext-op}}     @>{\sim}>>  \on{HC}(Z)  \\
@A{(\Upsilon_{\on{HC}})^{\on{ext-op}}}AA       @VV{\Psi_{\on{HC}}}V   \\      
\on{HC}^{\IndCoh}(Z)^{\on{ext-op}} @<{\sim}<< \on{HC}^{\IndCoh}(Z)
\endCD
\end{equation}
(in the sense that the composition map starting from any corner is canonically 
isomorphic to the identity map).
\end{prop}

\ssec{Proof of \propref{p:two isoms for HH}}

\sssec{}  

We will need the following general construction. 

\medskip

Let $\bC_1$ and $\bC_2$ be two DG categories, and let $\Phi:\bC_1\to \bC_2$ be a 
continuous functor. Assume that $\bC_1$ is compactly generated, and assume that 
$\Phi|_{\bC_1^c}$ is fully faithful. 

\medskip

We claim that in this case there exists a unique map of $\BE_2$-algebras
$$\Phi_{\on{HC}}:\on{HC}(\bC_2)\to \on{HC}(\bC_1),$$
compatible with $\Phi$, i.e., the functor $\Phi$ 
intertwines the action of $\on{HC}(\bC_2)^{\on{int-op}}$ on $\bC_1$, induced 
by $\Phi_{\on{HC}}$, and the tautological action of $\on{HC}(\bC_2)^{\on{int-op}}$ on $\bC_2$.

\medskip

Indeed, consider the object 
$$(\Phi,\on{Funct}_{\on{cont}}(\bC_1,\bC_2))\in (\StinftyCat_{\on{cont}})_{\Vect\!/}.$$
It is acted on the right and on the left by 
$$(\on{Id}_{\bC_i},\on{Funct}_{\on{cont}}(\bC_i,\bC_i)),\quad i=1,2,$$
respectively. 

\medskip

Applying the functor \eqref{e:from cat to alg}, we obtain the object
$$\CMaps_{\on{Funct}_{\on{cont}}(\bC_1,\bC_2)}(\Phi,\Phi)\in \BE_1\Alg,$$
equipped with an action of 
$\on{HC}(\bC_2)\in \BE_1\Alg(\BE_1\Alg)$, 
as well as a \emph{right} commuting action of $\on{HC}(\bC_1)\in \BE_1\Alg(\BE_1\Alg)$.

\medskip

Furthermore, the assumption on $\Phi$ implies that the action on the unit
defines an isomorphism of right modules over $\on{HC}(\bC_1)$: 
$$\on{HC}(\bC_1)\to \CMaps_{\on{Funct}_{\on{cont}}(\bC_1,\bC_2)}(\Phi,\Phi).$$

Since
$$\on{HC}(\bC_1)\to \on{End}_{\on{HC}(\bC_1)^{\on{ext-op}}\mod}(\on{HC}(\bC_1))$$
is an isomorphism of $\BE_2$-algebras, the left action of $\on{HC}(\bC_2)$ defines the desired homomorphism
$\Phi_{\on{HC}}$.

\medskip

The following assertion results from the construction:

\begin{lem}  \label{l:two isoms for HH}
Assume that in the above situation the category $\bC_2$ is also compactly generated, and that
the dual functor $\Phi^\vee:\bC_2^\vee\to \bC_1^\vee$ is such that $\Phi^\vee|_{(\bC_2^\vee)^c}$ 
is also fully faithful. Then the resulting diagram 
$$
\CD
\on{HC}(\bC_1)^{\on{ext-op}}    @>{\sim}>>  \on{HC}(\bC^\vee_1) \\
@A{\Phi_{\on{HC}}^{\on{ext-op}}}AA   @VV{\Phi^\vee_{\on{HC}}}V  \\
\on{HC}(\bC_2)^{\on{ext-op}}   @<{\sim}<<   \on{HC}(\bC^\vee_2)
\endCD
$$
commutes.  
\end{lem}

\begin{cor}  \label{c:two isoms for HH}
In the situation of \lemref{l:two isoms for HH}, the maps 
$\Phi_{\on{HC}}$ and $\Phi^\vee_{\on{HC}}$ are isomorphisms. 
\end{cor}

\begin{proof}[Proof of \propref{p:two isoms for HH}]

We apply \lemref{l:two isoms for HH} to 
$\bC_1:=\QCoh(Z)$, $\bC_2:=\IndCoh(Z)$, $\Phi:=\Upsilon_Z$, and so $\Phi^\vee=\Psi_Z$.
By definition, $\Psi|_{\Coh(Z)}$ is fully faithful. It also easy to see that $\Upsilon_Z|_{\QCoh(Z)^{\on{perf}}}$ is fully
faithful, as required.

\end{proof}

\ssec{Further remarks on the relation between $\on{HC}(Z)$ and $\on{HC}^{\IndCoh}(Z)$}

In what follows we will identify $\on{HC}(Z)$ and $\on{HC}^{\IndCoh}(Z)$, and unless specified
otherwise, we will do so using the isomorphism $\Upsilon_{\on{HC}}$.

\sssec{}

Assume that $Z$ is eventually coconnective. Recall that in this case the functor
$\Psi_Z$ admits a fully faithful left adjoint, denoted $\Xi_Z$. 

\medskip

In particular, by the first paragraph of the proof of \propref{p:two isoms for HH},
we obtain that there exists a unique homomorphism
$$\Xi_{\on{HC}}:\on{HC}^{\IndCoh}(Z)\to \on{HC}(Z),$$
compatible with $\Xi_Z$.

\medskip

Note, however, that the isomorphism $\Psi_{\on{HC}}$ is also compatible
with $\Xi_Z$, by adjunction. So, we obtain that $\Xi_{\on{HC}}$ provides
an explicit inverse to $\Psi_{\on{HC}}$.

\sssec{}

Assume now that $Z$ is eventually coconnective and Gorenstein (that is, $\omega_Z\in \Coh(Z)$
is a cohomologically shifted line bundle). For example, this is the case for quasi-smooth DG schemes.

\medskip

In this case, we can regard the tensor product by $\omega_Z$ as a self-equivalence
of the category $\QCoh(Z)$. We have a commutative diagram
$$
\CD
\QCoh(Z)  @>{\Xi_Z}>>  \IndCoh(Z)  \\
@A{\omega_Z\otimes -}AA      @AA{\on{Id}}A  \\
\QCoh(Z) @>{\Upsilon_Z}>> \IndCoh(Z).
\endCD
$$

Let $\omega_{\on{HC}}$ denote the automorphism of $\on{HC}(Z)$ compatible with the functor
$\omega_Z\otimes -$. Thus, we obtain:

\begin{lem} \label{l:Xi and Psi on HH}
We have
$$\Upsilon_{\on{HC}}\simeq \omega_{\on{HC}}\circ \Xi_{\on{HC}}$$
as isomorphisms $\on{HC}^{\IndCoh}(Z)\to \on{HC}(Z)$.
\end{lem}

\ssec{$\BE_2$-algebras arising from groupoids}   \label{ss:E2 and grpds}

\sssec{}  \label{sss:groupoids}

Let $Z$ be an affine DG scheme, and let
\begin{gather}
\xy
 (-15,0)*+{Z}="X";
(15,0)*+{Z}="Y";
(0,15)*+{\CG}="Z";
(0,-5)*+{Z}="W";
{\ar@{->}_{p_1} "Z";"X"};
{\ar@{->}^{p_2} "Z";"Y"};
{\ar@{->}^{\on{unit}} "W";"Z"};
\endxy
\end{gather}
be a groupoid acting on $Z$, where $\CG$ is itself an affine DG scheme.

\medskip

The category $\QCoh(\CG)$ acquires a natural monoidal structure
via the convolution product, and as such it acts on $\QCoh(Z)$.
The unit object in $\QCoh(\CG)$ is
\[\on{unit}_*(\CO_Z)\in \QCoh(\CG).\]
Hence, its endomorphism algebra
$$\CA_\CG:=\CMaps_{\QCoh(\CG)}(\on{unit}_*(\CO_Z),\on{unit}_*(\CO_Z))$$
is naturally an $\BE_2$-algebra. 

\sssec{}

The category $\IndCoh(\CG)$ also acquires a natural monoidal structure, 
and as such it acts on $\IndCoh(Z)$, where we use the functor $f^!$ for pullback, 
$f^\IndCoh_*$ for pushforward and $\sotimes$ for tensor product, where
$$\CF_1\sotimes \CF_2:=\Delta^!(\CF_1\boxtimes \CF_2).$$

\medskip

The unit object in $\IndCoh(\CG)$ is
\[\on{unit}^\IndCoh_*(\omega_Z)\in \IndCoh(\CG),\]
and we let
$$\CA^\IndCoh_\CG=\CMaps_{\IndCoh(\CG)}(\on{unit}^\IndCoh_*(\omega_Z),\on{unit}^\IndCoh_*(\omega_Z))$$
denote the $\BE_2$-algebra of its endomorphisms.

\sssec{}

We have the following assertion:

\begin{prop}  \label{p:two versions of A}
There exists a canonical isomorphism of $\BE_2$-algebras.
$$\CA_\CG\to \CA^{\IndCoh}_\CG.$$
\end{prop} 

The proof will be given in \secref{ss:proof of two versions}.

\begin{rem}  \label{r:E2 and grpds}
As is the case with $\on{HC}(Z)$ vs. $\on{HC}^{\IndCoh}(Z)$, the assertion of \propref{p:two versions of A}
is \emph{not} essential for the contents of the paper. Namely, the object that we need to work with is
$\CA^{\IndCoh}_\CG$. So, the reader who is not interested in the proof of \propref{p:two versions of A}
can simply substitute $\CA^{\IndCoh}_\CG$ for $\CA_\CG$ for any occurrence of the latter in the main
body of the paper.
\end{rem}

\sssec{Hochschild cochains via groupoids}  \label{sss:HC via grpds}

Let $Z$ be a quasi-compact DG scheme. Consider the groupoid
$\CG=Z\times Z$; the unit section is the diagonal morphism
$$\Delta_Z:Z\to Z\times Z.$$ 
The resulting $\BE_2$-algebra $\CA_{Z\times Z}$ (resp., $\CA^\IndCoh_{Z\times Z}$) identifies with
$\on{HC}(Z)$ (resp., $\on{HC}^{\IndCoh}(Z)$), using the identifications of the monoidal categories
$$\QCoh(Z\times Z)\to \on{Funct}_{\on{cont}}(\QCoh(Z),\QCoh(Z))$$ and 
$$\IndCoh(Z\times Z)\to \on{Funct}_{\on{cont}}(\IndCoh(Z),\IndCoh(Z)),$$
respectively.

\begin{rem}
It will follow from the proof of \propref{p:two versions of A} that the resulting 
isomorphism
$$\on{HC}(Z)\to \on{HC}^{\IndCoh}(Z)$$
identifies with the one given by $\Upsilon_{\on{HC}}$. 
\end{rem}

\sssec{}

For an arbitrary groupoid $\CG$, the algebra $\CA_\CG$ (resp., $\CA^{\IndCoh}_\CG$) 
naturally maps to $\on{HC}(Z)$ (resp., $\on{HC}^{\IndCoh}(Z)$).

\medskip

This can be viewed as a corollary of the functoriality of the assignment
$\CG\rightsquigarrow\CA_\CG$. Namely, a homomorphism of groupoids $f:\CG_1\to \CG_2$ induces homomorphisms
of monoidal categories 
$$f_*:\QCoh(\CG_1)\to \QCoh(\CG_2) \text{ and }f^{\IndCoh}_*:\IndCoh(\CG_1)\to \IndCoh(\CG_2),$$
and hence of $\BE_2$-algebras 
$$\CA_{\CG_1}\to \CA_{\CG_2} \text{ and } \CA^{\IndCoh}_{\CG_1}\to \CA^{\IndCoh}_{\CG_2}.$$
We apply this to $\CG_1=\CG$ and $\CG_2=Z\times Z$ and use \secref{sss:HC via grpds}.

\medskip

Equivalently, the monoidal category $(\CA_\CG)^{\on{int-op}}\mod$ (resp., $(\CA^{\IndCoh}_\CG)^{\on{int-op}}\mod$) 
acts on the category $\QCoh(Z)$ (resp., $\IndCoh(Z)$) via 
the canonical action of $\QCoh(\CG)$ on $\QCoh(Z)$ (resp., of $\IndCoh(\CG)$ on $\IndCoh(Z)$)
by convolution. 

\begin{rem}
Again, from the construction of the isomorphism of \propref{p:two versions of A}, it will follow
that for a homomorphism of groupoids $f:\CG_1\to \CG_2$, the following diagram
of $\BE_2$-algebras naturally commutes:
$$
\CD
\CA_{\CG_1}   @>{\sim}>>   \CA^\IndCoh_{\CG_1} \\
@V{f_*}VV     @VV{f^\IndCoh_*}V   \\
\CA_{\CG_2}   @>{\sim}>>   \CA^\IndCoh_{\CG_2}.
\endCD
$$  
\end{rem}

\sssec{Relative Hochschild cochains}   \label{sss:rel hh}

Let $Z\to \CU$ be a morphism of affine DG schemes. Consider the groupoid
$$\CG_{Z/\CU}:=Z\underset{\CU}\times Z.$$

By a slight abuse of notation we will continue to denote by $\Delta_Z$ the diagonal map
$Z\to Z\underset{\CU}\times Z$, which is the unit of the above groupoid.

\medskip

The groupoid gives rise to $\BE_2$-algebras $\CA_{\CG_{Z/\CU}}$ and $\CA^{\IndCoh}_{\CG_{Z/\CU}}$ of endomorphisms
of the unit objects $\Delta_{Z,*}\CO_Z\in\QCoh(\CG_{Z/\CU})$ and $\Delta^{\IndCoh}_{Z,*}\omega_Z\in\IndCoh(\CG_{Z/\CU})$, respectively.
We put
\begin{align*}
\on{HC}(Z/CU)&:=\CA_{\CG_{Z/\CU}}\\
\on{HC}^{\IndCoh}(Z/\CU)&:=\CA^{\IndCoh}_{\CG_{Z/\CU}}
\end{align*} 
and refer to these $\BE_2$-algebras 
as ``the algebras of Hochschild cochains on $Z$ relative to $\CU$." 

\medskip

It is easy to see that $\on{HC}(Z/\CU)$ identifies with the $\BE_2$-algebra 
$$\on{HC}(\QCoh(Z))_{\QCoh(\CU)}$$
of endomorphisms of the
identity functor on $\QCoh(Z)$ as a DG category tensored
over $\QCoh(\CU)$, i.e., the $\BE_2$-algebra of endomorphisms of the unit in the monoidal
category
$$\on{Funct}_{\QCoh(\CU)}(\QCoh(Z),\QCoh(Z)).$$
(This is because the natural homomorphism from $\QCoh(\CG_{Z/\CU})$ to the 
above category is an equivalence.)

\medskip

By essentially repeating the proof of \propref{p:two isoms for HH}, one can construct an isomorphism of
$\BE_2$-algebras
\[\on{HC}(Z/\CU)=\on{HC}(\QCoh(Z))_{\QCoh(\CU)}\overset{\sim}\to \on{HC}(\IndCoh(Z))_{\QCoh(\CU)}.\]
Here $\on{HC}(\IndCoh(Z))_{\QCoh(\CU)}$ is the $\BE_2$-algebra 
of endomorphisms of the
identity functor on $\IndCoh(Z)$ as a DG category tensored
over $\QCoh(\CU)$, that is, the $\BE_2$-algebra of endomorphisms of the unit in the monoidal
category
$$\on{Funct}_{\QCoh(\CU)}(\IndCoh(Z),\IndCoh(Z)).$$

\medskip

The content of \propref{p:two versions of A} in this case is that the natural map 
$$\CA^{\IndCoh}_{\CG_{Z/\CU}}\to \on{Funct}_{\QCoh(\CU)}(\IndCoh(Z),\IndCoh(Z))$$
is an isomorphism.

\ssec{Proof of \propref{p:two versions of A}}  \label{ss:proof of two versions}

\sssec{}     \label{sss:monoidal non-cocompl}

We will need the following variant of the construction of the assignment $\bO\mapsto \CA_\bO$
of \secref{ss:E2 and monoidal}.

\medskip

Let $'\bO$ be a \emph{non-cocomplete} DG category equipped with a monoidal structure. 
The datum of such a monoidal structure is, by definition, equivalent to
a datum of monoidal structure on the cocomplete DG category $\bO:=\Ind({}'\bO)$ such that the monoidal operation 
$\bO\times \bO\to \bO$ sends $'\bO\times {}'\bO\subset \bO\times \bO$ to $'\bO\subset \bO$. 

\medskip

We set, by definition, $\CA_{'\bO}:=\CA_\bO$. The assignment $'\bO\to \CA_{'\bO}$ is clearly functorial. In addition,
it has the following property:

\medskip

Let $'\Phi:{}'\bO\to{}'\bO_1$ be a monoidal functor, where $\bO_1$ is another monoidal DG category, not necessarily cocomplete. Then 
$'\Phi$ induces a continuous monoidal functor \[\Phi:\bO\to\bO_1,\]
where $\bO_1=\Ind({}'\bO_1)$, and provides a homomorphism \[\CA_{'\bO}=\CA_\bO\to \CA_{\bO_1}=\CA_{{}'\bO_1}.\] Assume now
that the original functor $'\Phi$ is fully faithful (but $\Phi$ does not have to be). Then the above homomorphism
$\CA_{'\bO}\to \CA_{'\bO_1}$ is an isomorphism. In fact, it suffices to assume that $\Phi'$ induces an isomorphism 
\[\CMaps_{'\bO}({\bf 1}_{'\bO},{\bf 1}_{'\bO})\to\CMaps_{`\bO_1}({\bf 1}_{'\bO_1},{\bf 1}_{'\bO_1}).\]
Equivalently, $\Phi'$ needs to be fully faithful on the full subcategory
of $'\bO$ that is \emph{strongly generated} by the unit object ${\bf 1}_{'\bO}\in{}'\bO$ (that is, on the smallest full DG subcategory
of $'\bO$ containing ${\bf 1}_{'\bO}$).

\sssec{}

Consider the full, but not cocomplete, subcategory 
$$\QCoh(\CG)^{\on{bdd.Tor}} \subset \QCoh(\CG),$$
consisting of the objects whose Tor amplitude with respect to $p_2$ is bounded on the left. That is,
$\CF\in \QCoh(\CG)$ if and only if the functor
$$\CF'\mapsto \CF\otimes p_2^*(\CF),\quad \QCoh(Z)\to \QCoh(\CG)$$
is of \emph{bounded cohomological amplitude on the left}. For example,
the object $\on{unit}_*(\CO_Z)$ belongs to $\QCoh(\CG)^{\on{bdd.Tor}}$.
(In what follows, we could actually replace the category $\QCoh(\CG)^{\on{bdd.Tor}}$
by its full subcategory strongly generated by the unit object $\on{unit}_*(\CO_Z)$.) 

\medskip

The monoidal structure on $\QCoh(\CG)$ preserves $\QCoh(\CG)^{\on{bdd.Tor}}$, which therefore
acquires a structure of a monoidal (non-cocomplete) DG category. 
Similarly, the monoidal structure on $\IndCoh(\CG)$ preserves $\IndCoh(\CG)^+$, and therefore
$\IndCoh(\CG)^+$ acquires a structure of a monoidal (non-cocomplete) DG category. 

\medskip

By \secref{sss:monoidal non-cocompl}, it is sufficient to construct a monoidal functor 
\[\QCoh(\CG)^{\on{bdd.Tor}}\to\IndCoh(\CG)^+\]
and verify that it is fully faithful on the subcategory strongly generated by the unit object.

\sssec{}

The operation of convolution in $\QCoh(\CG)$ defines an action of the monoidal category
$\QCoh(\CG)^{\on{bdd.Tor}}$ on $\QCoh(\CG)^+$. 

\medskip

Recall that the functor $\Psi_\CG$ defines an equivalence $\IndCoh(\CG)^+\to \QCoh(\CG)^+$. 
Hence, we obtain a monoidal action of $\QCoh(\CG)^{\on{bdd.Tor}}$ on $\IndCoh(\CG)^+$. 

\begin{lem} \label{lm:leftandright}
The above action of $\QCoh(\CG)^{\on{bdd.Tor}}$ on $\IndCoh(\CG)^+$ commutes with
the right action of $\IndCoh(\CG)^+$ on itself that is induced by the monoidal structure on
$\IndCoh(\CG)$.
\end{lem}

\begin{proof}
Follows from \cite[Lemma 3.6.13]{IndCoh}.
\end{proof}

On the other hand, the monoidal category of endomorphisms of the right $\IndCoh(\CG)^+$-module $\IndCoh(\CG)^+$ is identified 
with $\IndCoh(\CG)^+$ via its left action on itself. Thus,
we obtain a monoidal functor
$$\QCoh(\CG)^{\on{bdd.Tor}}\to \IndCoh(\CG)^+.$$

\medskip

If one ignores the monoidal structure, the functor can be given explicitly by
\begin{equation}\label{e:from bdd to ind+}
\CF\mapsto \CF\otimes p_2^*(\Psi_Z(\omega_Z))\in \QCoh(\CG)^+\simeq \IndCoh(\CG)^+\qquad(\CF\in\QCoh(\CG)^{\on{bdd.Tor}}).
\end{equation}
The advantage of using Lemma~\ref{lm:leftandright} to construct the functor (rather than treating \eqref{e:from bdd to ind+} as its 
definition) is that the monoidal structure appears naturally.

\sssec{}

It remains to show that the map
\begin{equation} \label{e:from A to ind}
\CMaps_{\QCoh(\CG)}(\on{unit}_*(\CO_Z),\on{unit}_*(\CO_Z))\to
\CMaps_{\IndCoh(\CG)}(\on{unit}^\IndCoh_*(\omega_Z),\on{unit}^\IndCoh_*(\omega_Z))
\end{equation}
induced by \eqref{e:from bdd to ind+}
is an isomorphism. 

\sssec{}

More generally, let $f:Z'\to Z$ be a map between affine DG schemes, and consider the corresponding subcategory
$$\QCoh(Z')^{\on{bdd.Tor}/Z}\subset \QCoh(Z),$$
consisting of objects $\CF'\in \QCoh(Z')$, for which the functor
$$\CF\mapsto \CF'\otimes f^*(\CF),\quad \QCoh(Z)\to \QCoh(Z')$$
is of bounded cohomological amplitude on the left (i.e., these are objects whose Tor dimension over $Z$ is bounded on the left). 

\medskip

Consider the corresponding  functor $$\Phi:\QCoh(Z')^{\on{bdd.Tor}/Z}\to \IndCoh(Z'),$$
defined by
$$\Phi(\CF'):=\CF'\otimes f^*(\Psi_Z(\omega_Z))\in \QCoh(Z')^+\simeq \IndCoh(Z')^+\subset \IndCoh(Z').$$

We will claim:

\begin{prop}
The functor $\Phi$ induces an isomorphism
$$
\CMaps_{\QCoh(Z')}(\CF'_1,\CF'_2)\to
\CMaps_{\IndCoh(Z')}(\Phi(\CF'_1),\Phi(\CF'_2)),
$$
whenever $\CF'_2\in \QCoh(Z')^{\on{bdd.Tor}/Z}\cap \QCoh(Z')^-$.
\end{prop}

\begin{proof}

Since $Z'$ was assumed affine, we can take 
$\CF'_1=\CO_{Z'}$. Denote $\CF'_2$ by $\CF'$. By the projection formula, it suffices to show that the map
$$\Gamma(Z,f_*(\CF'))\to \CMaps_{\QCoh(Z)}(\Psi_Z(\omega_Z),f_*(\CF')\otimes \Psi_Z(\omega_Z))$$
is an isomorphism.  

\medskip

The assumption that $\CF'\in \QCoh(Z')^{\on{bdd.Tor}/Z}$ implies that the right-hand side identifies
with $\CMaps_{\IndCoh(Z)}(\omega_Z, \CF\otimes \omega_Z)$, where $\CF:=f_*(\CF')$, and
$\otimes$ denotes the monoidal action of $\QCoh(Z)$ on $\IndCoh(Z)$.

\medskip

Hence, the required assertion follows from the next general lemma:

\begin{lem}  \label{l:psych}
Let $Z$ be a DG scheme. Then the map
\begin{equation} \label{e:psych}
\Gamma(Z,\CF)\to \CMaps_{\IndCoh(Z)}(\omega_Z,\CF\otimes \omega_Z)
\end{equation}
is an isomorphism, provided that $\CF\in \QCoh(Z)^-$.
\end{lem}

\end{proof}

\begin{proof}[Proof of \lemref{l:psych}]

Recall that if $Z$ is eventually coconnective, \eqref{e:psych} is an isomorphism for any $\CF\in\QCoh(Z)$,
see \cite[Corollary 9.6.3]{IndCoh}. Let $i_k:Z_k\to Z$ be the $k$-coconnective truncation of $Z$. Then
\[\omega_Z=\underset{\longrightarrow}{colim}\, (i_k)^\IndCoh_*(\omega_{Z_k}).\]

Therefore, the right-hand side of \eqref{e:psych} can be evaluated as
\begin{multline*}
\Hom_{\IndCoh(Z)}(\omega_Z,\CF\otimes\omega_Z)\simeq\underset{\longleftarrow}{lim}
\Hom_{\IndCoh(Z_k)}(\omega_{Z_k},(i_k)^!(\CF\otimes\omega_Z))\\
\simeq\underset{\longleftarrow}{lim} \Hom_{\IndCoh(Z_k)}(\omega_{Z_k},(i_k)^*(\CF)\otimes\omega_{Z_k})
\simeq\underset{\longleftarrow}{lim} \Hom_{\QCoh(Z_k)}(\CO_{Z_k},(i_k)^*(\CF))\\
\simeq\underset{\longleftarrow}{lim} \Hom_{\QCoh(Z)}(\CO_Z,(i_k)_*\circ (i_k)^*(\CF)).
\end{multline*}
It remains to notice that $\CF$ maps isomorphically to the inverse limit
\[\underset{\longleftarrow}{lim} (i_k)_*\circ (i_k)^*(\CF),\]
which is true because $\CF\in\QCoh(Z)^-$: indeed, if $\CF\in \QCoh(Z)^{\leq 0}$, then the map
$$\CF\to  (i_k)_*\circ (i_k)^*(\CF)$$ induces an isomorphism in the cohomological degrees $\geq -k$.
\end{proof}

\section{Hochschild cohomology and Lie algebras}    \label{s:Lie}

In this appendix we connect $\BE_2$-algebras arising from groupoids as in \secref{s:Hoch}
to Lie algebras in the symmetric monoidal category $\IndCoh$. \footnote{As we will be dealing
with Lie algebras, the assumption that our ground field $k$ has characteristic $0$
is crucial.} 

\ssec{Lie algebras arising from group DG schemes}

In this subsection we quote (without proof) two facts (Propositions \ref{p:Lie 1} and \ref{p:Lie 2}) 
about the relationship between group DG schemes and Lie algebras. The full 
exposition appears in \cite[Part IV.3, Sect. 3]{Algebroids}. 

\sssec{}

Let $Z$ be an affine DG scheme, and let 
$$p:\CG\rightleftarrows Z:\on{unit}$$
be a group DG scheme over $Z$. We claim:

\begin{prop} \label{p:Lie 1}  
The object $\on{unit}^*(T^*(\CG/Z))\in \QCoh(Z)$ has a natural structure of Lie co-algebra in
the symmetric monoidal category $\QCoh(Z)$.
\end{prop}

\sssec{}

Recall that all DG schemes are assumed to be almost of finite type. We have 
$$T^*(\CG/Z)\in \QCoh(Z)^-,$$ and note that it has \emph{coherent cohomologies}. Serre
duality identifies this category with a full subcategory of $\QCoh(Z)^+\simeq \IndCoh(Z)^+$
consisting of objects with coherent cohomologies. 

\medskip

Hence, 
$$T^{\on{IndCoh}}(\CG/Z):=\bD_Z^{\on{Serre}}(T^*(\CG/Z))\in \IndCoh(Z)$$
acquires a structure of Lie algebra, where $\IndCoh(Z)$ is regarded as a symmetric monoidal
category under the $\sotimes$ tensor product. 

\sssec{}

Let $\IndCoh(\CG)_Z\subset\IndCoh(\CG)$ be the full subcategory of objects set-theoretically supported
on the unit section. The embedding  $\IndCoh(\CG)_Z\hookrightarrow\IndCoh(\CG)$ admits a right adjoint, which we denote by
\[\CF\mapsto \CF^\wedge:\IndCoh(\CG)\to\IndCoh(\CG)_Z.\]
Consider the
object 
$$\omega_{\CG}^\wedge\in \IndCoh(\CG),$$
and its direct image $p^\IndCoh_*(\omega_{\CG}^\wedge)\in \IndCoh(Z)$.

\medskip

The group structure on $\CG$ makes $p^\IndCoh_*(\omega_{\CG}^\wedge)$ into an associative algebra
object in $\IndCoh(Z)$ (where $\IndCoh(Z)$ is, as always, considered as a symmetric monoidal
category under the $\sotimes$ tensor product). We claim:

\begin{prop} \label{p:Lie 2}
There is a canonical isomorphism of associative algebras in $\IndCoh(Z)$
$$p^\IndCoh_*(\omega_{\CG}^\wedge)\simeq U(T^{\on{IndCoh}}(\CG/Z)).$$
\end{prop} 

\sssec{}  \label{sss:Lie when perfect}

Assume now that $\CG$ is such that $T^*(\CG/Z)$ belongs to $\QCoh(Z)^{\on{perf}}$. Denote by
$T(\CG/Z)\in \QCoh(Z)^{\on{perf}}$ its monoidal dual.  

\medskip

By \propref{p:Lie 1}, the object $T(\CG/Z)$
acquires a structure of Lie algebra \emph{in the symmetric monoidal category} $\QCoh(Z)$.

\medskip

Recall also that the functor 
$$\Upsilon_Z=-\otimes \omega_Z:\QCoh(Z)\to \IndCoh(Z)$$
is symmetric monoidal.  By construction, we have
$$T^{\on{IndCoh}}(\CG/Z)\simeq \Upsilon_Z(T(\CG/Z)),$$
as Lie algebras in $\IndCoh(Z)$, and hence 
$$U(T^{\on{IndCoh}}(\CG/Z))\simeq \Upsilon_Z\left(U_{\CO_Z}(T(\CG/Z))\right),$$
as associative algebras in $\IndCoh(Z)$.

\ssec{Lie algebras arising from groupoids}

In this subsection we will relate the $\BE_2$-algebra $\on{HC}^{\IndCoh}(Z)$ to the
tangent sheaf of $Z$. This relationship is convenient for the direct invariant definition 
of singular support, as is given in \secref{ss:sing}.

\sssec{}   

Consider now the following situation. Let $i:Z\to W$ be a proper map of affine DG schemes,
equipped with a retraction, i.e., a map $s:W\to Z$ and an isomorphism $s\circ i\sim \on{id}_Z$.

\medskip

Note that in this case the groupoid $\CG_{Z/W}:=Z\underset{W}\times Z$ over $Z$ is actually a group
DG scheme over $Z$. 

\medskip

Assume now that $T^*(Z/W)$ is perfect, and consider its monoidal dual $T(Z/W)\in \QCoh(Z)^{\on{perf}}$. We will prove:

\begin{prop} \label{p:groupoid 1} \hfill

\medskip

\noindent{\em(a)} The object $T(Z/W)\in \QCoh(Z)$ has a structure of Lie algebra.

\medskip

\noindent{\em(b)} There is a canonically defined homomorphism of DG associative algebras 
$$\Gamma\left(Z,U_{\CO_Z}(T(Z/W))\right)\to 
\CMaps_{\IndCoh(W)}(i^\IndCoh_*(\omega_Z),i^\IndCoh_*(\omega_Z)),$$
which is an isomorphism if $Z$ is eventually coconnective. 
\end{prop}

\sssec{}

Before proving \propref{p:groupoid 1} let us show how we are going to apply it. Suppose we are in situation of
\secref{sss:groupoids}; thus, $\CG$ is a groupoid acting on $Z$, where both $Z$ and $\CG$ are affine DG schemes.
Assume that the relative cotangent complex $T^*(\CG/Z)$  (with respect to, say, projection $p_1$) is perfect; let
$T(\CG/Z)\in \QCoh(\CG)^{\on{perf}}$ denote its monoidal dual. 

\medskip

From \propref{p:groupoid 1} we obtain:

\begin{cor} \label{c:groupoid Lie}  \hfill

\medskip

\noindent{\em(a)} The object $\on{unit}^*(T(\CG/Z))[-1]\in \QCoh(Z)$ has a natural structure of Lie algebra.

\medskip

\noindent{\em(b)} There exists a canonically defined homomorphism 
$$\Gamma\left(Z,U_{\CO_Z}(\on{unit}^*(T(\CG/Z))[-1])\right)\to \CA^{\IndCoh}_\CG,$$
which is an isomorphism if $Z$ is eventually coconnective. 

\end{cor}

\begin{proof}

Apply \propref{p:groupoid 1} to $W:=\CG$ with the map $i$ being the unit map $Z\to \CG$ and the retraction being,
say, the first projection $p_1:\CG\to Z$. It remains to use the canonical identification
$$\on{unit}^*(T(\CG/Z))[-1]\simeq T(Z/\CG).$$

\end{proof} 

\begin{rem}
Note that the structure of Lie algebra on $\on{unit}^*(T(\CG/Z)[-1]$ depends on the choice of the retraction
of the map $\on{unit}:Z\to \CG$. If we chose a different retraction, namely, $p_2$ instead of $p_1$,
the resulting Lie algebra structure would be different. In general there does not
exist an isomorphism between the two resulting Lie algebras that induces the identity map on the
underlying object of $\QCoh(Z)$. 
\end{rem}

\sssec{}  \label{sss:HH as univ env}

Consider the particular case when $\CG=Z\times Z$. Assume that the cotangent complex of $Z$ is perfect. 
We obtain the following basic identification:

\begin{cor}  \label{c:HH as univ env}  \hfill

\smallskip

\noindent{\em(a)} The object $T(Z)[-1]$ has a natural structure of Lie algebra in the symmetric monoidal category $\QCoh(Z)$.

\smallskip

\noindent{\em(b)} 
The associative DG algebra underlying the $\BE_2$-algebra $\on{HC}^{\IndCoh}(Z)$ receives a canonical 
homomorphism from $\Gamma(Z,U_{\CO_Z}(T(Z)[-1]))$. This homomorphism is an isomorphism
if $Z$ is eventually coconnective. 
\end{cor}

More generally, for a map $Z\to \CU$ whose relative cotangent complex $T^*(Z/\CU)$ is perfect, we obtain: 

\begin{cor}  \label{c:HH as univ env rel} \hfill

\smallskip

\noindent{\em(a)} The object $T(Z/\CU)[-1]$ has a natural structure of Lie algebra in the symmetric monoidal category $\QCoh(Z)$.

\smallskip

\noindent{\em(b)} The associative DG algebra underlying the $\BE_2$-algebra $\on{HC}^{\IndCoh}(Z/\CU)$ 
receives a canonical homomorphism from $\Gamma(Z,U_{\CO_Z}(T(Z/U)[-1]))$. This homomorphism is an isomorphism
if $Z$ is eventually coconnective. 
\end{cor}

\sssec{Proof of \propref{p:groupoid 1}}

Let us start with a map $i:Z\to W$. Consider the groupoid $\CG_{Z/W}=Z\underset{W}\times Z$:
\begin{gather}
\xy
 (-15,0)*+{Z}="X";
(15,0)*+{Z.}="Y";
(0,15)*+{\CG_{Z/W}}="Z";
(0,-5)*+{Z}="W";
{\ar@{->}_{p_1} "Z";"X"};
{\ar@{->}^{p_2} "Z";"Y"};
{\ar@{->}^{\on{unit}} "W";"Z"};
\endxy
\end{gather}

By proper base change (see also \cite[Part II.2, Sect. 5.3]{Algebroids}), the monad
$i^!\circ i_*^{\IndCoh}$ acting on $\IndCoh(Z)$ identifies with the monad $(p_2)_*^\IndCoh\circ p_1^!$.

\medskip

In particular, for $\CF\in \IndCoh(Z)$, the structure of associative DG algebra on
\begin{equation} \label{e:ident E1 prel}
\CMaps_{\IndCoh(Z)}(\CF,(p_2)_*^\IndCoh\circ p_1^!(\CF))
\end{equation} 
identifies canonically with
$$\CMaps_{\IndCoh(W)}(i^{\IndCoh}_*(\CF),i^{\IndCoh}_*(\CF)).$$

\medskip

In particular, we obtain an identification of associative DG algebras
\begin{equation} \label{e:ident E1}
\CMaps_{\IndCoh(Z)}(\omega_Z,(p_2)_*^\IndCoh\circ p_1^!(\omega_Z))\simeq
\CMaps_{\IndCoh(W)}(i^{\IndCoh}_*(\omega_Z),i^{\IndCoh}_*(\omega_Z)).
\end{equation} 

\medskip

Assume now that the map $i:Z\to W$ is equipped with a retraction. In this case
$\CG_{Z/W}$ is a group DG scheme over $Z$, and by \secref{sss:Lie when perfect}, the object
$$\on{unit}^*(T(\CG_{Z/W}/Z))\in \QCoh(Z)$$ acquires a Lie algebra structure.  Point (a) of
\propref{p:groupoid 1} follows now by noticing that 
$$\on{unit}^*(T(\CG_{Z/W}/Z))\simeq T(Z/W).$$

\medskip

Denote $\CG:=\CG_{Z/W}$ and let $p$ denote its projection to $Z$ (which is canonically identified
with both $p_1$ and $p_2$). To prove point (b), we note that the associative algebra structure on 
\eqref{e:ident E1 prel} comes from the structure of associative algebra on $\IndCoh(Z)$ on
$p^\IndCoh_*(\omega_\CG)$, induced by the group structure on $\CG$. 

\medskip

Note now that the existence of the retraction $s$ implies that the map $i:Z\to W$ is a closed embedding. 
Hence, the maps $p:\CG\to Z$ and $\on{unit}:Z\to \CG$ induce isomorphisms of the underlying reduced
classical schemes. In particular, $\omega_\CG\simeq \omega^\wedge_\CG$. Therefore, combining 
\propref{p:Lie 2} and \eqref{e:ident E1},  we obtain an identification of associative DG algebras
\begin{equation} \label{e:ident E2}
\CMaps_{\IndCoh(Z)}\left(\omega_Z,U(T^{\on{IndCoh}}(\CG/Z))\right)\simeq 
\CMaps_{\IndCoh(W)}(i^{\IndCoh}_*(\omega_Z),i^{\IndCoh}_*(\omega_Z)).
\end{equation} 

\medskip

Finally, the symmetric monoidal functor $\Upsilon_Z$ defines a homomorphism
\begin{multline} \label{e:ident E3}
\Gamma\left(Z,U(\on{unit}^*(T(\CG/Z))[-1])\right)\simeq
\CMaps_{\QCoh(Z)}\left(\CO_Z,U_{\CO_Z}(T(\CG/Z))\right)\to \\
\to \CMaps_{\IndCoh(Z)}\left(\omega_Z,U(T^{\on{IndCoh}}(\CG/Z))\right).
\end{multline} 

Composing \eqref{e:ident E2} and \eqref{e:ident E3}, we obtain the desired map
\begin{equation} \label{e:ident E4}
\Gamma\left(Z,U(\on{unit}^*(T(\CG/Z))[-1])\right) \to \CMaps_{\IndCoh(W)}(i^{\IndCoh}_*(\omega_Z),i^{\IndCoh}_*(\omega_Z)).
\end{equation} 

If $Z$ is eventually coconnective, the functor $\Upsilon_Z$ is fully faithful, which implies that \eqref{e:ident E3}
is an isomorphism. Hence, \eqref{e:ident E3} is an isomorphism as well. 

\qed

\ssec{Compatibility with duality}   

\sssec{}

Consider the homomorphism of $\BE_1$-algebras
$$\Gamma(Z,U_{\CO_Z}(T(Z)[-1]))\to \on{HC}^{\IndCoh}(Z),$$
and the resulting map in $\Vect$
\begin{equation} \label{e:vector fields to HH}
\Gamma(Z,T(Z))[-1]\to \on{HC}^{\IndCoh}(Z).
\end{equation}

\sssec{}

We claim:

\begin{lem}  \label{l:Serre abs}
The following diagram commutes
$$
\CD
\Gamma(Z,T(Z))[-1] @>{\text{\eqref{e:vector fields to HH}}}>>  \on{HC}^{\IndCoh}(Z) \\
@V{\xi\mapsto -\xi}VV   @VV{\text{\eqref{e:inv on HH Ind}}}V   \\ 
\Gamma(Z,T(Z))[-1] @>{\text{\eqref{e:vector fields to HH}}}>>  \on{HC}^{\IndCoh}(Z)^{\on{ext-op}}
\endCD
$$
in $\Vect$.
\end{lem}

\begin{proof}

The map \eqref{e:vector fields to HH} is the composition of the natural morphism
\begin{multline*}
\Gamma(Z,T(Z))[-1]\simeq \Gamma(Z,N_{Z/Z\times Z})[-1]\to \\
\to \CMaps_{\IndCoh(Z\times Z)}((\Delta_Z)^\IndCoh_*(\omega_Z),(\Delta_Z)^\IndCoh_*(\omega_Z)),
\end{multline*}
and the homomorphism of objects of $\Vect$ (in fact, $\BE_1$-algebras)
$$\CMaps_{\IndCoh(Z\times Z)}((\Delta_Z)^\IndCoh_*(\omega_Z),(\Delta_Z)^\IndCoh_*(\omega_Z))\simeq \on{HC}^{\IndCoh}(Z).$$

Unwinding the definitions, we obtain that the isomorphism \eqref{e:inv on HH Ind} corresponds to the automorphism of 
$\CMaps_{\IndCoh(Z\times Z)}((\Delta_Z)^\IndCoh_*(\omega_Z),(\Delta_Z)^\IndCoh_*(\omega_Z))$, induced by the transposition
involution $\sigma$ on $Z\times Z$. 

\medskip

The assertion of the lemma follows now from the commutativity of the diagram
$$
\CD
T(Z) @>>>  N_{Z/Z\times Z}  \\
@V{\xi\mapsto -\xi}VV    @VV{\sigma}V  \\
T(Z) @>>>  N_{Z/Z\times Z}.
\endCD
$$

\end{proof}

\sssec{}  \label{sss:proof of Serre abs}

We are now ready to give a more direct proof of \propref{p:Serre}:

\begin{proof}

The assertion is local, so we can assume that $Z$ is affine. Fix $\CF\in\Coh(Z)$. It suffices to show that the action maps
$$\Gamma(Z,T(Z))[-1]\otimes \CF\to \CF \text{ and }
\Gamma(Z,T(Z))[-1]\otimes \BD_{Z}^{\on{Serre}}(\CF)\to \BD_{Z}^{\on{Serre}}(\CF)$$
correspond to each other under the automorphism of $\Gamma(Z,T(Z))[-1]$ given by $\xi\mapsto-\xi$. (In this statement, $\Gamma(Z,T(Z))[-1]$ 
appears merely as an object of $\Vect$: we make no statement about compatibility with the Lie algebra structure.) 

\medskip

The required statement follows from \lemref{l:Serre abs}
and the following observation:

\medskip

Let $\bC$ be a dualizable DG category, and $\bc\in \bC$ a compact object. Let $\bc^\vee$ denote the
corresponding compact object of $\bC^\vee$. Then under the isomorphism of $\BE_2$-algebras
$$\on{HC}(\bC)^{\on{ext-op}}\simeq \on{HC}(\bC^\vee)$$
(see \secref{sss:E2 action on dual}), the diagram
$$
\CD
\on{HC}(\bC)^{\on{op}} @>>>  \CMaps_{\bC}(\bc,\bc)^{\on{op}}  \\
@VVV    @VVV   \\
\on{HC}(\bC^\vee) @>>>  \CMaps_{\bC^\vee}(\bc^\vee,\bc^\vee),
\endCD
$$
commutes, where the horizontal arrows use the \emph{external} forgetful functor
$\BE_2\Alg\to \BE_1\Alg$.

\end{proof}

\ssec{The case of group-schemes} 

\sssec{}

Suppose that in the setting of \secref{sss:groupoids}, $\CG$ is a group DG scheme
over $Z$. In addition, we continue to assume that the relative
cotangent complex $T^*(\CG/Z)$ is perfect. In this case, we claim:

\begin{prop}  \label{p:structure abelian}
The Lie algebra structure on $\on{unit}^*(T(\CG/Z))[-1]\in \QCoh(Z)$, given by \corref{c:groupoid Lie}(a), is canonically trivial.
\end{prop}

\begin{proof}

We will deduce the assertion of the proposition from the following lemma\footnote{This lemma
is probably well known; the proof is given in \cite[Part IV.2, Theorem 2.2.2]{Algebroids}.}:

\begin{lem} \label{l:loop Lie}
Let $L$ be a Lie algebra in a symmetric monoidal category $\bO$ over a field of characteristic $0$.
Then the loop object $\Omega(L)$, considered as a plain Lie algebra in $\bO$ is canonically abelian,
i.e., identifies with object $L[-1]\in \bO$ with the trivial Lie algebra structure.
\end{lem}

Namely, we claim that the Lie algebra $\on{unit}^*(T(\CG/Z))[-1]\in \QCoh(Z)$, given by \corref{c:groupoid Lie}(a), 
is canonically isomorphic to 
$$\Omega(\on{unit}^*(T(\CG/Z))),$$ where $\on{unit}^*(T(\CG/Z))$ is the Lie algebra of \secref{sss:Lie when perfect}.

\medskip

For the latter, it suffices to notice that we have an identification of group DG schemes over $Z$:
$$Z\underset{\CG}\times Z\simeq \Omega(\CG).$$

\end{proof}

\begin{rem}
We emphasize that the isomorphism between $\Omega(L)$ and the abelian Lie algebra $L[-1]$
given by \lemref{l:loop Lie} does not
respect the structure of group-objects in the category of Lie algebras.
\end{rem}

\sssec{}

Combining \propref{p:structure abelian} with \corref{c:groupoid Lie}(b), we obtain:

\begin{cor}  \label{c:HC as E_1 for groups} 
The associative DG algebra underlying the $\BE_2$-algebra $\CA^{\IndCoh}_\CG$ receives a canonically
defined homomorphism from $\Gamma(Z,\on{Sym}_{\CO_Z}(L_\CG[-1]))$. This homomorphism is
an isomorphism if $Z$ is eventually coconnective. 
\end{cor}

\begin{rem}  \label{r:not E2}
Let us note that each side in the homomorphism
$$\Gamma(Z,\on{Sym}_{\CO_Z}(L_\CG[-1]))\to \CA^{\IndCoh}_\CG$$
of \corref{c:HC as E_1 for groups} has a structure of $\BE_2$-algebra. However, this homomorphism does not
respect this structure: it is a homomorphism of mere $\BE_1$-algebras. 
\end{rem}

\section{Other approaches to singular support}  \label{s:rev}

\ssec{$\IndCoh(Z)$ via the category of singularities}

\sssec{}\label{sss:parallel gci}
 Let us assume that $Z$ is an affine DG scheme, which is a global complete intersection. That is,
$Z$ can be included in a Cartesian diagram
$$
\CD
Z  @>{\iota}>> \CU  \\ 
@VVV    @VVV \\
\on{pt}  @>>>  \CV,
\endCD
$$
where $\CU$ and $\CV$ are smooth. Moreover, we assume that $\CV$ is parallelized; this allows us to replace $\CV$ with
its tangent space at the fixed point. Thus, we will assume that $\CV=V$ is a finite-dimensional vector 
space.

\begin{rem} 
In fact, the construction of this section remain valid in the setting of \secref{ss:parallel stacks}: that is, we may replace $Z$ with the 
zero locus of a section of a vector bundle on a smooth stack.
However, one can use affine charts of a stack to deduce this more general case from the special case that we consider here.
\end{rem}

\sssec{} Consider the product $V^*\times\CU$. The map $s:\CU\to V$ defines a function
\[
V^*\times\CU\to\BA^1:(u,\phi)\mapsto\langle\phi,s(u)\rangle,
\]
which we still denote by $s$. Let 
\[H:=(V^*\times\CU)\underset{\BA^1}\times\on{pt}\] 
be the zero locus of $s$. Note that $H$ is a classical
scheme (and then a closed hypersurface in $V^*\times\CU$) unless $s$ vanishes identically
on a connected component of $\CU$. Clearly, $H$ is conical, that is, it 
carries a natural action of $\BG_m$ lifting its action on $V^*\times\CU$ by dilations.
Moreover, it is easy to see that the singular locus of $H$ is identified with $\Sing(Z)\subset V^*\times\CU$.
 
\medskip

Let 
\[\CH:=H/\BG_m\simeq\left((V^*/\BG_m)\times\CU\right)\underset{\BA^1/\BG_m}\times(\on{pt}/\BG_m)\] 
be the quotient stack.
The following theorem is due to M.~Umut Isik: 

\begin{thm}\label{t:I} There is a natural equivalence
\begin{equation}\label{e:I}
\IndCoh(Z)\simeq\IndCoh(\CH)/\QCoh(\CH),
\end{equation}
where $\QCoh(\CH)$ is viewed as a full subcategory of $\IndCoh(\CH)$ using the functor $\Xi_\CH$.
\end{thm}

\thmref{t:I} is a version of \cite[Theorem~3.6]{Isik}. In \cite{Isik}, it is assumed that $Z$ is
classical (that is, that the coordinates of the map $s$ form a regular sequence of functions), but
this assumption is not used in the proof. 

\begin{rem} Let $S$ be a quasi-compact DG scheme. The category $\IndCoh(S)/\QCoh(S)$ is compactly
generated by the quotient
\[\Coh(S)/\QCoh(S)^{\on{perf}}.\]
The category $\Coh(S)/\QCoh(S)^{\on{perf}}$ is known as the \emph{category of singularities} of $S$, 
introduced in \cite{Or1}. Then $\IndCoh(S)/\QCoh(S)$ identifies with the ind-completion of the category
of singularities. The category $\IndCoh(S)/\QCoh(S)$ was introduced in \cite{KrS} under the name
``stable derived category." (As a minor detail, both \cite{Or1} and \cite{KrS} work with Noetherian classical schemes.)

\sssec{}

Theorem~\ref{t:I} is about the ``stable derived category" $\IndCoh(\CH)/\QCoh(\CH)$ of the stack $\CH$.

\medskip

Note that both $\IndCoh(\CH)$ and $\QCoh(\CH)$ are compactly generated: the former because $\CH$ is QCA (see \cite{DrGa0}),
the latter because $\CH$ is a quotient of an affine DG scheme by a linear group, which is a perfect stack (see \cite{BFN}).
(Alternatively, the two categories are compactly generated by \corref{c:category comp gen stacks}, since $\CH$ is
a global complete intersection.) Therefore, $\IndCoh(\CH)/\QCoh(\CH)$ is equivalent to the ind-completion of the
``category of singularities"
\[\Coh(\CH)/\QCoh(\CH)^{\on{perf}}.\]
In fact, \cite[Theorem~3.6]{Isik} gives an equivalence between the categories of compact objects
\[\Coh(Z)\simeq\Coh(\CH)/\QCoh(\CH)^{\on{perf}},\]
rather than between their ind-completions, as in Theorem~\ref{t:I}. 
\end{rem}

\sssec{Singular support via category of singularities}
The category $\IndCoh(\CH)/\QCoh(\CH)$ is naturally tensored over the monoidal category $\QCoh(\CH)$.
Using the natural morphism
\[\CH\to(V^*/\BG_m)\times\CU,\]
we can consider $\IndCoh(\CH)/\QCoh(\CH)$ as a category tensored over $\QCoh((V^*/\BG_m)\times\CU)$. 
Recall from \secref{sss:global intersection parallelized} that the category $\IndCoh(Z)$ is tensored over the category 
$\QCoh((V^*/\BG_m)\times\CU)$ as well. It is not hard to check that \eqref{e:I} is an equivalence
of $\QCoh((V^*/\BG_m)\times\CU)$-modules. 

\medskip

In particular, fix $\CF\in\IndCoh(Z)$ and let $\CF'\in\IndCoh(\CH)/\QCoh(\CH)$ be its image under \eqref{e:I}. 
We claim that 
\begin{equation}\label{e:supp}
\on{SingSupp}(\CF)=\on{supp}(\CF').
\end{equation}
Note that 
\[\on{SingSupp}(\CF)\subset\Sing(Z)\subset V^*\times\CU,\]
while the support of $\CF'$ can be defined naively, as the minimal closed subset
of $\CH$ (that is, a conical Zariski-closed subset of $H\subset V^*\times U$) such that $\CF'$ restricts to zero on its
complement. It is clear that for any $\CF'\in\IndCoh(\CH)/\QCoh(\CH)$, its support is a conical Zariski-closed subset
of the singular locus of $H$ (recall that the singular locus of $H$ is identified with $\Sing(Z)$).

\ssec{The category of singularities of $Z$}

\sssec{}
Denote by 
\[(\IndCoh(\CH)/\QCoh(\CH))_{\{0\}}\subset \IndCoh(\CH)/\QCoh(\CH)\]
the full subcategory of objects of $\IndCoh(\CH)/\QCoh(\CH)$
supported on the zero-section $$\{0\}\times\CU\subset V^*\times\CU.$$
Under the equivalence \eqref{e:I}, it corresponds to the full subcategory
$\QCoh(Z)\subset\IndCoh(Z)$ (where we identify $\QCoh(Z)$ with its image under $\Xi_Z$).  
This claim is not hard to check directly, but
it also follows from \eqref{e:supp}: indeed, $\CF\in\IndCoh(Z)$ belongs to the essential image of $\QCoh(Z)$
if and only if its singular support is contained in the zero-section (\thmref{t:zero sect}).

\sssec{}

Therefore, \thmref{t:I} induces an equivalence between the quotients
\[\IndCoh(Z)/\QCoh(Z)\simeq(\IndCoh(\CH)/\QCoh(\CH))/(\IndCoh(\CH)/\QCoh(\CH))_{\{0\}}.\]
Set
\[\CH':=\CH-(\{0\}/\BG_m)\times\CU\subset\CH.\]
Note that $\CH'$ is a DG scheme rather than a stack (in fact, $\CH'$ is a classical scheme unless the map
$s:\CU\to V$ vanishes on a connected component of $\CU$).
We can identify
\[\IndCoh(\CH')/\QCoh(\CH')\simeq (\IndCoh(\CH)/\QCoh(\CH))/(\IndCoh(\CH)/\QCoh(\CH))_{\{0\}}\]
(cf. \cite[Proposition~6.9]{KrS}). Therefore, \thmref{t:I} implies the following equivalence, due to D.~Orlov.

\begin{thm}\label{t:O} There is a natural equivalence
\[\IndCoh(Z)/\QCoh(Z)\simeq\IndCoh(\CH')/\QCoh(\CH').\]
\end{thm}

\thmref{t:O} is a variant of \cite[Theorem~2.1]{Or2}. Some minor differences include that \cite{Or2} works with the
category of compact objects, and assumes that $Z$ is classical. Besides, the equivalence of \thmref{t:O} is constructed
in a different and more explicit way than the equivalence of \thmref{t:I}; in fact, while the introduction to \cite{Isik} 
mentions the similarity between the two results, it also states that the agreement between the two constructions is not
immediately clear.

\begin{rem}\label{r:support without zero}
 Fix $\CF\in\IndCoh(Z)$. Just like \thmref{t:I} can be used to determine $\on{SingSupp}(\CF)$ (using \eqref{e:supp}), \thmref{t:O} determines 
\[\on{SingSupp}(\CF)\cap(V^*-\{0\})\times\CU,\]
that is, the complement to the zero-section in $\on{SingSupp}(\CF)$. However, one can easily reconstruct the entire singular support, because
\[\on{SingSupp}(\CF)\cap\{0\}\times\CU=\{0\}\times\on{supp}(\CF).\]
\end{rem}

\sssec{} Let $Y$ be a conical Zariski-closed subset of $\Sing(Z)$ that contains the zero-section. 
Such subsets are in one-to-one correspondence with Zariski-closed subsets of the singular locus of $\CH'$: the correspondence sends $Y$ to
\[Y':=(Y-\{0\}\times\CU)/\BG_m\subset\CH'.\]

\medskip

Since $Y$ contains the zero-section, the corresponding full 
subcategory $\IndCoh_Y(Z)$ contains $\QCoh(Z)$. Therefore, we can consider the quotient $\IndCoh_Y(Z)/\QCoh(Z)$,
which embeds as a full subcategory into $\IndCoh(Z)/\QCoh(Z)$.

\medskip

G.~Stevenson provides a complete classification of localizing subcategories of the triangulated category
$\on{Ho}(\IndCoh(Z)/\QCoh(Z))$ in \cite[Corollary~10.5]{St} (a triangulated subcategory is localizing if it is closed under
arbitrary direct sums). Such subcategories are in one-to-one correspondence with subsets of $\CH'$.
Subcategories that are generated by objects that are compact in $\on{Ho}(\IndCoh(Z)/\QCoh(Z))$ correspond to
specialization-closed subsets. Under this correspondence, the category
$\on{Ho}(\IndCoh_Y(Z)/\QCoh(Z))$ corresponds to the subset $Y'\subset\CH'$. (This is almost obvious because \cite{St} 
uses Orlov's equivalence of \thmref{t:O} to study $\IndCoh(Z)/\QCoh(Z)$.)

\medskip

To summarize, we can also reconstruct the singular support of an object $\CF\in\IndCoh(Z)$ using the results of \cite{St}.
Technically, this only allows us to reconstruct the support of the image of $\CF$ in $\IndCoh(Z)/\QCoh(Z)$, that is,
the complement to the zero-section in $\on{SingSupp}(\CF)$; the entire $\on{SingSupp}$ can be recovered using
\remref{r:support without zero}.

\ssec{$\IndCoh(Z)$ as the coderived category}

Let us comment on the relation between the results in the main body of the present paper and the non-linear Koszul transform introduced by 
L.~Positselski in \cite{Pos}.

\sssec{} The key notion that we need from \cite{Pos} is that of coderived category of modules over a curved DG algebra.
To simplify the exposition, we do not give the definitions in maximal generality, and only work with curved DG algebras whose 
curvature is central. This class of curved DG algebra suffices for our purposes.

\medskip

Let $A$ be a DG algebra. Fix a central element $c\in A^{2}$ such that $d(c)=0$. We refer to the pair $(A,c)$ as a 
``curved DG algebra"; $c$ is called the curvature of $A$. 

\medskip

A (left) module over the curved DG algebra $(A,c)$ is by definition a graded vector space $M$ equipped with an action of
$A$ and a degree one map $d:M\to M$ that satisfies the Leibniz rule and the identity $d^2=c$. Modules over $(A,c)$
form a DG category, which we denote by $A\mod_c$. Consider the corresponding triangulated category $\on{Ho}(A\mod_c)$.

\begin{defn} The full subcategory of coacyclic modules 
\[\on{Acycl}^{co}(A\mod_c)\subset\on{Ho}(A\mod_c)\] 
is the subcategory generated by
the total complexes of exact sequences 
\[0\to M_1\to M_2\to M_3\to 0\]
of $(A,c)$-modules. The coderived category $D^{co}(A\mod_c)$ is defined to be the quotient 
\[D^{co}(A\mod_c):=\on{Ho}(A\mod_c)/\on{Acycl}^{co}(A\mod_c)\]
\end{defn}

\sssec{} Suppose $Z$ is as in \secref{sss:parallel gci}. We consider two curved DG-algebras: one is
\[A:=\Sym(V^*)\otimes\Gamma(\CU,\CO_\CU)\]
with differential given by $s\in\Gamma(\CU,\CO_\CU\otimes V)$ and curvature zero. The other is
\[B:=\Sym(V[-2])\otimes\Gamma(\CU,\CO_\CU)\]
with zero differential and curvature $s$.

\medskip

The following is a variant of a special case of \cite[Theorem~6.5a]{Pos} (see also \cite[Theorem~6.3a]{Pos})
\begin{thm}\label{t:P}
There is a natural equivalence between the coderived categories
\[D^{co}(A\mod_0)\simeq D^{co}(B\mod_s).\]
\end{thm}

\begin{rem} We state \thmref{t:P} with algebras on both sides of the equivalence. This is one point of difference
from \cite{Pos}, where the equivalence relates algebras and coalgebras. However, note that $A$ is free of
finite rank over $\Gamma(\CU,\CO_\CU)$, so it is easy to pass from modules over $A$ to comodules over the dual coalgebra.
Another point of difference with \cite{Pos} is that Theorem~\ref{t:P} is ``relative" the correspondence is linear over
the algebra $\Gamma(\CU,\CO_\CU)$.
\end{rem}

\sssec{} Let us explain the relation between \thmref{t:P} and \thmref{t:I}.

\medskip
 
Indeed, $Z=\on{Spec}(A)$, and it follows from \cite[Theorem~3.11.2]{Pos} that $D^{co}(A\mod_0)$ can be identified with 
$\on{Ho}(\IndCoh(Z))$. 
On the other hand,
$(B,s)$-modules are similar to ``equivariant matrix factorization" (cf. \cite[Example~3.11]{Pos}), and it is natural that
they can be used to study the ``equivariant category of singularities" $\IndCoh(\CH)/\QCoh(\CH)$. Finally, note that,
just as the equivalence of \thmref{t:P} is given by a (non-linear) Koszul transform, the equivalence of \thmref{t:I}
(constructed in \cite{Isik}) is derived using a (linear) Koszul transform, namely, the linear Koszul transform of I.~Mirkovi\'c and 
S.~Riche \cite{MR}.

\newcommand{\eprint}[1]{Preprint {\tt #1}}\newcommand{\available}[1]{Available
  from \url{#1}}

\end{document}